\documentclass[a4paper,11pt]{book}
\usepackage[latin1]{inputenc}
\usepackage{fancyhdr}
\usepackage[OT1]{fontenc}
\usepackage{amssymb,amsmath,amsfonts,amssymb,graphicx,color,bbm}
\usepackage{comment}
\usepackage{amssymb}

\usepackage[colorlinks=true,
linkcolor=blue,citecolor=blue]{hyperref}

\usepackage{amsmath,amsfonts,xcolor} 

\newtheorem{thm}{thm}[section]   
\newtheorem{theorem}[thm]{Theorem}   
\newtheorem{problem}[thm]{Problem}   
\newtheorem{corollary}[thm]{Corollary}   
\newtheorem{proposition}[thm]{Proposition}   
\newtheorem{lemma}[thm]{Lemma}  

\newtheorem{remark}[thm]{Remark}  

\newcommand{\RR}{\mathbb{R}}   
\newcommand{\di}{\displaystyle}   
\newcommand{\CC}{\mathbb{C}}      
    
\newcommand{\NN}{\mathbb{N}}
\newcommand{\LL}{\mathcal{L}}
\newcommand{\TT}{\mathcal{T}}
\newcommand{\QQ}{\mathcal{Q}}
\newcommand{\JJ}{\mathcal{J}}
\newcommand{\II}{\mathcal{I}}

\newcommand{\MM}{\mathcal{M}}
\newcommand{\KK}{\mathcal{K}}

\newcommand{\DD}{{\cal D}}
\newcommand{\EE}{{\cal E}}

\newcommand{\un}{\mathbf{1}}
\newcommand{\e}{\varepsilon}
\def\un{{\mathbf{1}}}
\def\K{<\!K\!>}

\numberwithin{equation}{chapter}
\numberwithin{equation}{section}
\begin{document}

\title{\Huge The dynamics of front propagation in nonlocal reaction-diffusion equations}
\author{
Jean-Michel Roquejoffre\footnote{
Institut de Math\'ematiques,
Universit\'e Toulouse 3 - Paul Sabatier, 118 Route de Narbonne, 31062 Toulouse Cedex,
France. Email: jean-michel.roquejoffre@math.univ-toulouse.fr}}
\date{December 4, 2023}
\maketitle
\
\vfill\eject
\vglue10cm
\hfill To Oxana, my wife
\vfill\eject
\
\vfill\eject
\vglue2.5cm
\begin{center}
{\Large\bf Acknowledgements}
\end{center}

\vglue.35cm
\noindent In the first place, I want to thank the University of St Petersburg, the Chebyshev Laboratory and the Hadamard Foundation for awarding me the 2020 Gabriel Lam\'e Chair. The research visit took place in the fall 2021 and was, despite covid restrictions, a truly exceptional time. My thanks go  to F. Bakharev, A. Baranov and S. Tikhomirov, for optimal scientific conditions, and to M. Lantsova and M. Rossomakhina  for invaluable administrative support.  I am also grateful to A. Malac and the French Embassy in Russia, for providing all the help they could.

\medskip
\noindent One of the duties of the chair was to give a course,  I chose  front propagation in nonlocal models for biological invasions as the general topic.  This is the occasion to thank   the students and colleagues who attended the course, in presence or online.  As I had not initially thought about writing lecture notes, I am grateful to B. Perthame for suggesting me to write them in the form of a book. Thanks also to  A. Stevens, for her support of the idea. 
The scope of the book  was  stabilised at the occasion of another course given at the Mathematical Institute of Università Roma I at the end of spring 2022. I wish to thank to I. Birindelli and L. Rossi them for their invitation, and also for attending the course and asking a wealth of relevant questions. Thank you also to  the students who attended, and who showed interest from beginning to end.  Some parts were also taught in a master course 
given at my home university, together with G. Faye. I acknowledge the numerous discussions that we had, as well as the active participation of the students.
 
 \medskip
 \noindent Writing this book turned out to be a longhaul task that would not have been possible, or would have at least been  much more difficult, without interaction with many colleagues. I want to single out Y. Petrova and M. Zhang, who read parts of it; this allowed me to correct sometimes important mistakes. I am also indebted to all those who provided references that were unknown to me,  pointed out at an occasional error, suggested a shorter proof, raised an issue I had not immediately thought about, shared their knowledge with me, invited me for a talk - thus obliging me to clarify the exposition of some ideas - or simply made a constructive comment. The contribution of H. Berestycki, J. Berestycki, N. Boutillon, A.-C. Chalmin, O. Diekmann, G. Faye, C. Graham, F. Hamel, P. Maillard, B. Mallein, S. Mirrahimi, M. Pain, L. Rossi, L. Ryzhik\ldots, to the raising of my competence on the subject, and, as a consequence, to the improvement of this work,  is gratefully acknowledged.
 
 \bigskip
 \noindent This work has received support from the ANR project ReaCh, Grant ANR-23-CE40-0023-01.

\vfill\eject

\vfill\eject

\tableofcontents

\chapter*{Abstract}
\noindent The question addressed here is the long time evolution of the solutions to a class of one-dimensional reaction-diffusion equations, in which the diffusion is given by an integral operator. The underlying motivation, discussed in the first chapter,  is the mathematical analysis of models for biological invasions. The model under study, while   simple looking, is of current use in real life situations. Interestingly, it arises in totally different contexts, such as the study of branching random walks in probability theory.

\noindent While the model under study has attracted a lot of attention, and while many partial results about the time asymptotic behaviour of its solutions have been proved over the last decades, some basic questions on the sharp asymptotics have remained unanswered. One ambition of this monograph is to close these gaps and to provide a complete and unified treatment of the equation. 

\noindent In some of the situations that we envisage, the level sets  organise themselves into an invasion front that is asymptotically linear in time, up to a correction that converges exponentially in time to a constant. In other situations, that constitute the main and newest part of the work, the correction is asymptotically logarithmic in time. Despite these apparent different behaviours, there is an underlying common way of thinking in the study of all these situations.

\noindent The ideas presented in the book apply to more elaborate systems modelling biological invasions or the spatial propagation of epidemics. The models themselves may be multidimensional, but they all have in common a mechanism imposing the propagation in a given direction; examples are presented the problems that conclude each chapter. These ideas should also be useful in the treatment of further models that we are not able envisage at the time being.
 \chapter{Introduction}\label{Intro}
\section{Orientation of the book}
\noindent The objective is to provide a complete and self-contained  study of the asymptotic behaviour, as $t\to+\infty$, positive solutions $u(t,x)$ of the following class of equations
\begin{equation}
\label{e1.1.1}
u_t+\int_\RR K(x-y)\bigl(u(t,x)-u(t,y)\bigl)~dy=f(u),\quad t>0,x\in\RR.
\end{equation}
The kernel $K$ is smooth, nonnegative, even and compactly supported. The function $f$ is nonnegative, positive in the interval $(0,1)$, with $f(0)=f(1)=0$.
The initial datum $u(0,.)$ is a small nonnegative compactly supported function. Of course there is an underlying motivation that will be explained in more detailed in the subsequent sections of this introduction.

\noindent While Model \eqref{e1.1.1} has been much studied in the literature, especially from the point of view of travelling waves,  and while
many interesting, yet partial  results about the time asymptotic behaviour of its solutions have been proved over the last 50 years, some basic questions about the sharp asymptotics have remained unanswered. One ambition of this monograph is to close those gaps and to provide a unified treatment of this equation. It is to be hoped that some of the methods displayed will be useful to the study of more general models.

\noindent  For the time being, let us try heuristically to understand what the large time behaviour of $u(t,x)$ could resemble. For this, let us think for one moment that we are integrating  Model \eqref{e1.1.1} numerically:  in other words, we are trying to compute a sequence of functions $\bigl(u^n(x)\bigl)_{n\in\NN}$, where $u^n(x)$ is an approximation of $u(n\Delta t,x)$. The time step $\Delta t$  should fulfill  conditions on 
which we are not going to insist. An idea, among many, is to use a splitting method.
At even iterations, integrate the "diffusion" equation $u_t+\di\int_\RR K(x-y)\bigl(u(t,x)-u(t,y)\bigl)~dy=0$:
\begin{equation}
\label{e1.1.2}
u^{2n}=u^{2n-1}+\Delta t\int_\RR K(x-y)\bigl(u^{2n-1}(x)-u^{2n-1}(y)\bigl)~dy.
\end{equation}
We will see later that this is really a diffusion equation.
 At  odd iterations, simulate the differential equation $\dot u=f(u)$:
\begin{equation}
\label{e1.1.3}
u^{2n+1}=u^{2n}+\Delta tf\bigl(u^{2n}\bigl).
\end{equation}

\noindent If one believes that \eqref{e1.1.2} is really a diffusion equation, the intuition is that $u^{2n}$ is a sort of average of $u^{2n-1}$ in the vicinity of each of its points. On the other hand, \eqref{e1.1.3} is a typical reaction equation, in the sense that it pushes its solution to 1. There is therefore a succession of steps where the sequence $(u^n)_n$ is alternatively spread and pushed to 1: this entails the spatial invasion of the whole line by the stable state $u\equiv 1$. Let us take this for granted and let us come back to the solution $u(t,x)$ of the continuous problem:  intuitively, at each time $t$, the function $u(t,x)$ tends to 0 as $\vert x\vert$ tends to infinity. Therefore, a transition develops between the zone $u(t,x)\sim 1$ (that invades the whole domain) and the (receding) zone $u\sim 0$. So, for $\theta\in(0,1)$, one may define the furthest point $x$ to the right such that $u(t,x)=\theta$; call it $X_\theta(t)$. If we believe the preceding heuristics, then $X_\theta(t)$ tends to $+\infty$ as $t\to+\infty$, and the whole object of the monograph is to understand this in quantitative, rigorous terms. This will lead to an exploration of various questions touching comparison arguments, analysis of travelling waves, spectral analysis, approximation of Dirichlet heat kernels... and eventually enable us to discover an asymptotic expansion of $X_\theta(t)$ with precision $o_{t\to+\infty}(t)$. All these questions will be studied in a self-contained fashion: the various tools needed will be introduced without assuming any previous knowledge of them (although a background in diffusion equations cannot harm), and will first be used to reprove known results. However, when it comes to studying issues not previously answered, they will be at our disposition.

\noindent The underlying motivations of such an endeavour, namely the modelling of biological invasions, are presented in Section \ref{s1.2}. More elaborate multi-dimensional models, whose study still falls within the scope of the argument developped in the present monograph, are presented in Section \ref{s1.4}. The general organisation of this work is presented in Section \ref{s1.3}. 

\section{Two instances of biological invasions}\label{s1.2}
\noindent  When a living species of any sort: animal, vegetal, bacterium or virus, invades a given habitat: piece of land, sea, river..., one speaks of a  biological invasion. While it may sound aggressive, the term often describes a slow and peaceful process, such as 
the dissemination of a tree species. Of course the aggressive  aspect is present in phenomena occurring on a faster time scale, with dangerous species, such as rodents, mosquitoes, hornets... the ambition of 
this section is to explain how one can model, in a simple way, how reproduction and displacements work together and mechanically produce the invasion. More precisely, two processes interacting together will be presented: 
invasions driven by reproduction $+$ displacements, and invasions driven by remote contaminations.
\subsection{Invasions driven by reproduction and displacements}\label{s1.2.1}
\noindent Consider a population for which we assume that a notion of density is relevant; this in itself is an important assumption, not always justified. Let us proceed anyway, and call $u(t,x)$ the density at time $t$ and point $x$. Let us model the two mechanisms.
\begin{itemize}
\item[---] {\it Displacements.} The modelling assumption is that the proportion of individuals that will move, at time $t$ and point $x$, from $x$ to point $y$, is a time independent function of $\vert x-y\vert$, that is, in the case of one space dimension, an even function of $x-y$. The travel time is neglected. As  the individuals cannot not move instantly to infinity, it is legitimate to assume $K$ to be compactly supported.
\item[---] {\it Reproduction.} Assume, here only, the population to be homogeneous: $u(t,x)\equiv u(t)$. The time variation of the density  $\dot u(t)$, is assumed to be a time independent function $f(u)$ of the density. It is natural to assume it proportional to $u$: the more individuals are present, the more reproduction takes place. It is also natural to assume the proportionality coefficient to be negative for large values of $u$: too large a density is not sustainable. There is therefore a threshold, called the {\it carrying capacity}, such that the reproduction rate is $>0$ below this threshold, and $<0$ above. Assuming the carrying capacity to be equal to 1, the simplest relevant nonlinearity is $f(u)=u(1-u)$.
\end{itemize}
\noindent Coming back to a density that is not homogeneous in $x$ anymore, we want to express the fact that the variation of $u(t,x)$ is equal to  (i) the number of individuals that have travelled to $x$ at time $t$, minus (ii) the number of individuals leaving the point $x$ at time $t$, plus (iii) the number of individuals created/destroyed by the reproduction process. Term (i) is obviously $\di\int_\RR K(x-y)u(t,y)~dy$ (we have to count the individuals arriving from all over the world); term (ii) is, not less obviously, $u(t,x)\di\int_\RR K(x-y)dy$, and (iii) is $f(u)$. All in all, this yields Equation \eqref{e1.1.1}.

\noindent There is a close link between the nonlocal operator $\mathcal{J}u:=\di\int_\RR K(x-y)\bigl(u(x)-u(y)\bigl)~dy$ and the more usual diffusion operator $-\partial_{xx}u$. To see it, take for the kernel $K$ an approximation of the identity $K_\e(x)=\di\frac1\e \rho(\di\frac{x}\e)$, where $\rho$ is a smooth nonnegative even function supported in $(-1,1)$. We have, for any smooth function $u(x)$:
$$
\int_\RR K_\e(x-y)\bigl(u(x)-u(y)\bigl)~dy=\int_\RR\rho(z)\bigl(u(x)-u(x+\e z)\bigl)~dz=-\frac{d\e^2}2\bigl(-u''(x)+O(\e)\bigl),
$$
the number $d>0$ being the second moment of $\rho$.
Plugging this identity into \eqref{e1.1.1}, choosing a reproduction term $f(u)$ of the form $\e^2 g(u)$, accelerating time as $\tau=\e^2 t$, and finally throwing away the terms of order $\e^3$ and more, yields the equation
\begin{equation}
\label{e1.2.4}
\partial_\tau u-\partial_{xx}u=g(u),\quad (\tau>0,x\in\RR).
\end{equation}
This makes an important case for calling the operator $\mathcal{J}$ a diffusion operator, it was indeed not so obvious to see why at first sight. This also explains why the results that we are going to prove in the sequel are closely related to those pertaining to \eqref{e1.2.4}, as this model is really a distinguished limit of \eqref{e1.1.1}.

\noindent The model function $f(u)=u(1-u)$ falls in the class of nonlinearities that are said to be of the {\it Fisher-KPP} type, as a tribute to the seminal contributions of Fisher \cite{Fisher}, who probably produced the first analytical study of \eqref{e1.2.4}, and 
the extraordinary contribution of Kolmogorov-Petrovskii-Piskunov \cite{KPP}, that provides a deep insight in the asymptotic behaviour of $X_\theta(t)$ (they prove that it is asymptotically linear) as well as the asymptotic profile of $u(t,x)$. This class of nonlinearities is characterised by the property $f(u)\leq f'(0)u$, that is, the reproduction rate is maximised when the population is at its lowest. By extension, any phenomenon or mathematical feature pertaining to such nonlinearities will be, throughout the book, baptised "KPP", or "Fisher-KPP", or "of the Fisher-KPP type".

\noindent Not all functions $f$ satisfiy,  however, $f(u)\leq f'(0)u$. A very small reproduction rate  when the population density is very small is also a perfectly reasonable assumption. This legitimates another class of nonlinearities, that we will call the {\it ZFK} type terms. This time, the acronym pays a tribute to the physicists Zeldovitch and Frank-Kamemetskii. Far from the modelling in ecology, they were interested in understanding flame propagation, and nonlocal diffusive effects were not the primary of their concerns. They 
 have introduced such nonlinearities - theirs is essentially maximal when $u$ is very close to $1$ - in order to account for the fact that a propagating flame  profile has various characteristic lengths.  For us, a good paradigm for such nonlinearities will be that of a function $f$ having a small slope at 0, but a large mass. All this, and much more, is explained in \cite{ZBLM}. All these nonlinearities and, by extension, all features pertaining to them will be called "ZFK", or "of the ZFK type". 
\subsection{Invasions driven by remote contaminations}\label{s1.2.4}
\noindent In other words, we want to describe, in stylised terms, how an epidemic wave forms and spreads.To have an even faint grasp of the early stages, statistical tools are indispensable. However, deterministic modelling can give a lot of qualitative information about the stages immediately after the start, and this is the line that I choose to develop.  One has, essentially, to account for two mechanisms: contamination and diffusion. 
 \begin{itemize}
\item [---] {\it Contamination.} The most basic   mechanism can be described  by the SI model: in a homogeneous population, the density of individuals exposed to an infection (The {\it Susceptibles}) and those that have caught it (the {\it Infected}; hence the acronym for the system) vary according to the following common sense rules: the number of individuals that become infected at time $t$ is proportional both to the susceptible population, and also the number of infected. The proportionality coefficient depends on the contaminating power of the virus, the receptibility of the susceptibles, the time or possible exterior events such as the introduction of a vaccine... thus it is a nontrivial task to model it with with some precision. We wisely declare it constant, equal to $\beta>0$. The number of new infected are discounted from the susceptible population. Individuals leave the sickness (either by healing or in a more tragic way) at a rate that it is wise to assume proportional to the number of infected, and that we once again wisely assume to be constant. Thus the system reads
 \begin{equation}
 \label{e1.2.2}
 \left\{
 \begin{array}{rll}
 \dot S=&-\beta SI,\ \ \  \ \ \ \ \ 
 \dot I=\beta SI-\gamma I\\
 S(0)=&S_0>0,\ I(0)=I_0>0\ \hbox{small}.
 \end{array}
 \right.
 \end{equation}
The cumulated density of infected 
$u(t)=\di\int_0^tI(s)~ds
$ 
solves the differential equation
 \begin{equation}
 \label{e1.2.1}
 \dot u=S_0(1-e^{-\beta u})-\gamma u+I_0,\quad
 u(0)=0.
 \end{equation}
This simple ODE can be easily analysed, let us do it. Assume, as there is no spatial effect, that 
$I_0$ does not depend on $x$, so that $u$ does not depend on $x$ anymore, that is $u(t,x)\equiv u(t)$. For commodity,  set $f(u)=S_0(1-e^{-u})-\alpha u$. It is a concave function, and an important quantity to assess the large time behaviour of $u(t)$ 
is the sign of $f'(0)=S_0\beta-\alpha$. This legitimates the definition of 
\begin{equation}
\label{e1.2.3}
R_0={S_0\beta}/\alpha,
\end{equation}
so that $f'(0)=\alpha(R_0-1)$. Therefore, if $R_0\leq1$, the only zero of $f$ is $u\equiv0$, so that the equation $f(u)+I_0=0$ has a unique zero which is asymptotically proportional to $I_0$ if $R_0<1$ and to $\sqrt I_0$ if $R_0=1$. If $R_0>1$, then $f$ has a positive zero that we denote by $u_\infty^*$, and the equation $f(u)+I_0=0$ has a unique zero $u_\infty(I_0)>u_\infty^*$. Now the analysis of \eqref{e1.2.2} reveals that 
$\di
\lim_{t\to+\infty}=u_\infty(I_0)
$
for all values of $R_0$. However, the important differences between the cases $R_0<1$ and $R_0>1$ lie in the fact that $u_\infty(I_0)$ is bounded away from zero by $u_\infty^*$, no matter how small $I_0$ is, when $R_0>1$. As opposed to this, 
we have 
$
\di\lim_{I_0\to0}u_\infty(I_0)=0
$
if $R_0\leq1$. The consequence of this fact can be seen on the so-called final size of the population, that is, the limit of the density of susceptibles as time goes to $+\infty$. 
\item[---] {\it Diffusion.} Assume now the initial density of infected $I_0$ to depend on $x$. The basic situation that we want to describe is what happens if $I_0$ is supported in a bounded region, that is, the infected are not only introduced by a tiny amount, but they are also very much localised in space. One easily sees that, without any further mechanism, there is absolutely no chance that  Model \eqref{e1.2.2} will entail any sort of propagation. One way to make up for this is to assume that
 an infected person is contagious within a nontrivial spatial range. A very crude way to model it is to declare that an individual located at point $x$ will contaminate the individuals located at $y$ in some proportion. It is almost impossible to model this proportion with sufficient precision, at il may depend on the time of the day/week/year, temperature, ability of individuals to move, contaminating power of the virus... So, a wise attitude is to declare this proportion to be a function of the distance $\vert x-y\vert$, and
  that it is all the smaller as $\vert x-y\vert$ is large. So we take it, once and for all, compactly supported. And so, here comes again our favourite kernel $K$ and we replace, in (\ref{e1.2.2}), the term $SI$ by $(K*I)S$. 
  Thus, \eqref{e1.2.2} becomes
   \begin{equation}
 \label{e1.2.5}
 \left\{
 \begin{array}{rll}
 \dot S=&-\beta S(K*I),\ \ \ \ \ \ \ \ \ 
 \dot I=\beta SI(K*I)-\gamma I\\
 S(0)=&S_0>0,\ I(0)=I_0>0\ \hbox{small, compactly supported}.
 \end{array}
 \right.
 \end{equation}
The cumulative density  solves a Fisher-KPP type equation: 
 \begin{equation}
\label{e1.2.6}
u_t=S_0(1-e^{-\beta K*u})-\alpha u+I_0(x).
\end{equation}
\noindent  
\end{itemize}
In each of the subsequent chapters \ref{Cauchy} to \ref{short_range}, there will be a problem developping the study of \eqref{e1.2.5} and \eqref{e1.2.6}. The model can therefore be considered as completely understood.
\section{Models in several space dimensions with directed propagation}\label{s1.4}
\noindent The book will be focussed on one-dimensional models. However, the ideas developped to answer the questions listed in Section \ref{s1.3} can be fruitful for the understanding of   models posed in multi-dimensional geometries but where front propagation occurs in only one spatial direction. Loosely speaking, there may be a coefficient in the equation, or a heterogeneity, that directs the propagation. While the book does not, strictly speaking,
 present the analysis of the class of systems \eqref{e1.3.1}, these models can be understood essentially completely with the aid of the methods displayed in it. Two of them are presented in this section, and they are far from exhausting the vast subject of propagation driven by a one-dimensional structure.  The author of these lines believe  that, with the aid of some tools specific to second-order parabolic PDEs such as the Harnack inequality or spectral decomposition theorems, they can be studied completely with the ideas developped in this book. The beginning of the study of one of them is proposed in Problem \ref{P2.120}, others are proposed in the problems of Section \ref{ZFK_short_range}.
 \subsubsection*{Models in cylinders or with confining mechanisms}
 \noindent Let $\Sigma$ be the   cylinder $\{(x,y)\in\RR\times\omega\}$, where $\omega$ is a smooth bounded open subset of $\RR^{N-1}$. Let $\alpha(y)$ be a smooth function in $\overline\omega$. Consider the model, with unknown $u(t,x,y)$:
 \begin{equation}
 \label{e1.3.7}
 \left\{
 \begin{array}{rll}
 u_t+\mathcal{J}_xu-\Delta_y u+\alpha(y)u_x=&f(u),\quad \bigl(t>0,\ (x,y)\in\Sigma\bigl)\\
 \partial_\nu u=&0\quad\bigl((x,y)\in\partial\omega)
 \end{array}
 \right.
 \end{equation}
 Here, $\mathcal{J}_x$ means that our usual operator $\mathcal{J}$ is applied only to the $x$ variable, that is, to the function $x\mapsto u(t,x,y)$ with the variables $t$ and $y$ frozen. The operator $\Delta_y$ is the Laplacian in the variable $y$. 
 The function   $f(u)$ is smooth, of the Fisher-KPP or ZFK  type. 
 
 \noindent System \eqref{e1.3.7} is a stylised model that describes the influence of an external, imposed flow field on the dissemination of a population of individuals in a habitat; think, for instance, of a fish population in a river or a colony of insects submitted to the wind. The Neumann condition $\partial_yu=0$ at the boundary means that no individual may leave the domain; it is not the only one: a Dirichlet condition $u(t,x,y)=0$ in $(x,y)\in\omega$ is also legitimate from the modelling point of view. It means that the boundary is lethal for the individuals. The longitudinal diffusion can of course be replaced by the standard diffusion $-\partial_{xx}$, in this case, another situation of which \eqref{e1.3.7} is a stylised representation is the propagation of a flame submitted to an external advection field.
 
 \noindent Here, the propagation is dictated by the fact that we operate in a cylinder.  So, the tools and arguments developped in the sequel of this book, together with some adjustments that take into account the presence of the variable $y$ are, to the author's opinion, sufficient to carry out a
  through study if Model \eqref{e1.3.7}.
 
 \medskip 
 \noindent A related model, this time in the whole space, that we take to be $\RR^2$ without real loss of generality, but with a confining mechanism, is the following: 
 \begin{equation}
 \label{e1.3.8}
  u_t+\mathcal{J}_xu-\partial_{yy}u+y^2u=f(u),\quad \bigl(t>0,\ (x,y)\in\RR^2\bigl).
  \end{equation}
Here again, the operator $\mathcal{J}_x$ could be replaced by $-\partial_{xx}$, or the whole diffusion operator could be modelled by $u-K*u$, the kernel $K$ being this time radial, and the convolution acting on both variables $x$ and $y$. The underlying situation is that of the dissemination of a population characterised not only by its spatial density, but by one or several phylogenetic traits. In the present situation, $x$ is the spatial variable and $y $ accounts for the trait. Mutations are accouned for by the Laplacian in $y$, 
moreover the term $y^2u$ models the fact that the trait $y=0$ is the preferred one, as it presents the minimal losses. This is why we talk about a confining mechanism: one can rather easily show that there is no propagation in the direction $y$. An interesting variant is when the preferred trait is not $y=0$, but rather a trait $\beta(x)$ that depends on the space. In such a case the equation for $u$ reads
 \begin{equation}
 \label{e1.3.9}
  u_t+\mathcal{J}_xu-\partial_{yy}u+\bigl(y-\beta(x)\bigl)^2u=f(u),\quad \bigl(t>0,\ (x,y)\in\RR^2\bigl).
  \end{equation}
While Model \eqref{e1.3.8} is amenable to the methods developped in this book, model \eqref{e1.3.9} certainly requires new ideas if one wishes to understand the precise asymptotic behabviour of the solutions, even 
in the seemingly innocent configuration where $\beta$ is a bounded smooth function. 
\subsubsection*{Propagation directed by a line of fast, or nonlocal diffusion}
\noindent Let $\Omega_L$ be the strip $\{(x,y)\in\mathbb R \times (-L,0)\}$.  Consider the system, with unknowns $\bigl(u(t,x),v(t,x,y)\bigl)$:
\begin{equation}
\label{e1.3.1}
\left\{
\begin{array}{rll}
u_t-\mathcal{J}u+\mu u-v(t,x,0)=&0\  \  \  (t>0,x\in\mathbb R)\\
v_t-d\Delta v=&f(v)\  \  \  (t>0,(x,y)\in\Omega_L)\\
\  \\
dv_y(t,x,0)+v(t,x,0)=&\mu u(t,x)\  \  \  (t>0,x\in\mathbb R)\\
v_y(t,x,-L)=&0\  \  \  (t>0,x\in\mathbb R)\\
\end{array}
\right.
\end{equation}
 The positive numbers $\mu$ and $d$, are given. The function $f(v)$ is once again smooth, of the Fisher-KPP or ZFK  type.  Mathematically,  this model corresponds to a reaction-diffusion process in the strip, exchanging individuals with the upper line $\{y=0\}$ having a diffusion of its own. The diffusion on the line may assume many other forms, as it may be given by a standard Gaussian diffusion $-D\partial_{xx}$, or the fractional Laplacian.  One may even envisage more fancy boundary conditions, as the exchanges can occur on intervals, or only at points of the upper line... the boundary conditions have of course to be changed accordingly. 
 
 \noindent The underlying modelling motivation is to describe, in mathematical terms, how a line of fast diffusion (e.g. a road, a river, a railroad...) may drive a biological invasion.  And, as a matter of fact, that a line, or network of transportantion can enhance a biological invasion is a well-documented fact: one may for instance think of the pine processionary moth, an insect that thrives in pine forests and may cause health issues. It has invaded the whole Europe from the south. Climate change has long been assumed to be the main responsible of this fact, until scientists realised that the invasion speed was faster than what could be expected. The excellent European transportation network is now thought to have greatly facilitated the process. Another instance, concerning directly the homeland of the author of these lines, is the asian hornet. It is believed to have been introduced in France at the beginning of the century,  near the city of Bordeaux,  by a cargo ship coming from Asia. Its proliferation, tracked by entomologists over the years (and remarkably shown on interactive maps on the Museum d'Histoire Naturelle web site), shows that it rapidly penetrated the country by following the Garonne Valley, where a main railroad as well as an important road system run parallel to the river. From that, the dissemination occurred, efficiently but more slowly, in the other directions. It is therefore quite reasonable to postulate that the insect benefitted, at some stage, from all these lines of transportation.
 
 \noindent System \eqref{e1.3.1} is therefore a stylised model of the above situations. Its origin is, once again, discussed in the bibliographical notes. For the invasive species under considerations, the function $v(t,x,y)$ is the density of individuals in the domain $\Omega_L$, while $u(t,x)$ is the density of travelling individuals. One may by the way, take many more effects on the line: the individuals may experience a more severe mortality rate on the line, a situation that can be modelled by adding a term of the form $-\alpha u$ in the left handside of the first equation. Including transport by a flow, especially if the line models  a river, is relevant; the equation on the line becomes
 $$
 \partial_tu+\mathcal{J}u+q\partial_xu+\mu u-v(t,x,0)=0,
 $$
 the other equations being left unchanged. A seemingly innocent modification consists in replacing 
 the strip $\Omega_L$ by the lower half-plane. I claim that the situation is quite different, because there is a phenomenon of secondary propagation inside the domain that is outside the scope of the book. A class of models in the lower half plane that, 
 however, falls in the scope of the book, assumes that the line concentrates both fast diffusion and reproduction, while the lower half-plane is lethal. Assuming nonlocal diffusion on the line, the model reads
 \begin{equation}
\label{e1.3.6}
\left\{
\begin{array}{rll}
u_t-\mathcal{J}u+\mu u-v(t,x,0)=&f(u)\  \  \  (t>0,x\in\mathbb R)\\
v_t-d\Delta v+\alpha v=&0\  \  \  (t>0,(x,y)\in\Omega_L)\\
\ \\
dv_y(t,x,0)+v(t,x,0)=&\mu u(t,x)\  \  \  (t>0,x\in\mathbb R)\\
v_y(t,x,-L)=&0\  \  \  (t>0,x\in\mathbb R)\\
\end{array}
\right.
\end{equation}
As opposed to \eqref{e1.3.1} in the whole lower half plane, we are back to a globally one-dimensional propagation, and the model can be studied with the ideas presented for the basic nonlocal model \eqref{e1.1.1}. Here again, there is no real reason why the top boundary or the bottom boundary should be straight lines, an especially interesting question, outside once again of the scope of the book, is what happens when, for instance, the boundary $\{y=0\}$ becomes a curve $\{y=\varphi(x)\}$.  The function $\varphi$ may be well-behaved (for instance periodic or quasi-periodic) but nothing prevents a wilder behaviour: for instance there could be long bottlneck regions where the top boundary is close to the bottom one (that is, $\varphi(x)$ is verly close to $-L$ on very big regions) or, on the contrary, there could be very large regions where we have $\varphi\gg1$. All these situations require new ideas. 
 \section{Questions under study, organisation of the book}\label{s1.3}
 \noindent  We are back to Problem \eqref{e1.1.1}, and the basic issue is the following: when the initial density  is a function with bounded support, the unstable state 0 will be invaded (by the state $u\equiv 1$ if $\mu(x)\equiv 1$), and 
a transition will form between the region $\{u(t,x)\sim0\}$, and the region $\{u(t,x)\sim1\}$. The "speed of spreading" question is how fast this transition moves in each direction. In other words the issue is to find,
 a constant $c_*$ such that we have, uniformly on each compact in  $x$:
\begin{equation}
\label{e8.1}
\lim_{t\to+\infty}u\bigl(t,x\pm ct\bigl)=0\ \hbox{for all $c>c_*)$},\quad \liminf_{t\to+\infty}u(t,x\pm ct)>0\ \hbox{for all $c<c_*$}.
\end{equation}
And so, reintroducing the function $X_\theta(t)$, that is, the rightmost point $x$ such that $u(t,x)=\theta$, we have $\di\lim_{t\to+\infty}X_\theta(t)/t=c_*$. The next question concerns  the description of  this transition in terms of   asymptotic profiles, in other words: does $u(t,x)$ assume a universal profile if followed in a correct reference frame, that is, moving more or less as $\pm c_*t$
This is the issue of travelling waves. 

\medskip
\noindent Sharp asymptotics, that is, the prediction of the position $X_\theta(t)$ with arbitrary precision, at least to $o_{t\to+\infty}(1)$ precision, is the principal  goal of this monograph. We will encounter two types of behaviours: in one of them, there is an exponentially fast relaxation to linear propagation; in other words, there is $x_\infty\in \RR$ and $\omega>0$ such that 
$
X_\theta(t)=c_*t+x_\infty+O(e^{-\omega t}).
$

\noindent The other type of behaviour looks a little more strange: we will encounter situations where there is $\mu_*>0$ and $x_\infty\in\RR$ such that
$
X_\theta(t)=c_*t+x_\infty-\mu_*\mathrm{ln}~t+x_\infty+o_{t\to+\infty}(1).
$
Most certainly, the $o_{t\to+\infty}(1)$ does not decay exponentially fast. One of the objectives of this work is to sort out the cases when these behaviours occur, and to explain why.

\noindent  We say in Section In the sequel of this work, the attempt at relating the various findings on \eqref{e1.1.1} to the results pertaining to \eqref{e1.2.4} will be minimal, albeit not inexistent. There is a simple reason for that: the results that we are going to prove contain those pertaining to \eqref{e1.2.4}, and their proofs are quite a lot more involved. One notable exception is the Cauchy Problem: while its resolution does not take much effort in the nonlocal case, its resolution in the local case is more lengthy,  and is discussed in the problems section related to Chapter \ref{Cauchy}. There is no real reason to put it in the main body of this work, as it is quite a classical issue. On the other hand, the study of the nonlocal model takes much of its inspiration from the local one, and that   \eqref{e1.2.4} provides a lot of intuition to \eqref{e1.1.1}.  Sometimes, however, the intuition stops being helpful, and truly new arguments are needed. This is especially true for the sharp behaviour of the solutions of our problem.

\noindent One could question, in view of the underlying motivations, the relevance of devising sharp asymptotics in models for biological invasions or the spread of epidemics. Indeed, in these applications, the truly relevant quantity is the speed of spreading, something that does not require much efforts to figure out in Model \eqref{e1.1.1}. In real life situations, the uncertainty on the speed of spreading is so high that finding out the corrections seems to be a little help. True, but one should make one step aside. Once the applications have yielded the model, trying to understand all that it has to say is always a sound policy, and has a mathematical interest on its own. There is another case to make, more closely related to the applications: often,  especially in epidemiology, one needs to integrate numerically modes that are far from a mathematical reach, on very large time intervals. This needs specific schemes, and assessing their validity on large times is essential. For this, a good idea is to make sure that they reproduce, as faithfully as possible, the behaviour of simpler models, such as \eqref{e1.1.1}. In order to verify that, the first step is in turn to make sure that the behaviour of these simpler models is well understood, hence the usefulness of the theorems below.
\subsection*{Organisation of the book}
A {\it mode d'emploi} of this book, namely a presentation of its subsequent organisation, as well as an explanation of some expository choices, is provided below.
  \begin{itemize}
 \item[---]{\it Bibliographical references}. The book, while consisting for an important part of original research, reviews or uses results and concepts  due to other authors. In order not to interrupt the flow of exposition, I have chosen not to quote references in the body of the text.
However, each chapter comprises a section of bibliographical references, in which I endeavour to trace back the results that are not mine as far as I can, and where I also do my best to explain the relation of the 
theorems and methods displayed in the book to their general context. The reader is welcome to signal to me any important reference that I may have missed. I have made two exceptions to the above rule, with  the works of KPP and ZFK below. Besides the fundamental importance of these works, the acronyms are used everywhere in the text, so that the meaning has to be explained right away.  Existing, or potential  ramifications are also discussed. 
 More than once, I have chosen to quote PhD theses. While this choice is sometimes criticised, it is in fact totally legitimate. A PhD thesis is a scientific document that has been refereed at least as seriously as a journal paper, moreover those I quote are easily accessible online. In doing so, I also pay a tribute to hard working PhD students who have proved nontrivial theorems during their (sometimes difficult) years of initiation to research.
\item[---] {\it Open questions}. While this work goes reasonably far in the study of \eqref{e1.1.1}, it raises, in the end, more questions than answers. And so, in each chapter, a section on open problems comes after the bibliographical comments. I have restricted myself to issues directly pertaining to the situation studied in the book, that is, models in one space dimension (I allow an occasional excursion into models with directed one-dimensional propagation). Multi-dimensional issues, which are discussed at two places, and with no ambition of completeness, do not constitute  the central theme. This is a choice: the list of interesting questions would have, indeed,  probably doubled the size of the whole book.
These open problems are of two sorts. Some arise by relaxing an assumption in a seemingly harmless way, thus innocently bringing up a truly nontrivial issue. While I am not clueless about them, a definite answer is still quite far away. 
Others are a natural sequel to the questions studied in the chapter, they usually consist in bringing in  an additional element of complexity aimed at modelling new effects.  Once again, in these  cases, the material displayed in the chapter suggests an avenue of
 investigations which clearly requires other ideas than those displayed. Needless to say, I would be genuinely happy to see them solved, at least partly.
\item[---] {\it Problems.} Each chapter ends with a selection of problems, sometimes quite large. Some of them are really  exercises of application; sometimes they insist on, or make more precise, a particular point that was (willingly) elusive in the text; in which case they
should be regarded as a part of the chapter, that I did not develop myself in order to avoid repetitions. While they
 do not require anything else than a good understanding of the chapter, they often show that the scope of the considerations developped in this work go further than the basic model \eqref{e1.1.1}.
 Another large part of them, on the other hand, require actual thinking. They can be proposed as nontrivial master projects, and I hope that the students who choose to confront them will enjoy themselves.
 Some of them are published results by other authors (in such a case I give of course a reference) but the path that I propose is usually different from that of the original paper. They can also be nontrivial variants of these papers. Others are of my own creation.    While I have made no actual effort to distinguish between exercises of applications and the more elaborate problems, the former are, in general, located at the beginning of the selection.
 While I know (or, at least, I think I know) the answer to almost all,  I am not so sure that I can  answer a limited number of them. In such a case I  plainly confess it. 
\end{itemize}
\section{Bibliography, comments}\label{s1.6}
\noindent Let me dutifully put the  above principles to work.  An introductory book,  more complete than this introduction, is Roques \cite{Roques}. For the reader willing to grasp the breadth of the modelling in life science, the two-volume treatise of Murray \cite{Mu1}, \cite{Mu2} is highly recommanded. 
\subsubsection*{Invasions driven by reproduction and displacements (Section \ref{s1.2.1})}
\noindent The description of reproduction terms is quite rudimentary and I do not think that, given the scope of this book, it is necessary to go much further.  Let me only say that I have chosen to concentrate on positive nonlinearities $f(u)$, leaving aside an important class of nonlinear terms, namely the {\it bistable} terms. Such a term $f$ satisfies $f(0)=f(1)=0$ and $f'(0)<0$, $f'(1)=0$. From the modelling point of view it is perfectly legitimate, as the reproduction rate may be small, or even negative when the density os too small. This is called the {\it Allee effect,} and is discussed in \cite{Mu1}. A bistable nonlinearity has one or several zeroes inside the interval $(0,1)$, a feature that can generate complex dynamics. However, a for the significant subset of  those $f$s having only one zero inside $(0,1)$, the dynamics of a solution starting from a compactly supported initial data is essentially the same as that described in Chapter \ref{ZFK_short_range}, in the case when the bottom speed is strictly larger than the KPP speed.  

\noindent As for the modelling of displacements, only the essentials have been provided, with a brief discussion of the link between the general expression and the diffusion given by an elliptic operator. Let me stress, however, the following fact: while the resulting equation \eqref{e1.2.4} is quite simple looking, it describes actual natural phenomena in a surprisingly accurate way. It was for the first time put in action on a real situation by Skellam \cite{Sk}, in the analysis of the invasion of Central Europe by the muskrat in the first half of the previous century. This example is all the more valuable as it is a real life experience whose beginning is located precisely in space and time - that is, the moment when a few muskrats escaped from a farm where they were bred for their fur. In \cite{Sk}, the author has the idea of computing, from the data (whose accuracy is what it is), the evolution of the area occupied by the rodent, and finds out that its square root is, strikingly, linear in time. This behaviour is subsequently related to that of the Fisher-KPP equation.

\noindent For a much more detailed discussion, one 
should consult  the second volume of the treatise of Murray, that is, \cite{Mu2}.
\subsubsection*{Invasions driven by remote contaminations (Section \ref{s1.2.4})}
\noindent  The SIR model is a particular case of a much more elaborate class of models, elaborated by Kermack and McKendrick in \cite{KMK}, whose study is proposed as a problem. Contrary to that of the SIR model, which is elementary, the analysis of the homogeneous Kermack-McKendrick model is not entirely trivial and uses sub and super solutions. 
I have traced the SIR model with nonlocal contaminations back to Kendall \cite{Kend}, 
but there may be earlier occurrences. See, in any case, Murray  \cite{Mu1}, \cite{Mu2}. 

\noindent A well-known shortcoming of the SI  model is its inability to predict    documented dynamical phenomena, such as plateaus and rebounds. A recent  work of  Berestycki-Desjardins-Heintz-Oury \cite{BDHO}  proposes to structure the susceptible population by a parameter $a\in[0,1]$ representing the level of exposure to the risk. The number $S(t,a)$  of susceptibles with risk $a$, and the number of infected $I(t)$, evolves as
 \begin{equation}
 \label{e1.5.1}
 \left\{
 \begin{array}{rll}
  \partial_tS=\partial_{aa}S-\beta(a)SI,\ \ \ 
 \dot I=&\di\biggl[\di\int_0^1\beta(a)S(t,a)~da-\gamma(a)\biggl]I \quad \ \ (a\in(0,1))\\
 \partial_aS(t,0)=&\partial_aS(t,1)=0.
 \end{array}
 \right.
 \end{equation}
Numerical simulations  not only show the sought for plateaus and rebounds, but are in surprisingly good agreement with some covid-19 data in the South of France. On the other hand, 
the complete understanding of the dynamics of \eqref{e1.5.1} has not yet been achieved.

\noindent The first  diffusion mechanism, namely, Model \eqref{e1.2.5}, was  proposed  by Kendall, in a visionary answer to a paper in statistical epidemiology by Bartlett (J. Royal Stat. Soc. 1957).  Since then, a lot of more elaborate models have been proposed, especially SIR models on graphs or graph-like structures. There are so many of them that it is not possible to do justice to all of them. I just point out that, in all these more elaborate cases, a study of the spreading speed, 
not to speak about sharper asymptotics, still needs to be done.
\subsubsection*{Models of propagation directed by a line of fast, or nonlocal diffusion (Section \ref{s1.4})}
\noindent The first way that comes to the mind to ensure that the propagation of a front remains parallel to aline is to confine the front between walls parallel to this line. This is the situation studied around the beginning
 of the 1990's, in the framework of mathematical study of flame propagation theory. Model \eqref{e1.3.7}, with $\mathcal{J}$ replaced by $-\partial_{xx}$, represents the evolution of the temperature $u(t,x,y)$ of a flame propagating in the cylinder. The function $\alpha(y)$ represents an imposed flow field, and the term $f(u)$ is the reaction term. For a detailed study of the travelling waves with this particular diffusion, but with a different sort of nonlinearity $f$,I refer to Berestycki-Larrouturou-Lions \cite{BLL}. There, they introduce and make an extensive use of the sliding method, that will be reproduced in order to prove monotonicity and uniqueness of the wave profiles. The arguments developped there are generalised and very much sophisticated in Berestycki-Nirenberg \cite{BN-shear}. Model \eqref{e1.3.8}, which presents a more subtle way of confining the propagation to one direction, is proposed by Berestycki-Chapuisat \cite{BCh}, once again with $\mathcal{J}$ replaced by $-\partial_{xx}$. 
 
 \noindent Models of front propagation directed by a line of fast diffusion, of the type \eqref{e1.3.1}, were introduced and studied from the point of view of the spreading velocity by Berestycki, Rossi and the author of these lines in  \cite{BRR1},
  \cite{BRR2}, that they have popularised under the name "road-field model", for relatively obvious modelling reasons. Model \eqref{e1.3.6}, introduced by the same authors, presents a ramification of \eqref{e1.3.1} that models the spread on invasive plants in damaged environments.  In all these works, the diffusion is given by second order elliptic operators, as well as in the  beautiful study of travelling wave solutions to \eqref{e1.3.1} by L. Dietrich \cite{Diet}. The discrepancy of diffusion between line 
  and and interior of the domain generates new types of slow-fast dynamics, studied in detail by Dietrich and the author in \cite{DR}. The whole line of ideas is developped in \cite{BRR4} to model the spread of epidemics under the influence of transportation networks. Another model, pertaining to the class of SIR models on graphs, is proposed in Besse-Faye \cite{BF}.  
 
 \noindent  Berestycki-Hamel \cite{BH-book} give an extensive review of the general theory of reaction-diffusion equations.
\section{Problems}\label{s1.7}
\begin{problem}\label{P1.4}
Let  $K_0$ be $C^\infty$  an even, nonnegative, compactly supported function with mass 1. Give the limiting equation to $u_t+u-K_\e*u=0$ in the two following cases.
\begin{itemize}
\item[---] We have 
$
K_\e(x)=\di\frac1{\e}\biggl((1-\e)K_0(\frac{x}\e)+\e K_1(\frac{x}\e)\biggl),
$
where $K_1$  is a $C^\infty $ compactly supported function with zero mass. 
\item [---] We have
$\di
K_\e(x,z)=\frac1{\e}\biggl((1-\e) K_0(\frac{z}\e)+\e K_2(x,\frac{z}\e)\biggl),$
where $K_2: (x,z)\mapsto K_2(x,z)$ is such that, for all $x\in \RR$, the function $K_2(x,.)$ satisfies the above assumptions on $K_1$.
\end{itemize}
\end{problem}
\begin{problem}
\label{P1.11}
Let us revert to an even kernel $K_\e(x)$ that is an honest approximation of the identity, but assume this time that $\rho(x)\sim k/\vert x\vert^{1+2\alpha}$ as $x\to+\infty$, with $k>0$ and $0<\alpha<1$. Assume the function $\rho$ to have unit mass.
\begin{itemize}
\item[---] Show that the limiting equation to $u_t+u-K_\e*u=0$ is, after a proper rescaling in $t$: 
$$u_\tau+d(-\partial_{xx})^\alpha u=0,
$$ where $d>0$ is to be computed in terms of $k$ and $\alpha$, and $(-\partial_{xx})^\alpha$ is the fractional Laplacian of order $\alpha$. Its expression is given by $(-\partial_{xx})^\alpha u=c_\alpha\mathrm{P.V.}\di\int_\RR\di\frac{u(x)-u(y)}{\vert x-y\vert^{1+2\alpha}}~dy$, the constant $c_\alpha$ being chosen so that its Fourier symbol is $\vert \xi\vert^{1+2\alpha}$.
For more information, a good recent reference is Garofalo \cite{Garof}.
\item[---] Study the limiting cases $\alpha\to0$ and $\alpha\to1$.
\end{itemize}
\end{problem}
\begin{problem}
\label{P1.12}
For $L>0$, let $\mu(y)$ and $\nu(y)$ be smooth positive functions on $[-L,0]$. In the cylinder $\Omega_L=\{(x,y)\in\RR\times(-L,0)\}$, consider the system with unknowns $\bigl(u(t,x),v(t,x,y)\bigl)$:
\begin{equation}
\label{e1.6.8}
\left\{
\begin{array}{rll}
&\begin{array}{rll}
u_t+\mathcal{J}u+u\di\int_{-L}^0\mu(y)~dy=&\di\int_{-L}^0v(t,x,y)\nu(y)~dy\\
v_t-d\Delta v=&f(v)+\mu u-\nu v
\end{array}
\quad \bigl(t>0,(y,y)\in\Omega_L)\\
&v(t,-L,y)=0
\end{array}
\right.
\end{equation}
\begin{itemize}
\item[---] Interpret Model \eqref{e1.6.8} in the light of Model \eqref{e1.3.6}.
\item[---] If the function $\nu(y)$ has the form $\nu_\e(y)=\di\frac1\e\nu_0(\frac{x}\e)$, show that the limiting equation as $\e\to0$ has the form \eqref{e1.3.6}.
\end{itemize}
Model \eqref{e1.6.8}, as well as more general forms, is introduced and studied, from the point of view of the spreading velocity, by Pauthier \cite{Pau} (his original model is posed in the whole half-space). It also falls in the scope of the theory developped for Model \eqref{e1.1.1}.
\end{problem}
\begin{problem}
\label{P1.7}
If $f$ is of the KPP or ZFK type, one easily shows (do it) that the solution $u(t)$ of $\dot u=f(u)$, with $u(0)=u_0\in(0,1)$ converges to 1 as $t\to+\infty$. Let us introduce the following innocent modification: 
consider a function $f(t,u)$ which is smooth in both variables, 1-periodic in $t$, such that $f(t,0)=0$, $f_u(t,0)>0$ and such that $f(t,u)\leq-\alpha u$ for a given $\alpha>0$ and for $u$ larger than some large $U_0>0$. 
Consider the ODE
\begin{equation}
\label{e1.6.1}
\dot u=f(t,u).
\end{equation}
\begin{itemize}
\item[---] Let $u_0>0$ small be a Cauchy datum. Show that \eqref{e1.6.1}, with initial datum $u_0$, has a unique global in time solution $u(t)$ converges, as $t\to+\infty$, to a 1-periodic solution of \eqref{e1.6.1}.
(the sequence $\bigl(u(n)\bigl)_n$ is monotone).
\item[---] Show that two periodic solutions are strictly ordered.
\item[---] By considering the Cauchy Problem for \eqref{e1.6.1} with $u(0)$ small, show that \eqref{e1.6.1} has at least one positive periodic solution, let $u^+(t)$ be the minimal nontrivial periodic solution.
\item[---] Show that $\di\int_0^1f_u\bigl(t,u^+(t)\bigl)~dt\geq0$, and that it is $>0$ for some nonlinearities. 
\end{itemize}
\end{problem}
\begin{problem}
\label{P1.16}
Assume $f$ is a KPP nonlinearity that is $f_u(t,u)\leq f_u(t,0)$ for all $t\in\RR$ and all $u>0$. Anticipating a little on Chapter \ref{Cauchy}, Theorem \ref{t2.2.10}, show that $u^+$ is the only nontrivial periodic solution. To do this, one may compare $\sigma u$ and $v$, where $u$ and $v$ are periodic solutions to \eqref{e1.6.1}, and $\sigma\in[0,1]$.
\end{problem}
\begin{problem}\label{P1.15} The issue here is to study instance where \eqref{e1.6.1} has one, or several nontrivial periodic solutions.
\item[---] Investigate the uniqueness of $u^+$ if $f$ is $1/\e$-periodic in $t$, $\e>0$ small.
\item[---] What happens if $f$ is $\e$-periodic, $\e>0$ small?
\item[---] If $f(u)=u-u^2$ if $u\in[0,1]$ and 1-periodic in $u$, then of course every integer is a positive periodic solution to \eqref{e1.6.1}. Perturb $f$ in order to make it time dependent, while \eqref{e1.6.1} has an infinite number of
nonconstant periodic solutions.
\end{problem}
\begin{problem}
\label{P1.8}
In equation \eqref{e1.6.1}, keep all the already made assumptions on $f$, except the periodicity in $t$. Generalise the results of Problem \ref{P1.7} as much as you can. Find out what happens, for instance, when 
$f$ is almost constant in $t$  on time intervals of increasing size.
\end{problem}
\begin{problem} \label{P1.1} Consider the SI model \eqref{e1.2.2}.
\begin{enumerate}
\item [---] When $R_0\leq1$, expand $u_*(I_0)$  when $I_0\to0$.
\item [---] Show that there is $\alpha_*>0$ and $v_*>0$  such that $u(t)-u_*(I_0)\sim_{t\to+\infty} v_*e^{-\alpha_*t}.$
Study how $\alpha_*$ depends on $I_0$, especially when $R_0=1$.
\item  [---] Assume  $R_0>1$. Show that  $I(t)$ has a unique maximum point $I_*$; study how it depends on $I_0$, as well as the time $T_*$ at which it is reached.
\end{enumerate}
\end{problem}

\begin{problem}\label{P1.2} (Inspired from Murray \cite{Mu1}, where one may find many other similar models)
The SI model has an infinite number of sophistications, as there are indeed many sorts of infected individuals: the symptomatic ones, the pre-symptomatic ones, those who are symptomatic but not contagious, the contagious ones... 

\noindent This problem proposes a slightly more precise modelling, in the sense that it takes into account the density $E(t)$ of incubating individuals. They can be contagious, or not. The former case is of course the more worrisome in terms of public health. System \eqref{e1.2.2}  becomes
 \begin{equation}
 \label{e1.7.3}
 \left\{
 \begin{array}{rll}
 \dot S=&-\beta S(E+aI)\\
 \dot E=&\beta S(E+aI)-\alpha E\\
 \dot I=&\alpha E-\gamma I\\
 S(0)=&S_0>0,\ E(0)=E_0>0\ \hbox{small},\ I(0)=0.
 \end{array}
 \right.
 \end{equation}
The parameters $\alpha$, $\beta$, $\gamma$ and $a$ are $>0$.
\begin{enumerate}
\item [---]  Show that \eqref{e1.7.3} has a unique global solution $\bigl(S(t), E(t),I(t)\bigl)$.
 We introduce the cumulated densities
 $
 u(t)=\di\int_0^tE(s)~ds\ \hbox{and}\ v(t)=\di\int_0^tI(s)~ds.
 $
Show that  $\bigl(u(t),v(t)\bigl)$ solves the system
\begin{equation}
 \label{e1.7.4}
 \left\{
 \begin{array}{rll}
\dot u=&S_0(1-e^{-\beta(u+av)})-\alpha u+E_0.\\
\dot v=&\alpha u-\gamma v\\
u(0)=&v(0)=0.
\end{array}
\right.
\end{equation}
Show that $u(t)$ and $v(t)$ tend to limits $u_\infty(E_0)$ and $v_\infty(E_0)$ as $t\to+\infty$. {\rm Hint:} \eqref{e1.7.4} is a monotone system. An important  paper exploring this beautiful theory is Hirsch \cite{Hir})
\item [---] Propose a new quantity  $R_0$ such that $v_\infty(E_0)$ is larger than a fixed $v_*>0$ if $R_0$ is  $>1$.
\item [---] Suppose $R_0>1$ and $a>0$ small, that is, the most unfavourable case in terms of public health. When the number of infected reaches $\e S_0$, with $\e>0$ small, assume that there is a reaction of the authorities at time $T_*$. Assume that  it is so efficient that the transmission rate $\beta$ in \eqref{e1.7.3} falls to 0 for $t>T_*$.

\noindent Study the limits of $v(t)$ and  et $S(t)$ when $t\to+\infty$, as well as the maximal number of infected. You can also try to comment these results.
\end{enumerate}
\end{problem}
\begin{problem}\label{P1.14}
In the basic SI model, there is no real reason to assume that $\alpha$ and $\beta$ are constant in time. The transmission coefficient $\beta$ may vary seasonally: it is indeed well-known that the conditions for the propagation of a flu virus are much better (for the virus) during the winter than during the summer. Similarly, the population may acquire some resilience to the virus, so that $\alpha(t)$ decreases in time. So, it is a natural entreprise to investigate the time dependent system
\begin{equation}
\label{e1.7.12}
\dot S=-\beta(t)SI,\ \ \dot I=\beta(t)SI-\alpha(t)I.
\end{equation}
The function $\alpha$ and $\beta$ are smooth and bounded and positive.
\begin{itemize}
\item[---] Start from an initial condition $(S_0,I_0)$, both quantities being positive and $I_0$ small. Show the existence of $S_\infty\geq0$ such that $\di\lim_{t\to+\infty}\bigl(S(t),I(t)\bigl)=(S_\infty,0)$. Provide
an upper bound for $S_\infty$.
\item[---] Give examples of functions $\alpha(t)$ and $\beta(t)$ for which we (unfortunately) have $S_\infty=0$.
\item[---] Propose sufficient conditions (I am not   sure about a necessary one) of $S_0$ and the functions $\alpha$ and $\beta$ ensuring that $S_\infty$ is less than a  fraction of $S_0$, independent of the size of $I_0$.
\end{itemize}
\end{problem}
\begin{problem}\label{P1.3}
Variants attack! One way, among many others, to model this, is to declare that the surviving individuals become susceptible again.  
 \begin{enumerate}
\item [---] One proposes the following system:
 \begin{equation}
 \label{e1.7.5}
 \left\{
 \begin{array}{rll}
 \dot S=&-\beta SI+(\gamma-\e)I,\ 
 \dot I=\beta SI-\gamma I\\
 S(0)=&S_0>0,\ I(0)=I_0>0\ \hbox{small}.
 \end{array}
 \right.
 \end{equation}
  What does $\e I$ represent? Study (\ref{e1.7.5}) for all $\e\in[0,\gamma)$. What happens for small $\e$?
\item  [---] The government launches a vaccination campaign to try to get $\beta$ to drop as much as possible. Vaccine efficiency decreases, until a more adapted vaccine is found.
So, the constant $\beta$ is replaced by a function $\beta\phi(t)$. Assume $\phi$ periodic with period  $T$, discontinuous at all the times $nT$, $n\in\NN^*$, and affine between  $nT$ and $(n+1)T$. Assume
$
\phi(nT^-)=1,$ et $\phi(nT^+)=0$.  

\noindent Study the behaviour of $I(t)$ when  $t\to+\infty$, especially when $T$ is large.
\end{enumerate}
\end{problem}
\begin{problem}\label{P2.101} The full homogeneous Kermack-McKendrick model \cite{KMK} consists in considering that the transmission rate of an infected individual depends on the time at which he/she has been infected. So, this amounts to finding two functions
$S(t)$ (the density of susceptibles) and $I(t,\tau)$ (the density of individuals that have been infected for a duration exactly equal to $\tau$) satisfying
\begin{equation}
 \label{e1.7.6}
 \left\{
 \begin{array}{rll}
 \dot S=&-S\di\int_0^{+\infty}\beta(\tau)I(t,\tau)~d\tau\\
 \partial_tI+\partial_\tau I=&-\gamma(\tau)I\ \ (t>0,\tau>0)\\
 I(t,0)=&\di\int_0^{+\infty}\beta(\tau)I(t,\tau)~d\tau\\
 S(0)=&S_0>0,\ I(0,\tau)=I_0(\tau)\geq0\ \hbox{compactly supported}.
\end{array}
  \right.
  \end{equation}
Set  $f(u)=S_0(1-e^{- u}),$  the cumulated densities of infected 
$
u(t,\tau)=\di\int_0^tI(s,\tau)~ds
$
solve the system
 \begin{equation}
 \label{e2.7.1507}
 \left\{
 \begin{array}{rll}
  \partial_tu+\partial_\tau u=&-\gamma(\tau)u+I_0(\tau)\ \ (t>0,\tau>0)\\
u(t,0)=&f\biggl(\di\int_0^{+\infty}\beta(\tau)u(t,\tau)~d\tau\biggl),
\end{array}
\right.
\end{equation}
\begin{enumerate}
\item [---] Interpret this system.  How can one retrieve the SI system from it?
\item[---] Study the stationary solutions to \eqref{e2.7.1507}. Propose a new $R_0>0$ such that \eqref{e2.7.1507} has a nontrivial steady solution $u_*(\tau)>0$ if and only if $R_0>1$.
 Compute $u_\infty(\tau)$.
\item[---] Show that \eqref{e2.7.1507} has a global solution $u(t,\tau)$. 
\item[---] Assume $R_0>1$; show that $u(t,\tau)$ converges uniformly in $\tau$ to $u_*(\tau)$ as $t\to+\infty$. {\rm Hint:}  put a subsolution   below $u(1,\tau)$ (see Chapter \ref{Cauchy}).
Deduce  the large time behaviour of  $\bigl(S(t),I(t,\tau)\bigl)$. 
\end{enumerate}
\end{problem}
\begin{problem}
\label{P1.10}
The duration of infection is not, in practice, a continuous variable: you count for how long you have been infected by the flu (measles, covid, else...) in days, not in milliseconds.

\noindent Propose a discretisation of the Kermack-McKendrick model \eqref{e1.7.6} in terms of the duration of the infection; the resulting model should be an infinite system of differential equations, unless it is declared that the infection 
does not last more than a given number of time units. Study its large time behaviour.

\noindent Even if the transport equation in the $\tau$ variables is not the most sophisticated hyperbolic equation on earth, some care should be given to how the discretisation is conducted. Any good treatise in numerical analysis
studies the issue; as a tribute to the teachers that I had when I was a Master student, I quote Godlewski-Raviart \cite{GR}.
\begin{problem}
\label{P2.130}
Consider $\tau_0>0$. For  $\e>0$, consider a function $\gamma_\e(\tau)$ equal to a fixed function $\gamma_0(\tau)$ on $[0,\tau_0]$, and such that $\gamma_\e(\tau)\geq 1/\e$ for $\tau\geq\tau_0+\e$. Let $\bigl(S_\e(t),I_\e(t,\tau)\bigl)$  the solution of (\ref{e1.7.6}) with $\gamma=\gamma_\e$. Study and interpret its limit when   $\e\to0$.
\end{problem}
\begin{problem}\label{P2.109}
Study the system of Kermack and McKendrick \eqref{e2.7.1507} in the (sad) case $\gamma(\tau)\equiv0$. Do you know a concrete example? This case, by the way, falls into the class {\rm renewal equations,} see Perthame \cite{Pbook}, Chap. 3, for a survey of the subject.
\end{problem}

\end{problem}

\chapter{Cauchy Problem, related issues}\label{Cauchy}
\noindent The goal of this chapter is to derive the general properties of the solutions to the nonlocal problem that will occupy us throughout this work:
\begin{equation}
\label{e2.1.1}
\partial_tu+\int_{\RR}K(x-y)\bigl(u(t,x)-u(t,y)\bigl)~dy=f(u),\ \ t>0,\ x\in\RR.
\end{equation}
The kernel $K$ will always be smooth, even and nonnegative. It will be positive on   $(-1,1)$, and zero outside. One may, but this is not a requirement, ask that its total mass is 1. In such a case, equation \eqref{e2.1.1} becomes
\begin{equation}
\label{e2.1.100}
\partial_tu+u-K*u=f(u),\ \ t>0,\ x\in\RR.
\end{equation}
To alleviate the notations we will often denote the integral diffusion by
\begin{equation}
\label{e2.1.14}
\JJ u(x)=\int_{\RR}K(x-y)\bigl(u(t,x)-u(t,y)\bigl).
\end{equation}

\noindent The initial datum $u(0,x)=u_0(x)$ will be nonnegative, compactly supported. It will be useful to assume $u_0\leq1$. The function $f$ will be smooth, positive on $(0,1)$, with $f(0)=f(1)=0$, and negative everywhere
else. It will be useful, but by no means crucial, to assume that it is globally Lipschitz, and that all its subsequent derivatives are bounded. No other assumption will be made in this chapter.

\noindent In the whole chapter, the analogy with the reaction-diffusion equation with second order diffusion 
\begin{equation}
\label{e2.0.1000}
\partial_tu-\partial_{xx}u=f(u),\quad (t>0,\ x\in\RR)
\end{equation}
will be in the background. No systematic attempt will be made to relate our results to those pertaining to \eqref{e2.0.1000}, the reader having an experience with such equations can do it. The main feature that is lacking in the model under study, namely, \eqref{e2.1.1}, is the instant regularisation. This will be a low noise inconvenience in this chapter, but will turn out to be a source of important difficulties as the theory develops.

\subsubsection*{Organisation of the chapter}
\noindent Section \ref{s2.2} presents the most basic results and the comparison arguments that will be used freely throughout the book, especially the method of sub and super solutions. Section \ref{s2.25} is concerned with steady states, it is a first not totally trivial instance of what happens when no regularising effect comes to the help. The short section \ref{s2.3} presents a main character of the theory, that is, the linearised equation at $u=0$ and its waves, leading to the notion of critical speed. The main, and most original part of the chapter, section \ref{s2.4}, presents the asymptotics of the linearised equation with the critical speed as a drift. It will be very much used when it comes to the study of sharp asymptotics.
\section{The initial value problem}\label{s2.2}
\noindent As far as Equation \eqref{e2.1.1} is concerned, global existence and uniqueness to the Cauchy Problem for \eqref{e2.1.1}  resorts to nothing else than the Cauchy-Lipschitz theorem in the space, say, of all $C^k$ functions 
of $\RR$. The integer $k$ is chosen so that any differentiation that we will need is automatically justified.  In particular, if we assume $u_0$ to be a $C^k$ function, a $C^k$ solution will exist for all time, all derivatives being bounded by
exponentials, due to the Gronwall Lemma. 
\begin{remark}
\label{r21.1}
The smoothness of $u_0$ is not a crucial requirement, one may simply assume $u_0$ to be bounded. One then obtains a solution in $C^k(\RR_,L^\infty(\RR))$, for all $k\in\NN$.
\end{remark}
\subsection{The comparison principle}
\noindent The first specific insights that one may have in these systems is gained through comparison. The basic result is that two initially ordered initial data will give rise to two solutions that are ordered in the same fashion. 
\begin{proposition}
\label{p2.1.1}
Let $u_{10}\leq u_{20}$ be two initial data. For $i\in\{1,2\}$, let $u_i(t,x)$ the solution of \eqref{e2.1.1} emanating from $u_{i0}$. Then $u_1(t,x)\leq u_2(t,x)$.
\end{proposition}
\noindent{\sc Proof.} Set $v(t,x)=u_1(t,x)-u_2(t,x)$ and $a(t,x)=-\di\frac{f\bigl(u_1(t,x)\bigl)-f\bigl(u_2(t,x)\bigl)}{v(t,x)}$; the function $a$ is bounded on $[0,T]\times\RR$, for all $T>0$. 
An equation for $v$ is
\begin{equation}
\label{e2.1.1501}
v_t+\bigl(\K+a(t,x)\bigl)v=K*v,
\end{equation}
which, by Duhamel's principle,  integrates into
$$
v(t,x)=\mathrm{exp}\biggl(-\int_0^t\bigl(\K+a(s,x)\bigl)~ds\biggl)v_0(x)+\int_0^t\mathrm{exp}\biggl(-\int_s^t\bigl(\K+a(s',x)\bigl)~ds'\biggl)K*v(s,x)~ds.
$$
As $v_0\leq0$ and $K\geq0$, taking the positive parts yields
$$
v^+(t,x)\leq\int_0^t\mathrm{exp}\biggl(-\int_s^t\bigl(\K+a(s',x)\bigl)~ds'\biggl)K*v^+(s,x)~ds,
$$
so that 
$$
\Vert v^+(t,.)\Vert_\infty\leq \K\int_0^t\mathrm{exp}\biggl(-\int_s^t\bigl(\K+a(s',x)\bigl)~ds'\biggl)\Vert v^+(s,.)\Vert_\infty~ds,
$$
which eventually leads to $\Vert v^+(t,.)\Vert_\infty=0$. This yields $u_1(t,x)\leq u_2(t,x)$. \hfill$\Box$

\noindent As we look deeper into Model \eqref{e2.1.1}, we will need, as is common in diffusion equations, to compare solutions in spatial regions that are not the whole line. This is why we state the following maximum principle in domains with boundaries; as their proofs do not stem directly  from that of Proposition \ref{p2.2.1} we devote some care in stating and proving them. As everything happens at the level of linear equations, we study solutions of 
\begin{equation}
\label{e2.1.1400}
v_t+\mathcal{J}v+a(t,x)v\leq0,
\end{equation}
in various regions of $\RR_+\times\RR$, and where the notation $\mathcal{J}$ is given by \eqref{e2.1.14}.
\begin{proposition}
\label{p2.1.200}
Let $v(t,x)$ be a Lipschitz solution of \eqref{e2.1.1400} in the following two situations
\begin{enumerate}
\item We have  $X_1(t)<x< X_2(t)$, where $X_i$, $i\in\{1,2\}$, are Lipschitz functions of $t$.
Moreover,  we have $v(t,x)\leq0$ if $X_1(t)-1\leq x\leq X_1(t)$, $v(t,x)\leq0$ if $X_2(t)\leq x\leq X_2(t)+1$,  and $v(0,x)\leq0$ for $X_1(0)\leq x\leq X_2(0)$.
\item We have $x\geq X_1(t)$, where $X_1$ is a Lipschitz function of $t$, and there exists a family of bounded functions $(V_T)_{T>0}$ such that, for all $T>0$ we have
\begin{equation}
\label{e2.1.1401}
\sup_{0\leq t\leq T}\vert v(t,x)\vert\leq V_T(x).
\end{equation}
Moreover, $v(t,x)\leq0$ if $\leq X_1(t)-1\leq x\leq X_1(t)$,and $v(0,x)\leq0$ for $x\geq X_1(0)$.
\end{enumerate}
Then $v(t,x)\leq0$ for $X_1(t)\leq x\leq X_2(t)$ in Situation 1, and $v(t,x)\leq0$ for  $x\geq X_1(t)$ in Situation 2.
\end{proposition}
\noindent{\sc Proof.} Let us first examine Situation 1. For a continuous function $u(x)$, let $\mathrm{sgn}^+u(x)$ the function that is equal to 1 if $u(x)>0$, and $0$ if $u(x)\leq0$.
Write \eqref{e2.1.1400} as
$$
v_t+\bigl(1+a(t,x)\bigl)v\leq K*v,
$$
and multiply the equation by $\mathrm{sgn}^+v$, this yields
$$
v^+_t+\bigl(1+a(t,x)\bigl)v^+\leq (K*v)\mathrm{sgn}^+\leq K*v^+.
$$
For each $t>0$, this equation is integrated on $[X_1(t),X_2(t)]$, as $v$ is Lipschitz and $v^+\bigl(t,X_i(t)\bigl)=0$ we have (see Rudin \cite{Rud}, Chap. 7)
$$
\frac{d}{dt}\int_{X_1(t)}^{X_2(t)}v^+(t,x)~dx\leq-\int_{X_1(t)}^{X_2(t)}\biggl(\bigl(1+a(t,x)v^+(t,x)\bigl)+K*v^+(t,x)\biggl)~dx,
$$
taking into account that $v^+(t,x)=0$ if $x\notin[X_1(t),X_2(t)]$ there is a locally bounded function $M(t)$ such that 
$$
\frac{d}{dt}\int_{X_1(t)}^{X_2(t)}v^+(t,x)~dx\leq M(t)\int_{X_1(t)}^{X_2(t)} v^+(t,x)~dx,
$$
which entails $v^+\equiv0$. Situation 2 is treated in the same way, except that we have to use \eqref{e2.1.1401} in order to be able (see \cite{Rud} once again) to make the integral and the time derivation commute. \hfill$\Box$

\noindent Proposition \ref{p2.1.1} offers an important analogy with the weak maximum principle for parabolic equations. The analogy does not stop here, as we also have strong comparison principles. While they are much less involved than the corresponding results for parabolic equations, they will turn out to be extremely useful.
\begin{remark}
\label{r2.2.400}
Proposition \ref{p2.1.1} only needs the initial functions $u_{i}$ to be bounded.
\end{remark}
\noindent Just as in elliptic equations, some forms of comparison apply to solutions of steady equations. Let us state a generic result that will be needed in the refined study of travelling waves at infinity.  We make an easy start with the following statement.
\begin{proposition}
\label{p2.2.3}
Let $u(x)$ be a nonnegative continuous  function satisfying, at each of its zeroes $x_0$: $K*u(x_0)=0.$
Then $u\equiv0$ or $u>0$ on $\RR$.
\end{proposition}
\noindent{\sc Proof.}  By the commutativity of the convolution product ,we have, at every point $x$ such that $K*u(x)=0$:
$\di\int_{\RR}K(y)u(x_0-y)~dy=0.
$
As both $u$ and $K$ are nonnegative, this proves that $u\equiv0$ in an interval centred at $x$, thus proving that the zero set of $u$ vanishes is open in $\RR$. As it is closed, the zero set of $u$  is the whole line as soon as it is non void. \hfill$\Box$

\noindent This implies the useful, and easily proved proposition is the following.
\begin{proposition}
\label{p2.2.1}
Let $u(x)$ be a nonnegative $C^1$ function satisfying the integro-differential inequality
\begin{equation}
\label{e2.2.1}
\mathcal{J}u-cu'-a(x)u\geq0,\ \ \ x\in\RR,
\end{equation}
where $a(x)$ is a bounded continuous function.
 Then, either $u(x)>0$ everywhere, or $u(x)\equiv0$.
\end{proposition}
\noindent{\sc Proof.} Let us consider $x_0$ a zero of $u$, we have $K*u(x_0)\leq0$, so, because of the nonnegativity of both actors, $K*u(x_0)=0$. 
So, $u(x)=0$. \hfill$\Box$

\noindent Similarly to the strong maximum principle for parabolic equations, we have a strong maximum principle for the Cauchy Problem.
\begin{theorem}
\label{t2.1.1500}
Let $u_{10}\leq u_{20}$ be two smooth initial data. For $i\in\{1,2\}$, let $u_i(t,x)$ the solution of \eqref{e2.1.1} emanating from $u_{i0}$. If $u_{10}$ and $u_{20}$ do not coincide everywhere, then $u_1(t,x)<u_2(t,x)$.
\end{theorem}
\noindent{\sc Proof.} Set again $v(t,x)=u_1(t,x)-u_2(t,x)$. Assume the existence of $t_0>0$ such that the zero set of $v(t_0,.)$ is non void; let $x_0$ be such a point. We have $\partial_tv(t_0,x_0)\geq0$, so that we have $K*v(t_0,.)\geq0$ at $x=x_0$. As $v(t,x)\leq0$ we have $K*v(t_0,.)=0$ at $x=x_0$. This implies   $v(t_0,.)\equiv0$. Solving equation \eqref{e2.1.1501} for $v$ backwards in $t$ yields $v(0,.)\equiv0$, a contradiction. \hfill$\Box$
\subsection{Sub and super solutions}
\noindent At this stage, it may be useful to state a few facts about sub and super-solutions, that are very much inspired from the  notions pertaining 
second order elliptic or parabolic equations. While they are elementary, they will turn out to be helpful in the proofs of estimates that go beyond the standard ones. For $c\in\RR$, a super-solution (resp. sub-solution) to the equation 
\begin{equation}
\label{e2.2.15}
\partial_tu+\JJ u-c\partial_xu=f(u)
\end{equation}
is a locally Lipschitz function $\bar u(t,x)$ (resp. $\underline u(t,x)$)  that satisfies \eqref{e2.2.15} with $\geq0$ (resp. $\leq0$) instead of $=0$. A similar definition applies for a super-solution or sub-solution of the steady version of 
\eqref{e2.2.15}, that is, with $\partial_t=0$. The definition easily extends to sub- or super-solutions of \eqref{e2.2.15} on parts of the real line, with the caveat that one should know the function a little outside. For instance, a sub-(resp. super-) solution of \eqref{e2.2.5} on $[0,T]\times (a,b)$ should be defined on $[0,T]\times(a-1,b+1)$. If the support of the integral kernel was not compact (e.g. $\RR$), this would impose the function to be defined on $[0,T]\times\RR$.

\noindent The following proposition is a typical example of a useful result, whose proof is obvious. The reader may be interested in writing down its analogue for elliptic or parabolic equations.
\begin{proposition}
\label{p2.2.2}
Let $a$ be a real number, and let us consider a function $\underline u^+(t,x)$ defined on $[0,T]\times[a-1,+\infty)$ that is a sub-solution to \eqref{e2.2.15} on $[0,T]\times(a,+\infty)$.
Let $\underline u^-(t,x)$ be a function defined on $[0,T]\times[a-1,a]$ such that 
$$\underline u^-(t,x)\geq \underline u^+(t,x)\ \ \ \hbox{for $t\in[0,T]$ and $a-1\leq x\leq a$.}
$$
Then, $\underline u(t,x)$ defined by
$$
\underline u(t,x)=
\left\{
\begin{array}{rll}
&\underline u^-(t,x)\ \hbox{if $a-1\leq x\leq a$}\\
&\underline u^+(t,x)\ \hbox{if $x\geq a$}
\end{array}
\right.
$$
is a sub-solution to \eqref{e2.2.5} on $[0,T]\times[a,+\infty)$. A similar statement holds for a super-solution.
\end{proposition}

\noindent It is also true that the successive derivatives of $u(t,x)$ are uniformly bounded in time. This, however, is not a triviality anymore. It will be obtained as a corollary of much more precise estimates on $u(t,x)$ for large times, that will also 
stem from comparison.

\noindent We end this  section with the following proposition, that will be quite helpful when it comes to studying the large time behaviour of problems of the form \eqref{e2.2.15}. While the statement is standard, the lack of uniform bounds for the
derivatives forces to write the argument with  care.
\begin{proposition}
\label{p2.2.4}
Assume $c\neq0$.
\begin{enumerate}
\item Let $\underline u(x)$ be a nonnegative Lipschitz continuous sub-solution to \eqref{e2.2.15}, that is in addition $\leq1$. Let $v(t,x)$ be the solution of \eqref{e2.2.15} starting from $\underline u$. Then $\partial_tv(t,x)\geq0$ and $v(t,.)$ converges,
locally uniformly, to a smooth solution $u_\infty(x)$ of \eqref{e2.2.15} that is above $\underline u$.

\item Let $\bar  u(x)$ be a nonnegative Lipschitz continuous super-solution to \eqref{e2.2.15}. Let $v(t,x)$ be the solution of \eqref{e2.2.15} starting from $\bar  u$. Then $\partial_tv(t,x)\leq0$ and $v(t,.)$ converges,
locally uniformly, to a smooth solution $u_\infty(x)$ of \eqref{e2.2.15} that is below $\bar  u$.
\end{enumerate}
\end{proposition}
\noindent{\sc Proof.} We will prove the first statement, the second being proved in the same way. Set $w(t,x)=\partial_tv(t,x)$; as $\underline u(x)$ is Lipschitz continuous it is a
 well defined bounded function that satisfies  
$$
\partial_tw-c\partial_xw+\mathcal{J}w=f'(u)w,\ \ w(0,.)\geq0.
$$
Thus, $w(t,x)\geq0$ and the function $t\mapsto v(t,x)$ is a bounded nondecreasing function and, as such, converges in the pointwise sense to a bounded function $u_\infty(x)$. Obviously we have
$
\di\lim_{t\to+\infty}\mathcal{J}v(t,.)=\mathcal{J}u_\infty,
$
in the pointwise sense. Moreover, we have
$
\di\ \lim_{t\to+\infty}\partial_xv(t,.)=\partial_xu_\infty
$
in the sense of distributions. Defining $\tilde v_t(s,x)=v(t+s,x)$ we also have
$
\di\lim_{t\to+\infty}\partial_s\tilde v_t(.,.)=\partial_tu_\infty\equiv0,
$
 in the sense of distributions. Finally, thanks to Lebesgue's dominated convergence theorem, we have, still in the sense of distributions:
$
\di\lim_{t\to+\infty}f\bigl(v(t,.)\bigl)=f(u_\infty).
$
So we have  obtained a  solution $u_\infty(x)$ of
$
c\partial_xu_\infty=\mathcal{J}u_\infty-f(u_\infty).
$ in the distributional sense.
As $c$ is nonzero, $u_\infty$ is Lipschitz continuous, hence smooth by a bootstrap argument. As every $v(t,.)$ is Lipschitz continuous, and the convergence of $v(t,.)$ to $u_\infty$ is monotone,
Dini's Theorem ensures that it is locally uniform. \hfill$\Box$

\begin{remark}\label{r2.2.4}
If $c=0$, the property $\partial_tu\geq0$ 
(resp. $\partial_tu\leq0$) still holds.
\end{remark}
\section{Steady  states}\label{s2.25}
\noindent A steady state to Equation \eqref{e2.1.1} is a nonnegative bounded function $u(x)$ such that 
\begin{equation}
\label{e2.2.20}
\mathcal{J}u=f(u).
\end{equation}
We note that $u(x)$ is not necessarily continuous.  As steady states are expected to be the final state of the solution $u(t,x)$ of \eqref{e2.1.1} after the passage of the front, it is useful to classify them. Obvious solutions are $u\equiv0$ and $u\equiv1$. And in fact, they are the only ones.
\begin{theorem}
\label{t2.2.10}
Assume the function $u\mapsto {f(u)}/u$ to be nonincreasing. Then any bounded, nonnegative solution of \eqref{e2.2.20} is either $u\equiv0$ or $u\equiv1$.
\end{theorem}
Let us make the following benign computation, that will reappear later as the theory develops. Set
\begin{equation}
\label{e2.2.21}
e_\lambda(x)=\cos\lambda x,\ \ \lambda\in\RR.
\end{equation}
The evenness of $K$ yields
\begin{equation}
\label{e2.2.2221}
\mathcal{J}e_\lambda=\omega_0(\lambda)e_\lambda=\bigl(<\!y^2K\!>\lambda^2/2+O_{\lambda\to0}(\lambda^3)\bigl)e_\lambda, \quad\omega_0(\lambda)=\di\int_\RR K(y)(1-\cos\lambda y)~dy.
\end{equation}

\noindent{\sc Proof of Theorem \ref{t2.2.10}.}  When the diffusion operator is given by the Laplacian, every solution to \eqref{e2.2.20} is smooth. However, as the integral diffusion has no special regularising effect, this fact has to be proved with bare hands. What makes things work here is that, for any $\lambda>0$, the equation
\begin{equation}
\label{e2.2.22}
u-f(u)=\lambda
\end{equation}
has a unique solution. If   $u>0$ is a root of \eqref{e2.2.22}, we have 
$
{f(u)}/u=1-\lambda/{u}.
$
As $u\mapsto f(u)/u$ is nonincreasing, we have $f'(u)\leq f(u)/u$, so that $f'(u)<1$. The implicit functions theorem applied to $u\mapsto u-f(u)$ implies the uniqueness of the root, that we call $u_\lambda$. Moreover the function $\lambda\mapsto u_\lambda$ is smooth. 
Consider now a nontrivial steady solution $u(x)$, and $x_0\in\RR$. If $K*u(x_0)=0$, then $u\equiv0$ in an interval around $x_0$, so that $u$ is obviously smooth in a neighbourhood of  $x_0$. To be quite precise, what we really have  is $u=0$ almost everywhere in a neighbourhood of $x_0$, so that we may choose its continuous representative, that is, the zero function.
If  $K*u(x_0)>0$, then $K*u(x)$ is a smooth function that is bounded away from 0 in an interval $I_0$ centred at $x_0$.  The function $u(x)$ satisfies \eqref{e2.2.22}  with $\lambda(x):=K*u(x)$ in $I_0$, so that it is smooth by the smoothness of $\lambda\mapsto u_\lambda$ and that of $K*u$. 

\noindent Let us prove that a steady state $u(x)$ to Equation \eqref{e2.1.1} is below 1. As $u$ is bounded, we have $tu\leq 1$ if $t>0$ is small, so that there is a maximal $t$, called $t_{max}$, such that $tu\leq1$ if $t\leq t_{max}$. Assume $t_{max}<1$ and set $v=1-t_{max}u$, we have, setting $a_m(x)=-\di\frac{f\bigl(t_{max}u(x)\bigl)}{1-t_{max}u(x)}$:
$$
\begin{array}{rll}
\di\int_\RR K(x-y)\bigl(v(x)-v(y)\bigl)~dy=&a_m(x)v(x)+f\bigl(t_{max}u(x)\bigl)-t_{\max}f(u(x))\\
=&a_m(x)v(x)+\biggl(\di\frac{f\bigl(t_{max}u(x)\bigl)}{t_{max}u(x)}-\di\frac{f(u(x))}{u(x)}\biggl)t_{max}u(x)\\
\geq& a_m(x)v(x)\ \hbox{because $u\mapsto\di\frac{f(u)}u$ is nonincreasing}.
\end{array}
$$
If there is a contact point, that is, a point $x_0$ such that $t_{max}u(x_0)=1$, then we have $K*v(x_0)=0$, so that $t_{max}u\equiv1$ by Proposition \ref{p2.2.3}. This is impossible, as $t_{max}<1$ and $f$ does not have any zero beyond $u=1$.  If there is no contact point, there is a sequence $(x_n)_n$ going to infinity such that 
$
\di \lim_{n\to+\infty}u(x_n)=1.
$
As $u'$ is a bounded function, the sequence $(u_n)_n$, given by
$
u_n(x)=u(x+x_n)
$
is equicontinuous, therefore converges, up to a subsequence, to a steady state of  \eqref{e2.1.1} which additionally  has a contact point at $x=0$. This is once again a contradiction, so that we have $u\leq1$.

\noindent We would like to prove, symmetrically, that $u\geq1$, so that we have eventually $u\equiv1$. This, however, requires $u$ to be bounded away from zero, so that we have $u\geq t$ for $t$ small. This is what we are going to prove now, which will put an end to the proof of the theorem.  The computation \eqref{e2.2.2221} yields that, for $\lambda>0$ small enough, we have 
$
\di\mathcal{J}e_\lambda\leq\frac{f'(0)}2e_\lambda,
$
 and 
$
\di e_\lambda\geq0\ \hbox{on}\ [-\pi/{2\lambda},\pi/2\lambda]$.
Let us define $u_0>0$ such that $f(u)\geq {f'(0)u}/2$ for $u\leq u_0$, and let us choose $\delta_0>0$ small enough so that 
$$\delta_0\max_{\RR}e_\lambda\leq u_0.
$$
Thus, $\delta_0 e_\lambda\un_{(-\pi/2\lambda,\psi/2\lambda)}$ is, by virtue of Proposition \ref{p2.2.2}, a subsolution to \eqref{e2.2.20}. As it is compactly supported, we may assume, at the 
expense of decreasing $\delta_0$, that it is below $u$. Now, consider the set of all its translations
$
\{\delta_0e_\lambda(.+\alpha),\ \alpha\in\RR\}.
$
Then, either they are all below $u$, and we have $u\geq\delta_0$. Or such is not the case, and there is $\alpha_0\in\RR$ such that $\delta_0e_\lambda(.+\alpha_0)$ is below $u$, while having a contact point $x_0$ with $u$. We have therefore 
$$
K*\bigl(u-\delta_0e_\lambda(.+\alpha_0)\bigl)(x_0)=0.
$$
Application of Propostion \ref{p2.2.3} yields $u\equiv\delta_0e_\lambda-\lambda(.+\alpha_0)$, an impossibility once again. \hfill$\Box$

\noindent Once this result is in hand, the next step is to try to apply it to infer the large time behaviour to the initial value problem \eqref{e2.1.1}. Let us see at once that this is not the most difficult task on earth, and that it is in fact not limites to Fisher-KPP type nonlinearities.
\begin{theorem}
\label{t2.2.251}
Assume $f$ to be merely smooth, positive on $(0,1)$, with $f(0)=0$ and $f'(0)>0$. Let $u(t,x)$ be the solution of \eqref{e2.1.1} with a smooth nonnegative nonzero compactly supported initial datum. Then
$$
\lim_{t\to+\infty}u(t,x)=1,
$$
uniformly on each compact subset of $\RR$.
\end{theorem}
\noindent{\sc Proof.} If $\delta$ is small enough and $\lambda>1$, there is $\lambda_0>0$ and $\delta_0>0$ such that $f(u)\geq\lambda_0u$ for $0\leq u\leq\delta_0$. We have just seen that 
$$
\underline u_{\delta,\lambda}(x)=\left\{
\begin{array}{rll}
\delta e_\lambda(x)\ &\hbox{if $-\pi/2\lambda\leq x\leq\pi/2\lambda$}\\
0 &\hbox{everywhere else}
\end{array}
\right.
$$
 is a subsolution to the 
steady equation $\mathcal{J}u=f(u)$, provided that $\lambda\leq\lambda_0$, $\delta\leq\delta_0$ and $\lambda$ small enough so that $\pi/2\lambda>1$. Because of the strong comparison principle Theorem \ref{t2.1.1500}, we have $u(1,x)>0$ everywhere, so that we have
$\underline u_{\delta,\lambda}(x)\leq u(1,x)$ even if it means restricting $\delta$ further.

\noindent Consider $\underline f$ a concave function such that $\underline f>0$ on $(0,\delta_0)$, $\underline f'(0)=\delta_0$, and $\underline f(\delta_0)=0$. Clearly, with the chosen values of $\delta$ and $\lambda$, we have $\underline f(u)\leq f(u)$ and the function $\underline u_{\delta,\lambda}$ is a subsolution to the steady equation $\mathcal{J}u=\underline f(u)$. 

\noindent Let $\underline v(t,x)$ be the solution of 
\begin{equation}
\label{e2.2.340}
\begin{array}{rll}
\partial_tv+\mathcal{J}v=&\underline f(v)\quad (t>1,\ x\in\RR)\\
v(1,x)=&\underline v_{\delta,\lambda}(x).
\end{array}
\end{equation}
From Remark \ref{r2.2.4}, it is time increasing. As $\delta\leq\delta_0$ we have $\underline v(t,x)\leq\delta_0$. Therefore $\underline v(t,x)$ converges, pointwise to a function $\underline v_\infty(x)$. Let us show that $\underline v_\infty\equiv\delta_0$, for this we need to know that $v_\infty$ solves the stationary equation $\mathcal{J}u=\underline f(u)$; as Ascoli's theorem is not available here we need a little more care. We write the equation for $\underline v$ as
$
\underline v_t+\K \underline v=K*\underline v+f(\underline v),
$
so that, by the Duhamel formula and the change of variable $s\mapsto t-s$, we have
$$
\underline v(t,x)=e^{-t\K}\underline v(1,x)+\int_1^te^{-s\K}\biggl(K*\underline v(t-s,x)+\underline f\bigl(\underline v(t-s,x)\bigl)\biggl)~ds.
$$
As $\di\lim_{t\to+\infty}\underline v(t-s,x)=\underline v_\infty(x)$ pointwise in $s$ and $x$, Lebesgue's dominated convergence theorem implies, sending $t$ to infinity,  that
$
\K\underline v_\infty(x)=K*\underline v_\infty(x)+\underline f\bigl(u_\infty(x)\bigl).
$
Thus $\underline v_\infty$ is a steady solution. From Theorem \ref{t2.2.10}, we have $\underline v_\infty\equiv\delta_0$.

\noindent Let us now drop $\underline f$ and consider the solution $\underline u(t,x)$ of 
\begin{equation}
\label{e2.2.341}
\partial_tu+\mathcal{J}u=\underline f(u)\ (t>1,\ x\in\RR),\quad
u(1,x)=\underline  u_{\delta,\lambda}(x)
\end{equation}
Arguing as above, we discover a function $\underline u_\infty(x)$ such that $\di\lim_{t\to+\infty}\underline u(t,x)=\underline u_\infty(x)$, pointwise in $x$. As $\underline f\leq f$, the additional information that we have gained is $\underline u_\infty(x)\geq\delta_0$ everywhere.  

\noindent Let $\underline U(t)$ be the solution of the ODE $\dot{\underline U}=f(\underline U)$, with the initial datum $U(0)=\delta_0.$
From Remark \ref{r2.2.400}, we may apply the comparison principle Proposition \ref{p2.1.1}, that is, $u_\gamma(t,x)\geq\underline U(t)$. As $\di\lim_{t\to+\infty}U(t)=1$,
we have $\underline u_\infty=1$. This shows that our original solution converges to 1 for large times, pointwise in $x$. As the limit is continuous, and as $u(t,x)\geq \underline u(t,x)$, which converges monotonically to 1, the convergence is uniform on every compact set by Dini's Theorem. \hfill$\Box$

\section{Linear waves}\label{s2.3}
\noindent Linear waves are special solutions of \eqref{e2.1.1}, linearised around the steady state $u\equiv0$, that is 
\begin{equation}
\label{e2.2.2}
\partial_tv+\JJ v=f'(0)v,\ \ t\in\RR,\ x\in\RR.
\end{equation}
They have the form $v(t,x)=\phi_\lambda(x-ct)$, with $\phi_\lambda(x)=e^{-\lambda x}$.
The function $\phi_\lambda$ solves 
\begin{equation}
\label{e2.2.5}
\JJ\phi-c\phi'=f'(0)\phi,\ \ x\in\RR.
\end{equation}
When $c>0$ the wave is said to propagate rightwards, in the opposite case it is said to propagate leftwards. The situation being perfectly symmetric we will look for
rightwards propagating waves.

\noindent  Let us define the function
\begin{equation}
\label{e2.3.41}
D_c(\lambda)=2\int_{0}^1\bigl(\cosh(\lambda x)-1\bigl)K(x)~dx-c\lambda+f'(0).
\end{equation}
We will have many encounters with this function. Plugging  the expression  $\phi_\lambda(x)=e^{-\lambda x}$ into \eqref{e2.2.2} yields the equation
$D_c(\lambda)=0.
$
We have
\begin{equation}
\label{e2.3.42}
D_c'(\lambda)=2\di\int_0^1 xK(x)\sinh(\lambda x)~dx-c,\ \ D_c''(\lambda)=2\int_{0}^1x^2K(x)(\cosh(\lambda x)~dx.
\end{equation}
From the uniform strict convexity of $u\mapsto\cosh u$, there exists a critical $c_*>0$ such that  $D_c(\lambda)=0$   has two positive solutions $\lambda_-(c)<\lambda_+(c)$ if 
$c>c_*$, no solution if $c<c_*$, and exactly one, that we call $\lambda_*$ if $c=c_*$. If $c>c_*$ we have $D_c'(\lambda_-(c))<0<D_c'(\lambda_+(c))$, while 
 we have $
D_{c_*}'(\lambda_*)=0,\ D_{c_*}''(\lambda_*)>0$ if $c=c_*$.
We also notice that $xe^{-\lambda_* x}$  solves equation \eqref{e2.2.5}.

\noindent It will be useful to look for complex solutions of  $D_c(\lambda)=0$  for $c$ slightly below $c_*$. In this range, we write  $D_c(\lambda)=0$  under the form
$
\di\frac{D_{c_*}(\lambda_*)}2(\lambda-\lambda_*)^2=-\lambda(c_*-c)+O(\lambda-\lambda_*)^3,
$
that is,
$$
\lambda-\lambda_*=\pm i\sqrt{\frac{2(c_*-c)}{D''_{c_*}(\lambda_*)}}+O(c_*-c)^{3/4}.
$$
We write, therefore
$
\lambda=\lambda_*(c)\pm i\omega_*(c),\ \ \lambda_*(c)\in\RR,\ \omega_*(c)\in\RR,
$
with 
\begin{equation}
\label{e2.2.28}
\lambda_*(c)=\lambda_*+O(c_*-c)^{3/4},\ \ \omega_*(c)=\sqrt{\frac{2(c_*-c)}{D''_{c_*}(\lambda_*)}}+O(c_*-c)^{3/4}.
\end{equation}
The corresponding linear wave $\phi_\lambda$ may be taken as
\begin{equation}
\label{e2.2.29}
\phi_\lambda(x)=e^{-\lambda_*(c)x}\cos\omega_*(c)x.
\end{equation}
It will also be useful to examine the complex linear waves of 
$
\JJ\phi-c_*\phi'-f'(0)\phi=\mu\phi,
$
$\mu>0$ small. We have
$
D_{c_*}'''(\lambda_*)=2\di\int_0^{+\infty} x^3\sinh(\lambda_*x)~dx>0.
$
Therefore,  if $\phi(x)$ is sought for under the form $e^{-\lambda x}$, we have
$$
\frac12D_{c_*}''(\lambda_*)\bigl(\lambda-\lambda_*\bigl)^2+\frac16D_{c_*}'''(\lambda_*)\bigl(\lambda-\lambda_*\bigl)^3+O\bigl(\lambda-\lambda_*\bigl)^4+\mu=0,
$$
that is
\begin{equation}
\label{e2.2.27}
\lambda_\pm(\mu)=\lambda_*+\frac{D_{c_*}'''(\lambda_*)\mu}{3\bigl(D_{c_*}''(\lambda_*)\bigl)^2}\pm i\sqrt{\frac{2\mu}{D_{c_*}''(\lambda_*)}}+O(\mu^{3/2}).
\end{equation}
the decay exponent   slightly increases with $\mu$, this will turn out to be important in Chapter \ref{ZFK_short_range}.

\noindent Let us finally recall   the computation \eqref{e2.2.22}, for $c=0$. Setting this time
$
e_\lambda(x)=\cos~\lambda x$,
we have $\JJ e_\lambda=\omega_0(\lambda)e_\lambda,$
with 
$\omega_0(\lambda)= \lambda^2<\!y^2K\!>+O_{\lambda\to0}(\lambda^3).
$
\section{A heat kernel asymptotics}\label{s2.4}
\noindent The goal of this section is to understand in more depth the Cauchy problem for equation \eqref{e2.1.1}, linearised at 0, at the critical speed. In other words we wish to estimate the solutions $w(t,x)$ of 
\begin{equation}
\label{e2.4.10}
\partial_tw+\mathcal{J}w-c_*\partial_xw=f'(0)w\ (t>0,x\in\RR), \quad w(0,x)=w_0(x)
\end{equation}
where $w_0$ is a smooth function that decays at infinity faster than $e^{-\lambda_*x}$. More precisely, we will ask  that the function  $x\mapsto e^{\lambda_*x}w_0$ is in every $H^m(\RR)$. Setting $v(t,x)=e^{\lambda_*x}w(t,x)$ and $v_0(x)=e^{\lambda_*x}w_0(x)$, and using $D_{c_*}(\lambda_*)=D_{c_*}'(\lambda_*)=0$, equation \eqref{e2.4.10} reduces to computing $e^{-t\II_*}v_0$, with
\begin{equation}
\label{e2.4.6}
\II_*v(x)=-\int_{\RR}K_*(x-y)\bigl(v(y)-v(x)-(y-x)v_x(x)\bigl)~dy,\ \hbox{with}\ K_*(x)=e^{\lambda_*x}K(x).
\end{equation}
 Note that the new kernel $K_*$ is no longer symmetric. Also note that $\II_*v\geq0$ if $v$ is concave.
\subsection{Two identities for $\II_*$}\label{s2.2.1}
The first one that we wish to mention is the analogue of Kato's inequality for the laplacian: for all $v\in W^{2,p}(\RR^N)$ we have
$-\Delta \vert v\vert\leq -\mathrm{sgn}(v)\Delta v$. Similary to that we have, for all Lipschitz function $v$:
\begin{equation}
\label{e2.5.108}
\II_*\vert v\vert\leq\mathrm{sgn}(v)\II_* v.
\end{equation}
We have indeed
$$\begin{array}{rll}
\mathrm{sgn}(v)\II_*v=&\vert v\vert\di\int K_*(x-y)~dy+\partial_x\vert v\vert\int(y-x)K_*(x-y)~dy-\mathrm{sgn}(v)\int K_*(x-y)v(y)~dy\\
\geq&\vert v\vert\di\int K_*(x-y)~dy+\partial_x\vert v\vert\int(y-x)K_*(x-y)~dy-\int K_*(x-y)\vert v(y)\vert dy
=\II_*\vert v\vert.
\end{array}
$$
If now $\mathrm{sgn}^+(v)$ equals 1 if $v>0$, and 0 if $v\leq0$, the identity
${v^+}=(\vert v\vert+v)/2=\mathrm{sgn}^+(v)v$
yields
\begin{equation}
\label{e2.5.132}
\II_*v^+\leq\mathrm{sgn}^+(v)\II_* v.
\end{equation}
We mention the following property, valid for two Lipschitz functions $u(x)$ and $v(x)$:
\begin{equation}
\label{e2.4.50}
\II_*(uv)(x)=u(x)\II_*v(x)+v(x)\II_*u(x)-\int_\RR K_*(x-y)\bigl(u(x)-u(y)\bigl)\bigl(v(x)-v(y)\bigl)~dy.
\end{equation}
The proof is trivial and left to the reader, who
should convince him/herself that this identity is analogue to:
$(-\Delta u)v=u(-\Delta v)-2\nabla u.\nabla v+u(-\Delta v),$ the operator $\II_*$ playing the role of $(-\Delta)$.
\subsection{The main estimate}
\noindent The operator $\II_*$ does not look, at first sight, to have any interesting properties. In fact it even looks boring. One has to go beyond this first impression to discover the 
\begin{theorem}
\label{t2.4.2}
Set 
$d_*=1/2\di\int_\RR z^2e^{-\lambda_*z}K(z)~dz.
$
There is a nonnegative kernel $G_*(t,z)$, defined for $t\geq1$ and $z\in\RR$ and bounded on its domain of definition, such that we have, for all $\gamma\in(0,1/6)$, for all $t\geq1$ and $x\in\RR$:
\begin{equation}
\label{e2.4.11}
\Vert e^{-t\II_*}v_0-{G_*(t,.)}*v_0\Vert_{L^\infty(\RR)}\lesssim e^{-t^{1-2\gamma}}\Vert v_0\Vert_{H^1(\RR)}.
\end{equation}
If we choose
$ \delta>6\gamma,$  we have the following estimates.
\begin{itemize}
\item[---] If $\vert z\vert\leq t^{1/2+\delta}$, then $G_*(t,z)=\di\bigl(1+O(1/{t^{1-\delta}})\bigl){e^{-z^2/4d_*t}}/{\sqrt{4\pi d_*t}}.$
\item[---] If $\vert z\vert\geq t^{1/2+\delta}$, then 
$G_*(t,z)\lesssim e^{-t^\delta}.$
\end{itemize}
\end{theorem}
\noindent{\sc Proof.} As Problem \eqref{e2.4.10} is homogeneous in $x$, we gladly resort to the Fourier transform.  Let $\hat P_*(t,\xi)$ the Fourier transform in $x$ of $e^{-t\II_*}v_0$.
We have, with  the notation \eqref{e2.4.6}:
$$
\begin{array}{rll}
\hat P_*(t,\xi)=&\hat v_0(\xi)\mathrm{exp}\biggl(-t\bigl(\di\int_\RR K_*(z)~dz-i\xi\di\int_\RR zK_*(z)~dz-\hat K_*(\xi)\bigl)\biggl)\\
=&\hat v_0(\xi)\mathrm{exp}\biggl(-t\bigl(\hat K_*(0)+\xi\hat K_*'(0)-\hat K_*(\xi)\bigl)\biggl).
\end{array}
$$
Let us examine the phase in the exponential away from $\xi=0$. As $K_*$ is compactly supported, $\hat K_*$ belongs to the Schwartz class, so that there is 
$M_0>0$ such that , for $\vert\xi\vert\geq M_0$ we have:
$$
\vert\hat K_*(\xi)\vert\leq\frac12\di\int_\RR K_*={\hat K_*(0)}/2
$$
For $\vert\xi\vert\leq M_0$, we start from
$$
\mathrm{Re}~\bigl(\hat K_*(0)-\hat K_*(\xi)\bigl)=\int_{-1}^1\bigl(1-\cos(z\xi)\bigl)K_*(z)~dz.
$$
Pick $\gamma\in(0,1/2)$, and 
$t_\gamma=\bigl(\pi/4\bigl)^{-1/\gamma}.
$
For $t\geq t_\gamma$ and $\vert\xi\vert\in(t^{-\gamma},M]$, we use $1-\cos u\lesssim -u^2$ for $u$ close to $t^{-\gamma}$, together with the fact that 
$u\mapsto\cos u$ is less than 1, while not identically equal to 1. We obtain:
$$\int_{-1}^1\bigl(1-\cos(z\xi)\bigl)K_*(z)~dz\gtrsim t^{-2\gamma}.
$$
Putting everything together we obtain, for $\vert\xi\vert\geq t^{-\gamma}$:
\begin{equation}
\label{e2.4.17}
\mathrm{Re}~\bigl(\hat K_*(0)-\hat K_*(\xi)\bigl)\gtrsim t^{-2\gamma}.
\end{equation}
This calls for the decomposition
$$
\hat P_*(t,\xi)=\bigl(\un_{\vert\xi\lvert\leq t^{-\gamma}}+\un_{\vert\xi\lvert\geq t^{-\gamma}}\bigl)P_*(t,\xi):=\hat Q_*(t,\xi)+\hat R_*(t,\xi).
$$
From \eqref{e2.4.17} we have 
$\vert\hat R(t,\xi)\vert\lesssim e^{-t^{1-2\gamma}}.$
Therefore, if $R_*(t)$ denotes  the conjugate Fourier transform of $\hat R_*(t,.)$,
 we have, by Plancherel's equality:
$$
\Vert R_*(t)v_0\Vert_{L^\infty(\RR_x)} \lesssim e^{-t^{1-2\gamma}}\Vert\hat v_0\Vert_{L^1_\xi}\lesssim e^{-t^{1-2\gamma}}\Vert \vert\xi v_0\Vert_{L^2_\xi}=e^{-t^{1-2\gamma}}\Vert v_0\Vert_{H^1_x}.
$$
It remains to study $Q_*(t,x)$, the inverse Fourier transform of $\hat Q_*(t,\xi)$ .  Let us write
$$
\hat K_*(0)+\xi\hat K_*'(0)-\hat K_*(\xi)=-d_*\xi^2-\frac{\hat K_*'''(0)\xi^3}6-\frac{\xi^4}6\int_0^1(1-\sigma)^3\hat K_*^{(4)}(\sigma\xi)~d\sigma,
$$
so that we have, setting $\zeta=\sqrt t\xi$ :
$$
 Q_*(t,x)=\frac1{2\pi\sqrt t}\di\int_{-t^{1/2-\gamma}}^{t^{1/2-\gamma}}\int_\RR \mathrm{exp}\biggl(\frac{i\zeta (x-y)}{\sqrt t}-d_*\zeta^2-\frac{\hat K_*'''(0)\zeta^3}{6\sqrt t}-\frac{\zeta^4}{6 t}\int_0^1(1-\sigma)^3\hat K_*^{(4)}(\frac{\sigma\zeta}{\sqrt t})~d\sigma\bigl)\biggl) v_0(y)~dy.
 $$
By Fubini's Theorem,  we gave $Q_*(t,.)=G_*(t,.)*v_0$, with
$$
G_*(t,z)=\frac1{2\pi\sqrt t}\di\int_{-t^{1/2-\gamma}}^{t^{1/2-\gamma}} e^{\Phi_*(t,z,\zeta)}d\zeta,
$$
with
\begin{equation}
\label{e2.4.20}
\Phi_*(t,z,\zeta)= i \frac{z\zeta}{\sqrt t}-d_*\zeta^2-\frac{\hat K_*'''(0)\zeta^3}{6\sqrt t}-\frac{\zeta^4}{6 t}\int_0^1(1-\sigma)^3\hat K_*^{(4)}(\frac{\sigma\zeta}{\sqrt t})~d\sigma.
\end{equation}
Inspired by the computation of the Fourier transform of the Gaussian, we  want to change the line of integration. Define the line
\begin{equation}
\label{e2.4.23}
\Gamma_*=\{\zeta=\eta+\frac{iz}{2d_*\sqrt t},\ \ -t^\gamma\leq\eta\leq t^\gamma\},
\end{equation}
If $\hat K_*$ was a quadratic polynomial, as in the Gaussian case, the effect of this shift would be to make $\Phi_*$ real, equal to $-d_*\eta^2-\di\frac{z^2}{4d_*t}$ and this would finish the computation.
However there is a slight glitch, that is, $\hat K_*$ is not a quadratic polynomial.  And so, if we wish to mimick the classical computation, we have to make sure that the non-quadratic part of $\Phi_*$ is actually negligible. As a consequence, we will
have to worry about the size of $z$, and this will also constrain $\gamma$. The constraints are therefore:
$$
\frac1{\sqrt t}(\eta+\frac{iz}{2d_*\sqrt t})^3\ll 1,\ \frac1{t}(\eta+\frac{iz}{2d_*\sqrt t})^4\ll 1,\ \ \ \hbox{with}\ \vert\eta\vert\leq t^{\gamma}.
$$
It is enough to impose ${\vert\eta\vert}/{t^{1/6}}\ll1,$ and ${\vert z\vert^4}/{t^{3}}\ll1$.
As $\vert \eta\vert\leq t^\gamma$, we choose, in order to fulfill these two conditions:
${\vert z\vert}/{\sqrt t}\leq t^\delta,$
with $\delta\in(0,1/4)$, and possibly smaller if this suits us.  
Let us deal at once with $z$ outside the range ${\vert z\vert}/{\sqrt t}\leq t^\delta,$. If $z\geq0$, define  the vertical segments
$$
\Gamma_\pm=\{\zeta=\pm t^{\gamma}+i\eta,\ \ 0\leq\eta\leq 1\}.
$$
and the horizontal line
$$
\Gamma_*=\{\zeta=\eta+i,\ \ -t^\gamma\leq\eta\leq t^\gamma\}.
$$
If $z\leq0$ we replace $\eta+i$ by $\eta-i$, and $\pm t^{\gamma}+i\eta$ by $\pm t^{\gamma}-i\eta$.  We have
$$
\di\int_{-t^{1/2-\gamma}}^{t^{1/2-\gamma}} e^{\Phi_*(t,z,\zeta)}d\zeta=\int_{\gamma_-\cup\Gamma_*\cup\gamma_+}e^{\Phi_*(t,z,\zeta)}d\zeta.
$$
For $z\geq0$ (the same would of course also work for $z\leq0$), the following holds: on $\Gamma_\pm$, we have 
$$
\Phi_*(t,z,\zeta)=\pm\di{izt^{-\gamma}}-\frac{\eta z}{\sqrt t}-d_*(\pm t^{1/2-\gamma}+i\eta)^2+O(t^{\gamma-1/2}),\
$$ so that $\mathrm{Re}~\Phi_*(t,z,\zeta)\leq-d_*t^{1-2\gamma}(1+O(1)).$ On $\Gamma_*$  we have 
$$
\Phi_*(t,z,\zeta)=-\di\frac{z}{\sqrt t}+\frac{i\eta z}{\sqrt t}-d_*(\eta+i)^2+O(t^{\gamma-1/2}),
$$
{so that} $
\mathrm{Re}~\Phi_*(t,z,\zeta)\leq-\di\frac{z}{\sqrt t}+O(1)\leq-t^{\delta}+O(1).$
Thus we have
$
\di\biggl\vert\di\int_{-t^{1/2-\gamma}}^{t^{1/2-\gamma}} e^{\Phi_*(t,z,\zeta)}~d\zeta\biggl\vert\lesssim e^{-t^{\delta}},
$
and this  proves $\vert G_*(t,z)\vert \lesssim e^{-t^\delta}$.

\noindent Finally, let us deal with $z$ in the range ${\vert z\vert}/{\sqrt t}\leq t^\delta$. This time, if $z\geq0$, define  the vertical segments
$$
\Gamma_\pm=\{\pm t^{\gamma}+i\eta,\ \ 0\leq\eta\leq\frac{z}{2d_*\sqrt t} \}.
$$
while the horizontal line $\Gamma_*$ is defined by \eqref{e2.4.23}.
This time, the definition does not change if $z\leq0$.  
On $\Gamma_\pm$ we have, as $\delta<1/2-\gamma$:
$$
\mathrm{Re}~\Phi_*(t,z,\zeta)=\di-\frac{z\eta}{\sqrt t}-d_*(t^{1-2\gamma}-t^{2\delta})+O(t^{\gamma-1/2})
\leq-d_*t^{1-2\gamma}(1+o_{t\to+\infty}(1)).
$$
This  also accounts for the contribution of $\di\int_{\Gamma_-\cup\Gamma_+}e^{\Phi_*(t,z,\zeta)}dz$  on the new vertical segments $\Gamma_\pm$ in the area $\vert z\vert\geq t^{1/2+\delta}$. And so, what is left to us now is 
the Gaussian-like part. Recalling that we have this time $\zeta=\eta+\di\frac{iz}{2d_*\sqrt {t}}$ we decompose $\Phi_*(t,z,\zeta))=-d_*\eta^2-\di\frac{z^2}{4d_*t}-\psi_*(t,z,\eta).
$
From \eqref{e2.4.20}, and remembering that $\vert z\vert\leq t^{1/2+\delta}$, we estimate $\psi_*$ as
$$
\begin{array}{rll}
\psi_*(t,z,\eta)=&\di\frac{\hat K_*'''(0)\zeta^3}{6\sqrt t}+\frac{\zeta^4}{6 t}\int_0^1(1-\sigma)^3\hat K_*^{(4)}(\frac{\sigma\zeta}{\sqrt t})~d\sigma\\
=&\di\frac{\hat K_*'''(0)}{6\sqrt t}\bigl(\eta+\frac{iz}{2d_*\sqrt t}\bigl)^3+O(t^{4\gamma-1})=\di\frac{\hat K_*'''(0)\eta^3}{6\sqrt t}+O(t^{4\gamma-1}),
\end{array}
$$
where we have taken the worst estimate for $z$ and $\eta$, according to  the fact that $\vert\eta\vert\leq t^{1/2-\gamma}$.  Also notice that we have 
$\psi_*(t,z,\eta)$ is an $o_{t\to+\infty}(1)$ uniformly in $z\in(-t^{1/2+\delta},t^{1/2+\delta})$ and $\eta\in(-t^{\gamma},t^{\gamma})$, as it decays at least like $t^{4\gamma-1}$.  
With this final reduction, we may  compute $G_*$ in the zone $\vert z\vert\leq t^{1/2+\delta}$. We remember that, in the original integral  expanding the exponential as:
$$
\begin{array}{rll}
G_*(t,z)=&\di\frac1{2\pi\sqrt t}\di\int_{-t^{1/2-\gamma}}^{t^{1/2-\gamma}}\mathrm{exp}\biggl(-d_*\eta^2-\frac{z^2}{4d_*\sqrt t}-\psi_*(t,z,\eta)\biggl)~d\eta\\
=&\di \di\frac1{2\pi\sqrt t}\int_{-t^{1/2-\gamma}}^{t^{1/2-\gamma}}\biggl(1+\frac{\hat K_*'''(0)\eta^3}{6\sqrt t}+O(\psi_*(t,z,\eta))^2\biggl)e^{-d_*\eta^2-z^2/4d_*t}d\eta\\
=&\di\frac{e^{-z^2/4d_*t}}{2\pi\sqrt t}\int_{-t^{1/2-\gamma}}^{t^{1/2-\gamma}}\biggl(1+O(\psi_*(t,z,\eta))^2\biggl)e^{-d_*\eta^2}d\eta\ \ \hbox{by the oddness of $\eta\mapsto\eta^3$}\\
=&\di\bigl(1+O(\frac1{t^{1-6\gamma}})\bigl)\frac{e^{-z^2/4d_*t}}{\sqrt{4\pi d_* t}}.
\end{array}
$$
As $\gamma$ and $\delta$ are as small as we wish, and as we have chosen $\gamma<\delta/6$, the proof is complete. \hfill$\square$
\begin{remark}\label{r2.4.1}
One could use less derivatives of $v_0$, and resort to the $L^1$ norm of $\hat v$ in the estimation of $R_*$. The Plancherel formula would  then give
$
e(-t\II_*)v_0=G_*(t,.)*v_0+R_*(t)v_0,
$
where $R_*(t)$ is a linear continuous from $L^2(\RR)$ to $L^2(\RR)$, with $\Vert R_*(t)\Vert_{L^2(\RR)\to L^2(\RR}=O(e^{-t^{1-2\gamma}})$. This allows for the following estimate, similar to that of the classical heat equation:
$$\Vert e^{-t\II_*}v_0\Vert_{L^2(\RR)}\lesssim\frac1{t^{1/4}}\Vert v_0\Vert_{L^2(\RR)}.
$$
\end{remark}
\begin{remark}
\label{r2.4.3}
Estimating the deviation between $e^{-t\II_*}v_0$ and $G_*(t,.)*v_0$ by the $H^1$ norm of $v_0$ is not optimal. One could indeed do it with the $H^\sigma$ norm of $v_0$, for all $\sigma\in(1/2,1)$. It is, however, not easier to manipulate a nonlocal norm,  and, as we will see in Chapter \ref{short_range}, does not fundamentally makes things easier.
\end{remark}
\begin{remark}\label{r2.4.2}
One may wonder why such a universal behaviour arises for the heat kernel. In fact, Theorem \ref{t2.4.2} has nothing surprising: in spirit, its proof is the same as that of the classical Central Limit Theorem, which 
essentially assumes, in the landscape,  a characteristic function that one can expand near 0. It is all the less surprising that the full model \eqref{e1.1.1} has a probabilistic interpretation, at least for some very special functions $f$.
\end{remark}
\subsection{The solution emanating from a well-spread datum}
\noindent The objective of this section is to devise specific estimates for $e^{-\tau\II_*}v_s$,  where $v_s$ is an initial datum that is spread over the characteristic length $\sqrt s\gg1$. Of interest to us are times $\tau$ of the form $0\leq\tau\leq s^\e$, with $\e\ll1$. Anticipating on the main part of this book, the  reason for our interest in such an exercise concerns the sharp asymptotic behaviour of the solutons of the Fisher-KPP equation. 
 Theorem \ref{t2.4.2} below will be sufficient to devise bounds on the solution  at the door of the diffusive zone, that will be of the correct order of magnitude in time.  That will be enough to locate the level sets of $v$ with $O_{t\to+\infty}(1)$ precision. However, in order to reach the next order $o_{t\to+\infty}(1)$, we will need to update our observations at various  large times, and Theorem \ref{t2.4.2} will not be able to handle properly these updates on times ranges slightly larger than $s$, say, a tiny power of $s$.  
 
\noindent  Another way to see it, seemingly unrelated but that has in fact deep links with the above issue, is that, in the diffusive variables, the integral operator has no particular properties of regularisation. In any case, here is how we make up for the gap. 
 \begin{theorem}
\label{t2.4.4}
Consider a function $v_0$ belonging to $H^m(\RR)$, for all integer $m$. Consider $s>0$, a real number having the possibility of being very large. Set
 \begin{equation}
\label{e2.4.70}
v_s(y)=v_0({y}/{\sqrt s}).
\end{equation}
Pick $\omega\in(0,10^{-10})$. Then, there is a (possibly very large) integer $m_\omega$ such that, for $\tau\in(0,s^{1/2-\omega})$ we have:
\begin{equation}
\label{e2.4.71}
\big\vert e^{-\tau\II_*}v_s(x)-e^{\tau d_*\partial_{xx}}v_s(x)\big\vert\lesssim\frac{\tau}{s^{3/2-4\omega}} \bigl(\Vert v_0\Vert_{L^1}+\Vert v_0\Vert_{H^{m_\omega}}\bigl).
\end{equation}
\end{theorem}
\noindent{\sc Proof.} The strategy for proving the theorem is essentially the same as that for Theorem \ref{t2.4.2}: a Fourier integral will be subject to various changes of variables that will allow an asymptotic expansion of the phase, followed by a passage in the complex plane. There will, however, be an important difference: the passage to complex variables will have, in order to keep the computation meaningful, to occur at a later stage of the computations, thus leading to a lesser precision than that provided by Theorem \ref{t2.4.2}.  We set
$
e^{-\tau\II_*}v_s(x)=Q_*(s,\tau;x)+R_*(s,\tau;x),
$
with
$$
\begin{array}{rll}
Q_*(s,\tau;x)=&\di\frac1{2\pi}\int_{\vert\xi\vert\leq s^{\omega -1/2}}e^{ix\xi-\tau\bigl(\hat K_*(0)+\xi\hat K_*(0)-\hat K(\xi)\bigl)}\hat v_s(\xi)~d\xi\\
=&\di\frac1{2\pi}\int_{\vert\xi\vert\leq s^{\omega -1/2}}e^{i(x-y)\xi-\tau\bigl(\hat K_*(0)+\xi\hat K_*(0)-\hat K(\xi)\bigl)}v_s(y)~d\xi dy\\
=&\di\frac{1}{2\pi}\di\int_{\vert\zeta\vert\leq s^{\omega}}e^{i(\frac{x}{\sqrt s}-z)\zeta-\tau\bigl(\hat K_*(0)+\frac{\zeta}{\sqrt s}\hat K_*(0)-\hat K(\frac{\zeta}{\sqrt s})\bigl)}v_0(z)~dzd\zeta.
\end{array}
$$
We have not written explicitely the dependence in $v_0$ or $v_s$. The integral $R_*$ is readily estimated, using  the classical identity $\hat v_s(\xi)=\sqrt s\hat v_0(\sqrt s\xi)$,  the fact that $v_0$ is in every $H^m$, and that the real part of $\hat K_*(0)-\hat K_*(\xi)$ is nonnegative:
$$
\begin{array}{rll}
R_*(s,\tau;x)\leq&\di\frac{\sqrt s}{2\pi}\di\int_{\vert\xi\vert\geq s^{\omega -1/2}}\hat v_0(\sqrt s\xi)~d\xi\\
=&\di\frac{1}{2\pi}\di\int_{\vert\zeta\vert\geq s^{\omega}}\hat v_0(\zeta)~d\zeta\\
\leq&\di\frac{1}{2\pi}\biggl(\int_{\vert\zeta\vert\geq s^{\omega}}\frac{d\zeta}{1+\zeta^{2m}}\biggl)^{1/2}\biggl(\int_{\RR}\bigl(1+\zeta^{2m}\bigl)\vert\hat v_0(\zeta)\vert^2d\zeta\biggl)^{1/2}\\
\lesssim &s^{-\omega  m}\Vert v_0\Vert_{H^m}.
\end{array}
$$
Choosing $m$ large enough yields estimate \eqref{e2.4.71} for $R_*$. The integral $Q_*$ is a little more involved, which should be no surprise in view of Theorem \ref{t2.4.2}. Let us set
$$
Q_*(s,\tau;x)=\frac1{2\pi}\int_{\RR}G_*(s,\tau;\frac{x}{\sqrt s}-y)v_0(y)~dy,
$$
with, due to the change of variable $\zeta=\sqrt s\xi$:
$$
G_*(s,\tau;z)=\int_{-s^\omega }^{s^\omega }\mathrm{exp}\biggl(iz\zeta-\tau\bigl(\hat K_*(0)+\frac\zeta{\sqrt s}\hat K_*(0)-\hat K_*(\frac\zeta{\sqrt s})\bigl)\biggl)~d\zeta.
$$
Let $\Phi_*(s,\tau;z,\zeta)$ be the phase in the exponential; as in Theorem \ref{t2.4.2} we write
$$
\Phi_*(s,\tau;x,\zeta)=iz\zeta-\frac{\tau d_*\zeta^2}s-\psi_*(s,\tau;z\zeta),
$$
with, this time:
$$
\psi_*(s,\tau;z,\zeta)=\frac{\tau\hat K_*'''(0)\zeta^3}{6 s^{3/2}}+\frac{\tau\zeta^4}{6 s^2}\int_0^1(1-\sigma)^3\hat K_*^{(4)}(\frac{\sigma\zeta}{\sqrt s})~d\sigma.
$$
This entails, as $\zeta$ is in the range $(-s^\omega ,s^\omega )$:
$$
\begin{array}{rll}
e^{\Phi_*(s,\tau;z,\zeta)}=&\di\biggl(1-\psi(s,\tau;z,\zeta)+O\bigl(\psi(s,\tau;z,\zeta)\bigl)^2\biggl)\mathrm{exp}\bigl(iz\zeta-\frac{\tau d_*\zeta^2}s\bigl)\\
=&\di\biggl(1-\frac{\tau\hat K_*'''(0)\zeta^3}{6 s^{3/2}}\biggl)\mathrm{exp}\bigl(iz\zeta-\frac{\tau d_*\zeta^2}s\bigl)+O(\di\frac{\tau}{s^{2-4\omega }}).
\end{array}
$$
This is where we give up the possibility of moving the integration path, at the expense of a less precise error bound. The reward is that it enables an easy computation of $G_*$. Notice, that, in the second line below, that the difference between
the integral over $(-s^\omega ,s^\omega )$ is exponentially small in $s^\omega $, so that it is absorbed in the $O({\tau}/{s^{2-4\omega }})$. Once we have realised this, we may shift the integration path from $\RR$ to $\RR+\di\frac{isz}{2d_*\tau}$ with our mind at peace, and write the following equality:
$$
G_*(s,\tau,z)=\di\int_{-s^\omega }^{s^\omega }\biggl(1-\frac{\tau\hat K_*'''(0)\zeta^3}{6 s^{3/2}}\biggl)\mathrm{exp}\bigl(iz\zeta-\frac{\tau d_*\zeta^2}s\bigl)~d\zeta+O(\di\frac{\tau}{s^{2-4\omega }}).
$$
We simply estimate the $\zeta^3$ integral as:
$$
\biggl\vert\frac{1}{s^{3/2}}\di\int_{-s^\omega }^{s^\omega }\zeta^3\mathrm{exp}\bigl(iz\zeta-\frac{\tau d_*\zeta^2}s\bigl)~d\zeta\biggl\vert\lesssim \frac{\tau}{s^{3/2-4\omega }}.
$$
Therefore we have
$$
\begin{array}{rll}
G_*(s,\tau,z)=&\di\int_{-s^\omega }^{s^\omega }\mathrm{exp}\bigl(iz\zeta-\frac{\tau d_*\zeta^2}s\bigl)~d\zeta+O(\di\frac{\tau}{s^{3/2-4\omega }})
=\di\int_{\RR}\mathrm{exp}\bigl(iz\zeta-\frac{\tau d_*\zeta^2}s\bigl)~d\zeta+O(\di\frac{\tau}{s^{3/2-4\omega }})\\
=&\di\sqrt{\frac{\pi s}{d_*\tau}}e^{-\frac{sz^2}{4d_*\tau}}+O(\di\frac{\tau}{s^{3/2-4\omega }}).
\end{array}
$$

\noindent Convolution of $G_*(s,\tau,.)$ with $v_s$, and (by now usual) changes of variables yield the result. \hfill$\Box$
\subsection{The behaviour of $e^{-t\II_*}v_0$ away from the diffusive range}\label{s2.4.4}
\noindent As is usually the case when it comes to computing a heat kernel, Theorem \ref{t2.4.2} gives precise information about it behaviour slightly beyond the range of similarity, that is, in the present context, for 
$\di\frac{x}{\sqrt t}\geq t^\delta$, for $\delta>0$ small. In order to use it in comparisons, one needs to know how it behaves for $\vert x\vert$ very large. Here, a precise behaviour is not needed, what is requested is an estimate that will beat the 
algebraic expressions in $t$.
\begin{proposition}
\label{p2.4.1}
If $v_0$ is compactly supported, then, for all $A>0$, there exists $B>0$ such that, if ${x}/{\sqrt{t+1}}\geq B$, we have
$\big\vert e^{-t\II_*}v_0(x)\big\vert\leq \Vert v_0\Vert_\infty e^{-\frac{Ax}{\sqrt{t+1}}}.$
\end{proposition}
\noindent{\sc Proof.} We begin by rephrasing Section \ref{s2.3}, with  the following innocent computation. For $k\in\RR$, and $v(x)=e^{kx}$ we have, by Taylor's formula:
$$
\begin{array}{rll}
\II_*v(x)=&\di\int_\RR K_*(x-y)(x-y)^2\int_0^1\sigma v''\bigl(\sigma(y-x)\bigl)~d\sigma dy
=\di\frac{kv(x)}2\int_\RR zK_*(z)\bigl(e^{kz}-1\bigl)~dz\\
=&d_*k^2\bigl(1+\omega(k)\bigl)v(x),
\end{array}
$$
where $\omega(k)$ is a real analytic function such that $\omega(0)=0$. Setting $\xi={x}/{\sqrt{t+1}}$, we have 
$$
\begin{array}{rll}
\bigl(\partial_t+\II_*)\bigl)e^{-\frac{Ax}{\sqrt{t+1}}}=&\di\frac{A}{t+1}\biggl(\frac\xi2-Ad_*\bigl(1+\omega(\di\frac{A}{\sqrt{t+1}})\bigl)\biggl)e^{-\frac{x}{\sqrt{t+1}}}\\
\geq&0\ \ \hbox{if $\xi\geq Ad_*(1+\Vert\omega\Vert_{L^\infty([0,1])}$}).
\end{array}
$$
As $e^{-t\II_*}v_0(x)\leq\Vert v_0\Vert_\infty$, the maximum principle entails the estimate. \hfill$\Box$

\section{Bibliography, comments, open questions}\label{s2.106}
\subsection*{The initial value problem, steady solutions (Sections \ref{s2.2} and \ref{s2.25})}
\noindent These sections can hardly be seen as original material. I have written it as a gentle introduction, and also to underline how simple the study of the Cauchy Problem can be when the diffusion is an integral operator. 
 A lot is known on the Cauchy Problem for intego-differential problems of the form \eqref{e2.1.1} because comparison is available. A minor point
 should, however, be pointed out here. People having a background in parabolic equations, such as, for example, the humble author of these lines, may find the following fact a little disturbing: the solution of the Dirichlet problem needs not be continuous at the boundary. The reader may have indeed noticed that some care has been taken in the writing of inequalities for $x$ in Proposition \ref{p2.1.200}. Two minutes of not so intense thinking will convince him/her that this situation is perfectly normal, and that what {\rm is} abnormal is the situation of the  Dirichlet fractional Laplacian. There, the continuity of the solution and, even, the H\"older continuity of the solution across the boundary is a consequence of the singularity of the kernel. See Ros Oton-Serra \cite{ROS} for an elaborate discussion of this fact. The phenomenon is of course even more remarkable in the case of the Dirichlet classical Laplacian, where Lipschitz continuity across the boundary is the rule. Therefore, given that our operator displays no singularity, discontinuity is anything but surprising.

\noindent The reader is advised to consult the Habilitation dissertation of Coville \cite{CovHDR}, where all the above  issues, and  more, such as Harnack inequalities, are discussed and placed in their historical context. While it is not so standard nowadays to refer to such documents,  they are accessible, useful and refereed sources of information that it would be a mistake to overlook.

\noindent Unfortunately, this nice comparison feature that we have encountered since the beginning of this section does not always hold. A first instance is when competition is present: the simplest instance concerns models of the form 
$$
n_t=R(x,I(t))n,\ \ \ I(t)=\int_\RR K(x)n(t,x)~dx.
$$
This model represents the evolution of the density $n(t,x)$ of a living species at time $t$ and with trait $x$. Individuals with various traits compete for a resource and this is represented by the term $I(t)$. The competition $R$ is positive in a limited range of $x$ , and negative everywhere else. This seemingly innocent model give rise to large time concentration phenomena, see Perthame \cite{Pbook}. Of course much richer models exist. 

\noindent The study of the steady states owes much to ideas pertaining to nonlinear elliptic equations. For concave (or concave-like) nonlinearities, the idea of comparing a solution to a multiple of the other is due to Berestycki \cite{HB}, and can be seen as a generalisation of (parts of )the Krein-Rutman Theorem. In unbounded domains, the role of uniform bounds from below is first noticed, and analysed, in Beretsycki-Hamel-Roques \cite{BHR}, still in the context of elliptic equations. Their interesting argument, that we borrow for our cause, complete, in a nontrivial way, a general argument of \cite{HB} for elliptic equations with concave-like nonlinearities.

\noindent That an invasion occurs in models of this type was probably first understood by Kendall in his extraordinary text \cite{Kend}.  There, he not only introduces the SI model with nonlocal contaminations, but also proves the invasion. He calls this property the {\it Pandemic threshold Theorem}, that is, for a basic reproduction number strictly larger than 1, the final size of the susceptible population is diminished by a factor that is nontriavially less than 1, uniformly in every bounded region, and irrespective of the initial size of the infected population. Theorem \ref{t2.2.251} is known, in the context of parabolic equations, as the {\it hair trigger effect}. This terminology was introduced in the pioneering paper of Aronson-Weinberger \cite{AW}, where the phenomenon is investigated in several space dimensions.

\noindent There are of course many variants to the Fisher-KPP equation. One of them aims at modelling the evolution of a species within which competition occurs. One way to do it is to introduce a nonlocality in the absorption term \cite{BNPR}, which leads to 
\begin{equation}
\label{e2.8.4}
u_t-u_{xx}=u(1-\phi*u),
\end{equation}
where $\phi$ is a convolution kernel similar to our favourite $K$. The convolution with $\phi$ accounts for the fact that saturation at a given point  may occur from incursion of individuals from the neighbouring sites. This seemingly benign modification causes the maximum principle  to fail and, therefore, the model displays effects different from those described so far. Still, the arguments that we have laid out so far allow a good grasp of this sort of questions, and Problems \ref{P2.103}-\ref{P2.105} show how we can deal with the Cauchy Problem for them. Simply, one has to work  more. 
\subsection*{Asymptotics for the heat kernel (Section \ref{s2.4})}
\noindent The mathematical treatment of  linear equations of the type  
$v_t+\II_*v=0
$ {\it via} the Fourier transform is not new, see for  instance in Bates-Chen \cite{BC} for equations of the form $u_t+\mathcal{J}u=0$. What is rediscovered is in fact the Central Limit Theorem for which, I believe, it is not necessary to give a precise reference. Its ingredient is precisely that the symbol of the involved diffusion behaves like $\xi^2$ in the vicinity of 0. When this is not so, the heat kernel has a different behaviour, see for instance Chasseigne et al. \cite{CCR}, Garofalo \cite{Garof} or, much less recently, Polya \cite{Pol} (who would have believed it? The issue arises from a question in number theory).  What seems more original is the special form of the operator $\II_*$ that has occupied us in the best part of the chapter, as well as the lengthy study of what happens for a well-spread initial datum. This is perhaps what makes the most important difference between the nonlocal model \eqref{e2.1.1} and the standard Fisher-KPP equation:
in the former case, no instant regularisation occurs, which forces this careful study in self-similar variables at intermediate times.

\noindent Heat kernel  estimates outside the standard diffusive context  is an important question for which I have chosen not to give an extensive bibliography, that would be outside the scope of the book. Sharp estimates, such as those presented here, are not so often available. This is why I mention the beautiful work of Coulombel-Faye \cite{CouFa}, which estimates convolution powers of operators that are (but not only) discretisations of diffusive operators. They retrieve and extend previous results that were initially proved with probabilistic tools. 
\subsection*{Open questions}
\noindent We have just seen that it is possible to make inroads in the questions addressed in this chapter when 
comparison is available. They turn out to be much more difficult as soon as it fails. Consider, for instance, the system with two unknowns $\bigl(u(t,x),v(t,x)\bigl)$:
\begin{equation}
\label{e2.7.800}
u_t+\mathcal{J}_{K_1}u=vf(u),\ \ 
v_t+\mathcal{J}_{K_2}v=-vf(u)
\end{equation}
The assumptions and notations are the following. The diffusion operators $\mathcal{J}_{K_i}$ write, as usual, 
$
\mathcal{J}_{K_i}u(x)=-<\!K_i\!>u(x)+K_i*u(x),
$
where $K_i$ are smooth kernels. They can also be singular at $0$ which may, in this case, render the problem less difficult because there would probably be regularisation. In any case they have finite second moments and can be assumed to be compactly supported. One may suppose one of them to be the Laplacian. The function $f$ is smooth, nonnegative. One may also, if it helps, assume that it is globally Lipschitz, but I do not believe that it is drastically simplifies  the question under study.

\noindent  Start from a nonnegative compactly supported $u(0,.)$ and, say, $v(0,.)\equiv1$. The issue, which at first sight may, deceptively,  look like a $X\!X^{th}$ century question on a toy model, is global existence of a solution to \eqref{e2.7.800}, and, more crucially, boundeness of $u$. One may indeed observe that $v(t,x)$ is trivially bounded by 1. One may also observe that the hypothesis $K_1\equiv K_2$ transforms the system in  to an equation that looks pretty much like \eqref{e2.2.1}. When the diffusions are both Laplacians,  a recent development is the preprint of La, Ryzhik and the author \cite{LRR}; while it goes one step further in the study of this question it certainly does not close it. The rich story of this question is depicted in \cite{LRR}.

\noindent  Any nontrivial element of answer would have implications on, at least, models in epidemiology. As explained in Problem \ref{P2.102}, there is no reason why the population of susceptibles is forbidden to move. A relevant system would resemble something of the following sort: suppose that the susceptibles diffuse like $-d_S\partial_{xx}$, whereas the infected, less prone to displacements, would 
diffuse like $-d_I\partial_{xx}$, with $d_I<d_S$, or even not at all ($d_I=0$). On the other hand, the infected can infect individuals around them, so that the resulting system is
\begin{equation}
\label{e2.7.801}
S_t-d_S\partial_{xx}S=-\beta(K*I)S,\ \ 
I_t-d_I\partial_{xx}I=\beta(K*I)S-\alpha I.
\end{equation} 
Introducing these additional diffusions would certainly, besides transforming a well-mastered system into something much more open, induce important new effects. 

\medskip
\noindent A question of a more academic type, but certainly not less interesting, is a feature which seems to me intimately linked to second order one-dimensional parabolic equation, and whose extension to nonlocal models seems problematic to me. The nonincrease of the lap number, or, in other words, the fact that the zero set of the solution of a linear parabolic equation of the form
$$
v_t-v_{xx}+a(t,x)v=0\quad (t>0,x\in(a,b))
$$
is, provided that $u(t,a)$ and $u(t,b)$ are never zero: (i) finite as soon as $t>0$, (ii) nonincreeasing in number. The pioneering author is Matano \cite{Mat}; this principle has been generalised and used a lot of times in one-dimensional equations. One can in fact trace it back to the KPP paper, which may be thought of being at the beginning of this line of research (among many other things). Does a principle of this form hold, and what is a correct and sound formulation is, is something I would be quite happy to know.
 \section{Problems}\label{s2.107}
\begin{problem}\label{P2.104}
How does the strong maximum principle get modified if the support of $K$ does not comprise the point $x=0$?  What would, by the way, be an interpretation of this situation? How can the conclusions of Proposition \ref{p2.2.3} be improved if $K>0$ on $\RR$?\end{problem}
\begin{problem} \label{P2.112} 
 Consider two continuous functions $u_l(x)$ and $u_r(x)$, and the Dirichlet problem
\begin{equation}
\label{e2.7.1560}
\left\{
\begin{array}{rll}
u_t+\JJ u=&f(u)\quad (t>0,\ a\leq x\leq b)\\
u(t,x)=&u_l(x)\quad(t>0,\ a-1\leq x<a)\\
u(t,x)=&u_r(x)\quad(t>0,\ b<x\leq b+1)
\end{array}
\right.
\end{equation}
Show that \eqref{e2.7.1560}, supplemented with your favourite Lipschitz initial datum, has a unique solution, such that $\partial_xu(t,.)$ is bounded in every set of the form $[0,T]\times[a,b]$. Find examples of functions $u_l$ and $u_r$ such that $u(t,x)$ is discontinuous at $x=a$ or $x=b$. Conversely, find sufficient conditions for continuity. How far the study can be pushed is not entirely clear to me.
\end{problem}
\begin{problem}\label{P2.110}
Are the exponentials the only solutions to \eqref{e2.2.5}, for $c\geq c_*$?
\end{problem}
\begin{problem}\label{P2.115}
Assume $f$ to be a Fisher-KPP type term, with, if needed, that $f'(u)\leq f'(0)$. Consider the equation, posed on the interval $(-a,b)$:
\begin{equation}
\label{e2.7.10}
\JJ u=f(u),\quad -a\leq x\leq b.
\end{equation}
The issue is to understand it with various values of $a$ and $b$, and various sets of boundary conditions.
\begin{itemize}
\item[---] Intuition may be borrowed from standard diffusion equation $-u''=f(u)$. So, assume for one moment that $\JJ$ has been replaced by the standard diffusion $-\partial_{xx}$. Thus, Problem \eqref{e2.7.10} is a boundary value problem for an ODE, whose solutions can be computed (semi)-explicitely. Examine what happens in the following cases: (i) $a=0$, $b=+\infty$ and $u(0)=0$, $u(+\infty)=1$; (ii). $u(a)=0$ and $u(b)=1$; (iii). $b=a$ and $u(\pm a)=0$.
\item[---] Assume $a=0$ and $b=+\infty$. Under the conditions $u(x)=0$ for $-1\leq x<0$, and $u(+\infty)=1$, prove the existence of a unique solution $u_0$. Show that it is {\rm discontinuous} at $x=0$ and that it is increasing in $x$. Inspiration can be taken from the sliding method of Berestycki-Nirenberg \cite{BN-shear}, up to the fact that some care should be given to the boundary conditions. Some clues may also be found in the PhD thesis of Coville \cite{Cov1}.

\noindent If you like elliptic equations, do you see an analogy with the Hopf Lemma?
\item[---] Assume $a>0$. Under the boundary conditions $u(x)=1$ if $-a-1\leq x<a$, and $u(x)=0$ if $b<x\leq b+1$, show the existence of a unique solution $u_{a,b}$ to \eqref{e2.7.10}. Show that it is decreasing in $x$, and {\rm discontinuous} at $x=\pm a$.  If $b-a>0$ is large enough, let $x_{a,b}$ be the unique $x$ such that $u_{a,b}(x)=1/2$. Show that $(x_{a,b})_{a,b}$ is bounded, and study the limit of $u_{a,b}$ as $b-a\to+\infty$.

\item[---] Assume $a>0$ and $b=a$. Under the boundary condition  $u(x)=0$ if $a<\vert x\vert<a+1$, show the existence of $a_*$ such that \eqref{e2.7.10} has no nontrivial solution if $a\leq a_*$, and a unique positive solution $v_a$ if $a>a_*$. This of course has to do with the stability of the zero solution with respect to the evolution problem in $(-a,a)$. For the uniqueness of $v_a$, one should try to adapt the argument of Berestycki \cite{HB} and examine how the treatment of the boundary condition gets modified.
\end{itemize}
\end{problem}
\begin{problem}\label{P2.170}
\noindent As explained in Chapter \ref{Intro}, the book will not devote any special part to the reaction-diffusion equation with standard diffusion 
\begin{equation}
\label{e2.7.20}
u_t-\partial_{xx}u=f(u),\quad t>0,x\in\RR,
\end{equation}
as the investigations for this model are facilitated by a number of features. One thing, however, requires more care, and it is the Cauchy Problem. It is quite a classical subject, already treated in great detail in the KPP paper \cite{KPP}. For less specific nonlinearities, one may consult the book of Henry \cite{Hy}. The reader interested in the fascinating parabolic machinery is invited to consult the monumental masterpiece of Ladyzhenskaya-Solonnikov-Ural'tseva \cite{LSU}. I quote it here not only because it is appropriate, for someone who has spent time at the {\rm alma mater} of the authors, to quote them. It is also appropriate because essentially everything one needs to know about parabolic equations is there.

\noindent In the sequence of problems ranging from this one (which sets the framework) to Problem \ref{P2.174}, the reader is taken, in a pedestrian way, into the proof of the following theorem:
\begin{theorem}
\label{t2.7.1}
{\rm Consider a Cauchy datum $u_0(x)$ for \eqref{e2.7.20}, which is continuous and with values in $[0,1]$. Then, if $f$ is smooth, the Cauchy Problem for \eqref{e2.7.20} has a unique solution $u(t,x)$ which is, in addition, in $C^\infty(\RR_+^*\times\RR)$,
and which satisfies $0<u(t,x)<1$.}
\end{theorem}
Most of the effort will be devoted to 
\begin{equation}
\label{e2.7.21}
v_t-\partial_{xx}v=g(t,x),\quad t>0,x\in\RR,
\end{equation}
where $g$ is a bounded function of $\RR_+\times\RR$. Let $v_0$ be a bounded Cauchy datum for \eqref{e2.7.21}. For short we will sometimes denote, for a function $v(t,x)$, by $v(t)$ the function $v(t,.)$. Consider $t>0$, and $v(t)$ a solution of \eqref{e2.7.21} which is once differentiable in $t$ and twice differentiable in $x$ over $\RR_+^*\times\RR$. Fix $t>0$, by differentiating the function $s\mapsto e^{(t-s)\partial_{xx}}v(s)$, show te Duhamel formula:
\begin{equation}
\label{e2.7.22}
v(t)=e^{t\partial_{xx}}v_0+\int_0^te^{(t-s)\partial_{xx}}g(s)~ds.
\end{equation}
We call the function $v(t,x)$ a {\rm weak solution} (other authors would call it mild, in order to distinguish it from distributional solutions). The main issue is to show that it is indeed an honest classical solution to \eqref{e2.7.21}, with all the prescribed derivatives needed.
\end{problem}
\begin{problem}\label{P2.171} ($C^{1,1/2}$ estimates for the weak solutions of \eqref{e2.7.21}) Consider $v(t)$ given by \eqref{e2.7.22}, with $g$ merely bounded. For simplicity (the problem is of course not here), assume $v_0=0$.
\begin{itemize}
\item[---] Show, by direct differentiation, that $\partial_xv$ is bounded on every set of the form $[0,T]\times\RR$.
\item[---] If $w_0(x)$ is bounded, show that, for every $\delta\in(0,1)$, we have $\Vert\partial_x\bigl(e^{\delta\partial_{xx}}w_0\bigl)\Vert_\infty\lesssim \Vert w_0\Vert_\infty/\sqrt\delta$.
\item[---] Fix $t>0$ and $t'\in[t,t+1]$. Show that $v(t)-v(t')=\di\int_0^t\bigl(I-e^{(t'-t)\partial_{xx}}\bigl)e^{(t-s)\partial_{xx}}g(s)~ds+O(t'-t)$. Deduce from the above that $\vert v(t)-v(t')\vert\lesssim \sqrt{t'-t}$.
\item [---] Conclude that, for all $T>0$ and for all $(t,t',x,x')\in[0,T]^2\times\RR^2$, we have 
\begin{equation}
\label{e2.7.23}
\vert v(t,x)-v(t',x')\vert\lesssim \bigl(\sqrt{\vert t-t'\vert}+\vert x-x'\vert\bigl).
\end{equation}
\end{itemize}
\end{problem}
\begin{problem}\label{P2.172} ($C^{2,1}$ differentiability). Fix $t>0$. For $\e\in(0,t)$, let $v_\e(t)=\di\int_0^{t-\e}e^{(t-s)\partial_{xx}}g(s)~ds.$ Let $G(t,z)=e^{-z^2/4t}/\sqrt{4\pi t}$ be the 
Gaussian solving $\partial G-\partial_{xx}G=0$, with $G(0,.)=\delta_{x=0}$.
\begin{itemize}
\item[---] Show that $\partial_{xx}v_\e(t,x)=\di\int_0^{t-\e}\int_\RR\partial_{zz}\bigl(G(t-s,x-y)\bigl)g(s,y)~dyds.$
\item[---] Deduce that $\partial_{xx}v_\e(t,x)=\di\int_0^{t-\e}\int_\RR\partial_{zz}\bigl(G(t-s,x-y)\bigl)\bigl(g(s,y)-g(t,x)\bigl)~dyds.$ {\rm Hint:} the integral of the derivative of a $C^1$ integrable function on $\RR$ is zero.
\item[---] Assume now that $g$ satisfies an inequality of the form \eqref{e2.7.23}. Show that $(\partial_{xx}v_\e)_\e$ converges uniformly, as $\e\to0$, on every set of the form $[\delta,T]\times\RR$. Deduce that 
$v$ is twice differentiable in $x$ and once differentiable in $t$.
\item[---] (Not needed for the sequel, this involves cutting the integral into various pieces, or to consult \cite{LSU}) Show that $\partial_tv$ and $\partial_{xx}v$ satisfy an inequality of the type \eqref{e2.7.23}.
\end{itemize}
\end{problem}
\begin{problem}\label{P2.173} (Weak and strong maximum principle, in a nutshell). Consider a solution $v(t,x)$ of the linear equation
\begin{equation}
\label{e2.7.24}
v_t-\partial_{xx}v=a(t,x)v\quad{t>0,x\in\RR},
\end{equation}
where $a(t,x)$ is bounded continuous on $[0,T]\times\RR$. We assume $v(0,x)\leq0$ for all $x\in\RR$.
\begin{itemize}
\item[---] Assume first that $a(t,x)\geq0$. From the Duhamel formula (with $g(t,x)=a(t,x)v(t,x)$), prove an inequality of the form 
$\Vert v^+(t)\Vert_\infty\leq C_T\di\int_0^t\Vert v^+(s)\Vert_\infty~ds$, and deduce from this that $v(t,x)\leq0$ for $(t,x)\in[0,T]\times\RR$.
\item[---] If $a$ is not assumed to be nonnegative anymore, show that the result still  holds by considering $w(t)=e^{\Lambda t}v(t)$, $\Lambda$ well chosen.
\item[---] If $a(t,x)$ is assumed again to be nonnegative, and $v(0)$ not identically zero, 
prove that $v(t,x)<0$ by noticing $v(t)\leq e^{t\partial_{xx}}v(0)$. Deduce the result when $a$ is not assumed nonnegative anymore.

\noindent I am grateful to B. Mallein for pointing out this argument, that bypasses the traditional (and more elaborate) strong maximum principle, at least in this simple setting.
\end{itemize}
\end{problem}
\begin{problem}\label{P2.174} (Theorem \ref{t2.7.1}, at last) Prove the result by a fixed point argument, complete the smoothness argument by induction, and show that Problem \ref{P2.173} implies $0<u(t,x)<1$.
\end{problem}

\begin{problem}\label{P2.103}
Study the equation $D_c(\lambda)=0$, with unkn0wn $\lambda\in\CC$, in the full range $c\in[0,c_*)$.
\end{problem}
\begin{problem}
\label{P2.111}
Figure out a sharp asymptotics of $e^{-t\II_*}v_0-G_*(t,.)*v_0$ as $t\to+\infty$.
\end{problem}
\begin{problem}\label{P2.144}
Consider a $C^1$ function $f$ such that $f(0)=f'(0)=0$ and $f(1)=0$. Are there sufficient conditions on the size of $f$ near 0 so that the hair trigger effect holds?  Aronson-Weinberger \cite{AW} discuss it in detail when the diffusion is given by the Laplacian.
\end{problem}
\begin{problem} \label{P2.133}  (The Berestycki-Nadin-Perthame-Ryzhik \cite{BNPR} model for competing species).  Consider Model \eqref{e2.8.4} and, {\rm Une fois n'est pas coutume}, we start with an honest diffusion given by the Laplacian. Let $u_0$ be our favourite initial datum, and $u(t,x)$ the solution, that is for the moment local in time.
\begin{itemize}
\item[---] Show  that $u(t,x)\geq0$, then that $u(t,x)\leq e^te^{t\partial_{xx}}u_0(x)$. Conclude that $u(t,x)$ is global.
\item[---] Assume that, at some time $t_0$, $u$ reaches a global maximum $M>1$ at $x_0$. Show that the measure of $[x_0-1/2,x_0+1/2]\cap\{u(t_0,.)\geq\sqrt M\}$ does not exceed $1/\sqrt M$.
\item[---] Using $u(t,x)\leq e^{t-t_0}e^{(t-t_0)\partial_{xx}}u(t_0,.)$, show that
$
u(t_0+1,x)\lesssim \sqrt M+M^{-1/2}\di\sum_{n\in\mathbb{Z}}e^{-n^2/4}.
$
\item[---] Conclude that $u(t,.)$ is bounded and that $\di\liminf_{t\to+\infty}u(t,x)>0$, uniformly on compact sets.
\end{itemize}
This result is due to Hamel-Ryzhik \cite{HRy}. The main issue here is the global bound, and the path proposed here is essentially that followed in \cite{HRy}. Why is  the last statement not formulated more precisely?
\end{problem}
\begin{problem} \label{P2.154} In \eqref{e2.8.4}, there is no real reason why the diffusion should be represented by a Laplacian, let us see what happens for the diffusion $\JJ$. Consider, therefore, a convolution kernel $K$ satisfying the usual assumptions, and a competition kernel $\phi$. Setting 
$\JJ u=u-K*u$, we consider the system
\begin{equation}
\label{e2.8.5}
u_t+\JJ u=u(1-\phi*u),
\end{equation}
still with the initial datum $u_0$.   Assume that there is $\beta>0$ for which $K(x)\leq\beta\phi(x)$. Using the fact that $u(t,x)$ solves $u_t\leq\phi*u(\beta-u)$, show the existence of $M\geq\beta$ such that $u(t,x)\leq M$. Show an analogue of the hair trigger effect.
\end{problem}
\begin{problem} \label{P2.105} When the support of $K$ is larger than that of $\phi$, the above argument that shows the global boundedness of $u$ is not valid anymore. Show that 
$e^{-t\JJ}$ has a heat kernel type behaviour, that is, a behaviour similar to that of Theorem \ref{t2.4.2}. Invoke the argument of Problem \ref{P2.103} and cross your fingers that it is sufficient.
\end{problem}
\begin{problem}\label{P2.108}
Consider Kendall's model for the spatial spread of epidemics:
\begin{equation}
\label{e2.7.1500}
\partial_tI+\alpha I=SK*I,\ \
\partial_tS=-SK*I\quad(t>0,\ x\in\RR).
\end{equation}
The quantity $\beta>0$ is still the mass of $K$.  We assume a uniform initial susceptible density: $S(0,x,y)\equiv S_0>0$, and that  $I(0,x)=I_0(x)$ compactly supported. 
\begin{enumerate} 
\item[---] Show that the cumulative numbers of infected individuals $u(t,x)$ solves 
\begin{equation}
\label{e2.7.1501}
\partial_t u+\alpha u=S_0(1-e^{-K*u(t,.)})+I_0(x).
\end{equation}
\item[---] Show that \eqref{e2.7.1501}, starting from the initial datum $u(0,x)\equiv0$, has a unique solution $u(t,x)$. 
\item [---] For $I_0\in\RR_+^*$, let $u_*(I_0)$ be the unique solution of 
$\alpha u=S_0(1-e^{-\beta u})+I_0.
$.
Let $u_*$ denote, when it exists, the unique positive solution  with $I_0=0$.

\noindent Show that $0\leq u(t,x)\leq u_*(\Vert I_0\Vert_\infty)$.
\item[---] Show that the steady problem
\begin{equation}
\label{e2.7.1512}
\alpha u=S_0(1-e^{-K*u(t,.)})+I_0(x)
\end{equation}
has a unique solution $u_\infty[I_0]$. 
\item[---] Introduce $R_0=S_0\beta/\alpha$. If $R_0<1$, show that $u_\infty[I_0]$ tends to 0 exponentially fast. If $R_0>1$, show that $u_\infty[I_0]$ tends to $u_*$, also with an exponential rate. In both cases, estimate the rate with the data of the problem.
\item[---] Show Kendall's {\rm Pandemic thereshold theorem:} if $R_0>1$ we have $\di\liminf_{t\to+\infty}u(t,x)\geq u_*$, uniformly on every compact subset of $\RR$.
\end{enumerate}
\end{problem}
\begin{problem} \label{P2.102}  We consider a situation when infected and susceptible individuals move. There is, however, no reason why the diffusions should be the same. We propose the model
$$ \left\{
 \begin{array}{rll}
 \partial_tS-d\partial_{xx}S=&-\beta SI\\
 \partial_tI+I-K*I=&\beta SI-\gamma I\\
 S(0)=&S_0(x)>0,\ I(0)=I_0(x)\geq0,\ \hbox{small and compactly supported}.
 \end{array}
 \right.
$$
Show the existence of a global solution $\bigl(S(t,x),I(t,x)\bigl)$. Of course there is no reason to be limited to these to sorts of diffusions, you are welcome to test other sets of diffusions.  Notice that here, it does not seem possible to use the cumulative densities.

\noindent  Is $I(t,x)$ uniformly globally bounded?  I am not so sure to know the answer to this question.
\end{problem}
\begin{problem}
\label{P2.50}
This is the first problem of a series of five (together with Problems \ref{P2.51} and \ref{P2.52} below, Problems \ref{P3.52} in Chapter \ref{TW}), and Problem \ref{P4.50} in Chapter \ref{short_range}) aimed at elucidating the sharp asymptotic behaviour of the solutions to the time-dependent version of the basic model \eqref{e1.1.1}, that is
\begin{equation}
\label{e2.7.50}
u_t-K*u+u=f(t,u)\quad(t>0,x\in\RR)
\end{equation}
where the kernel is, for commodity, assumed to have unit mass. The function $f$ is as in Problem \ref{P1.7}, that is, smooth, 1-periodic in $t$, and such that 
$f(t,0)=0$, $f_u(t,0)>0$, $f(t,u)\leq-\alpha u$ for $u$ large. For commodity, $f$ will be assumed to be of the KPP type, in the strongest form, that is $u\mapsto f_u(t,u)$ 
strictly decreasing. This is not a limitation, as the ideas developped in the monograph would carry over to ZFK type terms. 
Set, in any case:
\begin{equation}
\label{e2.7.52}
m=\di\int_0^1f_u(t,0)~dt,\quad D_c(\lambda)=2\di\int_0^1\bigl(\cosh(\lambda x)-1\bigl)~dx-c\lambda+m.
\end{equation}
This problem considers the linear waves, that is, solutions $\psi(t,x)$ to 
\begin{equation}
\label{e2.7.51}
\psi_t-K*\psi+\psi=f_u(t,0)\psi,
\end{equation}
sought for under the form $\psi(t,x)=a(t)e^{-\lambda(x-ct)}$. Define $c_*$ as the lowest $c>0$ such that the equation $D_c(\lambda)=0$ has positive solutions. 

\noindent Show that linear waves under the prescribed form exist if and only if $c\geq c_K$. Are all solutions to \eqref{e2.7.51} that are time periodic in a reference frame of speed $c$ linear waves?
\end{problem}
\begin{problem}
\label{P2.51}
The goal of this problem is to prove a hair trigger, or pandemic threshold theorem, as the reader prefers, for Problem \eqref{e2.7.50}. Let $u^+(t)$ be the minimal positive periodic solution to $\dot u=f(t,u)$.
\begin{itemize}
\item[---] Pick $\underline m\in(0,m)$ and let $\psi(x)=e^{-\lambda x}\cos(\omega x)$ be a (complex) linear 
wave with zero speed with $f_u(t,0)$ replaced by $\underline m$, with $\omega>1$. Set $\underline\psi(x)=\e\psi(x)\un_{-\pi/2\omega,\pi/2\omega}(x)$, assume $\e>0$ to be small enough to ensure $\underline\psi<u^+$. 
\begin{itemize}
\item[--] Show that one can restrict $\e$ even further to make $\underline\psi$ a subsolution to \eqref{e2.7.50}.
\item[--] Let $\underline u(t,x)$ solve \eqref{e2.7.50} with initial datum $\underline\psi$. Adapt Theorem \ref{t2.2.251} to show that $\di\lim_{t\to+\infty}\bigl(\underline u(t,x)-u^\infty(t,x)\bigl)=0$, uniformly on every compact in $x$, where $u^\infty(t,x)$ is a solution of \eqref{e2.7.50}
that is 1-periodic in $t$.
\end{itemize}
\item[---] Consider a positive (but not necessarily bounded away from 0) time-periodic solution $u(t,x)\leq u^+(t)$ of \eqref{e2.7.50} 
\begin{itemize}
\item[--] Show that $\di\int_0^1f(t,u)/u~dt\leq 1$, and hence that $\di\int_0^1f_u(t,u)~dt<1$.
\item[--] Writing the equation for $u$ under the form $u_t+u-f(t,u)=K*u$, deduce that $x\mapsto u(t,x)$ is continuous (this is really a stability argument). 
\item[--] Show that $u$ is smooth in $t$ and $x$.
\item[--] Conclude that $u(t,x)\equiv u^+(t)$.
\end{itemize}
\item[---] Let $u_0\leq\di\min_{t\in[0,1]}u^+(t)$ be an initial datum for \eqref{e2.7.50},and $u(t,x)$ the corresponding solution. Show that $\di\lim_{t\to+\infty}\bigl(\underline u(t,x)-u^+(t)\bigl)=0$, uniformly on every compact in $x$, 
\end{itemize}
\end{problem}
\begin{problem}
\label{P2.52}
Let the operator $\II_*$ be defined by \eqref{e2.4.6}, and  consider the linear Cauchy Problem
\begin{equation}
\label{e2.7.54}
v_t+\II_*v=\bigl(f_u(t,0)-m\bigl)v,\ v(0,x)=v_0(x) \hbox{compactly supported}
\end{equation}
so that $w(t,x)=e^{\lambda_*x}v(t,x)$ solves $w_t+w-K*w-c_*\partial_xw=f_u(t,0)w$ with $w(0,x)=e^{\lambda_*x}v(t,x)$.

\noindent Show that the analysis and conclusions of of Theorems \ref{t2.4.2} and \ref{t2.4.4} carry over to \eqref{e2.7.54}, up to a harmless factor $\mathrm{exp}\bigl(\di\int_0^tf_u(s,0)~ds\bigl).$
\end{problem}
\begin{problem}
\label{P2.120}
\noindent  (The Berestycki-Chapuisat model  \cite{BCh} with nonlocal longitudinal diffusion). The following three problems set the stage for the study of this model, which is subsequently developped in Chapter \ref{ZFK_short_range}. Let us consider our favourite kernel $K$, with mass 1, and consider the issue of the spreading of a population structures by a phenotypic trait. The model that is proposed below is  variant of Berestycki-Chapuisat \cite{BCh}; the density $u(t,x,y)$ is assumed to solve the equation
\begin{equation}
\label{e2.7.15}
\partial_tu+u-K*u-d\partial_{yy}u+\alpha y^2u=f(u).
\end{equation}
In \cite{BCh}, the diffusion operator is the Laplacian in both variables $x$ and $y$. The parameter $\alpha$ is positive, and the function $f$ is the usual honest KPP or ZFK type term. The present problem studies the linear waves for \eqref{e2.7.15}, the subsequent one presents an analogue of the hair trigger effect, or, f the reader prefers, the pandemic threshold theorem.
\begin{itemize}
\item[---] Show that the eigenvalue problem 
$-\psi''+y^2\psi=\mu\psi,
$ posed, for instance, in $L^2(\RR)$,
has an increasing sequence $(\mu_n)_{n\in\NN}$, with $\mu_{n+1}=\mu_n+1/2$. Compute $\mu_1$.

\noindent {\rm Hint:} as it is a classical question, consult Reed-Simon \cite{RS}, or or the last chapters of Brezis \cite{Bre}. Alternatively (this is more instructive) do it with bare hands: write down the solution of the heat equation $u_\tau-u_{\zeta\zeta}$ in the self-similar variables $t=\mathrm{ln}~\tau$, $y=\zeta/\sqrt\tau$, then set $u(t,y)=e^{\lambda t-y^2/8}v(t,y)$, for a well-chosen $\lambda$.
\item[---] Deduce, from a suitable scaling  in $y$, that  the eigenvalue problem 
$-\psi''+\alpha y^2\psi=\mu\psi,
$ posed in $L^2(\RR)$
has an increasing sequence of eigenvalues $\bigl(\mu_n(\alpha)\bigl)_n$. 
\item [---] A linear wave to \eqref{e2.7.15} is a solution of the linear equation $\partial_tu+u-K*u-d\partial_{yy}u+\bigl(\alpha y^2u-f'(0)\bigl)u=0$, of the form $e^{-\lambda(x-ct)}\psi(y)$. Show that such a solution is a positive linear wave if and only if $(c,\lambda)$ satisfies the equation
\begin{equation}
\label{e2.7.16}
2\int_0^1K(x)\bigl(\cosh(\lambda x)-1\bigl)~dx+f'(0)-\mu_1(\alpha)=c\lambda.
\end{equation}
\item[---] Study all the linear waves, irrespective of their possible positivity.
\item[---] Let $\bar\alpha$ the only $\alpha>0$ such that $\mu_1(\alpha)=f'(0)$. If $\alpha<\bar\alpha$, show the existence of a bottom velocity $c_*(\alpha)>0$ for \eqref{e2.7.16}.
\item[---] If $\JJ$ is replaced by $-\partial_{xx}$, show that $c_*=2\sqrt{f'(0)-\mu_1(\alpha)}$.
\end{itemize}
\end{problem}
\begin{problem}
\label{P2.122}
Assume $f'(0)>\mu_1(\alpha)$ and denote, for short, $\mu_1,\ldots,\mu_n,\ldots$ the eigenvalues of $-\partial_{yy}+\alpha y^2$. Let $\II_*$ have the notation \eqref{e2.4.6}. Consider the linear equation
\begin{equation}
\label{e2.7.30}
v_t+\II_*v-\partial_{yy}v+\alpha y^2v=0\ \ \  \bigl(t>0,(x,y)\in\RR^2\bigl),\quad v(0,x,y)=v_0(x,y)
\end{equation}
with $v_0$ compactly supported in $\RR^2$. Noticing that the operators $\II_*$ and $-\partial_{yy}+\alpha y^2$ commute, show  that the solution $v(t,x,y)$ of \eqref{e2.7.30} behaves as in Theorem \ref{t2.4.2}, that is: there is $d_*>0$ and
a nonnegative kernel $G_*(t,z)$, defined for $t\geq1$ and $z\in\RR$ and bounded on its domain of definition, such that we have, for all $\gamma\in(0,1/6)$, for all $t\geq1$ and $x\in\RR$:
$$
\Vert e^{-t\II_*}v_0-{G_*(t,.)}*_x<\!e_1^*,v_0\!>e_1\Vert_{L^\infty(\RR^2)}\lesssim e^{-t^{1-2\gamma}}\Vert v_0\Vert_{H^1(\RR)}.
$$
Moreover, if
$ \delta>6\gamma,$  then
\begin{itemize}
\item[---] If $\vert z\vert\leq t^{1/2+\delta}$, then $G_*(t,z)=\di\bigl(1+O(\frac1{t^{1-\delta}})\bigl)\frac{e^{-z^2/4d_*t}}{\sqrt{4\pi d_*t}}.$
\item[---] If $\vert z\vert\geq t^{1/2+\delta}$, then 
$G_*(t,z)\lesssim e^{-t^\delta}.$
\end{itemize}
Show the analogue of Theorem \ref{t2.4.4} for initial data that are well-spread in $x$.
\end{problem}
\begin{problem}
\label{P2.121} Things are now ready for a proof of the hair trigger effect. Let $\bar\alpha$ be defined as in Problem \ref{P2.120}.
\begin{itemize}
\item[---] Show the existence of a unique solution $v(t,y)$ to the Cauchy Problem for 
\begin{equation}
\label{e2.7.18}
\partial_tv-\partial_{yy}v+\alpha y^2v=0\ (t>0,y\in\RR)\quad v(0,y)=v_0(y).
\end{equation}
Setting $L=-\partial_{yy}+\alpha y^2$, show that $e^{-tL}u_0>0$ if $u_0\geq 0$ and is not identically zero. 
\item[---] Show the existence of a unique solution to the Cauchy Problem \eqref{e2.7.15}, and that the strong comparison principle holds. A good way to tackle  it is to look for a solution of \eqref{e2.7.15} in the space $C\bigl(\RR_+\times\RR_x,C(\RR_y)\bigl)$.
\item[---] Consider $u_{10}(x,y)\leq u_{20}(x,y)$ two Cauchy data for \eqref{e2.7.15}. Show that we have, for the corresponding solutions $u_i(t,x,y)$: $u_1(t,x,y)<u_2(t,x,y)$.
\item[---] Show the existence of a large constant $M>0$ such that, if $\alpha>M$, we have $\di\lim_{t\to+\infty}u(t,x,y)=0$ uniformly on $\RR^2$.
\item[---] Show the existence of $\bar\alpha>0$ such that \eqref{e2.7.15} has a unique nontrivial steady solution $U(y)$ if and only if $\alpha<\bar\alpha$. 
\item[---] If $\alpha\in[0,\bar\alpha)$, show that $u(t,x,y)$ converges, locally uniformly, to $U(y)$.
\item[---] What happens if $\alpha\geq\bar\alpha$?
\item[---] Study the limits $\alpha\to0$ and $\alpha\to\bar\alpha$.
\item [---] Scale $\alpha$ and $f$ in such a way that, when $K$ is the usual approximation of the identity, that is, $K(x)=\e^{-2}\rho(x/\e)$, one retrieves the Berestycki-Chapuisat model in limit $\e\to0$.
\end{itemize}
\end{problem}

\chapter{Travelling waves}\label{TW}
\noindent Travelling waves are special propagating solutions of the nonlinear equation \eqref{e2.1.1}. They play the same role as the linear waves in equation \eqref{e2.2.2}. They have the form $u(t,x)=\varphi(x-ct)$
with the same notion of leftwards or rightwards propagation as for the linear waves. As they are supposed to be transition profiles between the two steady states of the equation, one requires them to connect 0 at one infinity to 1 at the other infinity. So, for the rightwards propagating waves, the problem is
\begin{equation}
\label{e3.7.1}
\left\{
\begin{array}{rll}
\JJ\varphi-c\varphi'=&f(\varphi)\\
\varphi(-\infty)=1,\ \ \varphi(+\infty)=&0.
\end{array}
\right.
\end{equation}
In the whole chapter we will take, by commodity: $\K=1$. The deviation of $f$ to its initial slope will be denoted $g(u)$, that is:
\begin{equation}
\label{e2*.1.11}
f(u)=f'(0)u-g(u).
\end{equation}
 It will be smooth, nonnegative, not identically 0 with
\begin{equation}
\label{e4.0.340}
g(0)=g'(0)=0,\ \ \ \lim_{u\to+\infty}\frac{g(u)}u=+\infty.
\end{equation}
We will see that the results will be strongly influenced by the sign of $g$.

\noindent  As we will need to observe the solutions of \eqref{e2.1.1} in reference frames moving with various speeds, we set, for any $c\geq c_*$:
\begin{equation}
\label{e2.2.6}
\mathcal{N}_cu=\JJ u-c\partial_xu-f(u),
\end{equation}
and we denote by $\LL_c$ the linear operator
\begin{equation}
\label{e2.2.7}
\LL_c\phi=\JJ\phi-c\partial_x\phi-f'(0)\phi=\mathcal{N}_cu-g(u).
\end{equation}
We denote it by $\LL_*$ if $c=c_*$.  

\noindent Let $c_K$ be the Fisher-KPP speed, that is, the least $c>0$ such that a linear wave $\phi_\lambda(x-ct):=e^{-\lambda(x-ct)}$ solves $v_t+v-K*v-f'(0)v=0$. In other words, it is the least $c>0$ such that the equation with unknown $\lambda>0$
\begin{equation}
\label{e3.7.4}
D_c(\lambda):=-c\lambda+f'(0)+2\int_0^{+\infty}K(y)\bigl(\cosh(\lambda y)-1\bigl) dy=0
\end{equation}
has a solution. If $c=c_K$, \eqref{e3.7.4} has a unique solution $\lambda_K$, if $c>c_K$ there are two solutions $\lambda^-(c)<\lambda^+(c)$. For $\lambda\in\RR$ we have
\begin{equation}
\label{e3.7.6}
\phi_\lambda-K*\phi_\lambda-c\phi_\lambda'=-D_c(\lambda)\phi_\lambda.
\end{equation}
An important issue in this chapter will be to decide what the bottom propagation velocity of a travelling wave will be. Sometimes it will be $c_K$, sometimes it will be strictly larger. In any case, we will call $c_*$ the bottom velocity. If it is $c_K$, we will use the notation $c_*$ to refer to both the bottom linear wave speed or travelling wave speed. In the case that the bottom speed is larger than $c_K$, it will stilll be denoted by $c_*$ but the distinction with $c_K$ will be clearly made.

\noindent In the whole chapter, an extensive use will be made of the following corollary of Proposition \ref{p2.2.4}. As it is of independent interest we state it as a theorem.
\begin{theorem}
\label{t3.0.0}
Consider $c>0$ for which there is $\underline\varphi\leq\overline\varphi$ satisfying
$$
\mathcal{J}\underline\varphi-c\underline\varphi'-f(\underline\varphi)\leq0\leq\mathcal{J}\overline\varphi-c\overline\varphi'-f(\overline\varphi);
$$
in other words $\underline\varphi$ (resp. $\overline\varphi$) is a sub-(resp. super-) solution of the equation. There is a solution $\varphi$ of the equation such that $\underline\varphi\leq\varphi\leq\overline\varphi$.
\end{theorem}
\subsection*{Organisation of the chapter}
\noindent The chapter starts with the study of the Fisher-KPP waves; existence, provided by the sub-super solutions argument of Theorem \ref{t3.0.0}, allows for a smooth start. The asymptotic behaviour is rather more involved, the strategy followed here is to give a representation of the solutions {\it via} a Fourier argument. This will allow, when one pushes it to its end, a representation by means of an integral equation opening the door to a detailed study. Once this is done, attention is turned on the ZFK waves, to which three sections are devoted.  Section \ref{s3.102} presents some generalities, from which it appears that the bottom velocity $c_*$ of these waves may, or may not be, the Fisher-KPP speed $c_K$. Section \ref{s3.103} studies the case when $c_*>c_K$, while Section \ref{s3.105} studies those with bottom speed $c_*=c_K$. A discussion on their decay at infinity comes along with these considerations.
\section{The Fisher-KPP travelling waves} \label{s3.101}
Set $c_*:=c_K$.  The main result of the section is summarised in the following
\begin{theorem}
\label{t2.3.1}
For a given $c>0$, equation \eqref{e3.7.1} has a nonnegative solution $\varphi_c(x)$ if and only if $c\geq c_*$. Moreover, $\varphi_c$ is a $C^\infty$ function, with $\varphi'<0$. It is also unique modulo translations.

\noindent Denote by $\varphi_c(x)$  the solution to \eqref{e3.7.1} such that $\varphi(0)=1/2$.  If $c>c_*$,  there exists $k_c>0$ and $\delta>0$ such that
\begin{equation}\label{e2.3.3}
\varphi_c(x)=k_ce^{-\lambda_-(c)x}+O(e^{-(\lambda_-(c)+\delta)x})\ \ \ \hbox{as $x\to+\infty$.}
\end{equation}
If $c=c_*$, there is $k_*\in\RR$ such that
\begin{equation}\label{e2.3.4}
\varphi_{c_*}(x)=(x+k_*)e^{-\lambda_*x}+O(e^{-(\lambda_*+\delta)x})\ \ \ \hbox{as $x\to+\infty$.}
\end{equation}
\end{theorem}
\noindent Consider $\varphi$ a travelling wave solution of \eqref{e2.1.1}, and $u_0$ any nonnegative, compactly supported function that is below $\varphi$. If $u(t,x)$ is the solution of \eqref{e2.1.1} emanating from $u_0$, we have $u(t,x)\leq\varphi(x-ct)$. This proves that $c\geq c_*$, so that  we will only need to worry about speeds that will be larger than $c_*$.  Existence of travelling waves will follow from a rather classical sub/super-solution argument.

\noindent Before delving into the proof of Theorem \ref{t2.3.1}, let us notice that it is easy to study the asymptotic behaviour of the bottom wave when the diffusion is given by the Laplacian. Consider indeed the problem 
$$
-u''-cu'=u-u^2,\quad, u(-\infty)=1,\ u(+\infty)=0.
$$
We have trivially $c_*=2$, and $\lambda_*=1$.  So, in order to prove \eqref{e2.3.4} it suffices to notice that the function $v(x)=e^xu(x)$ solves $-v''+e^{-x}v^2=0$, with $v(x)\sim e^x$ as $x\to-\infty$. Integrating twice from $-\infty$ to $x$ and intgrating by parts yields $v(x)=\di\int_{-\infty}^x(x-y)e^{-y}v^2(y)~dy$. The classical linearisation theorems for ODEs yield that $u(x)=O(e^{-(1-\e)x})$ as $x\to+\infty$, for all $\e>0$, so that the integrals $\di\int_\RR e^{-y}v^2~dy$ 
and $\di\int_\RR ye^{-y}v^2~dy$ converge. This is enough to yield \eqref{e2.3.4}.
\subsection{Existence}\label{s2.1}

\noindent Consider first $c>c_*$, and 
\begin{equation}
\label{e2.3.32}
\bar\varphi(x)=\inf(1,e^{-\lambda_-(c)x}).
\end{equation}
It is a bounded, Lipschitz super-solution to $\mathcal{N}_c\varphi=0$.  Then, for all $\e>0$, and all $\delta>0$, consider
\begin{equation}
\label{e2.3.38}
\underline\varphi(x)=\e\bigl(e^{-\lambda_-(c)x}-e^{-(\lambda_-(c)+\delta)x}\bigl)_+.
\end{equation}
Suppressing the dependence of $\lambda_\pm$ on $c$, we have $\underline\varphi(x)=\e\bigl(e^{-\lambda_-x}-e^{-(\lambda_-+\delta)x}\bigl)$ for $x>0$. We choose
$\delta<\min\bigl(\lambda_+-\lambda_-,\lambda_-\bigl);
$
 we have, on $\RR_+$:
$$
\begin{array}{rll}
\mathcal{N}_c\underline\varphi(x)=&\e D_c(\lambda_-)e^{-(\lambda_-+\delta)x}+O\biggl(\e^2\bigl(e^{-\lambda_-x}-e^{-(\lambda_-+\delta)x}\bigl)^2\biggl)\\
=&\e e^{-(\lambda_-+\delta)x}\biggl((\delta\bigl(D_c'(\lambda_-)+O(\delta)\bigl) +O\bigl(\e e^{-(\lambda_--\delta)x}\bigl)\biggl)\\
\leq&\delta\e e^{-(\lambda_-+\delta)x}\bigl(D_c'(\lambda_-)+O(\delta+\e)\bigl)\ \ \hbox{because of the choice of $\delta$.}
\end{array}
$$
As $D_c'(\lambda_-)<0$, we have $\mathcal{N}_c\underline\varphi(x)\leq0$ on $\RR_+$ as soon as $\delta$ and $\e$ are small enough. Restrict $\e$ further if necessary, and
apply Proposition \ref{p2.2.2} with $a=0$: we have just constructed a subsolution $\underline\varphi$ below a nonincreasing supersolution $\bar \varphi$. This grants the existence of a nonincreasing, nonzero solution $\varphi$ of 
$\mathcal{N}_c\varphi=0
$
on the whole line. Obviously it goes to 0 at $+\infty$; as it is nonincreasing and bounded it goes to a finite limite $l\in(0,1]$ at $-\infty$. Sending $x$ to $-\infty$   shows that
$
\lim_{x\to-\infty}\varphi'(x)=-{f(l)}/c,
$
which entails $l=1$ and ends the construction of a solution to the full problem  \eqref{e3.7.1} for $c>c_*$.

\noindent We could repeat the process for $c=c_*$, but we prefer to use a computation free argument, based on compactness. Consider $\varphi_c$ a solution of \eqref{e3.7.1} such that $\varphi_c(0)=1/2$. As $c\to c_*$, the family is bounded in $C^{1,1}(\RR)$,
so that a subsequence converges to a solution $\varphi_{c_*}$ of \eqref{e3.7.1} with $c=c_*$, that is both nonzero by the normalisation condition, and nonincreasing. The same argument as above shows that $\varphi_{c_*}$
 goes to 0 at $+\infty$, and to 1 at $-\infty$, which painlessly ends the construction of travelling wave solutions. A more involved task is to study the behaviour of any solution of \eqref{e3.7.1} at infinity, this is the goal of the next section.
\subsection{Asymptotic behaviour at infinity and consequences}\label{s2.2.2}

\noindent A few easy facts should first be noticed, one of them being that a solution
$\varphi(x)$ to \eqref{e3.7.1}  converges to 1 exponentially fast as $x\to-\infty$. Indeed, for $x\ll-1$, the equation for $\varphi$ entails 
$$
\JJ(1-\varphi)-c(1-\varphi)'+{\vert f'(1)\vert}(1-\varphi)/2\leq0,
$$
while $1-\varphi$ is a bounded function; one may therefore bound $1-\varphi$ by an exponential that decays as $x\to-\infty$. Let us also notice that all the derivatives of $\varphi$ are bounded. Finally, we observe that $\varphi\in L^1(\RR_+)$. Indeed, if $M$ is a positive number, integration of \eqref{e3.7.1} on $[M,+\infty)$ implies that
$$
\begin{array}{rll}
f'(0)\di\int_M^{+\infty}\varphi-\di\int_M^{+\infty}g(\varphi)
=&c\varphi(M)+\di\int_M^{+\infty}\bigl(\int_{\RR}K(z)(\varphi(x)-\varphi(x-z))~dz\bigl)~dx\\
=&c\varphi(M)+\di\int_M^{+\infty}\varphi(x)~dx-\di\int_\RR\varphi(y)\bigl(\int_{M}^{+\infty}K(x-y)~dx\bigl)~dy\\
=&c\varphi(M)+\di\int_M^{+\infty}\varphi(x)~dx-\di\int_{M-1}^{+\infty}\varphi(y)\bigl(\int_{M-2}^{M_2}K(x-ydx\bigl)~dy\\
=&c\varphi(M)+\di\int_{M-1}^{M+1}\varphi(x)~dx\\
=&O_{M\to+\infty}(1).
\end{array}
$$
By definition of $\varphi$, we have $g(\varphi)=o_{M\to+\infty}(\varphi)$. This implies the finiteness of the integral of $\varphi$. These two properties pass over to the successive derivatives of $\varphi$.

\noindent The main body of the proof of Theorem \ref{t2.3.1} will be, similarly to the proof of the asymptotic behaviour of the waves with standard diffusion, to write down an integral equation for $\varphi$. For this we will rely on a Fourier transform argument. In what follows, the Fourier transform of  an $L^1$ function $u(x)$ will be denoted by 
\begin{equation}
\label{e3.7.15}
\hat u(\xi)=\int_{\RR}e^{-ix.\xi}u(x)~dx;
\end{equation}
if $\hat u\in L^1(\RR)$ we also have
\begin{equation}
\label{e3.7.16}
u(x)=\frac1{2\pi}\int_{\RR}e^{i\xi.x}\hat u(\xi)~d\xi:=\frac{\tilde{u}}{2\pi}.
\end{equation}
Hence, $\tilde u$ is the conjugate Fourier transform of $\hat u$. Definition \eqref{e3.7.15} and reconstruction formula \eqref{e3.7.16} extend, by duality, to tempered distributions. 

\noindent Define
\begin{equation}
\label{e2.3.6}
\DD_c(\xi)=2\int_{\RR_+}K(z)\bigl(1-\cos(z\xi)\bigl)~dz-ic\xi-f'(0).
\end{equation}
The function $\DD_c$ is an entire function of $\CC$, that is nonzero on $\RR$, and whose modulus grows indefinitely on $\RR$ in a linear fashion. On the imaginary axis, we have $\DD_c(i\lambda)=-D_c(\lambda)$; there,
the roots  of  $\DD_c$ are $i\lambda_\pm(c)$ if $c>c_*$, and $i\lambda_*$ if $c=c_*$.

\noindent When the diffusion is given by the Laplacian, we have 
$\DD_c(\xi)=\xi^2-ic\xi-f'(0),
$
and its zeroes are easily studied. Such does not seem to be the case for the function $\DD_c$ given by \eqref{e2.3.6}, and we will need an intermediate lemma.
\begin{lemma}
\label{l2.3.1}
Consider $c\geq c_*$, and a real exponent $\lambda\leq\lambda_-(c)$ (resp. $\lambda\leq\lambda_*$ if $c=c_*$). Consider $k$ distinct real numbers $\xi_1,\ldots \xi_k$, and the function
\begin{equation}
\label{e3.7.1700}
\mathcal{P}(x)=\sum_{j=1}^k\alpha_je^{i\xi_jx},
\end{equation}
where the $\alpha_j$'s are complex numbers. Assume that $u(x)=e^{-\lambda x}\mathcal{P}(x)$ is a solution of the linear equation \eqref{e2.2.5}. 

\noindent Then, if the $\alpha_j$'s are not all zero, the function $\mathcal{P}$ cannot have a constant sign, unless we have $\lambda=\lambda_-(c)$ (resp. $\lambda=\lambda_*$ if $c=c_*$), $k=1$ and $\xi_1=0$.
\end{lemma}
\noindent{\sc Proof.} Assume the contrary, we may always suppose that  $\mathcal{P}(x)\geq0$. The assumption on $u(x)=e^{-\lambda x}\mathcal{P}(x)$ implies that $\mathcal{P}$ solves the following integral equation:
\begin{equation}
\label{e3.7.18}
\mathcal{J}_\lambda\mathcal{P}-c\mathcal{P}'+\DD_c(i\lambda)\mathcal{P}=0.
\end{equation}
We have denoted
\begin{equation}
\label{e3.7.1900}
\mathcal{J}_\lambda v(x)=\int_{\RR}e^{\lambda(x-y)}K(x-y)\bigl(v(x)-v(y)\bigl)~dy,
\end{equation}
and we recall that $\DD_c(i\lambda)=-D_c(\lambda)$ is a real number. We claim the existence of $\delta>0$ such that $\mathcal{P}(x)\geq\delta$. Assume the contrary: if $\mathcal{P}(x_0)=0$ for some $x_0$ in $\RR$, the strong comparison principle for \eqref{e3.7.18}
implies that $\mathcal{P}$ is identically zero. This means that all $\alpha_j$'s are zero, by the linear independence of the complex exponentials. If, for a sequence $(x_n)_n$ going to $+\infty$ or $-\infty$ we have
$
\di\lim_{n\to+\infty}\mathcal{P}(x_n)=0,
$
we set
$
\mathcal{P}_n(x)=\mathcal{P}(x+x_n)=\di\sum_{j=1}^k\beta_j^ne^{i\xi_jx},$ with  $\beta_j^n=\alpha_je^{i\xi_jx_n}.
$
Consider a subsequence $(n_p)_p$ such that each sequence $(\beta_j^{n_p})_p$ converges to some $\beta_j^\infty$, and set
$
\mathcal{P}_\infty(x)=\di\sum_{j=1}^k\beta_j^\infty e^{i\xi_jx},
$
The function $\mathcal{P}_\infty$ is still a nonnegative solution of \eqref{e3.7.18}, and vanishes at $x=0$, therefore it is identically zero. As $\vert\alpha_j\vert=\vert\beta_j^\infty\vert$, we have $\alpha_j=0$ for all $j$, once again.

\noindent  The function $\underline u(x)=e^{-\lambda x}$ solves
$
\JJ\underline u-c\underline u'-f'(0)\underline u\leq0;$
set 
$\delta_0=\inf \mathcal{P}(x),
$ 
there is a contact point between $u$ and $\delta_0\underline u$, either at finite distance, or at infinity. In the first case, we readily contradict the strong comparison principle for equation \eqref{e2.2.5}; in the second case one just has to consider a sequence
$(x_n)_n$ such that $\bigl(u(x_n)-\delta_0\underline u(x_n)\bigl)_n$, and argue by compactness on the sequence
$$
v_n(x)=e^{\lambda x_n}\bigl(u(x+x_n)-\delta_0\underline u(x+x_n)\bigl).
$$
The limiting function $v(x)\geq0$ has a minimum at $x=0$, solves the equation \eqref{e2.2.5}, therefore contradicts the strong comparison principle. \hfill$\square$

\noindent Before going to the actual proof of Theorem \ref{t2.3.1}, let us state a Gronwall type lemma, of independent interest.
\begin{lemma}
\label{l2.3.2}
Let $\psi(x)$ be a uniformly continuous nonnegative function in $L^\infty(\RR)\cap L^1(\RR_+)$. Assume the existence of $0<\lambda_0<\bar\lambda$ and two nonnegative functions $e_1(x)$ and $e_2(x)$, both in $L^2(\RR)$, such that
\begin{equation}
\label{e3.7.11}
\psi(x)\leq \int_{-\infty}^xe^{-\bar\lambda(x-y)}e_1(x-y)\psi^2(y)~dy+\int_x^{+\infty}e^{\lambda_0(x-y)}e_2(x-y)\psi^2(y)~dy.
\end{equation}
Then we have, for
any $\lambda<\bar\lambda$:
\begin{equation}
\label{e3.7.12}
\psi(x)=O_{x\to+\infty}(e^{-\lambda x}).
\end{equation}
\end{lemma}
\noindent{\sc Proof.} Let us first reduce inequality \eqref{e3.7.11} to a slightly more amenable form, as the presence of the functions $e_1$ and $e_2$, which are only in $L^2(\RR)$, may look bothering. We have, by Cauchy-Schwartz inequality:
\begin{equation}
\label{e2.3.30}
\int_{\RR_-}e^{-\bar\lambda(x-y)}e_1(x-y)\psi^2(y)~dy\leq e^{-\bar\lambda x}\Vert e_1\Vert_{L^2(\RR_-)}\Vert e^{\bar\lambda y}\psi^2\Vert_{L^2(\RR_-)},
\end{equation}
so that \eqref{e3.7.11} becomes
\begin{equation}
\label{e3.7.13}
\psi(x)\leq c_1e^{-\bar\lambda x}+ \int_0^xe^{-\bar\lambda(x-y)}e_1(x-y)\psi^2(y)~dy+\int_x^{+\infty}e^{\lambda_0(x-y)}e_2(x-y)\psi^2(y)~dy.
\end{equation}
Estimate \eqref{e3.7.12} will now result from a bootstrap argument. Consider $\e\in(0,\lambda_0)$, that will be chosen suitably small, and $\lambda=\lambda_0-\e$. The function $\psi$ being uniformly continuous and integrable, 
it goes to 0 as $x\to+\infty$ and we may assume that $\psi(x)\leq\e$ on $\RR_+$. Hence, \eqref{e3.7.13} implies the following weaker inequality:
\begin{equation}
\label{e3.7.14}
\psi(x)\leq c_1e^{-\lambda x}+ \e\int_0^{+\infty}e^{-\lambda_0\vert x-y\vert}e(x-y)\psi(y)~dy,
\end{equation}
with $e(z)=\max\bigl(e_1(z),e_2(z)\bigl)$.  We have
$
\di\int_0^{+\infty}e^{-\lambda_0\vert x-y\vert}e(x-y)e^{-\lambda y}dy\lesssim\frac{\Vert e\Vert_{L^2}}{\sqrt{\lambda_0-\lambda}}=\Vert e\Vert_{L^2}\e^{-1/2},
$
so, we choose $\e$ so that $c_1\Vert e\Vert_{L^2}\sqrt\e<1$; iteration of \eqref{e3.7.14} yields
$
\psi(x)\leq e^{-\lambda x}\di\sum_{n=0}^{+\infty}(c_1\Vert e\Vert_{L^2}\sqrt\e)^n$, that is, an $O(e^{-\lambda x}).
$
Thus, the constant $\e>0$ being as small as we wish, estimate \eqref{e3.7.12} holds for all $\lambda\in(0,\lambda_0)$. Fix such a $\lambda$, and plug this newly found inequality into \eqref{e3.7.13}, we obtain
$
\psi(x)\lesssim e^{-\min(2\lambda,\bar\lambda)x}.
$
This operation can be iterated: for all integer $k$, we obtain
$
\psi(x)\lesssim e^{-\min(2^k\lambda,\bar\lambda)x},
$
with a possible degeneracy if $2^{k_0}\lambda=\bar\lambda$. The exponent $\lambda$ being chosen arbitrarily in $(0,\lambda_0)$, estimate \eqref{e3.7.14} is proved. \hfill$\square$

\noindent{\sc Proof of \eqref{e2.3.3} and \eqref{e2.3.4}.}  The strategy that we are going to follow is to derive a representation formula for the travelling wave $\varphi$, by using the full force of the equation satisfied by 
$\varphi$:
\begin{equation}
\label{e2.3.100}
\LL_c\varphi=g(\varphi),
\end{equation}
the operator $\LL_c$ being given by \eqref{e2.2.7}.
This formula, together with the aid of Lemma \ref{l2.3.2}, will give us the dominant term in the asymptotic expansion of $\varphi$. We will then conclude by using the positivity of $\varphi$.

\noindent The Fourier transform of $\varphi$, denoted by  $\hat\phi(\xi)$, satisfies
$
\DD_c(\xi)\hat\phi(\xi)=-\widehat{g(\varphi)}(\xi).
$
As $g(\varphi)$ does not belong to $L^1(\RR)$, it is better to approximate it by a sequence of $L^1$ functions. For all $\e>0$, let us define
$$
g_\e(x)=\un_{[-1/\e,+\infty)}(x)g(\varphi(x)),
$$
the sequence $(g_\e)_\e$ converges to $g(\varphi)$ in the sense of tempered distribution, so that we have
\begin{equation}
\label{e2.3.7}
\varphi=-\frac1{2\pi}\di\lim_{\delta\to0}\widetilde{\frac1{\DD_c}}*g_\delta.
\end{equation}
the limit being understood in the tempered distributions sense. The function $\di\frac1{\DD_c(\xi)}$, while being defined on the whole real line, is  in $L^2(\RR)$, but not in $L^1(\RR)$. The usual way to write it explicitely
 it is to take the distributional  limit, as $\delta\to0$, of the inverse Fourier transform of $\di\frac{e^{-\delta\xi^2}}{\DD_c(\xi)}\hat{g_\e}$, then to take the limit as $\e\to0$.
So, everything boils down to the following: consider $a(x)\in L^1(\RR)$, we wish to study
$$
\EE_\delta[a](x)=\int_\RR\frac{e^{ix\xi-\delta\xi^2}}{\DD_c(\xi)}\hat a(\xi)~d\xi.
$$
By Fubini's theorem, we have
$$
\EE_\delta[a](x)=\int_\RR a(y)\biggl(\int_\RR\frac{e^{i(x-y)\xi-\delta\xi^2}}{\DD_c(\xi)}d\xi\biggl)~dy.
$$
For any $\lambda>0$, let us define the strip of the complex plane
\begin{equation}
\label{e2.3.9}
\Sigma_\lambda=\{\xi=\xi_1+i\xi_2: \ \vert\xi_2\vert\leq\lambda\}.
\end{equation}
 As $\DD_c$ is nonzero on $\RR$, and as $\vert\DD_c(\xi)\vert$ grows to infinity as $\vert\xi\vert\to\infty$, there is $\lambda_0>0$ such that this property is preserved in $\Sigma_{\lambda_0}$.  The discussion will now follow the sign of $x-y$ in the integral inside the expression of $\EE_\delta[a]$.
 
\noindent Let $\EE^+_\delta[a](x)$ (resp.  $\EE^-_\delta[a](x)$) correspond to the chunk of $\EE_\delta[a](x)$ for which  $x-y\geq0$ (resp. $x-y\leq0$).  In order to compute $\EE_\delta^-$ we move the integration domain from $\RR$ to $\RR-i\lambda_0$, so that we have
 $$
\EE_\delta^-[a](x)=e^{\lambda_0x}\int^{+\infty}_x e^{-\lambda_0y}a(y)\biggl(\int_\RR\frac{e^{i(x-y)\xi-\delta(\xi+i\lambda_0)^2}}{\DD_c(\xi-i\lambda_0)}d\xi\biggl)~dx.
$$
As $e^{-\lambda_0y}a(y)\in L^1([x,+\infty)])$, the above expression makes sense, and we have
\begin{equation}
\label{e2.3.8}
\lim_{\e\to0}\bigl(\lim_{\delta\to0}\EE_\delta^-[g_\e](x\bigl))=e^{\lambda_0x}\int^{+\infty}_x e^{-\lambda_0y}g(\varphi(y))e_-(x-y)~dy:=\EE^-[g(\varphi)](x),
\end{equation}
with
$$
 e_-(z)=\int_\RR\frac{e^{iz\xi}}{\DD_c(\xi-i\lambda_0)}d\xi\in L^2(\RR).
 $$
The real number $\lambda_0$ will be fixed once and for all in the expression of $\EE^-$, as this expression will essentially act as a perturbation. In
order to compute $\EE_\delta^+$ we move the integration domain from $\RR$ to $\RR+i\lambda$,  but this time we will choose $\lambda$ more carefully. We 
have, in a similar way as in the computation of $\EE^-$:
\begin{equation}
\label{e2.3.12}
\lim_{\e\to0}\bigl(\lim_{\delta\to0}\EE_\delta^+[g_\e](x\bigl))=e^{-\lambda x}\int_{-\infty}^x e^{\lambda y}g(\varphi(y))e_+(\lambda,x-y)~dy:=\EE^+[\lambda,g(\varphi)](x),
\end{equation}
with
\begin{equation}
\label{e2.3.15}
 e_+(\lambda,z)=\int_\RR\frac{e^{iz\xi}}{\DD_c(\xi+i\lambda)}d\xi\in L^2(\RR).
\end{equation}
Consequently, the integral equation \eqref{e2.3.7} giving $\varphi$ now becomes
\begin{equation}
\label{e2.3.10}
\varphi(x)=-\frac1{2\pi}\biggl(\EE_-[g(\varphi)](x)+\EE_+[\lambda,g(\varphi)](x)\biggl),
\end{equation} 
the integral operators $\EE_\pm$ being given by \eqref{e2.3.8} and \eqref{e2.3.12}, the expression being valid as long as the function $\xi\mapsto \DD_c(\xi+i\lambda)$ does not vanish on $\RR$.

\noindent The issue is now to understand what to make out of \eqref{e2.3.10}, in other words, how to make the exponents $\lambda_\pm(c)$ or $\lambda_*$ appear. Let $\bar\lambda$ be the minimal $\lambda$ such that 
$\xi\mapsto \DD_c(\xi+i\lambda)$ does not vanish on $\RR$; examination of the expression \eqref{e2.3.6} for $\DD_c$ reveals that its zero set is bounded in every horizontal strip of the form $\Sigma_\Lambda$ with $\Lambda>0$. Therefore,
$\DD_c$ has a nontrivial zero set on the line $\RR+i\bar\lambda$; because $\DD_c$ is analytic this zero set is a finite collection of isolated points 
\begin{equation}
\label{e2.3.16}
\mathcal{Z}_{\bar\lambda}=\{\bar \zeta_j=\bar\xi_j+i\bar\lambda,\ 1\leq j\leq\bar n\},
\end{equation}
and there is $\e_0>0$ such that the rest of the zero set in the upper complex half plane is outside $\Sigma_{\bar\lambda+i\e_0}$. Consider now the function $e_+(\lambda,z)$ entering in the expression of $\EE_+[\lambda,g(\varphi)]$. We have
$$
e_+(\bar\lambda,z)=\sum_{j=1}^{\bar n}\int_{\bar\Gamma_j}\frac{e^{iz\xi}}{\DD_c(\xi+i\bar\lambda)}d\xi+e_+(\bar\lambda+i\e_0,z),
$$
where $\bar\Gamma_j$ is a circle of arbitrarily small radius enclosing $\bar\zeta_j$. Let $k_j$ be the multiplicity of the zero $\bar\zeta_{j}$, the residue theorem yields, with the usual convention $0!=1$:
$$
e_+(\bar\lambda,z)=2i\pi\sum_{j=1}^{\bar n}\frac{i^{k_j-1}z^{k_j-1}e^{iz\bar\xi_j}}{(k_j-1)!\DD_c^{(k_j)}(\bar\zeta_j)}+e_+(\bar\lambda+i\e_0,z),
$$
so that we have a  more explicit integral equation for $\EE_+[\bar\lambda,g(\varphi)]$. 

\noindent   From Lemma \ref{l2.3.2}, we
obtain the dominant term in the asymptotic expansion of $\varphi$. Let us indeed see why. The main remark is that, if we choose $\lambda=2\bar\lambda/3$, all integrals of the form
$$
\int_\RR \vert y\vert^ke^{\bar\lambda y}g(\varphi(y))~dy
$$
will be finite, as $g(\varphi(y))=O_{y\to+\infty}(\varphi(y))^2=O_{y\to+\infty}(e^{-4\bar\lambda y/3})$ and that we will have
$$
\int_{-\infty}^x \vert y\vert^ke^{\bar\lambda y}g(\varphi(y))~dy=\int_\RR \vert y\vert^ke^{\bar\lambda y}g(\varphi(y))~dy+O_{x\to+\infty}(e^{-\bar\lambda x/6}).
$$
Then, we may consider the trigonometric polynomial
\begin{equation}
\label{e2.3.18}
\mathcal{P}_+[\bar\lambda,g(\varphi)](x)= \sum_{j=1}^{\bar n}\frac{i^{k_j-1}}{(k_j-1)!\DD_c^{(k_j)}(\bar\zeta_j)}\int_{\RR} e^{\bar\lambda y}(x-y)^{k_j-1}e^{i(x-y)\bar\xi_j}g(\varphi(y))~dy.
\end{equation}
We may write $\EE_+[\bar\lambda,g(\varphi)](x)$ as:
\begin{equation}
\label{e2.3.17}
\EE_+[\bar\lambda,g(\varphi)](x)=2i\pi e^{-\bar\lambda x}\biggl(\mathcal{P}_+[\bar\lambda,g(\varphi)](x)+\EE_+[\bar\lambda+i\e_0,g(\varphi)]+O_{x\to+\infty}(e^{-\bar\lambda x/6})\biggl),
\end{equation}
and inspection of $\EE_+[\bar\lambda+i\e_0,g(\varphi)]$ reveals that it is at most an $O_{x\to+\infty}(e^{-(\lambda+\e_0/2)x})$, where we have chosen (something that we may always assume) $\e_0<1/6$, and where $\e_0/2$ stands instead 
of $\e_0$ in order to take into account some possible polynomial degeneracies. 
Identity \eqref{e2.3.8} for $\EE_-[\bar\lambda,g(\varphi)]$, then implies
 \begin{equation}
\label{e2.3.19}
\varphi(x)=-ie^{-\bar\lambda x}\mathcal{P}_+[\bar\lambda,g(\varphi)](x)+O(e^{-(\bar\lambda+\e_0/2)x}),
\end{equation}
which is the sought for dominant term. The term of highest degree of $-i\mathcal{P}_+$, that we denote $\mathcal{P}$, has the form \eqref{e3.7.17} with $k=\bar n$, and $u(x)=e^{-\bar\lambda x}\mathcal{P}(x)$ solves the linear equation \eqref{e2.2.5}. As $\varphi$ is positive, $u(x)$ should also be positive,
as the dominant term in the expansion of $\varphi$. Lemma \ref{l2.3.1} then  implies $\bar n=1$, and $\xi_1=0$. This implies in turn $\bar\lambda=\lambda_-(c)$. As no ambiguity is possible here, we delete  the dependence in $c$, and we finally obtain:
\begin{equation}
\label{e2.3.31}
\varphi(x)=-\frac{i^{k_1}e^{-\lambda_- x}}{\DD_c^{(k_1)}(i\lambda_-)}\int_{\RR} e^{\lambda_- y}(x-y)^{k_1-1}g(\varphi(y))~dy+O(e^{-(\lambda_-+\e_0/2)x}).
\end{equation}
If $c>c_*$, we have $k_1=1$, and
$
D_c'(\lambda_-)=-i\DD_c'(i\lambda_-),
$
with $D_c'(\lambda)$ given by \eqref{e2.3.41}. As $\lambda_-$ is the minimal root of $\DD_c(i\lambda)=0$, we have $D_c'(\lambda_-)<0$, so that:
\begin{equation}
\label{e2.3.310}
\varphi(x)=\frac{e^{-\lambda_- x}}{-D_c'(\lambda_-)}\int_{\RR} e^{\lambda_- y}g(\varphi(y))~dy+O(e^{-(\lambda_-+\e_0/2)x}),
\end{equation}
which proves \eqref{e2.3.3}. If $c=c_*$, then $k_1=2$ and we have this time
$$
D_{c_*}''(\lambda_*)=\DD_{c_*}''(i\lambda_*),
$$
with $D_{c_*}''(\lambda_*)>0$ given by \eqref{e2.3.41} once again. Equation \eqref{e2.3.31} yields
\begin{equation}
\label{e2.3.3100}
\varphi(x)=\frac{e^{-\lambda_*x}}{D_{c_*}''(\lambda_*)}\int_{\RR} e^{\lambda_* y}(x-y)g(\varphi(y))~dy+O(e^{-(\lambda_*+\e_0/2)x}).
\end{equation}
This proves \eqref{e2.3.4}.  \hfill$\square$

\noindent From then on, one may infer the uniqueness of the wave profile in the case $f'(1)<0$, and let us see how it works on the wave with speed $c>c_*$, the discussion being the same for the wave with bottom speed.

\noindent For $c>c_*$, let $\varphi_c(x)$ be the wave constructed in the previous paragraph, and let $\psi$ be another wave solution. We use a sliding argument to prove that $\psi$ is nothing else than a translate of $\varphi_c$, the main step being to prove that a possibly very negative translate of $\varphi_c$ is above $\psi$. So, let us consider the family of translations $\varphi_c(x-\tau)$, for $\tau>0$. Because of estimate \eqref{e2.3.3} of Theorem \ref{t2.3.1}, there is $x_+>0$ and $\tau_0>0$ such that, for all $\tau\geq\tau_0$ and all $x\geq x_0$ we have $\varphi_c(x-\tau)\geq\psi(x)$.

\noindent We claim that, possibly by enlarging $\tau_0$, then $\varphi_c(x-\tau_0)\geq\psi(x)$ for very negative $x$. Notice indeed the existence of $\theta\in(0,1)$ such that $f'(u)<f'(1)/2$ if $u\in(\theta,1)$. By assumption, there is $x_-<0$ such that $\varphi_c(x-\tau_0)$ and $\psi(x)$ are above $\theta$ if $x\geq x_-$. From Proposition \ref{p2.2.1} (strong comparison principle) we have  $\psi(x)<1$ if $x\in[x_--1,x_-]$. Therefore we may increase $\tau_0$ such that we have
$$
\phi_x(x-\tau_0)\geq\psi(x)\ \hbox{if $x_--1\leq x\leq x_-$.}
$$
As $\varphi_c'\leq0$, we have $f'(\varphi_c(x-\tau))\leq f'(1)/2$ if $x\leq x_-$ and $\tau\geq\tau_0$. We also have by construction, $f'(\psi(x))\leq f'(1)/2$ for $x\leq x_-$. As a consequence, for $\tau\geq\tau_0$, the function
\begin{equation}
\label{e2.3.50}
v_{\tau}(x)=\varphi_c(x-\tau)-\psi(x)
\end{equation}
solves
$$
\JJ v_\tau-cv_\tau'+b(x)v_\tau=0\ \hbox{for $x\leq x_--1$},\quad v_\tau(x)\geq 0\ \hbox{for $x_--1\leq x\leq x_-$.}
$$
with $b(x)=\di\frac{f(\varphi_c(x-\tau))-f(\psi(x))}{\varphi_c(x-\tau)-\psi(x)}\leq0$ by construction. This entails $v_\tau(x)\geq0$ for $\tau\geq\tau_0$ and $x\leq x_--1$. And so, as $\varphi_c'\leq0$, it suffices to increase $\tau_0$ even more to have $v_\tau(x)$ everywhere.

\noindent One may now conclude to uniqueness, as one may now define $\tau_m$, the least $\tau\in\RR$ such that $v_\tau(x)\geq0$; we notice that $\tau_m$ is trivially finite. The discussion  organises around the behaviour of $\varphi_c(x)$ and $\psi(x)$ as $x\to+\infty$. Let us first assume that $\psi(x)$ and $\varphi_c(x-\tau_m)$ are equivalent as $x\to+\infty$, that is, according to \eqref{e2.3.19}:
$$
\int_{\RR} e^{-\lambda_- y}g(\varphi(_cy-\tau))~dy=\int_{\RR} e^{-\lambda_- y}g(\psi(y))~dy.
$$
Let us now use the assumption $g'\geq0$. Then we have $g(\varphi_c(x-\tau_m))=g(\psi(x))$, and, in particular, $\varphi_c(x-\tau_m)=\psi(x)$ on a half line going to $+\infty$. So,  $x_0\in\RR$ is the least $\bar x$ such that $\varphi_c(x-\tau_m)=\psi(x)$ 
for $x\geq\bar x$, let us write the equation for $v_{\tau_m}$ at $x=x_0$, the function $v_\tau$ being given by \eqref{e2.3.50}:
$$
\int_{\RR}K(x_0-y)v_{\tau_m}(y)~dy=v_{\tau_m}(x_0)-cv_{\tau_m}'(x_0)-f(\varphi_c(x_0-\tau_m))+f(\psi(x_0))=0.
$$
As $K$ is supported in $(-1,1)$, we have $v_{\tau_m}(x)=0$ if $x_0-1\leq x\leq x_0$, contradicting the minimality of $x_0$. This implies $\varphi_c(x-\tau_m)=\psi(x)$ for all $x\in\RR$. If now $\psi(x)$ and $\varphi_c(x-\tau_m)$ are not equivalent as $x\to+\infty$, 
we may find $\e_0>0$ such that $v_{\tau_m-\e}(x)\geq0$ for large $x$. Defining $x_-$ as above, we may restrict $\e_0$ such that this inequality holds for $x\geq x_-$ and $\tau\in[\tau_m-\e_0,\tau]$. The maximum principle, applied as above to $v_{\tau}$ on $(-\infty,x_--1]$ implies $v_{\tau}(x)\geq0$ for $\tau\in[\tau-\e_0,\tau_m]$. This contradicts the minimality of $\tau_m$, hence he assumption. All in all, we have proved that $\varphi_c(x-\tau_m)$ and $\psi(x)$ coincide, which finishes the proof of Theorem \ref{t2.3.1}.
\section{The ZFK travelling waves: generalities}\label{s3.102}
\noindent  We will see later that, much more than in the case of a KPP type nonlinearity, travelling waves will be at the heart of the propagation mechanism. So, they deserve a special attention.
\begin{theorem}\label{t3.7.1}
There is $c_*\geq c_K$ such that \eqref{e3.7.1} has solutions if and only if $c\geq c_*$. If $\varphi_c$ is a solution to \eqref{e3.7.1}, any other solution of \eqref{e3.7.1} is a  translate of $\varphi_c$.
\end{theorem}
The section is organised in a very classical fashion: first we prove existence, then uniqueness.
\subsection{General existence results}
\noindent As a warm up, let us prove the existence of a solution to \eqref{e3.7.1} for large enough $c$ by the method of sub and super solutions. A subsolution is rather easily found: because of Assumption \eqref{e3.7.14}, there is 
a concave function $\underline f\leq f$ on $[0,1]$, with $\underline f'(0)=f'(0)$, $f>0$ on $(0,1)$ and $f'(1)<0$. 
The problem
\begin{equation}
\label{e3.7.8}
\varphi-K*\varphi-c\varphi'=\underline f(\varphi)\ (x\in\RR),\ \ \ \ \ \varphi(-\infty)=1,\ \varphi(+\infty)=0
\end{equation}
has a solution $\underline\varphi_c$ as soon as $c\geq c_K$, which we assume from now on, with strict inequality. We also have $\varphi_c'<0$. A supersolution is found by putting two supersolutions together: for $\delta>0$ let us set
$\bar \varphi^+=\delta\phi_{\lambda_+(c)/2},$
and let us write $f(u)\leq f'(0)u+qu^2$ for $0\leq u\leq1$. Then we have, for $x\geq0$:
$$
\di\frac1\delta\mathcal{N}_c\bar \varphi^+\geq\bigl(-D_c(\lambda)-q\delta\varphi_{\lambda_+(c)/2}\bigl)\varphi_{\lambda_+(c)/2}\\
\geq\bigl(-D_c(\lambda)-q\delta\bigl)\varphi_{\lambda_+(c)/2}.
$$
Due to the uniform strict convexity of $D_c$, and as
 $D_c(\lambda)<0$ on $(\lambda_-(c),\lambda_+(c))$ we easily obtain that
$
\di \lim_{c\to+\infty}D_c\bigl({\lambda_+(c)}/2\bigl)=-\infty.
$
 Choose $\bar c>c_K$ such that $\di D_{\bar c}\bigl({\lambda_+(\bar c)}/2\bigl)\leq-1$, and choose from then on $c\geq\bar c$. Then
 $\mathcal{N}_c\bar \varphi^+\geq0$ on $\RR_+$ as soon as 
$\delta\leq\delta_0:=1/q.
$
Let $\delta_0$ be so chosen.

\noindent The second piece of supersolution, $\bar \varphi^-$,  is hardly more difficult to find. Let us set
$M_0=\Vert f\Vert_{L^\infty([0,1])},$
then the function $\bar\varphi^-(x)=\delta-{M_0x}/c$
is obviously a supersolution, as $\mathcal{J}\bar\varphi^-=0$ due to the evenness of $K$. We have $\bar\varphi^-(0)=\bar\varphi^+(0)$, additionally we have
\begin{equation}
\label{e3.7.17}
\frac{d\bar\varphi^-}{dx}(0)=-M_0/c,\ \ \ \frac{d\bar\varphi^+}{dx}(0)=-{\delta_0\lambda_+(c)}/2,
\end{equation}
as $\di\lim_{c\to+\infty}\lambda_+(c)=+\infty$ we choose $c_0>\bar c$ such that 
$$
\frac{d\bar\varphi^-}{dx}(0)<\frac{d\bar\varphi^+}{dx}(0),\ \ \ \ \bar\varphi^+\leq\bar\varphi^-\ \hbox{on $[0,1]$.}
$$
From then on we fix $c=c_0$, the function
\begin{equation}
\label{e3.7.1901}
\bar\varphi(x)=\left\{
\begin{array}{rll}
&\bar\varphi^-(x)\ \ \hbox{if $x\leq0$}\\
&\min\bigl(1,\bar\varphi^+(x)\bigl)\ \ \hbox{if $x\geq0$.}
\end{array}
\right.
\end{equation}
satisfies $\mathcal{N}_{c_0}\bar\varphi\geq0$. Can one conclude? Not so fast, an unfortunate feature of these two functions is that they cannot be ordered, as $\bar\varphi^+$ decays faster 
than $\bar\varphi^-$ at positive infinity. The regrettable feature of this computation is that it seems difficult - at least the author was not able to achieve it - to do something similar by replacing $\phi_{\lambda_+(c)/2}$ by
$\phi_{\lambda_-(c)}$, that is, a function which is kind enough to have the same decay at infinity as $\underline \varphi$. This is why we have to take a more indirect path.

\noindent We start from the fact that, while $\bar\varphi$ is not always above $\underline\varphi$, it is above $\underline\varphi$ on a significant region. Consider $a>0$, that will soon be made to grow indefinitely, and let $\tau_a>0$ be chosen such that 
$\bar\varphi(a-\tau_a+1)=\underline\varphi(a+1)$. Then we have $\bar\varphi(x)\geq\underline\varphi(x)$ for all $x\leq a+1$. If $a>0$ is large enough, this inequality is certainly true for $x\in[a,a+1]$. To see that it holds everywhere to the left of $a$, we slide:
translate $\bar\varphi(.-\tau_a)$ sufficiently to the right so that it is indeed above $\underline\varphi$ to the left of $a$, then slide back until the inequality no longer holds. We can only be return to $\bar\varphi(.-\tau_a)$, as we would hit a contact point if the 
sliding process was blocked before. Consider the problem 
\begin{equation}
\label{e3.7.20}
\mathcal{N}_c\varphi=0\ \ (x\in[-a,a]),\ \ 
\varphi=\underline\varphi\ \ (x\in[-a-1,-a)\cup(a,a+1]).
\end{equation}
As $\underline\varphi\leq\bar\varphi(.-\tau_a)$ on $[-a-1,a+1]$, Problem \eqref{e3.7.20} has a solution $\varphi_a(x)\in[\underline\varphi(x),\bar\varphi(x-\tau_a)]$, which is additionally nonincreasing in $x$. 

\noindent Consider $\varphi_a(0)$, it is bounded away from 0 as $a\to+\infty$, as $\varphi_a\geq\underline\varphi$. If $\varphi_a(0)$ is bounded away from 1, we win: the equation $\mathcal{N}_c\varphi_a=0$ implies, because $c>0$, that $\varphi_a'$ is uniformly bounded, so that a subsequence $(\varphi_{a_n})_n$ converges in $C^1_{loc}(\RR)$. If $\varphi_\infty$ is the limit, we have $\mathcal{N}_c\varphi_\infty=0$, and $\varphi_\infty$ cannot be 0 or 1 as $\varphi_\infty(0)$ cannot be 0 or 1. Since $\varphi_\infty'<0$, then $\varphi_\infty$ has the correct limits at $\pm\infty$.

\noindent Assume the contrary.  For every $\tau>0$ let $\varphi_{\tau,a}$ solve \eqref{e3.7.20}, with $\underline\varphi(x)$ replaced by $\underline\varphi(x+\tau)$. Obviously, $\underline\varphi(x+\tau)\leq\bar\varphi(x-\tau_a)$. We have $\di \lim_{\tau\to+\infty}\varphi_{\tau,a}=0,$
uniformly on $[-a,a]$. Indeed, we have $\varphi_{\tau,a}'<0$ and $(\varphi_{\tau,a})_a$ goes to 0 uniformly on $[-a-1,-a]$. Therefore, there is $\tilde\tau_a$ such that $\varphi_{a,\tilde\tau_a}(0)=1/2$. Repeating the above argument produces a solution to 
\eqref{e3.7.8}, which shows the existence of a solution for at least one large value of $c$. Let us notice that the argument that we have just developped shows that, if  \eqref{e3.7.8} has a solution for a value of $c$, it has a solution for all values above: we just have to repeat the sub and super solution argument for all $c'>c$, and use $\varphi_c$ as a super-solution. This allows us to define $c_*$ as the least $c\geq c_K$ such that \eqref{e3.7.8} has a solution; notice that the above argument indeed shows that it is really a minimum and not an infimum. 
\subsection{Uniqueness and asymptotic behaviour}
\noindent A last general fact is the  uniqueness, the cornerstone of the argument being the asymptotic behaviour as $x\to+\infty$. This will allow a sliding argument that will not be detailed further. The asymptotic behaviour is given by the following theorem, which says that, rather expectedly, the decay of the wave at infinity is given by the positive exponential solutions of the linear equation. 
\begin{theorem}
\label{t3.7.500}
Consider $c>c_K$ and $\varphi_c$ a solution of \eqref{e3.7.1} such that $\varphi_c(0)=1/2$. There exist $\lambda\in\{\lambda_-(c),\lambda_+(c)\}$ and $\alpha_c>0$ such that 
\begin{equation}
\label{e3.7.21}
\varphi_c(x)=\alpha_ce^{-\lambda x}+O\bigl(e^{-(\lambda+\delta)x}\bigl)\ \hbox{as $x\to+\infty$,}
\end{equation}
for some $\delta>0$.

\noindent If there is a solution $\varphi_K$ to \eqref{e3.7.1}  with bottom speed $c_*=c_K$, and if $\varphi_K$ is normalised such that $\varphi_K(0)=1/2$, there is $\alpha_K\geq0$ and $\beta_K\in\RR$ such that 
\begin{equation}
\label{e3.7.2100}
\varphi_K(x)=(\alpha_Kx+\beta_K)e^{-\lambda x}+O\bigl(e^{-(\lambda_K+\delta)x}\bigl)\ \hbox{as $x\to+\infty$.}
\end{equation}
If $\alpha_K=0$, then $\beta_K>0$.
\end{theorem}
For this we need the following two results. The first one  completes Lemma \ref{l2.3.1}:
\begin{lemma}
\label{l3.7.1}
Consider $c\geq c_K$, and a real exponent $\lambda\geq\lambda_-(c)$. Consider $k$ distinct real numbers $\xi_1,\ldots \xi_k$, and the function
\begin{equation}
\label{e3.7.22}
\mathcal{P}(x)=\sum_{j=1}^k\alpha_je^{i\xi_jx},
\end{equation}
where the $\alpha_j$'s are complex numbers. Assume that $u(x)=e^{-\lambda x}\mathcal{P}(x)$ is a solution of the linear equation. 

\noindent Then, if the $\alpha_j$'s are not all zero, the function $\mathcal{P}$ cannot have a constant sign, unless we have $\lambda=\lambda_-(c)$ or $\lambda=\lambda_+(c)$, $k=1$ and $\xi_1=0$.
\end{lemma}
\noindent{\sc Proof.} If $\lambda=\lambda_-(c)$ we just apply Lemma \ref{l2.3.1}. If $\lambda\geq\lambda_+(c)$, we follow the argument of Lemma \ref{l2.3.1}. If $\lambda\in(\lambda_-(c),\lambda_+(c)]$, we follow he same lines, up to the fact that now, we set  $\bar u(x)=e^{-\lambda x}\mathcal{P}(x)$, and we have
$$
\mathcal{J}\bar u-c\bar u'-f'(0)\bar u\geq0.
$$
We compare $\bar u$ to a small multiple of $u$, and arrive at the same conclusion. \hfill$\Box$ 

\noindent The second one says that a wave  cannot decay faster than exponentially for $x\to+\infty$.
\begin{lemma}
\label{l3.7.2}
Still when $c>c_K$, a solution $\varphi_c$ of \eqref{e3.7.1} satisfies, for some constant $q>0$:
\begin{equation}
\label{e3.7.25}
\varphi_c(x)\geq qe^{-\lambda_+(c)x},\ \hbox{for $x\geq0$}
\end{equation}
\end{lemma}
\noindent{\sc Proof.} Recall the existence of $\theta>0$ such that $f(u)\geq f'(0)u$ for $u\in[0,\theta]$. Assume $\varphi_c(0)=\theta$, and consider $q\in(0,\theta)$ such that $\varphi_c(x)\geq q$ if $0\leq x\leq1$. Consider the 
exponentials
$
\phi_{\lambda}=e^{-\lambda x}.
$
Pick $\delta\in(0,q)$ small, and consider 
\begin{equation}
\label{e3.7.26}
\underline\varphi_\delta(x)=q\bigl(\phi_{\lambda_+(c)}(x)-\delta\phi_{\lambda_+(\bar c)}(x)\bigl).
\end{equation}
As $\lambda_+(\bar c)\leq\lambda_+(c)$, there is $x_\delta>0$ such that $\underline\varphi_\delta(x)>0$ if $x\in(0,x_\delta)$ and $\underline\varphi_\delta(x)<0$ if $x>x_\delta$. If $\delta$ is small enough we have $x_\delta>1$. Consider $\delta$ so chosen; then as have
$\bigl(\mathcal{J}- c\partial_x-f'(0)\bigl)\phi_{\lambda_+(\bar c)}\geq0$ because $\bar c\leq c$, we have
$$
\mathcal{J}\underline\varphi_\delta^+-c\partial_x\underline\varphi_\delta^+\geq f(\underline\varphi_\delta^+).
$$
As $\mathrm{supp}~\underline\varphi_\delta^+\cap\RR_+$ is compact, we may use sliding to conclude: for $\tau>0$ small, we have $\tau\underline\varphi^+_\delta\geq\varphi_c$ on $\RR_+$, increasing $\tau$ until 
the inequality no longer holds  we conclude that we should have $\tau=1$, otherwise there should be a contact point between $\tau\underline\varphi_\delta^+$ and $\varphi_c$. So, for all $\delta$ small enough we have $\underline\varphi_\delta^+\leq\varphi_c$. Sending $\delta$ to 0 implies the sought for inequality \eqref{e3.7.25}.\hfill$\Box$

\noindent{\sc Proof of Theorem \ref{t3.7.500} (sketch).} It is simply a matter of following the proof of equalities \eqref{e2.3.3} and \eqref{e2.3.4} in Section \ref{s2.2.2} of the preceding chapter. We still decompose $f$ as
$$
f(u)=f'(0)u-g(u),\ \ \ g(u)=O_{u\to0}(u^2);
$$
this time, however, the function $g$ does not need to have a constant sign between 0 and 1. Repeating the proof  of  Formula \eqref{e2.3.31} in Chapter \ref{short_range}, we arrive at the following identity:
\begin{equation}
\label{e3.7.23}
\varphi_c(x)=\frac{e^{-\lambda_-(c) x}}{-D_c'(\lambda_-(c))}\int_{\RR} e^{\lambda_- (c)y}g(\varphi_c(y))~dy+O(e^{-(\lambda_-(c)+\delta)x}),
\end{equation}
for some $\delta>0$. However, we are not entirely done yet: as $g(u)$ changes sign, it may well be that the integral $\di\int_{\RR} e^{\lambda_- (c)y}g(\varphi_c(y))~dy$ is zero. Assume indeed that this is the case, venturing a little further in the proof of equalities 
\eqref{e2.3.3} and \eqref{e2.3.4} reveals the existence of $\lambda>\lambda_-(c)$, a real number $\alpha_c$, an integer $k$, and a trigonometric sum $\mathcal{P}$ of the form \eqref{e3.7.22}, such that we have, as $x\to+\infty$:
\begin{equation}
\label{e3.7.24}
\varphi_c(x)=\alpha_cx^ke^{-\lambda x}\mathcal{P}(x)+O(e^{-(\lambda+\delta)x}).
\end{equation}
so that $\mathcal{P}>0$.  From Lemma \ref{l3.7.1},  $\lambda=\lambda_+(c)$ and $k=0$ if $c>c_K$, $k\in\{0,1\}$ if $c=c_K$. \hfill$\Box$
\begin{remark}
\label{r3.2.1}
The asymptotic behaviours \eqref{e3.7.21} and \eqref{e3.7.2100} are also valid for $\varphi_{c_*}$, with $\alpha_Kx+\beta_K)$ replaced by $\bigl(-\lambda_K(\alpha_Kx+\beta_K)+\alpha_K\bigl)$ if $c_*=c_K$, and $\alpha_*$ replaced by $-\alpha_*\lambda^+(c_*)$ if $c_*>c_K$.
\end{remark}
\section{ZFK waves with speed $c_*>c_K$}\label{s3.103}
\noindent What will occupy us until the end of this chapter is the following grave question: how do we relate the bottom wave speed $c_*$ to the KPP speed $c_K$? The answer, in a nutshell, is that it depends on how much the function $f$ deviates from a good natured KPP type reaction term, in other words, whether $g$ takes positive values, how large and how often. A first clue is that, if $f'(0)$ is small and $g$ large, something will happen.
\begin{theorem}\label{t3.7.2}
The set of nonlinearities $f$ for which we have $c_*>c_K$ is nonempty. In addition, if $c_*>c_K$ and $\varphi_{c_*}$ the solution of \eqref{e3.7.1} with $c=c_*$ and $\varphi_{c_*}(0)=1/2$, there exists $\alpha_*>0$ and $\delta>0$ such that we have,
as $x\to+\infty$:
$$
\varphi_{c_*}(x)=\alpha_*e^{-\lambda_+(c_*)x}\bigl(1+O(e^{-\delta x})\bigl).
$$
\end{theorem}
In other words, the wave with bottom speed chooses the maximal decay at infinity, in contrast with what happens for the Fisher-KPP type nonlinearities. 

\noindent  For small $\e>0$,  $\theta\in(0,1)$ and $\theta_1\in(\theta,1)$, let $f_\e$ be the nonlinearity defined by
\begin{equation}
\label{e3.7.27}
f_\e(u)=
\left\{
\begin{array}{rll}
&\e u\ \ \ \hbox{for $0\leq u\leq\theta$}\\
&\theta_1-u\ \ \ \hbox{for $\theta<u<1$.}
\end{array}
\right.
\end{equation}
Let aside the fact that $f_\e$ is discontinuous at $u=\theta$, expanding $\lambda\mapsto2\di\int_{\RR_+}\bigl(\cosh(\lambda y)-1\bigl)K(y)~dy$ around 0 reveals that the Fisher-KPP velocity, denoted by $c_{K,\e}$, is given by
\begin{equation}
\label{e3.7.28}
c_{K,\e}=2\sqrt{d\e}\bigl(1+o_{\e\to0}(1)\bigl),\ \ \ \ d=\int_{\RR_+}y^2K(y)~dy.
\end{equation}
For $\e>0$ and $c>c_{K,\e}$ we will construct a subsolution $\underline\varphi_\e$ to 
\begin{equation}
\label{e3.7.29}
\mathcal{J}\varphi-c\varphi'=f_\e(\varphi),\ \ \ \  \varphi(-\infty)=\theta_1,\ \mathrm{supp}\varphi\cap\RR_+\ \hbox{compact.}
\end{equation}
This will show the following:
\begin{proposition}\label{p3.7.1}
Consider $\e>0$ such that \eqref{e3.7.29} has a subsolution. Take any reaction term $f$ as described at the beginning of this section, and assume, in addition, that $f\geq f_\e$ and $f'(0)=\e$. If $c_*$ is the bottom wave speed we have
$c_*>c_{K,\e}$ as soon as $\e$ is small enough.
\end{proposition}
The proof of the proposition, once the subsolution is known, relies once again on a sliding argument: consider $f$ as in the theorem, pick $c>c_{K,\e}$ and assume the existence of a solution $\varphi_c$ to \eqref{e3.7.1}. Then, for a correct translation of 
$\varphi_c$, we have $\varphi_c\geq\underline\varphi_\e$, with a contact point, an impossibility that should be by now standard. This entails $c_*>c_{K,\e}$, and also shows that, given the large choice for $f$, the situation $c_*>c_K$ occurs.

\noindent Let us therefore prove the existence of a subsolution to \eqref{e3.7.29}, for $c>c_{K,\e}$ small enough. As usual it will be made up with two pieces: one on $\RR_+$, one on $\RR_-$. On $\RR_+$, we take our inspiration from
Lemma \ref{l3.7.2}. Pick a small $\delta>0$, and consider 
\begin{equation}
\label{e3.7.30}
\underline\varphi_{\delta,\e}(x)=\theta\bigl(\phi_{\lambda_+(c)}(x)-\delta\phi_{\lambda_+(c_{K,\e})}(x)\bigl).
\end{equation}
As $f(u)\geq\e u$ on $[0,\theta]$ we have $\mathcal{J}\underline\varphi_{\delta,\e}-c\underline\varphi_{\delta,\e}'\leq f(\underline\varphi_{\delta,\e})$ on $(1,+\infty)$. The second piece is devised, according to the definition \eqref{e3.7.27}
of $f_\e$, by finding a solution to 
$$
\mathcal{J}\varphi-c\varphi'=\theta_1-\varphi,\ \ \ \varphi(-\infty)=\theta_1,\ \varphi(0)=\theta.
$$
A solution is readily given by
$
\varphi_1(x)=\theta_1-\bigl(\theta_1-\theta\bigl)e^{\lambda_1(c)x},
$
where $\lambda_1$ solves the equation
$$
\di 2\int_{\RR}\bigl(\cosh(\lambda y)-1\bigl)K(y)~dy=-c\lambda+1.
$$
Indeed, still due to the evenness and  strict convexity of the left handside, with respect to $\lambda$, the equation has a unique positive root, which is the sought for $\lambda_1(c)$. We have $\lim_{c\to0}\lambda_1(c)>0, $ and so, we effortlessly construct our subsolution as
$$
\underline\varphi_\e(x)=\left\{
\begin{array}{rll}
&\underline\varphi_{\delta,\e}(x)\ \ \hbox{on $\RR_+$}\\
&\varphi_1(x)\ \ \ \hbox{on $\RR_+$.}
\end{array}
\right.
$$
The only item to check is $\varphi_1(x)\geq\varphi_{\delta,\e}(x)$ for $-1\leq x\leq0$. However we have ${\varphi_1}'(0)<0$, while, as $\e\to0$ we have, from \eqref{e3.7.28}: 
$\di\lim_{\e\to0}\varphi_{\delta,\e}'(0)=0$. This grants that $\varphi_1(x)$ remains above $\underline\varphi_{\delta,\e}(x)$ on a large interval of $\RR_-$, thus the full subsolution property.

\noindent Let us now concentrate on a solution $\varphi_{c_*}$ of \eqref{e3.7.1}, and assume that $c>c_*$. We wish to show that $\varphi_{c_*}$ chooses the fastest exponential decay at $+\infty$, 
thus ending the proof of Theorem \ref{t3.7.2}. This fact is a rather easy consequence of the following proposition, which will occupy us until the end of this section:
\begin{proposition}
\label{p3.7.2} 
Consider $c_0>0$ such that \eqref{e3.7.1} has a solution $\varphi_{c_0}$ for which we have $\di\lim_{x\to+\infty}e^{\lambda_-(c_0)x}\varphi_{c_0}(x)>0$. Then there is $\e_0>0$ such that the wave problem \eqref{e3.7.1}
has solutions for all $c\in[c_0-\e_0,c_0]$.
\end{proposition}
\noindent As we only need to work for $c$ close to $c_0$, a natural idea is to look for a solution  $\varphi_c$ of \eqref{e3.7.1} as a perturbation of $\varphi_{c_0}$, that is, $\varphi_{c_0-\e}=\varphi_{c_0}+\e\psi$. Doing this without further precautions, however, would not lead us too far, as such a form does not take into account the game that is played by the solutions at infinity. It is, therefore, a good moment to introduce weigthed space that we will encounter in various places of this  chapter. For all $r>0$,  let $w_r(x)$ be the weight function
\begin{equation}
\label{e3.7.36}
w_r(x)=1+e^{rx}.
\end{equation}
Let $B_{w_r,0}$ be the space of all uniformly continuous functions $v(x)$ such that $w_rv$ is a bounded, uniformly continuous function on $\RR$; for $v\in B_{w_r,0}$ we define
$
\Vert v\Vert_{r,0}=\Vert w_rv\Vert_\infty.
$
For any integer $k\geq1$, the set Let $B_{w_r,k}$ will be  the space of all uniformly continuous functions $v(x)$ such that, for all $j\in\{0,\ldots,k\}$ the $j^{th}$ derivative of $v$, denoted by $v^{(j)}$, is in
$B_{w_r,0}$. We will work with these functional spaces until the end of the section, with various ranges of $r$.

\noindent Before proving Proposition \ref{p3.7.2}, let us set 
\begin{equation}
\label{e3.7.43}
\mathcal{M}_c u=\mathcal{J} u-c u'-f'(\varphi_c) u,
\end{equation}
where $\varphi_c$ is a wave with speed $c$; we wisely notice that $\mathcal{M}_c$ is defined only if a travelling wave $\varphi_c$ exists. Then, we introduce an operator that will be useful at various places. Let $\gamma(x)$ denote, this time, a smooth nonnegative nonincreasing function that is 0 on the positive real line and 1 on 
$(-\infty,-1]$. For every $C^1$ function $u(x)$ we set
\begin{equation}
\label{e3.7.430}
\mathcal{M}^0_{c_0}u=\mathcal{J}u-cu'-\gamma(x)f'(1)u-\bigl(1-\gamma(x)\bigl)f'(0)u,
\end{equation}
and $\mathcal{K}_{c_0}=\mathcal{M}_{c_0}-\mathcal{M}^0_{c_0}$. The cornerstone of the argument is the continuous invertibility of 
 $\mathcal{M}^0_{c_0}$, as it will be used later we formulate it as an important lemma.
 \begin{lemma}
 \label{l3.7.10}
  For all $r\in(\lambda_-(c_0),\lambda_+(c_0))$, $\mathcal{M}^0_{c_0}$ is an isomorphism from $B_{w_r,0}$ to $B_{w_r,1}$.
 \end{lemma}
\noindent{\sc Proof.} Consider $f$ in $B_{w_r,0}$ and let us try to solve
\begin{equation}
\label{e3.7.45}
\mathcal{M}^0_{c_0}v=f.
\end{equation}
As $\vert f(x)\vert\leq\bar f(x)=e^r\Vert f\Vert_{w_r,0}\inf(1,e^{-r(x+1)})$, a solution to \eqref{e3.7.45} will be found as soon as we find a super-solution to the equation with $f$ replaced by $\bar f$. For $x\leq -1$, a super-solution is given by
$
\bar u(x)=\alpha$, with $\alpha\geq{e^r}/{\vert f'(1)\vert}.
$
For $x\geq-1$, we have
$$
\begin{array}{rll}
e^{rx}\mathcal{M}^0_{c_0}e^{-rx}=&2\di\int_0^{+\infty}\bigl(1-\cosh(ry)\bigl)K(y)~dy+cr-\gamma(x)f'(1)-\bigl(1-\gamma(x)\bigl)f'(0)\\
=&D_{c_0}(r)-\gamma(x)\bigl(f'(1)-f'(0)\bigl)
\geq-D_{c_0}(r).
\end{array}
$$
Recall that $D_c$ is given by \eqref{e2.3.41} and that, as $r\in\lambda_-(c_0),\lambda_+(c_0)$, we have $D_{c_0}(r_0)<0.$ As a consequence, a super-solution for $x\geq-1$ is given by
$
\bar u(x)=\alpha e^{-rx},$ with $\alpha\geq\di\frac{e^r\Vert f\Vert_{w_r,0}}{-D_{c_0(r)}}.
$
As a consequence, $\bar u(x)=\alpha\inf(1,e^{-r(x+1)})$ is a super-solution to \eqref{e3.7.45}, which entails the existence of a solution $u$ to \eqref{e3.7.45} that satisfies, in addition:
$\Vert u\Vert_{w_r,0}\lesssim \Vert f\Vert_{w_r,0}.$ This inequality can be converted into a $B_{w_r,1}$ estimate for $u$, as  $cu'$ is estimated by $u$ from the equation.

\noindent As for uniqueness, let us investigate \eqref{e3.7.45} with $f\equiv0$. Note that, because we are looking for $u\in B_{w_r,0}$, we have $u(x)\leq M\inf(1,e^{-r'x})$, for some $M>0$, and for all $r'<r$. So, pick any $r'\in(\lambda_-(c_0),r)$ and
 consider the least $M$ such that this inequality is true. In this case, the blocking only occurs through a contact point between $u$ and $M\inf(1,e^{-r'x})$, an impossibility unless $M=0$. repeating the argument for $-u$ yields $u\equiv0$, hence the lemma. 
 \hfill$\Box$

\noindent{\sc Proof of Proposition \ref{p3.7.2}.} Let us fix the translation in $\varphi_{c_0}$: for $\lambda>0$,  let  $\phi_\lambda(x)$ be the exponential given by \eqref{e3.7.6}; we will see that this alleviates the notations, even if it does not look obvious at the moment. We normalise $\varphi_{c_0}$ as
$\varphi_{c_0}(x)=\phi_{\lambda_-(c_0)}(x)\bigl(1+O(e^{-\delta x})\bigl),
$
something that we may write by virtue of Theorem \ref{t3.7.5} and the assumption that the proposition is false. For $M>0$ large, let $\gamma_M(x)$ be the usual cut-off function that is zero for $x\geq2M$, 1 for $x\leq M$ and nonnegative nondecreasing, with $\Vert\gamma_M'\Vert\leq1/M$. Proposition \ref{p3.7.2} will be proved by constructing a travelling wave $\varphi_{c_0-\e}$, $\e>0$ small, under the form
$$
\varphi_{c_0-\e}(x)=\gamma_M(x)\varphi_{c_0}(x)+\bigl(1-\gamma_M(x)\bigl)\phi_{\lambda_-(c_0-\e)}(x)+\psi(x),
$$
with $\psi$ small in the norm $B_{w_r,1}$ with
$2r=\min\bigl(2\lambda_-(c_0),{\lambda_+(c_0)+\lambda_-(c_0)}/2\bigl).$
For notational brevity, we will sometimes set $c=c_0-\e$.
\subsubsection* {Evaluating $\mathcal{N}_c\varphi_{c_0-\e}$} 
\noindent As the computation is somewhat tedious, it is not irrelevant to introduce an additional notation. If $u(x)$ is a $C^1$ function, we 
denote $\mathcal{O}(\e,x,u,u')$ any function $g_\e(x,u,u')$ that satisfies
\begin{equation}
\label{e3.7.41}
\left\{
\begin{array}{rll}
\vert g_\e\vert+\vert\partial_xg\vert
\lesssim&\bigl(\phi_{2\lambda_-(c_0)}(x)+u^2+(u')^2+\e x\un_{(-\infty,2M]}(x)\phi_{\lambda_-(c_0)}(x)\bigl)H(x-M)\\
\vert \partial_ug\vert+\vert \partial_{u'}g\vert
\lesssim&\bigl(\phi_{2\lambda_-(c_0)}(x)+\vert u\vert +\vert u'\vert+\e x\un_{[M,2M]}(x)\phi_{\lambda_-(c_0)}(x)\bigl)H(x-M),
\end{array}
\right.
\end{equation}
where $H$ is the Heaviside function.
Let us note that, even after a (normal) slightly negative reaction at the view of such a complexity, this notation is in fact expected and natural. We have
$$
\begin{array}{rll}
f(\varphi_{c_0-\e})=&f\bigl(\varphi_{c_0}+(1-\gamma_M)(\phi_{\lambda_-(c_0-\e)}-\varphi_{c_0})+\psi\bigl)\\
=&f(\varphi_{c_0})+(1-\gamma_M)(\phi_{\lambda_-(c_0-\e)}-\varphi_{c_0})f'(\varphi_{c_0})+f'(\varphi_{c_0})\psi+\mathcal{O}(\e,x,\psi,\psi').
\end{array}
$$
Writing 
$$
\begin{array}{rll}
f(\varphi_{c_0})=&\gamma_Mf(\varphi_{c_0})+(1-\gamma_M)f(\varphi_{c_0})
=\gamma_Mf(\varphi_{c_0})+(1-\gamma_M)f'(0)\varphi_{c_0}+O(\phi_{2\lambda_-(c_0)})\\
=&\gamma_Mf(\varphi_{c_0})+(1-\gamma_M)f'(\varphi_{c_0})\varphi_{c_0}+O(\phi_{2\lambda_-(c_0)}),
\end{array}
$$
one finally obtains
$$
f(\varphi_{c_0-\e})=\gamma_Mf(\varphi_{c_0})+\e f'(\varphi_{c_0})\psi+(1-\gamma_M)f'(0)\phi_{\lambda_-(c_0-\e)}+\mathcal{O}(\e,x,\psi,\psi').
$$
Playing the same game with $\varphi_{c_0-\e}'$, we have
$$
(c_0-\e)\varphi_{c_0-\e}'=c\psi'+c_0\gamma_M\varphi_{c_0}'+(1-\gamma_M)c\phi_{\lambda_-(c)}'+\mathcal{O}(\e,x,\psi,\psi');
$$
here we have put the term $\e\gamma_M(x)\varphi_{c_0}'(x)$ in the basket $\mathcal{O}(\e,x,\psi,\psi')$. Finally we have
$$
\mathcal{J}\varphi_{c_0-\e}=\gamma_M\mathcal{J}\varphi_{c_0}+(1-\gamma_M)\mathcal{J}\phi_{\lambda_-(c)}+\mathcal{J}\psi+\mathcal{O}(\e,x,\psi,\psi').
$$
We use $\mathcal{N}_{c_0}\varphi_{c_0}=0$ and 
$\LL_{c_0-\e}\phi_{\lambda_-(c_0-\e)}=0$, so that:
\begin{equation}
\label{e3.7.44}
\mathcal{N}_{c_0-\e}\biggl(\gamma_M\varphi_{c_0}+\bigl(1-\gamma_M\bigl)\phi_{\lambda_-(c_0-\e)}+\psi\biggl)=\mathcal{M}_{c_0}\psi+\mathcal{O}(\e,x,\psi,\psi').
\end{equation}
In other words, we have cleared the way to a simple application of the Banach fixed point theorem, provided that we can invert $\mathcal{M}_{c_0}$. This is what we are going to see next.
\subsubsection*{Inverting $\mathcal{M}_{c_0}$} 
\noindent  Write 
$
\mathcal{M}_{c_0}=\mathcal{M}^0_{c_0}+\mathcal{K}_{c_0}=\bigl(\mathcal{I}_{B_{w_r,0}}+\mathcal{K}_{c_0}(\mathcal{M}^0_{c_0})^{-1}\bigl)\mathcal{M}^0_{c_0}.
$
As the multiplication of a regularising operator by a smooth function tending to 0 at infinity, the operator $\mathcal{K}_{c_0}(\mathcal{M}^0_{c_0})^{-1}$ is compact in $B_{w_r,0}$. This leaves the injectivity 
of $\mathcal{M}_{c_0}$ as the only thing left to check for the invertibility of $\mathcal{M}_{c_0}$. As the argument that follows will be invoked a few times under various disguises until the end of the section, it is worth going through it in detail once, and refer to it as the "Linear Sliding Argument".

\noindent We want to prove that, if $u\in B_{w_r,1}$ solves $\mathcal{M}_{c_0}u=0$, it is automatically zero. We first notice that, because $u\in B_{w_r,0}$ and 
$r>\lambda_-(c_0)$, and because $\varphi_{c_0}$ (hence, $\varphi_{c_0'}$ by 
Remark \ref{r3.2.1}) decays exactly as $e^{-\lambda^+(c_*)x}$, then, for every $M>0$,  there exists $k_M>0$ such that $u\leq k_0\varphi_{c_0}$ on $[-M,+\infty)$.  Pick $M>0$ such that 
$f'\bigl(\varphi_{c_0}(x)\bigl)\leq{f'(1)}/2\ \ \hbox{for}\ x\leq -M.
$
As $u$ and $\varphi_{c_0}'$ solve the same equation, which has a positive zero order coefficient on $(-\infty,-M]$, the maximum principle implies $u\leq -k_M\varphi_{c_0}'$ on $(-\infty,-M)$. Let $k_{min}$ be the least $k$ such that
 $u\leq -k\varphi_{c_0}'$; assume $k_{min}>0$. As the strong maximum principle precludes any contact point, and because $u$ decays faster than $-\varphi_{c_0}'$ as $x\to+\infty$, we may lower $k_{min}$ a little and assert the existence of $\delta>0$ such that 
 $u\leq-(k_{min}-\delta)\varphi_{c_0}'$ on $[-M,+\infty)$. By the maximum principle, this inequality is also valid on $(-\infty,-M]$, so that it holds everywhere. This contradicts the minimality of $k_{min}$, so $k_{min}\leq0$. Arguing 
 similarly for $-u$ entails $u\equiv0$. 
 
 \subsubsection*{Conclusion} 
 \noindent From \eqref{e3.7.44}, solving $\mathcal{N}_{c_0-\e}\varphi=0$ amounts to solving
 $\mathcal{M}_{c_0}\psi=\mathcal{O}(\e,x,\psi,\psi')$ with $\Vert\psi\Vert_{w_r,1}$ small.
This amounts to
$
\psi=\mathcal{M}_{c_0}^{-1}\mathcal{O}(\e,x,\psi,\psi'),
$ $\Vert\psi\Vert_{w_r,1}$ small. This is 
an equation that, given the properties \eqref{e3.7.41} of $\mathcal{O}(\e,x,\psi,\psi')$, is solved by the Banach Fixed Point Theorem. \hfill$\Box$
\section{ZFK waves with speed $c_*=c_K$ and minimal decay at infinity}\label{s3.105}
\noindent \noindent Quoting one of the heroes in the monumental masterpiece \cite{MG} of Roger Martin du Gard, a scientific researcher is someone who works for a life time on five legged sheep, only to realise, at the end, that four  legged 
sheep also exist. It is legitimate to wonder whether we are venturing now in this territory.

\noindent There is, however, no need to think a lot about where one should look for such a behaviour: perturbing slightly a Fisher-KPP nonlinearity will indeed yield the result. 
Here is how we can proceed. Consider $\theta\in(0,1)$ and $g(u)$ any nonnegative $C^1$ function such that $g\equiv0$ on $[0,\theta]$.  Consider $f$ our favourite Fisher-KPP nonlinearity. For small $\e>0$ we investigate the perturbed problem
\begin{equation}
\label{e3.5.500}
\JJ\varphi-c_*\varphi'=\bigl(1+\e g(\varphi)\bigl)f(\varphi)\ (x\in\RR),\quad \varphi(-\infty)=1,\ \varphi(+\infty)=0.
\end{equation}
We will sometimes set $f_\e(u)=\bigl(1+\e g(u)\bigl)f(u)$, we evidently have $f_0=f$. The result is 
\begin{theorem}\label{t3.7.100}
There is $\e_g>0$ such that Problem \eqref{e3.5.500} has a solution, denoted by $\varphi_{c_*,\e}$, such that  we have, at positive infinity:
\begin{equation}
\label{e3.5.501}
\varphi_{c_*,\e}(x)=(x+k_*)e^{-\lambda_Kx}\bigl(1+O(e^{-\delta x})\bigl).
\end{equation}
\end{theorem}
A solution  obtained by such a kind of perturbation is  out of  the realm of five legged sheeps.  

 \noindent The theorem, however, does not exclude the existence of waves that decay exactly as $e^{-\lambda_K(c_K)x}$. This is indeed a non void possibility. It occurs much less often than \eqref{e3.5.501}, but it does occur. We leave it for the bibliographical comments, and for a problem. 
 
 \noindent{\sc Proof of Theorem \ref{t3.7.100}.} Let  $\varphi_{c_*}$ denote the solution of the honest Fisher-KPP travelling wave problem, that is, \eqref{e3.5.500} with $\e=0$, such that expansion \eqref{e3.5.501} holds for $\varphi_{c_*}$. We already notice that $\varphi_{c_*}$ is a sub-solution to \eqref{e3.5.500}; therefore, the proof of the theorem boils down to producing a super-solution to \eqref{e3.5.500} that is above $\varphi_{c_*}$.
 
 \noindent In this scope, consider $\delta\in(0,k_*/100)$ and define $M>10(\vert k_*\vert+1)$ such that we have, for $x\geq M$:
 \begin{equation}
\label{e3.5.502}
(x+k_*-\delta)e^{-\lambda_Kx}\leq\varphi_{c_*}(x)\leq(x+k_*+\delta)e^{-\lambda_Kx}.
\end{equation}
We claim that the problem
\begin{equation}
\label{e3.5.503}
\JJ\varphi-c_*\varphi'=f_\e(\varphi)\ (x\leq M),\quad \varphi(-\infty)=1,\ \varphi(x)=\varphi_{c_*}(x)
\end{equation}
has a unique solution, that we call $\varphi_{c_*,M,\e}$. Not too much imagination is required to figure that out, it suffices to note that 1 is a super-solution, while $\varphi_{c_*}$ is still a sub-solution that is below 1. Uniqueness and monotonicity in $x$ are obtained by the now routine sliding argument. 

\noindent The challenge is now to extend $\varphi_{c_*,M}$ past $x=M$. Letting $\e\to0$ and uniqueness imply that $\varphi_{M,c_*,\e}$ and $\varphi_{c_*}$ are $o_{\e\to0}(1)$ close, so that the behaviour \eqref{e3.5.502} can be turned into
 \begin{equation}
\label{e3.5.504}
(x+k_*-2\delta)e^{-\lambda_Kx}\leq\varphi_{c_*,M,\e}(x)\leq(x+k_*+2\delta)e^{-\lambda_Kx},\quad M/2\leq x\leq M,
\end{equation}
provided that $\e>0$ is sufficiently small.
Consider now  
$
\psi_{\alpha,\beta}(x)=(x+k_*+\alpha+\beta)e^{-\lambda_K(x+\alpha)},
$
 with $\alpha>0$ and $\beta+k_*>0$. Each function $\psi_{\alpha,\beta}$ solve the linear problem $\JJ\psi-c_*\psi'=f'(0)\psi$; given  the fact that $g(u)=0$ if $u\leq a$, it is a super-solution to $\JJ\psi-c_*\psi'=f_\e(\psi)$ as soon as it falls below $a$. As $\alpha>0$, we have $\di\lim_{x\to+\infty}\psi_{\alpha,\beta}(x)/\varphi_{c_*}(x)=0$, in particular, $\psi_{\alpha,\beta}$ is below $\varphi_{c_*}$ for large values of $x$. So, we are going to choose $\alpha$ and $\beta$ such that the two graphs have an intersection point somewhere between $M/2$ and $M$, taking the infimum will produce the sought for super-solution. Define
$$
x_{\alpha,\beta}=-k_*+(\alpha+\beta)\frac{e^{-\lambda_*\alpha}}{1-e^{-\lambda_*\alpha}},
$$  
an easy computation tells us that $
\psi_{\alpha,\beta}(x)\geq(x+k_*)e^{-\lambda_Kx}$ as soon as $x\leq x_{\alpha,\beta}$ and  $
\psi_{\alpha,\beta}(x)\leq(x+k_*)e^{-\lambda_Kx}$ as soon as $x\leq x_{\alpha,\beta}$. Pick any $\alpha\in(0,1)$ and $\beta:=\beta_M>0$ such that $x_{\alpha,\beta}=2M/3$. Given the asymptotic behaviour \eqref{e3.5.502}, have $\psi_{\alpha,\beta_M}(x)\leq\varphi_{c_*,M,\e}$ for $x\geq 2M/3+O(\delta)$, and $\psi_{\alpha,\beta_M}(x)\geq\varphi_{c_*,M,\e}$ for $x\leq 2M/3+O(\delta)$. Therefore, 
$\overline \varphi_{c_*,M}(x)=\inf\bigl(\varphi_{c_*,M},\psi_{\alpha,\beta_M}(x)\bigl)$ is the sought for super-solution. A suitable translation to the right pushes $\overline \varphi_{c_*,M}$ above $\varphi_{c_*}$, which ends the construction of the super-soloution, hence the proof of the theorem. \hfill$\Box$
 
 \begin{remark}\label{r3.4.20}
 A close examination of the proof of Theorem \ref{t3.7.100} reveals that it does not use much of the structure of the nonlinearity, and that it can easily be pushed into the following statement: consider $\e>0$, not necessarily close to 0, such that the problem \eqref{e3.5.500} has a solution $\varphi_{c_*,\e}$ satisfying \eqref{e3.5.501} for $\e=\e_0$. There is $\e_1>\e_0$ such that it still has a solution satisfying the expansion \eqref{e3.5.501}.
 \end{remark}
 \noindent Remark \ref{r3.4.20} entails the following slightly more general corollary, that will be useful in Section \ref{ZFK_short_range}.
 \begin{corollary}\label{c3.4.1}
 Let $f$ be a ZFK term. Consider a nonnegative $C^2$ function $g$ such that $g\equiv0$ on $[0,\theta_0]$. For a small $\e\in\RR$, set $f_\e(u)=f(u)$, and assume that $f_\e(1+\e)=0$, $f_\e>0$ on $(0,1+\e)$. Assume that Problem \eqref{e3.5.500} with $\e>0$ has a solution with minimal decay.
 
 \noindent Consider Problem \eqref{e3.5.500}, with the nonlinearity $f_\e$.  If $\vert\e\vert$ is small enough, then, for $c=c_*$, Problem \eqref{e3.5.500}, with $f=f_\e$ and $\varphi(-\infty)=1+\e$ instead of 1, has a solution $\varphi_{c_*}^\e$ with minimal decay at $+\infty$. In addition, there is  $\omega(\e)=o_{\e\to0}(1)$, such that
 $
 \bigl\vert \varphi_{c_*}(x)-\varphi_{c_*}^\e(x)\bigl\vert\lesssim \omega(\e) \min\bigl(1,\vert x\vert e^{-\lambda_*x})\bigl). 
 $
 \end{corollary}
\section{Bibliographical elements and comments}\label{s3.80}
\subsection*{Fisher-KPP waves (Section \ref{s3.101})}
\noindent The waves that we have studied at length  bear a lot of resemblance with those satisfying the original Fisher-KPP equation, about which I like to refer to Problem 82 of Arnold's trivium 
\cite{Arn}.  The ingredient that opens many doors in the study of travelling waves with classical diffusion is of course the possibility of phase portraits. When the diffusion is given by an integral operator, they
are not available anymore; however, another source of inspiration comes from the qualitative study of elliptic equations in cylinders, pioneered by Berestycki and Nirenberg in the 1990's. The arguments to show uniqueness and monotonicity are directly borrowed from them, and I have chosen to quote the paper  \cite{BN-shear}, as it covers almost all the aspects of their theory. The argument showing that the set of admissible velocities is an interval is also  due to them, in \cite{BN-shear}. A lot of their ideas is put to work on the specific setting of the nonlocal model by
Coville in his PhD thesis \cite{Cov1}, as well as  sharp maximum principle and qualitative properties by the same author in \cite{Cov2} (monotonicity), \cite{Cov3}. 

\noindent The precise asymptotic behaviour of the wave at infinity, entailing uniqueness up to translations, is due to Carr and Chmaj \cite{CC}, following earlier results by Diekmann \cite{Diekmann-78}. The asymptotic behaviour for bistable nonlocal waves is presented in Zhao-Ruan \cite{ZR}. 
The proof that we present differs from them, it is based on the Fourier transform rather than the Laplace transform in \cite{CC}. As it makes a minmal use of the structure of the nonlinearity, it allows  generalisations not covered in \cite{CC} to ZFK nonlinearities and, 
more importantly,  the multi-dimensional models displayed in Section \ref{s1.4}. Examples are presented in the problems below. 
\subsection*{ZFK travelling waves and their bottom speeds (Sections \ref{s3.102} to  \ref{s3.105})}
\noindent As said in the introduction, the nonlinearities that are studied in section  are sometimes known as ZFK nonlinearities, the acronym being this time Zeldovich-Frank Kamenetskii. The scientific background of these source terms is quite different from ecology, as they
were devised to study the different spatial scales in models for flame propagation. For the physics of this fascinating subject we refer to the book   \cite{ZBLM}, which has been a rich source of inspiration to many fundamental works that we simply cannot undertake to cite here. A glimpse of the basic asymptotic phenomena at work is sketched in Problem \ref{P3.20}, and the exploration of further mathematical questions is proposed in the problems at the end of this section.

\noindent With this precise nonlinearity, and even if this does not look that obvious at first sight, the model belongs to the big family of models with bistable nonlinearities,  rather than to the family of Fisher-KPP problems. The first study is due to  Bates-Fife-Ren-Wang \cite{BFRW}, and has seen many developments, it is difficult to cite them all. It is relevant to quote   Besse-Capel-Faye-Fouilh\'e \cite{BCFF}, in particular for the study of the phenomenon of pinning for nonlocal diffusion, something that was known for discrete models.

\noindent The result on the asymptotic behaviour $x\to+\infty$ of the wave $\varphi_c$, with $c>c_*$ is, at least for part of it, claimed in \cite{Lv}, and announced to follow the Carr-Chmaj \cite{CC} line. The proof that we propose is, once again, different, and has the potential to extend to some multi-dimensional configurations. The case $f'(0)=0$ is proposed as a problem. An interesting question is the transition from KPP to ZFK. A qualitative discussion relevant to the modelling in combustion  may be found in Clavin-Searby \cite{ClaS}. Eventually,  one finally encounters a creature that has been looming over this chapter, namely, a five legged sheep.  In other words, it is a wave with minimal speed $c_K$, yet having the maximal decay at infinity. An, Henderson, Ryzhik \cite{AHR2} baptise it in a different way, as they call it a "pushmi-pullyu" front, in reference to an animal described in detail by H. Lofting in his childrens' series featuring Dr Dolittle. In \cite{AHR2} this transition, and much more, is discussed in a wealth of details when the diffusion is given by $-\partial_{xx}$. The study of the transition at the level of travelling waves  is also proposed in Problems \ref{P3.9} and \ref{P3.271}. 

\section{Problems}\label{s3.9}
\begin{problem}\label{P3.00}
Consider the travelling wave problem
\begin{equation}
\label{e3.7.2}
\JJ\varphi-c\varphi'=f(\varphi)\ (x\in\RR),\quad \varphi(-\infty)=1,\varphi(+\infty)=0
\end{equation}
An alternative argument to Theorem \ref{t3.0.0}, due from Sattinger \cite{Sat}, is the following.
Suppose the existence of $\underline\varphi\leq\bar\varphi$, respectively a Lipschitz bounded sub-solution and a Lipschitz bounded super-solution to \eqref{e3.7.2}; one  constructs inductively a sequence $(\varphi_n)_{n\geq1}$ as
\begin{equation}
\label{e2.3.37}
\JJ\varphi_{n+1}-c\varphi_{n+1}'+L\varphi_{n+1}-f(\varphi_{n+1})=L\varphi_n,\ \ \ \varphi_1=\bar \varphi.
\end{equation}
where $L$ is a large constant.
 Work out the argument in detail.
\end{problem}
\begin{problem}\label{P2.2.10}
In the original  paper work of Kolmogorov, Petrovskii and Piskunov \cite{KPP}, we have $f(u)=u(1-u)^2$ instead of $u-u^2$. The KPP assumption is, by the way, well grounded in the biological setting.
To what extent do the results of this chapter extend to this case? 
\end{problem}
\begin{problem}\label{P3.81} This problem proposes a quick proof of the asymptotic behaviour of the  non-critical Fisher-KPP waves with standard diffusion, that is, of the solutions $\varphi_c(x)$ of
$
-\varphi''-c\varphi'=\varphi-\varphi^2
$ on $\RR$, $c>2$, with $\varphi(+\infty)=0$. 
Set $\lambda_\pm(c)=\di\frac{c\pm\sqrt{c^2-4}}2$, and
$
G_c(x)=\di\frac{e^{-\lambda_-(c)x}-e^{-\lambda_+(c)x}}{\lambda_+(c)-\lambda_-(c)}.
$
\noindent Show that $\varphi_c=G_c*\varphi_c^2$, and deduce from this integral equation the behaviour of $\varphi_c$ as $x\to+\infty$. How is the proof modified when $c=c_*$?
\end{problem}
\begin{problem}
\label{P3.71}
\noindent Work out completely the proof of Corollary \ref{c3.4.1}. Is $\omega(\e)$ of the order $\e$?.
\end{problem}
\begin{problem}\label{P3.5100} The question is here to study the travelling wave solutions to the SIR model with nonlocal contaminations. Recall the model:
$$
\partial_tI+\alpha I=SK*I\,\ \ \ \partial_tS=-SK*I,
$$
where $K$ has mass $\beta$. 
Assume, as usual, $S(0,x)$ to be a constant $S_0$ and $R_0:S_0\beta/\alpha>1$. Consider, as always, $u_*>0$ the unique positive root of $\alpha u=S_0(1-e^{\beta u})$. We look for a solution $(c,\phi,\psi)$ of
\begin{equation}
\label{e3.8.20}
\left\{
\begin{array}{rll}
&\begin{array}{rll}
-c\varphi'+\alpha\varphi=&\psi K*\varphi\\
-c\psi'=&-\psi K*\varphi
\end{array}
(x\in\RR)\\
&\varphi(\pm\infty)=0,\ \ \psi(-\infty)=S_0e^{-\beta u_*},\ \psi(+\infty)=S_0.
\end{array}
\right.
\end{equation}
\noindent The equivalent quantity of the cumulated density of infected is
$\phi(x)=\di\frac1c\int^{+\infty}_x\varphi(y)~dy.$ It solves
\begin{equation}
\label{e3.8.21}
-c\phi'=S_0(1-e^{-K*\phi})-\alpha\phi.
\end{equation}
\begin{itemize}
\item [---] Show that the conditions at infinity for $\phi$ are $\phi(-\infty)=u_*$, $\phi(+\infty)=0$.
\item [---] Let $c_*$ be the least linear wave speed for  
$
v_t+S_0\bigl(K*v-\beta v\bigl)=\alpha(R_0-1)v.
$
Show that a solution to \eqref{e3.8.21} exists if and only if $c\geq c_*$.
\item[---] Study the corresponding solutions $(\varphi,\psi)$ of \eqref{e3.8.20}. Show that $\psi'>0$. How many maxima does $\varphi$ have? Intuition suggests one, I am not so sure it is that easy.
\item[---] Study how $(\varphi,\psi)$ behaves as $R_0\to1$.
\end{itemize}
\end{problem}
\begin{problem}
\label{P3.8.31}
Instructed by the examples of Fisher-KPP and Kendall, propose assumptions on $f$ so ensuring that the problem
$$
u_t=f(u,K*u)
$$
has travelling waves for which the ranges of speeds is a semi-infinite interval. This is specifically investigated in Schumacher \cite{Sch}.
\end{problem}
\begin{problem}\label{P3.50}
From the identity \eqref{e3.7.23}, show that a ZFK wave $\varphi_c(x)$ decays exactly like $e^{-\lambda_-(c)x}$ as $x\to+\infty$, at least for large $c$. From this consideration, propose an interpretation of why the minimal speed can be larger than $c_K$.
\end{problem}
\begin{problem}\label{P3.25}
\noindent If $f$ is zero on some nontrivial range near 0, that is, $f\equiv0$ on $[0,\theta]$ and $f>0$ on $(\theta,1)$, with $0<\theta<1$, show that there is only one possible wave speed. Conversely, if $f$ is zero on the range $[0,\theta]$, approximate it by a sequence of ZFK functions 
$(f_\e)_{\e>0}$
with $f'(0)=\e$. Investigate what happens to the wave with bottom speed, then to those with higher speeds.
\end{problem}
\begin{problem}\label{P3.9}
Consider a function $a(x)$ such that $a(x)\sim_{x\to+\infty}a_0/x$, and that is positive for large negative $x$. Assume that, for a given $c>0$, there is a unique (up to multiplication) solution $v_c(x)$ to the equation
$\mathcal{K}_cv=0$, where $\mathcal{K}_c$ is the linear operator $\JJ-c\partial_x-a(x)$. Assume the existence of $\alpha>0$ such that $v_c(x)\sim_{x\to+\infty}1/x^\alpha$.
\begin{itemize}
\item[---] Show that $\mathcal{K}_c(1/x^r)\geq0$ for large $x>0$.
\item[---] For $r>0$, let $w_r(x)$ be this time a smooth positive function such that $w_r(x)=1$ for $x\leq0$, and $w_r(x)=1/x^r$ for $x\geq1$.  Define the   spaces $B_{w_r,k}$ as in Section \ref{s3.103}. Show that, for $r>\alpha$ large enough, the operator $\mathcal{K}_c$ is invertible from $B_{w_r,0}$ to $B_{w_r,1}$.
\item [---] What assumptions on $v_c$ can be removed, and proved instead?
\end{itemize}
\end{problem}
\begin{problem}\label{P3.8} (A degenerate version of ZFK). 
Assume $f(u)=u^{1+\delta}\bigl((1+o_{u\to 0}(1)\bigl),$
 with $\delta>0$. 
 \begin{itemize}
 \item[---] Prove the existence of $c*>0$ such that there is a wave $\varphi_c$ with speed $c$, for all $c\geq c_*$.
 \item[---] Walk again the path leading to Formula \eqref{e2.3.10}, but this time taking $\lambda_0=0$, to show that $\varphi_c$ satisfies
 \begin{equation}
 \label{e3.8.100}
 \varphi_c(x)=\frac1c\int_x^{+\infty}f\bigl(\varphi_c(y)\bigl)~dy+O(e^{-\lambda x}),
 \end{equation}
 for some $\lambda>0$. This, by the way, is easy to prove if the diffusion is given by $-\partial_{xx}$ and is another manifestation of the deep analogies between all these diffusions.
 \item[---] Show that a solution $\varphi(x)$ of \eqref{e3.8.100} has two possible behaviours: either $\varphi(x)\sim_{x\to+\infty}k/x^{1/\delta}$ (Type 1), or $\varphi(x)$ decays exponentially fast (Type 2).
 \item [---] Show that $\varphi_c$ is of Type 1 for large $c$.
 \item[---] Assume $\varphi_c$ to be of Type 1. Show, using Problem \ref{P3.9}, that there is $\e>0$ such that there are waves of velocity $c'\in[c-\e,c)$.
\item[---] Show that $\varphi_{c_*}$ is of Type 2. The result is claimed in \cite{ZLW}; I am not entirely sure that that the first step, Theorem 3.1, does not need more justifications than what the authors have written.
\item[---] Study the behaviour of the waves with higher speeds in the ZFK model when $f'(0)$ goes to 0.
\end{itemize}
\end{problem}
\begin{problem}\label{P3.271} Start from the standard Fisher-KPP nonlinearity $f(u)=u-u^2$ and devise a new $C^1$ nonlinearity $f_\sigma$, $0\leq \sigma\leq1$, as follows: set $f_1(u)=u-u^2$, keep $f_\sigma(u)=u-u^2$ on $[1/2,1]$, while setting $f_\sigma(u)=\sigma u$ for $0\leq u\leq 1/3$.  Thus, for $\sigma>0$ small enough we are in the framework of Theorem \ref{t3.7.100}.

\noindent Show the existence of a threshold $\sigma_0\in(0,1)$ above which we have $c_*=2\sqrt \sigma$, and below which $c_*>2\sqrt \sigma$. This phenomenon is sometimes referred to as "transition from KPP to ZFK". 
\end{problem}
\begin{problem}
\label{P3.8.30} (Transition from KPP to ZFK, continued). Another example is due to Hadeler-Rothe \cite{HR}: they consider  
$$
-u''-cu'=f_a(u),\ \ u(-\infty)=1,\ u(+\infty)=0,
$$
 with $f_a(u)=u(1+au)(1-u)$. The transition from KPP to ZFK occurs at $a=2$.
Find out the transition value of $a$, or, at least, an asymptotic expansion of the transition value,   when  the diffusion is given by a kernel $K$ which is the approximation of identity $K(x)=\di\frac1\e\rho(\frac{x}\e)$,
with $\e>0$ small. 
\end{problem}
\begin{problem}\label{P3.27}
Study how the size of the support of $K$ influences the minimal wave speed.
\end{problem}
\begin{problem}\label{P3.500}
\noindent Study the velocity of the travelling wave solutions when the kernel $K$ is not symmetric anymore.
\end{problem}
\begin{problem}\label{P3.60}
\noindent Study the  asymptotic behaviour of a (KPP or ZFK) travelling wave solutions in the double limit $x\to+\infty$ and $c\to c_*$. I am not so sure that the nature of the nonlinearity influences things that much.
\end{problem}
\begin{problem}\label{P3.20} Still in the realm of ZFK nonlinearities, consider, for $\e>0$, the function
$$
f_\e(u)=\frac1{\e^2}(1-u)\mathrm{exp}\biggl(\frac{u-1}\e\biggl),
$$
that we modify in an $\e^2$ neighbourhood of 0 so as to make it a $C^1$ function with $f_\e(0)=0$. Study the travelling wave solutions as $\e\to0$. It may be useful to start with the case of a second order diffusion, studied in detail in Berestycki, Nicolaenko, Scheurer \cite{BNS}. In order to have an interesting limit, the kernel $K$ should resemble a (perhaps large or small) approximation of the identity.
\end{problem}

\begin{problem}
\label{P3.7.32} Study the sign of  $\JJ(x^\alpha)$, for $x>0$ large and $\alpha>0$. Use the result to figure out a super-solution to
$$
\JJ \varphi-c_*\varphi'=\bigl(1+\e g(u)\bigl)f(\varphi),\ (x\in\RR),\quad \varphi(-\infty)=1,\ \varphi(+\infty)=0
$$
for $\e>0$ small, when we only assume $g$ to be $C^1$ and nonnegative, with $g'(0)>0$.
\end{problem}
\begin{problem}
\label{P3.7.33}
Let $f$ be our favourite ZFK nonlinearity: $f(0)=f(1)=0$, $f'(1)<0<f'(0)$. Let $c_K$ the Fisher-KPP speed. Assume that  \eqref{e3.7.2}
has a solution $\varphi_K$ for $c=c_K$, having the minimal decay. Find out the largest possible class of perturbations of $f$ such that this property persists.
\end{problem}
\begin{problem}\label{P3.51}
\noindent Assume the kernel $K$ to have a radically different behaviour at infinity, that is
$$
K(z)\sim_{\vert z\vert\to+\infty}{k}/{\vert z\vert^{1+2\alpha}}.
$$
Show that \eqref{e1.1.1} does not have a travelling wave solution for $\alpha\leq1/2$, but has travelling wave solutions in the range $\alpha>1/2$. Study their asymptotic behaviour at infinity in this last case.

\noindent Let $K_\e(x)$ be a smooth approximation of $\un_{[-1/\e,1/\e]}(x)K(x)$, supported, say, in $(-2/\e,2/\e)$. If $\varphi_{*}^\e$ is the wave with bottom speed $c_*^\e$, study the behaviour of $c_*^\e$ and $\varphi_{*}^\e$ as $\e\to0$. The most interesting case is of course $\alpha\leq1/2$, but the other case is quite instructive.

\noindent Some elements in the case of the fractional Laplacian are found in  Coville-Gui-Zhao \cite{CGZ}.
\end{problem}
\begin{problem}
\label{P3.70}
Study the existence of travelling waves for the discrete version of the KPP or ZFK problem, that is, one looks for a collection of functions $\bigl(u_i(t))_{i\in\mathbb{Z}}\bigl)$ solving the system 
$$
\dot u_i=2u_i-u_{i+1}-u_{i-1}+f(u_i).
$$
In particular, show that there are cases when the bottom speed is the Fisher-KPP speed, while it is strictly larger in other instances.

\noindent One may do it directly (hints may be found in Guo-Hamel \cite{GuoH} and Chen-Fu-Guo \cite{CFG}), or, even more interestingly, try to treat this system as a limiting case of the nonlocal problem \eqref{e1.1.1}.
\end{problem}
\begin{problem}
\label{P3.52}
Consider the model \eqref{e2.7.50}, that first appears in Problem \ref{P2.50}. Travelling waves are solutions of the form $\varphi(t,x-c_*t)$, with $\varphi(t,\xi)$ 1-periodic in $t$. the equation for $\varphi$ is thus
\begin{equation}
\label{e3.8.30}
\partial_t\varphi+\varphi-K*\varphi-c\varphi_x=f(t,\varphi),\ \ \varphi(t,+\infty)=0.
\end{equation}
The function $f$ is a Fisher-KPP type nonlinearity, the precise assumptions are presented in Problem \ref{P2.50}. In this problem we consider solutions of \eqref{e3.8.30} that are below the minimal positive solution $u^+(t)$ of $\dot u=f(t,u)$. 
\begin{itemize}
\item[---] Show that \eqref{e3.8.30} has a solution $\varphi_c$ for $c\geq c_*$, and that $\di\lim_{x\to-\infty}\bigl(\varphi(t,x)-u^+(t)\bigl)=0$.
\item[---] Let $\mathcal{D}_c(\xi)$ be defined as \eqref{e2.3.6}, the term $f'(0)$ being replaced by $m$, the average of $f_u(t,0)$. We also define $\phi_c(x)=\varphi_c(0,x)=\varphi_c(1,x)=\ldots=\varphi_c(n,x)\ldots$ Define $g(t,u)=f_u(t,0)u-f(t,u)$; the function $g$ is nonnegative in the vicinity of $u=0$ and (at least) quadratic.
\begin{itemize}
\item[--] With the notations of Section \ref{s2.2.2}, show that 
\begin{equation}
\label{e3.8.31}
\hat\phi_c(\xi)=\frac1{1-e^{-\mathcal{D}_c(\xi)}}\int_0^1\widehat{g\bigl(t,\varphi_c(t,.)\bigl)}(\xi)~dt.
\end{equation}
\item[--] Denote by $\varphi_c(x)$  the solution $\varphi$ to \eqref{e3.7.1} such that $\phi(0)=1/2$.  Let $a_c^\pm(t)e^{-\lambda_\pm(c)(x-ct)}$ be linear waves, with $\lambda_\pm(c)$ replaced by $\lambda_*$ if $c=c_*$. Adapt the strategy
 of Theorem \ref{t2.3.1} to prove the following. If $c>c_*$,  there  is  $k_c(t)>0$, 1-periodic, and $\delta>0$ such that
\begin{equation}\label{e2.3.3}
\varphi_c(t,x)=k_c(t)e^{-\lambda_-(c)x}+O(e^{-(\lambda_-(c)+\delta)x})\ \ \ \hbox{as $x\to+\infty$.}
\end{equation}
If $c=c_*$, there ais $k_*(t)>0$ and $l_*(t)$, both 1-periodic in $t$, such that
\begin{equation}\label{e2.3.4}
\varphi_{c_*}(x)=\bigl(k_*(t)x+l_*(t)\bigl)e^{-\lambda_*x}+O(e^{-(\lambda_*+\delta)x})\ \ \ \hbox{as $x\to+\infty$.}
\end{equation}
\item[--] Show that $\varphi_c$ is the unique solution of \eqref{e3.8.30} up to translations in $x$, and that $\partial_x\varphi_c<0$.
\end{itemize}
\end{itemize}
When $\JJ$ is replaced by the Laplacian, the results are proved by Nadin-Rossi \cite{NaRo}. 
\end{problem}
\begin{problem}
\label{P3.53}
Consider Problem \eqref{e3.8.30} with a ZFK reaction term. Show, by perturbing a constant ZFK term by small 1-periodic nonlinearities,  that we have $c_*>c_K$ in many cases.
\end{problem}
\begin{problem}\label{P3.40} The last three problems of this section are devoted to the analysis of the travelling waves of the Berestycki-Chapuisat model \cite{BCh}, with nonlocal
 longitudinal diffusion and KPP nonlinearity. So, for this problem and the next one, we are given a smooth function $f(u)=f'(0)u-g(u)$, with $f'(0)u>0$, $g(u)\geq 0$, $g(u)/u\to-\infty$ as $u\to+\infty$, and $g(u)=O_{u\to0}(u^2)$. Consider $\alpha>0$ such that the first eignevalue $\mu_1(\alpha)$ of 
$-\partial_{yy}+\alpha y^2$ satisfies $\mu_1(\alpha)<f'(0)$. All this grants the existence of a unique solution $U(y)$ to the problem
\begin{equation}
\label{e3.7.10}
-U''+\alpha y^2U=f(u),\quad U(\pm\infty)=0.
\end{equation}
The issue is to find a couple $(c,\varphi)$ solving
\begin{equation}
\label{e3.7.110}
\begin{array}{rll}
\JJ\varphi-c\varphi_x-\partial_{yy}\varphi+y^2\varphi=&f(u)\quad(x,y)\in\RR^2\\
\varphi(-\infty,y)=U(y),\ \varphi(+\infty,y)=&0.
\end{array}
\end{equation}
The convergence to the limits as $x\to\pm\infty$ is requested to be uniform in $y$. 

\noindent In the sequel, as $\alpha$ will not be made to vary, the dependence of the various quantities on $\alpha$ will be omitted. In particular the eigenvalues of $-\partial_{yy}+\alpha y^2$ will be denoted by $\mu_1,\mu_2,\ldots\mu_n,\ldots$ and the eigenfunctions $e_1,e_2,\ldots e_n,\ldots$. Because $\mu_1<f'(0)$ the equation $D_c(\lambda)=-\mu_1$ has roots only if $c$ is larger than the bottom speed $c_K$, that will also be called $c_*$ in these two problems.
\begin{itemize}
\item[---] Assume $c>c_*$. Pick $\lambda\in\bigl(\lambda_-(c),\lambda_+(c)\bigl)$, and show that $\underline\varphi(x,y)=\e\bigl(e^{-\lambda_-(c)x}-e^{-\lambda x}\bigl)^+e_1(y)$ is a subsolution to \eqref{e3.7.110} as soon as $\e>0$ is small enough.
Show that $Me^{-\lambda_-(c)x}e_1(y)$ is a super-solution, and conclude that, for all $c>c_*$, Problem \eqref{e3.7.110} has a solution $\varphi_c$.  Extend the result to $c=c_*$.
\item[---] Show that any solution of the PDE satisfying $\varphi(+\infty,y)=0$ is such that $\varphi(x,y)\leq U(y)$.
\item[---] Can one relax the uniform convergence requirement at $-\infty$ into a pointwise convergence or, even, by nothing except a feeble statement of the sort "$\varphi>0$ in $\RR^2$"? (sometimes the question is nontrivial and the answer is "no", I do not, however, believe that we are in such a configuration).
\end{itemize}
In the general case, that is, the ZFK case, things can be a little more elaborate and I am not so sure that I know about all the possible cases. Consider a nonlinearity  $f(u)=f'(0)u+g(u)$, where $g\geq0$ on an interval of the form $[0,\theta_0]$, $g(u)=O_{u\to0}(u^2)$ and $g(u)/u\to-\infty$ as $u\to+\infty$. The following specific features are worth a serious look:
\begin{itemize}
\item[---] If $f'(0)>\mu_1$,  the theory developped for the ZFK waves should work {\rm mutatis mutandis}.
\item[---] If $f'(0)<\mu_1$, travelling waves may not exist, but also may exist. For instance, if $g$ is huge on, say, the range $[1/4,3/4]$, the equation $-U''+\bigl(\alpha y^2-f'(0)\bigl)U=g(U)$ has at least one positive solution $U_+(y)$, while $u=0$ is a stable solution. The number of stable solutions is not that clear, but the existence of a travelling wave connecting 0 to the minimal stable solution is certainly true. In such a case the speed is unique.  See Berestycki-Nirenberg \cite{BN-shear} for related models in cylinders.
\end{itemize}
\end{problem}
\begin{problem}\label{P3.41}
Let $<\!e_n*,.\!>e_n$ be the projection on the eigenvector $e_n$ of $-\partial_{yy}+\alpha y^2$. Let $\varphi_c(x,y)$ be a travelling wave solution to \eqref{e3.7.110}. Set $\psi_n(x)=<\!e_n^*,\varphi_c(x,.)\!>$. As $f$ can, at least for the general results,
of the KPP or the ZFK type, the distinction will have to be made between $c_*$ and $c_K$. 
\begin{itemize}
\item[---] Starting from the equations $\JJ\psi_n-c\psi_n'=\bigl(f'(0)-\mu_n\bigl)\psi_n=<\!e_n^*,g(\varphi_c)\!>$, and analysing them in the Fourier variables,  show the existence of $\lambda\in[0,\lambda_-(c)]$, two integers $p$ and $q$, real numbers $\xi_1,\ldots\xi_q$, $\alpha_1,\ldots\alpha_q$ 
 such that 
 $$
 \varphi_c(x,y)=x^pe^{-\lambda x}e_1(y)\sum_{j=1}^q\alpha_je^{i\xi_jx}+O\bigl(e^{-(\lambda+\delta)x-y^2/8}\bigl).
 $$
 \item[---] From the positivity of $\varphi_c$, infer that $p=0$ and $q=1$ and $\lambda=\lambda_-(c)$ if $c>c_*$.
 \item [---] If $f$ is of the KPP type, then and $p=1$ and $q=2$ if $c=c_K=c_*$. If $f$ is of the ZFK type and $c_*>c_K$, then $\lambda=\lambda_+(c_*)$ if $c=c_*$. 
 \item[---] Show the uniqueness up to translations for \eqref{e3.7.110}, as well as the monotonicity of $\varphi_c$.
\end{itemize}
\end{problem}

\chapter{Sharp Fisher-KPP spreading}\label{short_range}
This chapter is devoted to the large time behaviour of the solution $u(t,x)$ to the 
equation
\begin{equation}
\label{e4.1.1}
\partial_tu+\int_{\RR}K(x-y)\bigl(u(t,x)-u(t,y)\bigl)~dy=f(u),\ \ t>0,\ x\in\RR.
\end{equation}
The initial datum $u(0,x)=u_0(x)$ will be nonnegative, compactly supported. The function $f$ will be positive on $(0,1)$ with $f(0)=f(1)=0$.  The major assumption will be
\begin{equation}
\label{e4.1.2}
f(u)\leq f'(0)u\ \ \ \hbox{for all $u\in(0,1)$.}
\end{equation}
Recall that we will often denote the integral diffusion by
\begin{equation}
\label{e2*.1.14}
\JJ u(x)=\int_{\RR}K(x-y)\bigl(u(t,x)-u(t,y)\bigl).
\end{equation}
In this chapter, the function $g$, deviation of $f(u)$ from $f'(0)u$, and defined by \eqref{e2*.1.11}, will be nonnegative on $[0,1]$. We will also need, at some places, the following stronger property:
\begin{equation}
\label{e4.0.341}
g'(u)\geq {g(u)}/u\ \hbox{for $u\geq0$.}
\end{equation}
Property \eqref{e4.0.340} implies that $g(u)=O_{u\to0}(u^2)$. It also implies the existence of a least $u_+>0$ such that $f(u_+)=0$, we will without loss of generality assume that $u_+=1$.   Equation \eqref{e4.0.341} implies that the 
function $u\mapsto\di\frac{g(u)}u$ is nondecreasing. 
When is is convenient to us, we will also assume $f'(1)<0$. All these assumptions will be quite  useful at many places, they are probably not, however, the optimal ones.  Among the nonlinearities satisfying \eqref{e4.0.340}-\eqref{e4.0.341},
are the convex $g$s such that $\di\lim_{u\to+\infty}g'(u)=+\infty$, such as $g(u)=u^2$. However, Assumptions \eqref{e4.0.340}-\eqref{e4.0.341} are more general: we may, for instance, perturb the function $f(u)=u-u^2$ into a non-concave 
function in the following way. Modify indeed the function $u\mapsto u$ into the function $h(u)$ that is equal to $u$ on $[0,1]$, 1 on $[1,2]$ and $u-1$ on $[2,+\infty)$.  Then, the function $g(u)=uh(u)$ is not convex, and one may easily regularise it by convolution into a $C^\infty$ function. Many more examples may be constructed.

\noindent Recall that the kernel $K$ is even,  smooth, compactly supported. For convenience it will be assumed to have unit mass. 
In the sequel, when \eqref{e4.1.1} is written in the reference frame moving with speed $c*$ we will denote it  for short:
 \begin{equation}
\label{e4.2.8}
\partial_tu+\mathcal{N}_*u=0,
\end{equation}
or, equivalently,
\begin{equation}
\label{e4.2.9}
\partial_tu+\mathcal{L}_*u+g(u)=0.
\end{equation}

\noindent The issue is to follow, in the most precise fashion as possible, the level sets of $u$. More precisely, for a given $\theta>0$, we look for an asymptotic expansion, as $t\to+\infty$, of the quantity
\begin{equation}
\label{e4.0.342}
X_\theta(t)=\sup\{x\in\RR:\ u(t,x)=\theta\}.
\end{equation}
The main goal of this chapter will be to prove an asymptotic expansion for $X_\theta$, precise up to $o_{t\to+\infty}(1)$ terms. Let $c_K$ denote the bottom linear wave speed, in this chapter it will be denoted by $c_*$. Let $\lambda_*$ be the unique solution of the equation with unknown $\lambda$:  $D_{c_*}(\lambda)=0$; the function $D_c$ is defined by \eqref{e2.3.41}. The whole chapter is devoted to the proof of the following result:
\begin{theorem}
\label{t2.1.1}
There is $x_\infty(\theta)$, smooth in $\theta$, such that $X_\theta(t)=c_*t-\di\frac3{2\lambda_*}\mathrm{ln}~t+x_\infty(\theta)+o_{t\to+\infty}(1).$
\end{theorem}
\subsection*{Organisation of the chapter}
The first question that one can ask is the existence of an asymptotic spreading speed, namely, whether the quantity
$
\di\lim_{t\to+\infty}{X_\theta(t)}/t
$
exists. A related question is what the solution looks like well behind $X_\theta(t)$: all this is examined in Section \ref{s4.1}. We will see that the state that the solution leaves behind is simply the upper equilibrium solution $u\equiv 1$.  For Models of the form \eqref{e4.1.1} these last two questions are not so difficult, and we will come back for a moment to the realm of simple, yet important results. The sharp asymptotics, however, are much more involved. The original work starts in Section \ref{s4.25}. There, we construct various sorts of barriers, which will be put to work in the analysis of the solution in the far field. An important tool will be the heat kernel estimated in Chapter \ref{Cauchy},  as well as a new look at what complex linear waves can do.  The logarithmic delay up to $O_{t\to+\infty}(1)$ terms is derived in Sections \ref{s4.3}, and we  show in Section \ref{s4.4} that the deviation from the logarithmic delay is asymptotically constant. 
\section{The spreading velocity}\label{s4.1}
Let $u(t,x)$ be the solution of \eqref{e4.1.1} emanating from $u_0$. 
\begin{theorem}
\label{t2.2.2}
For all $\theta\in(0,1)$ we have $X_\theta(t)=c_*t+o_{t\to+\infty}(t).$
More precisely we have
\begin{enumerate}
\item For all $\e>0$, $\di\lim_{t\to+\infty}\sup_{\vert x\vert\geq(c_*+\e)t}u(t,x)=0,$
\item for all $\e\in(0,c_*)$,  $\di  \lim_{t\to+\infty}\sup_{\vert x\vert\geq(c_*-\e)t}\vert u(t,x)-1\vert=0.$
\end{enumerate}
\end{theorem}

\noindent{\sc Proof.} Consider $c>c_*$ and 
$
\bar \phi(x)=e^{-\lambda_-(c)x},\ \ \ \bar u(t,x)=\bar\phi(x-ct).
$
As $f(u)\leq f'(0)u$, any multiple of  $\bar u(t,x)$ is a super-solution to \eqref{e4.1.1}. So, if we choose any $c>c_*$ and $u(0,x)\leq M\varphi(x)$, with $M$ large enough, 
then $u(t,x)\leq M\bar u(t,x)$, which proves Item 1.  

\noindent Consider $c<c_*$, it is suffcient to choose $c=c_*-\e$, with $\e>0$ small. In the reference frame with speed $c$, the solution $u(t,x)$ of \eqref{e4.1.1} solves Therefore, in the reference frame moving to the right with speed $c$, equation \eqref{e4.1.1} becomes 
\begin{equation}
\label{e4.2.15}
u_t+\JJ u-cu_x=f(u) \quad (t>0,\ x\in\RR).
\end{equation} 
In the same spirit as in the proof of Theorem \ref{t2.2.10}, we are going to construct a family of compactly supported sub-solutions to \eqref{e4.2.15} of arbitrarily small size. Let $\lambda_*(c)+i\omega_*(\e)$ be a solution of  $D_c(\lambda)=0$ . We notice that the discussion of the complex linear waves in Section \ref{s2.3} perturbs if $f'(0)$ is replaced by $f'(0)-\alpha$, where $\alpha$ is a suitably small positive number. This generates a new family of solutions to  $D_c(\lambda)=0$  that will be denoted by $\lambda^\alpha=\lambda_*^\alpha(c)+i\omega_*^\alpha(\e)$. We fix such an $\alpha$ once and for all.

\noindent Let $\underline\phi_{\alpha}$ be defined by $\phi_{\lambda^\alpha}$ - that is, given by \eqref{e2.2.29} with $\lambda=\lambda^\alpha$ - if $\vert x\vert\leq\di\frac\pi{\omega_*^\alpha(\e)}$, and 0 everywhere else. As $\e>0$ is assumed to be small, we may always assume $\di\frac\pi{\omega_*^\alpha(\e)}>1$. Let $\theta_\alpha>0$ be such that 
$
f(u)\geq \bigl(f'(0)-\alpha\bigl)u\ \hbox{if $0\leq u\leq\theta_\alpha$.}
$
Consider $\delta_\alpha>0$ small enough so that, for all $\delta\in(0,\delta_\alpha]$ we have
$
\delta\underline\phi_\alpha\leq\theta_\alpha.
$
This makes $\delta\underline\phi_\alpha$ a sub-solution to \eqref{e4.2.15} which can be, at the expense of a last restriction on the size of $\delta$, assumed to be below $u(1,.)$. We call $\underline u$ this sub-solution.

\noindent We may now conclude. Let $\tilde u(t,x)$ the solution of \eqref{e4.2.15} starting from $\underline u$, obviously it is below $u(t,x)$.  From Proposition \ref{p2.2.4}, it converges to a steady solution $u_\infty$ 
of \eqref{e4.2.15} that is above $\underline u$. There is something else: translating $\underline u$ as in the proof of Theorem \ref{t2.2.10} for the steady solutions of $\mathcal{J}u=f(u)$, we realise that $u_\infty$ is in fact above the maximum $\mu_\alpha$ of $\underline u$. Therefore $u_\infty$ is above the solution of \eqref{e4.2.15} starting from $\mu_\alpha$, that is, the solution of the simple ODE
$
\dot u=f(u)$ with initial datum $u(0)=\mu_\alpha,
$
which tends to 1 as $t\to+\infty$. Therefore $u_\infty\geq1$, which implies, eventually, that 
$\di
\liminf_{t\to+\infty}u(t,x)\geq1.
$

\noindent So far, we have proved Item 2 of Theorem \ref{t2.2.2} for $\e>$ small, to push the result up to $c_*$ we argue as follows. Consider $\mu_0>0$ such that 
\begin{equation}
\label{e4.2.3000}
f(u)\geq\mu_0u\ \hbox{if}\ 0\leq u\leq3/4,
\end{equation}
and, remebering that $\JJ e_\lambda=\omega_0(\lambda)e_\lambda,$ we choose $\lambda_0>0$ such that $\omega_0(\lambda_0)<\mu_0/2$. Finally, define
\begin{equation}
\label{e4.2.3001}
\underline\phi_0(x)=\left\{
\begin{array}{rll}
\di\frac{e_{\lambda_0}(x)}2\ &\hbox{if}\ \vert x\vert\leq\di\frac\pi{2\sqrt{\lambda_0}}\\
0&\hbox{for all other $x$.}
\end{array}
\right.
\end{equation}
As should be by now a routine property, $\underline\phi_0$ is a steady subsolution to \eqref{e4.2.15} with $c=0$. Fix $\delta>0$  small such that Item 2 of Theorem \ref{t2.2.2} holds for $\e=\delta$. Accordingly, consider $s_0>0$ large such that, for all $s\geq s_0$:
\begin{equation}
\label{e4.2.3003}
u(s,x)\geq\frac12\ \hbox{for}\ \bigl\vert x-(c_*-\delta)s\bigl\vert\leq \frac\pi{\sqrt{\lambda_0}}.
\end{equation}
Consider now
$t\geq2s_0,$, with $s_0\leq \di s\leq\frac{c_*-2\delta}{c_*-\delta}t.$
Because of \eqref{e4.2.3003}, we have 
$$
u(s,x)\geq\underline\phi_0\bigl(x-(c_*-\delta)s\bigl)\ \hbox{for}\ \bigl\vert x-(c_*-\delta)s\bigl\vert\leq 1+\frac\pi{2\sqrt{\lambda_0}}.
$$
We have assumed implicitely that $\lambda_0$ was small enough to that $\omega_0(\lambda_0)$ was also small enough so that $\di\frac\pi{2\sqrt{\lambda_0}}>1.$ Running the Cauchy Problem Item 2 of Theorem \ref{t2.2.2} for $s'\in[s,\di\frac{c_*-2\delta}{c_*-\delta}t)$ we infer, by comparison:
$$
u\bigl(\frac{c_*-2\delta}{c_*-\delta}t,x\bigl)\geq\underline\phi_0\bigl(x-(c_*-\delta)s\bigl)\ \hbox{for}\ \bigl\vert x-(c_*-\delta)s\bigl\vert\leq 1+\frac\pi{2\sqrt{\lambda_0}}.
$$
As a consequence, we obtain
$$
u\bigl(\frac{c_*-2\delta}{c_*-\delta}t,x\bigl)\geq\frac12\ \hbox{if}\ (c_*-\delta)s_0\leq\vert x\vert \leq (c_*-\delta)t.
$$
In fact, this inequality is valid for all $x\leq(c_*-\delta)t$, as one may apply Theorem \ref{t2.2.10} for $\vert x\vert\leq(c_*-\delta)s_0$ for $s'\in\bigl[s_0,\di\frac{c_*-2\delta}{c_-\delta}t\bigl]$. So, we may conclude. Consider $q_0>0$ such that $f(u)\geq q_0(1-u)$ for $1/2\leq u\leq1$. Pick $\e\in(0,1/10)$ and $t_\e>0$ such that 
$$
u(t,x)\geq1-\e\ \hbox{if}\ \vert x\vert\in[c_*-\delta)t,(c_*-\delta)t+1].
$$
The function $v(s,x)=1-u(s,x)$ solves
\begin{equation}
\label{e4.2.3005}
\begin{array}{rll}
v_s+\JJ v+q_0v\leq&0\quad \bigl(\di\frac{c_*-2\delta}{c_*-\delta}t<s\leq t,\ X_1(t)<x<X_2(t)\bigl)\\
v(s,x)\leq&\e\quad (X_1(s)-1\leq x\leq X_1(s),\ X_2(s)\leq x\leq X_2(s)+1),
\end{array}
\end{equation}
with $X_1(s)=-(c_*-\delta)s$ and $X_2(s)=(c_*-\delta)s$. Let
$
\bar v(s)=\e+\di\mathrm{exp}\biggl(-q_0\bigl(s-\frac{c_*-2\delta}{c_*-\delta}t\bigl)\biggl)/2.
$
It is therefore a super-solution to \eqref{e4.2.3005} on it domain of definition, which is larger than $v$ at the initial time $s=\di\frac{c_*-2\delta}{c_*-\delta}t$ and $X_1(s)-1\leq x\leq X_1(s),\ X_2(s)\leq x\leq X_2(s)+1$. So, $v(s,x)\leq\bar v(s)$ for $s\geq \di\frac{c_*-2\delta}{c_*-\delta}t$ and $X_1(s)\leq x\leq X_2(s)$, by Proposition \ref{p2.2.4}. At time $s=t$ we have 
$\bar v(s)=\e+\mathrm{exp}\bigl({-q_0\delta t}/{(c_*-\delta)}\bigl),
$
and it now suffices to enlarge $t$ to finish the proof. \hfill$\Box$
%
\section{Dirichlet heat kernels, barriers}\label{s4.25}
Unsurprisingly, the operator $\LL_*$ appears again. Recall that we have
$
\LL_*u(x)=e^{-\lambda_*x}\II_*\bigl(e^{\lambda_*x}\bigl)(x),
$ for all function $u(x)$. 
Here $\lambda_*$ is the critical exponent related to $c_*$, and where we also recall 
\begin{equation}
\label{e4.4.0001}
\II_*v(x)=\int_{\RR}K_*(x-y)\bigl(v(x)+(y-x)v_x(x)-v(y)\bigl)~dy,
\end{equation}
where 
$K_*(x)=e^{\lambda_*x}K(x).$
We have already seen that $e^{-t\II_*}$ has a heat kernel like behaviour at large times. This insight will not, for a reason of scale variation the will be duly explained in Section \ref{s4.3} below,  be sufficient. This is why we undertake, before analysing the solutions of \eqref{e4.2.9} in more detail, to understand what a Dirichlet problem for $e^{-t\II_*}$ may look like.
In any case, independently of the context, it is  natural to ask what the solutions of an equation of the form
\begin{equation}
\label{e4.3.3000}
v_t+\II_* v=0\quad(t>0,\ x>0)\ \ \ \ \ 
v(t,x)=0\ (t>0,-1\leq x<0)
\end{equation}
will look like. Note the strict inequality for $x$, we do not, indeed, know whether a possible solution of
\eqref{e4.3.3000}, extended by 0 outside $\{x>0\}$, is continuous at $x=0$.  In fact, there is no explicit solution we can think of. Fortunately, we may take our inspiration from the classical heat equation: the Dirichlet problem can simply be solved by computing the solution of the heat equation on the whole line, emanating from the odd extension of the initial datum. This simple idea does not allow us to retrieve a true solution of \eqref{e4.3.3000}. However, the large time behaviour theorem \ref{t2.4.2} tells us that the solution of the Dirichlet heat equation is a faithful representation of the solution of \eqref{e4.3.3000} on a large part of the domain of integration. The goal of this section is to show that, in the parts that will be of interest to us, we may use the Dirichlet heat equation as a {\it bona fide} reresentation.
 \subsection{The Dirichlet heat equation}
The following material is exceedingly standard, it is nevertheless useful to  the crucial part of the analysis of $v(t,x)$, namely  in the area $x\sim t^\gamma$ with $\gamma\in(0,1/2)$. So, we decide to describe, in a few lines, what is 
going to be useful to us. Consider the solution $w(t,x)$ of
\begin{equation}
\label{e2*.5.100}
\partial_tw-d_*\partial_{xx}w=0\  (t>0,\ x>0)\ \ \ \ \
w(t,0)=0\ \ \ \ \
w(0,x)=v_0(x)\ \hbox{compactly supported}.
\end{equation}
Let $v_0^*$ be the odd extension of $v_0$ to $\RR$, in other words, the odd function coinciding with $v_0$ on $\RR_+$. Then $w(t,x)=e^{t\partial_{xx}}v_0^*(x)$, that is, the heat kernel of the whole line applied to $v_0^*$. The expression of $w(t,x)$ is therefore
$$
w(t,x)=\int_0^{+\infty}\frac{e^{-\frac{(x-y)^2}{4d_*t}}-e^{-\frac{(x+y)^2}{4d_*t}}}{\sqrt{4\pi d_*t}}v_0(y) dy.
$$
Inspection of this integral yields, setting $\eta={x}/{\sqrt{d_*t}}$, that 
$
w(t,x)\di\sim_{t\to+\infty}\eta e^{-\frac{\eta^2}4}/({2\sqrt\pi d_*t}).
$
The asymptotic slope of $w(t,x)$ is therefore the first moment of $v_0$ and we have, for every $\gamma\in(0,1/2)$:
\begin{equation}
\label{e2*.5.102}
\lim_{t\to+\infty}t^{3/2-\gamma}w(t,t^\gamma)=\frac1{2\sqrt\pi d_*^{3/2}t^{3/2}}\int_0^{+\infty}yv_0(y)~dy.
\end{equation}
All this is readily seen  by the change of variables $z={y}/{\sqrt t}$, and the fact that $v_0$ is compactly supported. This simple formula will be our holy grail for the solution $v(t,x)$ of the full Fisher-KPP problem, possibly at the expense of restricting $\gamma$ a little. When the initial datum is well spread, that is, $v_0(y)=w_0({y}/{\sqrt s})$, then $e^{td_*\partial_{xx}}v_0$ depends on how the relative values of $t$ and $s$, as is seen on the formula
$$
w(t,x)=\int_0^{+\infty}\frac{e^{-\frac{s}{4d_*t}(\frac{x}{\sqrt s}-\zeta)^2}-e^{-\frac{s}{4d_*t}(\frac{x}{\sqrt s}+\zeta)^2}}{\sqrt {4\pi d_*t}}\sqrt s w_0(\zeta)~d\zeta.
$$
If $s\gg t$, ${x}/{\sqrt{s}}\sim s^\gamma$, and down to quantities of order 1, we have
\begin{equation}
\label{e2*.5.105}
w(t,x)=e^{\frac{t}std_*\partial_{\eta\eta}}(w_{0,*}^+)(x/\sqrt s)\bigl(1+O(e^{-\frac{s^\gamma}t})\bigl).
\end{equation}
It indeed suffices to notice that the main contribution in the integral yielding to \eqref{e2*.5.102} comes from $\zeta\sim x/\sqrt s$. The integral may  then be seen as a regularisation of $v_0$ at the scale $\sqrt s$. When $t$ becomes significantly larger than $s$, $w(t,x)$ behaves as the solution of \eqref{e2*.5.100} and we have
\begin{equation}
\label{e2*.5.106}
w(t,x)\sim_{\substack{t\to+\infty\\x=O(\sqrt t)}}\di\frac{sx}{2\sqrt\pi d_*^{3/2}t^{3/2}}\di\int_0^{+\infty}\zeta w_0(\zeta)~d\zeta.
\end{equation}
If $w_0(\zeta)$ is proportional to $s^{-1}$, formula \eqref{e2*.5.106} is compatible with \eqref{e2*.5.102}.
\subsection{The action of $e^{-t\II_*}$ on odd extensions}
\noindent We are interested in pursuing the analogy between the Dirichlet heat kernel of $\II_*$ and that if $-\partial_{xx}$, not only for large times but for finite times. Formulated in such a general way, the entreprise does not look that promising. However, there is a type of data for which there is a hope to do it, and they are precisely the well spread ones, already encountered in Chapter \ref{Cauchy}.  
\begin{theorem}
\label{t4.3.3000}
Consider  $v_s(x)$ a $C^\infty$ function, such that all the derivatives of $v_s$ are bounded on $\RR_+$, independently of $s$. Assume the existence of a function $w_s(\eta)\in L^1(\RR_+)\cap L^\infty(\RR_+)$, such that we have, for all $\eta>0$, for some $A>1$ and small $\delta>0$:
\begin{equation}
\label{e2*.6.5}
(\eta-s^{\delta-\frac12})e^{-A\eta}\un_{[0,A]}(\eta)\lesssim w_s(\eta)\lesssim(\eta+s^{\delta-\frac12})e^{-A\eta},
\end{equation} 
and such that 
\begin{equation}
\label{e2*.6.6}
v_s(x)=s^{-1}w_s({x}/{\sqrt s}).
\end{equation}
Let $v^*_s$ be the odd extension of $v_s$ on $\RR_-$.  
Consider $\e>0$, $\e\gg s^{-1}$. 
\begin{enumerate}
\item {\bf (Behaviour near the origin)}
We have, for $\tau\geq0$:
\begin{equation}
\label{e2*.6.3}
{\e\sqrt s}/{(\tau+s)^{3/2}}\lesssim e^{-\tau \II_*}v^*_s(\e\sqrt s)\lesssim{\e\sqrt s}/{(\tau+s)^{3/2}}.
\end{equation}
Similarly we have, to the left:
\begin{equation}
\label{e2*.6.4}
-{\e\sqrt s}/{(\tau+s)^{3/2}}\lesssim e^{-\tau \II_*}v^*_s(-\e\sqrt s)\lesssim-{\e\sqrt s}/{(\tau+s)^{3/2}}.
\end{equation}
\item {\bf (Slope estimate).} For all $B>0$,and all $\gamma\in(0,1/2)$, there is $q_{\gamma,B}\in(0,1)$ such that, for all $s\geq1$, all $\tau>0$, the following estimates hold.
\begin{itemize}
\item[--] If $x\in((\tau+s)^{\gamma},B\sqrt{\tau+s})$, then 
$${x}/{(\tau+s)^{3/2}}\lesssim e^{-\tau\II_*}v^*_s(x)\lesssim{x}/(\tau+s)^{3/2}.$$
\item [--] If $0\leq x\leq (\tau+s)^\gamma$, then, for all $\delta>0$ we have
$$e^{-\tau\II_*}v^*_s(x)\lesssim\max\biggl({x}/(\tau+s)^{3/2},1/{(\tau+s)^{3/2-\delta}}\biggl).
$$
\end{itemize}
\end{enumerate}
\end{theorem}
\noindent{\sc Proof.} 
Let us see why \eqref{e2*.6.3} and \eqref{e2*.6.4} are true. In what follows, we will argue according to the relative sizes of $\tau$ and $s$. As soon as $\tau\geq s^\kappa$, for a small $\kappa>0$, we use Theorem \ref{t2.4.2} to claim
$
e^{-\tau \II_*}v^*_s(x)=e^{\tau d_*\partial_{xx}}v^*_s(x)+O\bigl(e^{-\tau ^\delta}\bigl).
$
So, $e^{-\tau d_*\partial_{xx}}$ is the dominant term, and we only need to invoke Formula \eqref{e2*.5.106}, with $t$ replaced by $\tau$. We then argue as follows:
\begin{itemize}
\item [--] If $\tau\geq\e s$, we may write, for $x=\e\sqrt s$:
$
e^{-{(x-y)^2}/{4d_*\tau}}-e^{-{(x+y)^2}/{4d_*\tau}}\simeq {\e\sqrt sy}/{d_*\tau},
$
and use the asymptotic formula \eqref{e2*.5.106} once again to conclude.  Some care should be given to the fact that the remainder term uses the $H^1$ norm of $v_s^*$; this is, however, easily circumvented:  the assumptions on $v_s^*$ and the integration by parts $\di\int_{\RR_+}(\partial_xv_s^*)^2=-\di\int_{\RR_+}v_s^*\partial_{xx}v_s^*+O(1)$, the $O(1)$ being independent of $s$, indeed imply that $\Vert\partial_xv_s^*\Vert_{H^1(\RR_+)}$ is estimated by $\sqrt s$. Hence the remainder term is estimated by $\sqrt se^{-\tau^{1-2\gamma}}=O( e^{-\tau^{1-3\gamma}})$.
\item[--] If $s^\delta\leq \tau\leq \e s$, the order of magnitude of the term $e^{-{(x+y)^2}/{4d_*\tau}}$ becomes negligible in front of that of $e^{-{(x-y)^2}/{4d_*\tau}}$. So we use, this time, formula \eqref{e2*.5.105}.
\item[--] In the regime $0\leq \tau\leq s^\delta$, we resort to Theorem \ref{t2.4.4} on the well-spread initial data, to the caveat that the smoothness of $v_s$ is an issue: indeed, the theorem involves the $H^m$ norm of $w_s$, $m$ possibly large,
and we do not have an estimate for $\Vert w_s\Vert_{H^m}$. We will see in our application that the only available estimate will be the insufficient   one $\Vert w_s\Vert_{H^m}\lesssim s^{\frac{m}2}$. To circumvent this, we use Inequality  \eqref{e2*.6.5}; in order to prove \eqref{e2*.6.3} we invoke the left inequality: let
$\underline w_s(\eta)$ be equal to a multiple of $(\eta-s^{\delta-1/2})e^{-A\eta}$ on $[\e^2,+\infty)$, to a multiple of $[\eta-s^{\delta-1/2})e^{A\eta}$, and regular on $\RR$, so that it is in the end below $w_s^*(\eta)$. Then, Theorem \ref{t2.4.4} yields
$$
e^{-\tau\II_*}v^*_s(x)\geq\di\frac1se^{-\tau\II_*}\underline w_s^*(\di\frac{x}{\sqrt s})\\
=\di\frac1se^{\tau d_*\partial_{xx}}\underline w_s^*(\di\frac{x}{\sqrt s})+O\bigl(\di\frac{\Vert \underline w_s^*\Vert_{L^1}+\Vert \underline w_s^*\Vert_{H^m}}{s^{\frac52-4\omega-\delta}}\bigl).
$$
Now, in the integral expression for $e^{-\tau\II_*}v^*_s(x)$, computed at $x=\e\sqrt s$, only the contribution of the interval $[0,\e\sqrt s]$ really matters, the contribution coming from the integration on $\RR_-$ being of the order $e^{-\e^2s^{1-\delta}}$. The latter, hence, is negligible in comparison of the integral on $\RR_+$, which is of the order $\e/{s}$: this entails \eqref{e2*.6.3}. For the upper bound we work with the upper bound in \eqref{e2*.6.5}.

\noindent For \eqref{e2*.6.4} we repeat the argument with, this time, the right inequality in \eqref{e2*.6.5}.
\end{itemize}

\noindent The proof of Point 2 of the theorem is essentially identical to the proof of Point 1, as we just have to use Theorem \ref{t2.4.2} for $\tau\geq s^\delta$ and Theorem \ref{t2.4.4} in the range $\tau\leq s^\delta$; only the place of $x$ changes in 
the former range. The arguments are, however, the same. \hfill$\Box$

\noindent Theorem \ref{t4.3.3000} entails a corollary that may look trivial at first sight, but that is probably false in general. Only the particular form of our initial datum makes it true, and this will be important when we construct barriers for the nonlinear problem \eqref{e4.2.9}.
\begin{corollary}
\label{c4.3.1}
Let $v_s$ be as in Theorem \ref{t4.3.3000}. Consider $\gamma\in(0,1/2)$ small. Then there is $s_\gamma>0$ such that, for $s\geq s_\gamma$, $x\geq (t+s)^\gamma$ and $t\geq0$ we have: $e^{-t\II_*}v_s^*(x)\geq0$.
\end{corollary}
\noindent{\sc Proof.} In Theorem \ref{t4.3.3000}, we choose $\delta\ll\gamma$ and $s\geq s_\gamma$ accordingly. Point 2 implies the existence of a band of the form 
$
\{(t+s)^{\gamma}-1\leq x\leq (t+s)^{\gamma}\},
$
on which $e^{-t\II_*}$ is nonnegative. As it is also nonnegative at time $t=0$ for $x\geq s^\gamma$, the comparison principle Proposition \ref{p2.1.200} applies with $X_1(t)=(t+s)^{\gamma}$. \hfill$\Box$

\noindent We end this section by a proposition quite similar to Proposition \ref{p2.4.1}, which essentially says that the $t^{-1}$ decay observed for $e^{-t\II_*}v_0^*$ in the range $x\sqrt t=O(1)$ is conserved for larger $x$, with a factor integrable in $x\sqrt t$. This is by no means surprising, as this behaviour is observed for the standard Dirichlet heat equation.
\begin{proposition}
\label{p4.2.3}
If $v_0$ is compactly supported, then, for all $A>0$, there exists $B>0$ such that, if ${x}/{\sqrt{t+1}}\geq B$, we have
$\big\vert e^{-t\II_*}v_0^*(x)\big\vert\leq \Vert v_0\Vert_{H^1} e^{-\frac{Ax}{\sqrt{t+1}}}/t.$
\end{proposition}
\noindent {\sc Proof.} We reproduce the computation of Proposition \ref{p2.4.1}, with this time the function $e^{-Ax/\sqrt{t+1}}/(t+1)$ as a candidate for a super-solution. We have
$$
\begin{array}{rll}
\bigl(\partial_t+\II_*)\bigl)\bigl(\di\frac{e^{-\frac{Ax}{\sqrt{t+1}}}}{t+1}\bigl)=&\di\frac{A}{t+1}\biggl(\frac\xi2\frac1{A(t+1)}-Ad_*\bigl(1+\omega(\di\frac{A}{\sqrt{t+1}})\bigl)\biggl)e^{-\frac{x}{\sqrt{t+1}}}\\
\geq&0\ \ \hbox{if $A>1$ and $\xi\geq 1+Ad_*(1+\Vert\omega\Vert_{L^\infty([0,1])}$}).
\end{array}
$$
So, let us declare $\bigl(\partial_t+\II_*)\bigl)\bigl(e^{-\frac{Ax}{\sqrt{t+1}}}/(t+1)\bigl)\geq0$ if $x/\sqrt{t+1}\geq B$. For $x\in[B\sqrt{t+1}-1,B\sqrt{t+1}]$, we use Theorem \ref{t2.4.2} and estimate \eqref{e2*.5.106}, which says exactly that $e^{td_*\partial_{xx}}v_0^*$ behaves like $t^{-1}$ when $x\sqrt{t+1}$ is bounded. An routine application of the maximum principle finishes the proof. \hfill$\square$
\subsection{Correcting barriers around the origin}
As the approximation of the solution $u(t,x)$ of the full problem \eqref{e4.2.9} by (an exponential multiple of) $e^{-t\II_*}$ become less and less precise as we get close to close to the origin,  we will need an additional set of barriers to make up for this effect. To devise it, we apply the computation of  Proposition \ref{p2.4.1}, Section \ref{s2.4.4} to the following situation: pick $A>1$ and $\alpha\in(0,1/2)$. 
As we have $\II_*v(x)=d_*k^2\bigl(1+\omega(k)\bigl)$, we take $k=i(1+t)^{-\alpha}$; there are two  real analytic functions $\omega_j(k)$  satisfying $\omega_j(0)=0$ such that:
$$
\begin{array}{rll}
&\di\bigl(\partial_t+\II_*\bigl)\bigl((1+t)^{-A}\cos(\frac{x}{(1+t)^\alpha})\bigl)\\
=&\di\mathrm{Re}~\biggl(\bigl(\partial_t+\II_*\bigl)((1+t)^{-A}e^{ix/(1+t)^\alpha})\biggl)\\
=&\di\mathrm{Re}~\biggl[\biggl(\frac{d_*}{(1+t)^{2\alpha}}\bigl(1+\omega_1((1+t)^{-\alpha})\bigl)-\frac{A}{1+t}-i\bigl(\frac{\alpha x}{(1+t)^{1+\alpha}}-\frac{d_*\omega_2((1+t)^{-\alpha})}{(1+t)^{2\alpha}}\bigl)\biggl)\frac{e^{ix/(1+t)^\alpha}}{(1+t)^{A}}\biggl]\\
=&(1+t)^{-A}\di\biggl(\frac{d_*}{(1+t)^{2\alpha}}\bigl(1+\omega_1((1+t)^{-\alpha})\bigl)-\frac{A}{1+t}\biggl)\cos(\frac{x}{(1+t)^\alpha})\\
&\di+(1+t)^{-A}\bigl(\frac{\alpha x}{(1+t)^{1+\alpha}}-\frac{d_*\omega_2((1+t)^{-\alpha})}{(1+t)^{2\alpha}}\bigl)\sin(\frac{x}{(1+t)^\alpha}).
\end{array}
$$
This computation, in view of the targetted expression \eqref{e2*.4.49}, tells us how we are going to control $v(t,x)$ in a (large) vicinity of $x=0$, and where the cosine comes from. We may push the computation a little more, as we are interested in 
a restriction of $\cos\bigl({x}/{(1+t)^\alpha}\bigl)$. Notice indeed that $\cos\bigl({x}/{(1+t)^\alpha}\bigl)$ vanishes at $x={3\pi(1+t)^\alpha}/{2}$, and remains nonnegative from there until $x={5\pi(1+t)^\alpha}/{2}$, that is, well outside the
support of $K_*$.  Proposition \ref{p2.2.2} implies:
\begin{proposition}
\label{p4.3.4004}
For $1\leq x\leq{3\pi(1+t)^\alpha}/2$, we have
$$
\begin{array}{rll}
&\di(1+t)^A\bigl(\partial_t+\II_*\bigl)\bigg(\un_{(-\infty,\frac{3\pi(1+t)^\alpha}2]}\bigl((1+t)^{-A}\cos(\frac{x}{(1+t)^\alpha})\bigl)\biggl)\\
\geq&\di\biggl(\frac{d_*}{(1+t)^{2\alpha}}\bigl(1+\omega_1(t)\bigl)-\frac{A}{1+t}\biggl)\cos(\frac{x}{(1+t)^\alpha}\bigl)\di+\bigl(\frac{\alpha x}{(1+t)^{1+\alpha}}-\frac{d_*\omega_2(t)}{(1+t)^{2\alpha}}\bigl)\sin(\frac{x}{(1+t)^\alpha}),
\end{array}
$$
where $\omega_i(t)$ denote two possibly different functions of $t$, going to 0 as $t\to+\infty$.

\end{proposition}
\section{The logarithmic delay estimate to $O(1)$ terms}\label{s4.3}
\noindent  The  idea that 0 is the most unstable value of the range of $u$, which had worked  well in finding out the spreading speed,
will turn out to be a useful guideline here again. There is, however, one aspect that will require some attention:  the precise form of $u(t,x)$ in the range where it is almost 0 will have to be scrutinised. So, in order to uncover the small values of $u$, we remove the exponential at which it is expected to decay:
$u(t,x)=e^{-\lambda_*x}v(t,x), 
$
so that equations \eqref{e4.2.8} or \eqref{e4.2.9} become
\begin{equation}
\label{e4.4.7}
\partial_tv+\II_*v+h(x,v)v=0,
\end{equation}
with the initial datum 
$v_0(x)=e^{\lambda_*x}u_0(x).$
We have  set $h(x,v)={e^{\lambda_*x}g(e^{-\lambda_*x}v)}/{v}$,
we have, for $x\geq0$ and $v$ in a bounded set:
$
h(x,v)\geq0,$ $h(x,v)\lesssim e^{-\lambda_*x}v.
$
The goal of this section is to prove 
\begin{theorem}
\label{t4.4.4}
Let $X_\theta(t)$ be given by \eqref{e4.0.342}. Then we have $\di X_\theta(t)=c_*t-\frac3{2\lambda_*}\mathrm{ln}~t+O_{t\to+\infty}(1)$.
\end{theorem}

\subsection{Analysis at $x\sim t^\gamma$}\label{s2.5.2}
\begin{theorem}
\label{t2.4.5}
For all $\gamma\in(0,1/2)$ we have $
1/{t^{3/2-\gamma}}\lesssim v(t,t^\gamma)\lesssim 1/{t^{3/2-\gamma}}.$
\end{theorem}
Another way to formulate the theorem is that, as $t\to+\infty$ and $x\in(t^\gamma,t^{1/2})$, the slope of $v(t,x)$, in the rescaled coordinate $\eta=x/\sqrt t$, is of the same order as ${1}/{t^{3/2}}$.
The main ingredient in the analysis of $v(t,x)$ for large $x$ is a comparison with the heat kernel computed in the preceding section.   Before we state the main lemma that will imply the theorem, let us prepare the initial datum. Consider a scale parameter $s>0$ that may be enlarged at request.
\begin{itemize}
\item[---] Let $\overline v_s(x)$ have the form $s^{-1}\overline w_0(x/\sqrt s)$, with $\underline w_0(\eta)=B\eta e^{-A\eta}$. Pick  $\delta>0$, a constant that will be as small as we wish. From Corollary \ref{c4.3.1}, we have $e^{-t\II_*}\overline v_s(x)\geq0$ if $x\geq s^\delta$.
\item[---] Let $\underline v_s(x)$ have the form $s^{-1}\underline w_0(x/\sqrt s)$, with $w_0(\eta)=B^{-1}\eta e^{-A\eta}\Gamma(\eta/A)$, where $\Gamma$ is smooth nonnegative nonincreasing, equal to 1 on $[0,1]$ and 0 on $[2,+\infty)$. For the same $\delta>0$, we have, still from Corollary \ref{c4.3.1}: $e^{-t\II_*}\overline v_s(x)\geq0$ if $x\geq s^\delta$ and $e^{-t\II_*}\overline v_s(x)\leq0$ if $x\leq -s^\delta$.
\item[---] Assume $v_0$ to be supported in $(1,+\infty)$; we always may assume this at the expense of a  translation.  If $B>0$ is large enough, we have $\underline v_s(x)\leq v_0(x)\leq\overline v_s(x)$.  
\end{itemize}
\begin{lemma}
\label{l2.4.1}
Fix $\alpha\in(1/3,1/2)$. There  are two bounded positive functions $\xi_0^\pm(t)>0$, with $\xi_0^+$ nondecreasing and $\xi_0^-$ nonincreasing, bounded and bounded away from 0,   such that
the following holds.
\begin{itemize}
 \item[--] For  $x\geq {3\pi t^\alpha}/2$, and $t\geq1$, we have
\begin{equation}
\label{e2*.4.47}
\xi_0^-(t)e^{-(t-1)\II_*}\underline  v_s^*(x)\leq v(t,x)\leq\xi_0^+(t)e^{-t\II_*}\overline v_s^*(x).
\end{equation}
\item[--] Fix $\beta\in(0,1/2)$. If $-1\leq x\leq{3\pi t^\alpha}/2$, , and $t\geq1$, we have
\begin{equation}
\label{e2*.4.45}
\xi_0^-(t)e^{-(t-1)\II_*}\underline  v_s^*(x)+O(1/{t^{3/2-\beta}})\leq v(t,x)\leq\xi_0^+(t)e^{-t\II_*}\overline v_s^*(x)+O(1/{t^{3/2-\beta}}).
\end{equation}
\end{itemize}
\end{lemma}
\subsubsection  {The upper barrier}
\noindent In an ideal world, the solution $v(t,x)$ of \eqref{e4.2.8} would only be controlled from above by the sole function $e^{-t\II_*}\overline v_s^*(x)$, which is indeed a super-solution. Alas, we do not very well know what
 the value of $v(t,0)$ is, and, in particular, if it is actually below $e^{-t\II_*}v_0$ on $[0,1]$. Thus we have to compare it to the barrier slightly to the left, say, at a distance  $t^\delta$. There, $v$ is  exponentially small, thus controllable by $e^{-t\II_*}v_0$. But then, assuming that we have extended $v_0$ in an odd fashion, $e^{-t\II_*}v_0$ is negative, thus to be supplemented by an additional ingredient in the barrier. This is the role of the cosine function. 
 
\noindent So, consider $\delta\in(0,1/2)$, that will be as small as needed, we will estimate $v(t,x)$ ahead of $x=-t^\delta$. 
The idea is to construct a barrier function of the form:
\begin{equation}
\label{e2*.4.49}
 \bar  v(t,x)=\xi_0^+(t)e^{-t\II_*}\overline v_s^*(x)+\di\frac1{(1+t)^A}\cos(\frac{x}{(1+t)^\alpha})\un_{(-\infty,\frac{3\pi(1+t)^\alpha}2]}(x),
\end{equation}
where   $\xi_0^+$ as well as the positive parameters $\alpha\in(\delta,1/2)$  and $A>0$ are
to be adjusted in the course of the investigation, that will be done for $t>0$ large enough. We have
$$
\bigl(\partial_t+\II_*\bigl)\bar  v=\dot\xi_0^+e^{-t\II_*}\overline v_s^*(x)+\bigl(\partial_t+\II_*\bigl)\bigg(\un_{(-\infty,\frac{3\pi(1+t)^\alpha}2]}\bigl((1+t)^{-A}\cos(\frac{x}{(1+t)^\alpha})\bigl)\biggl),
$$
and we note that $(\partial_t+\II_*)\overline v=0$ in the zone $\{x\geq3\pi(t+1)^\alpha/2$. So, three zones are to be singled out. As said before, we have $e^{-t\II_*}\overline v_s^*(x)\geq0$ for all $t\geq0$ and $x\geq s^\delta$.

\noindent -- {\it The case $x_0\leq x\leq{\pi(1+t)^\alpha}/4$}. We have
$\di\cos\bigl({x}{/(1+t)^\alpha}\bigl)\geq{\sqrt 2}/2,$ and, for all $A>0$ and $\alpha\in(0,1/2)$ there is $t_{A,\alpha}>t_0$ such that we have, for $t\geq t_{A,\alpha}$:
$$
\frac{d_*}{(1+t)^{2\alpha}}\bigl(1+\omega_1((1+t)^{-\alpha})\bigl)-\frac{A}{1+t}\geq\frac{d_*}{2(1+t)^{2\alpha}},
$$
and, by possibly by enlarging $t_{A,\alpha}$ we have, from Proposition \ref{p4.3.4004}:
\begin{equation}
\label{e2*.4.56}
\bigl(\partial_t+\II_*\bigl)\bigl((1+t)^{-A}\cos(\frac{x}{(1+t)^\alpha})\bigl)\geq\frac{d_*}{2(1+t)^{2\alpha+A}}.
\end{equation}
Indeed we have ${x}/{(1+t)^{1+\alpha}}\lesssim1/t$. 
So, if we now request $\dot\xi_0^+(t)\geq0$, we win.

\noindent -- {\it The case $-t^\delta\leq x\leq x_0$}. This time, $e^{-t\II_*}$ can be negative, and we have to rely even more on the cosine perturbation. Theorem \ref{t2.4.2} teaches us that $e^{-t\II_*}\overline v_s^*$ is, at worst, distant from
$e^{t\partial_{xx}}\overline v_s^*$ by a quantity that  decays exponentially in time. Therefore, estimate \eqref{e2*.5.102} teaches us in turn that we have, at worst:
$$
\vert\dot\xi_0^+(t)e^{-t\II_*}\overline v_s^*(x)\vert\gtrsim{\dot\xi_0^+(t)}/{(1+t)^{3/2-\delta}},
$$
while we still have \eqref{e2*.4.56}.  Therefore we should impose
\begin{equation}
\label{e2*.4.57}
1/{(1+t)^{2\alpha+A}}\gg{\dot\xi_0^+(t)}/{(1+t)^{3/2-\delta}}.
\end{equation}

\noindent -- {\it The case ${\pi(1+t)^\alpha}/4\leq x\leq{3\pi(1+t)^\alpha}/2$}. In this area, the cosine perturbation is negative. However, we may use the full force of $e^{-t\II_*}\overline v_s^*$, which is like $x/t^{3/2}$, with
$x\sim(1+t)^\alpha$. So,
we require this time
\begin{equation}
\label{e2*.4.58}
{\dot\xi_0^+(t)}/{(1+t)^{3/2-\alpha}}\gg1/{(1+t)^{2\alpha+A}}.
\end{equation}

\noindent -- {\it Summary.} From \eqref{e2*.4.57} and \eqref{e2*.4.58} we should achieve
\begin{equation}
\label{e2*.4.70}
1/{(1+t)^{3\alpha+A-3/2}}\ll\dot\xi_0^+(t)\ll1/{(1+t)^{2\alpha+\delta+A-3/2}},
\end{equation}
which is always possible as soon as $\delta<\alpha$, something we have assumed. It remains to choose $A$, we devise it to have: for all $x\in[0,1]$,
$
\bar  v(t,-t^\delta+x)\geq v(t,-t^\delta+x),
$ the last quantity being of the order $e^{-t^\delta}$. Therefore it is enough to have the cosine perturbation dominate $e^{-t\II_*}\overline v_s^*(-t^\delta)$ by a small algebraic order. This is achieved if we choose
$A=3/2-\beta$, with $\delta<\beta<2\beta<\alpha$, and
\begin{equation}
\label{e2*.4.600}
\dot\xi_0^+(t)=1/{(1+t)^{3\alpha-2\beta}}.
\end{equation}
If $\delta\in(0,1/2)$ and $\beta<\alpha$ are small enough, we may choose $\alpha>1/3$ so that $3\alpha-2\beta>1$ and $\dot\xi_0^+$ is integrable.

\noindent With this, we not only have $\bar  v(t,x)\geq0$ everywhere, but there exists $t_0\geq0$ such that   $\bigl(\partial_t+\II_*\bigl)\bar  v(x)\geq0$ for $t_0>0$, $x\geq -t^\delta$.
The only item that remains to be checked is that $\bar  v(t_0,x)$ is above $v_0(x)$ for $x\geq -1$. This may not be the case; however,  we may put a large multiple of $\bar  v(t_0,x)$ above $v_0(x)$. 
So, we have proved the right handsides of \eqref{e2*.4.47} and \eqref{e2*.4.45}.

\subsubsection{The lower barrier} 
\noindent In an ideal world once again, we would just try to control $v(t,x)$ by some multiple of $e^{-t\II_*}\underline v_s^*$ from below. Indeed, we do not care about $e^{-t\II_*}\underline v_s^*$ to be negative, 
so we could think of repeating the above argument, with this time a boundary control near $x_0$. There is one case where the world is ideal, namely when the function $g$, entering in the composition \eqref{e2*.1.11},
is zero on a small range of parameters $(0,\theta_0)$.  Consider indeed the function
$
\phi_*(x)=x^+e^{-\lambda_*x}$, it is a steady linear wave solution, positive on $\RR_+$, zero on $\RR_-$. Now, if $A>0$ is large enough, the equation $A\phi_*(x)=1$ has two roots $0<x_A^-<x_A^+$, and their distance is larger than 1, even if it means enlarging $A$. Fix such an $A$, and define
$
\bar  u(x)=
A\phi_*(x)\ \hbox{if $x\geq x_A^+$,}$ and $\bar u(x)=1\ \hbox{if $x\leq x_A^+$}$.
As $x_A^+-x_A^->1$, Proposition \ref{p2.2.2} applies, and $\bar  u$ is a super-solution to equations \eqref{e4.2.8} or \eqref{e4.2.9} with $c=c_*$, that can additionally be put above $u_0$. Consider $x_s$ such that $\bar  u(x)\leq \theta_0$ if $x\geq x_s$,
we may, at the expense of a translation, assume $x_s=-s^\delta-1$. Therefore, as $u(t,x)\leq\bar  u(x)$, we have $g(u(t,x))\equiv0$ if $x\geq-s^\delta-1$. The equation for $v$ is therefore the pure diffusive equation; as $e^{-(t-1)\II_*}\underline v_s(x)\leq0$ if $x\in[-s^\delta-1,-s^\delta$ we have $e^{-(t-1)\II_*}\underline v_s(x)\leq v(t,x)$ if $x\geq-s^\delta$ and $t\geq1$.

\noindent However, the  nonlinear term is usually positive. So, if we   hope to control
$v(t,x)$ by $e^{-t\II_*}v_0^*(x)$, we would be well advised to do it in an area where the nonlinear term is negligible. Hence the idea to estimate $v(t,x)$ ahead of $t^\delta$.

\noindent In this area, we use the fact that $v(t,x)\leq\bar  v(t,x)$, where $\bar  v(t,x)$ is the constructed super-solution given by \eqref{e2*.4.49}. We use the rough estimate that $\bar  v(t,x)$ is bounded by a constant, so that we have:
\begin{equation}
\label{e2*.4.59}
\hbox{for $x\geq t^\delta$,} \ \ \ \ \ h(x,v(t,x))\lesssim e^{-t^\delta}:=\bar h(t).
\end{equation}
Therefore, any  solution $\underline v(t,x)$ of $\partial_t \underline v+\II_*\underline v+\bar h(t)\underline v\leq0,
$
on a domain of the form $\{t\geq t_0,\ x\geq t^\delta\}$ and that  satisfies
$
\underline v(t_0,x)\geq0\ \hbox{for $x\geq0$},$ $v(t,x)\geq\underline v(t,x)\ \hbox{for $t^\delta\leq x\leq t^{\delta}+1$},
$
will be a barrier for $v$. Notice that, as $\bar h(t)$ is an integrable function, the change of functions $\underline w(t,x)=e^{-\int_0^t\bar h(s)~ds}\underline v(t,x)$ allows the search of a subsolution for the sole
integral equation $\partial_tw+\II_*w=0$, in other words the computations of the preceding section remain valid.

\noindent And so, a function of the form
\begin{equation}
\label{e2*.4.60}
\underline v(t,x)=\xi_0^-(t)e^{-t\II_*}v_0^*(x)-\di\frac1{(1+t)^{3/2-\beta}}\cos(\frac{x}{(1+t)^\alpha})\un_{(-\infty,\frac{3\pi(1+t)^\alpha}2]}(x),
\end{equation}
with this time 
\begin{equation}
\label{e2*.4.61}
\dot\xi_0^-(t)=-1/{(1+t)^{3\alpha-2\beta}}.
\end{equation}
is a good candidate. There are, however, one small catch, which  is that $\xi_0^-(t)$ should be positive for large times, something that is not guaranteed by \eqref{e2*.4.61}. It is, however, sufficient to modify the definition  \eqref{e2*.4.60}
by 
$$
\underline v(t,x)=\xi_0^-(t)e^{-t\II_*}v_0^*(x)-\di\frac\e{(1+t)^{3/2-\beta}}\cos(\frac{x}{(1+t)^\alpha})\un_{(-\infty,\frac{3\pi(1+t)^\alpha}2]}(x),
$$
$\e>0$ small. The function $\xi_0^-(t)$  is characterised by $\xi_0^-(0)=1$, and $\dot\xi_0^-$ given by \eqref{e2*.4.61}. 
These modifications ensure the lower barrier property.
\subsection{Retrieving the information to $x=O(1)$}\label{s2.5.3}
\noindent The situation is the following: around $x=0$, the solution $v(t,x)$ of \eqref{e4.2.9} is essentially of the order $1/{t^{3/2}}$. While this is an interesting information, this does not answer the question we are really interested in,
 that is, the location of a level set of $u(t,x)$ of the original equation \eqref{e4.1.1}. Obviously, trying to fish for a level set of $u(t,x)$ of fixed value at this place is bound to fail, and we have to look a little back in order to see 
 nontrivial values for $u(t,x)$. This is exactly the rationale of the logarithmic delay.

\noindent In order to retrieve the information to the back, we need a vehicle, and we have a beautiful one: the travelling wave solution. Let us recall that $\varphi_{c_*}$ is the wave with bottom speed. We could normalise it once and for all to $\varphi_{c_*}(0)=1/2$, but there is something smarter to do here. From Theorem \ref{t2.3.1}, we may normalise it according to its behaviour at $+\infty$, and ask:
\begin{equation}
\label{e2*.5.200}
\varphi_{c_*}(x)\sim_{x\to+\infty} (x+k_*)e^{-\lambda_*x}.
\end{equation}
We are sure that it is the only one: for any translate of $\varphi_{c_*}$, the  exponential  would find itself multiplied by a constant different from 1.  This being set once and for all, we are now going to look for 
two functions $x_-(t)<x_+(t)$ such that
\begin{equation}
\label{e2*.5.201}
\varphi_{c_*}(x-x_-(t))\leq u(t,x)\leq \varphi_{c_*}(x-x_+(t))\ \ \hbox{for $-t^\delta\leq x\leq t^\delta$.}
\end{equation}
As $x_\pm$ will be much smaller than $t^\delta$, this will locate the nontrivial level sets of $u$.

\noindent Consider $\delta>0$ small, $\beta\in(0,\delta)$ and $\alpha>1/3$, close to $1/3$. We start from the following consequence of Theorem \ref{t2.4.5}:
\begin{equation}
\label{e2*.5.203}
u(t,x)\leq {Cxe^{-\lambda_*x}}/{t^{3/2}}\ \ \hbox{for $x\in[t^\delta,t^\delta+1]$.}
\end{equation}
Let us devise $x_+(t)$ so that $\varphi_{c_*}(x-x_+(t))\geq u(t,x)$ for $x\in[t^\delta,t^\delta+1]$. Given the equivalent \eqref{e2*.5.201}, it is sufficient to impose
$$
\bigl(x-x_+(t)+k_*\bigl)e^{-\lambda_*(x-x_+(t))}\geq{Cxe^{-\lambda_*x}}/{t^{3/2}}\ \ \hbox{for $x\in[t^\delta,t^\delta+1]$,}
$$
for a constant $C$ that is possibly different from the one in \eqref{e2*.5.203}. An elementary computation yields
\begin{equation}
\label{e2*.5.204}
x_+(t)\geq-\frac3{2\lambda_*}\mathrm{ln}~t+K_+,
\end{equation}
with $K_+>0$ large enough. Similarly, a sufficient condition for a function $x_-(t)$ to satisfy the left handside of \eqref{e2*.5.201} is
\begin{equation}
\label{e2*.5.205}
x_-(t)\leq-\frac3{2\lambda_*}\mathrm{ln}~t-K_-,
\end{equation}
with, again, $K_->0$ large enough. The functions $x_\pm(t)$ are now fixed as in \eqref{e2*.5.204} and \eqref{e2*.5.205}, with equalities in the place of inequalities, and this is of course where we are going to look for the level sets of $u(t,x)$. The common point between $x_+$ and $x_-$ being the $\mathrm{ln}~t$ term, we make the change of reference frame
$
y=x+\mathrm{ln}~3t/{2\lambda_*}.
$
We rename $x$ the variable $y$; in the new reference frame the equation for $u$ is
\begin{equation}
\label{e2*.5.600}
u_t+\frac3{2\lambda_*t}u_x+\mathcal{L}_*u-u+g(u)=0,
\end{equation}
so that the equation for $v(t,x)=e^{\lambda_*x}u(t,x)$
is
$$
v_t+\frac3{2\lambda_*t}(v_x-v)+\mathcal{I}_*v+h(x,v)=0.
$$
We also set
$
\psi_*(x)=e^{\lambda_*x}\varphi_{c_*}(x),
$
and we study the differences
$w_\pm(t,x)=v(t,x)-\psi_*(x-K_\pm).
$
The domain of study would be, strictly speaking:
$
\{-t^\delta+\di\frac3{2\lambda_*}\mathrm{ln}t\leq x\leq t^\delta+\di\frac3{2\lambda_*}\mathrm{ln}t\}.
$
As $t$ is supposed to be large, we may keep the domain as $\{-t^\delta \leq x\leq t^\delta\}.$

\noindent We are going to prove that 
\begin{equation}
\label{e2*.5.215}
\limsup_{\substack {t\to+\infty\\
-t^\delta\leq x\leq t^\delta}}w_+(t,x)\leq0,\ \ \limsup_{\substack {t\to+\infty\\
-t^\delta\leq x\leq t^\delta}}w_-(t,x)\geq0.
\end{equation}
The equation for $w_\pm(t,x)$ is 
$
\partial_tw_\pm+\di\frac3{2\lambda_*t}(\partial_xw_\pm-w_\pm)+\mathcal{I}_*w_\pm+a(t,x)w_\pm=0,
$
with as usual,
$
a(t,x)={e^{\lambda_*x}\biggl(g\bigl(e^{-\lambda_*x}v(t,x)\bigl)-g\bigl(e^{-\lambda_*x}\psi_*(x-K_\pm)\bigl)\biggl)}\bigl/w_\pm(t,x).
$
Note that $a(t,x)\geq0$ because $g'\geq0$. Let us, for definiteness, examine the large time behaviour of $w_+(t,x)$. The right handside of \eqref{e2*.5.215} is equivalent to
$\di\lim_{\substack {t\to+\infty\\
-t^\delta\leq x\leq t^\delta}}w_+^+(t,x)=0,
$
where $w_+^+(t,x)$ is the positive part of $w_+(t,x)$. Let us prove it, we will in fact see that the above limit is uniform in $x\in[-t^\delta,t^\delta]$.
Using the analogue \eqref{e2.5.108} of Kato's inequality,  we easily derive an inequation for $\vert w_+^+(t,x)\vert$:
\begin{equation}
\label{e2*.5.216}
\begin{array}{rll}
&\begin{array}{rll}
\partial_tw_+^++\di\frac3{2\lambda_*t}(\partial_xw_+^+-w_+^+)+\mathcal{I}_*w_+^+
\leq&\di\frac3{2\lambda_*t}(\psi_*'+\psi_*)(x-K_+)\\
\lesssim&\di\frac1{t^{1-\delta}},
\end{array}\ \ \ \ (-t^\delta\leq x\leq t^\delta)\\
&w_+^+(t,x)=O(e^{-t^\delta})\ (-t^\delta\leq x\leq -t^\delta+1),\quad
w_+^+(t,x)=0,\ \ \ \ (-t^\delta\leq x\leq t^\delta).
\end{array}
\end{equation}
Note that the last line of \eqref{e2*.5.216} is perhaps the most crucial one: the region $\{t^\delta-1\leq x\leq t^\delta\}$ is indeed the zone where the long range behaviour of $v(t,x)$ and that where the finite range behaviour of $v(t,x)$ communicate.
A barrier for $w^+_+$ will be devised using, once again, the computation toolbox of Section \ref{s2.2.1} will be useful. We pick again $\alpha>1/3$, close to 
$1/3$, and we keep in mind that $\delta$ is small. For $t\geq1$ we set
$
\bar  w(t,x)=\cos\bigl({x}/{t^\alpha}\bigl)/t^A;
$
from Proposition \ref{p4.3.4004} again we have, for $t\geq1$ and $-t^\delta\leq x\leq t^\delta$:
$$
\bigl(\partial_t+\II_*\bigl)\bar  w(t,x)\gtrsim1/{t^{A+2\alpha}},
$$
while
$$
\bigl\vert\partial_x\bar  w(t,x)-\bar  w(t,x)\bigl\vert/t\lesssim1/{t^{A+\alpha+1}}.
$$
Therefore, $\bar  w(t,x)$ is a barrier as soon as $1/{t^{A+2\alpha}}$ dominates the right handside of \eqref{e2*.5.216}, that is, $1/{t^{1-\delta}}$. This is achieved as soon as $\delta$ is small and $\alpha$ slightly larger than $1/3$: in this setting there is a room for finding a suitable $A$; here, any $A<1-2\alpha-\delta$ will work. The barrier $\bar  w(t,x)$ tending to 0 as $t\to+\infty$, we may put a large multiple of $\bar  w(t,x)$ above $w^+_+(t,x)$, and this proves the first part of\eqref{e2*.5.215}. As similar considerations can be made for $w_-^-(t,x)$, we have proved \eqref{e2*.5.215}. This in turn proves the main result of the whole section.

\section{Gradient bounds}\label{s4.10}

\noindent As already mentionned, there is no smoothing mechanism in \eqref{e4.1.1}, so that we have to rely on the conservation of initial smoothness and explore a mechanism that prevents the inflation of its derivatives. 
Yet we have the 
\begin{theorem}
\label{t2.4.6}
Assume $u_0\in C^\infty(\RR)$. For all $m\in\mathbb{N^*}$, the quantity $\Vert \partial_x^mu(t,.)\Vert_{L^\infty(\RR)}$ is bounded.
\end{theorem}
The starting point is equation  for $u(t,x)$ \eqref{e2*.5.600} in the reference frame moving like $c_*t-3\mathrm{ln}~t/{2\lambda_*}$. The derivative $u_x(t,x)$ solves
\begin{equation}
\label{e2*.5.6010}
\partial_tu_x+u_x-K*u_x-(c_*-\frac3{2\lambda_*t})\partial_xu_x=u_x-g'(u)u_x.
\end{equation}
We define the function $v(t,x)$ as, this time, $u_x(t,x)=e^{-\lambda_*x}v(t,x);$
the function $\vert v(t,x)\vert$ solves
\begin{equation}
\label{e2*.5.602}
\partial_t\vert v\vert+\II_*\vert v\vert-\frac3{2t}(\partial_x\vert v\vert-\vert v\vert)+g'(u)\vert v\vert=0.
\end{equation}
We always may arrange the picture to be translated by a constant so that there is $q_*>0$ so that
\begin{equation}
\label{e2*.5.601}
g'\bigl(u(t,x)\bigl)\geq q_*\ \ \ \hbox{for $x\leq0$.}
\end{equation}
We are going to construct an upper barrier $\bar  v(t,x)$ for $\vert v(t,x)\vert$ for large times; as $\vert v(t,x)\vert$ is bonded for finite times by the sole virtue of the Gronwall  lemma, this will provide the sought for barrier.

\noindent The barrier function $\bar  v(t,x)$ will be the infimum of two super-solutions to \eqref{e2*.5.602}, that we call $\bar  v^+(t,x)$ and $\bar  v^-(t,x)$.

\subsubsection{Construction of $\bar  v^+$} 
\noindent It will be devised according to the pattern that has already proved to be useful in Section \ref{s2.5.2}, up to the fact that the new definition of $v(t,x)$ means the suppression of the $t^{-3/2}$ factor. Let $\Gamma(\eta)$ be a smooth nonnegative nondecreasing function that is zero on $[0,1]$ and equal to 1 on  an interval of the form $[\eta_0,+\infty)$. We set
\begin{equation}
\label{e2*.5.801}
\bar  v^+(t,x)=\xi(t)\bigl(t^{3/2}e^{-t\II_*}\bar  v_0^*\bigl)+q(t)\Gamma(\eta)e^{-b\eta}.
\end{equation}
We have denoted  $\eta={x}/{\sqrt t},$
in order to alleviate the notations, and the functions $\bar  v_0$, $\xi(t)$ and $q(t)$, as well as the constant $b>0$, are to be properly chosen. Unsurprisingly we will
impose $\dot\xi\geq0$. We have denoted by $\bar  v_0^*$ the odd extension of $\bar  v_0$ to the left of $x=0$. In order to alleviate the notations somehow we set
$S_*(t,x)=t^{3/2}e^{-t\II_*}\bar  v_0^*(x).$
We may now compute
$$
\begin{array}{rll}
&\di\biggl(\partial_t+\II_*+\frac3{2t}(\partial_x-1)\biggl)\bar  v^+(t,x)=\dot\xi(t)S_*+\dot q\Gamma e^{-b\eta}+q\Gamma(\partial_t+\II_*)e^{-b\eta}\\
+&\di\frac3{2t}\xi\partial_xS_*-q\Gamma'\di\frac{x}{2t^{3/2}}e^{-b\eta}+qe^{-b\eta}\II_*\Gamma(\eta)+\di\frac{3q}{2t}\partial_x\bigl(\Gamma(\eta)e^{-b\eta}\bigl)-\di\frac{3q\Gamma e^{-b\eta}}{2t}\\
-&q\di\int_{\RR}K_*(x-y)\bigl(\Gamma(\di\frac{x}{\sqrt t})-\Gamma(\di\frac{y}{\sqrt t})\bigl)\bigl(e^{-bx/\sqrt t}-e^{-by/\sqrt t}\bigl)~dy.
\end{array}
$$
This makes a string of nine terms, that we number from $T_1$ to $T_9$. In order to reach the super-solution property, we will prove that $T_1$ to $T_3$ prevail over the others, that will be treated as perturbations.
Fortunately, they are relatively easy to examine. As is now usual, three zones are to be investigated.

\noindent{\bf 1. The area $\eta\geq\eta_0+1$.} As  $\Gamma\equiv1$ here, we have, provided $\eta_0$ is large enough:
$$
T_1+T_2+T_3\geq\bigl(\dot q+qe^{b\eta}(\partial_t\II_*)e^{-b\eta}\bigl)e^{-b\eta}
\geq(\dot q+b\eta q/t)e^{-b\eta}.
$$
On the other hand, as $\Gamma'\equiv0$ we only have to worry about $T_4$, $T_6$, $T_7$, $T_8$. We have
$$
\vert T_4\vert\lesssim\vert \partial_xS_*(t,x)/t\vert\lesssim {x^2}/te^{-x^2/4d_*t}\lesssim{\eta^2 e^{-\eta^2/4d_*}} \ll {e^{-b\eta}}/t,
$$
as soon as $\eta_0$ is large enough. As for the remaining terms, we have
$
\vert T_6\vert+\vert T_7\vert+\vert T_8\vert\lesssim{qe^{-b\eta}}/t.
$
Therefore this last sum is absorbed by ${b\eta}q/t$, so that the function $q(t)$ should, in the end, satisfy 
$
\dot q+{b\eta_0}q/t\geq1/t,
$
that is, 
$q(t)\simeq1/{b\eta_0}$, $\dot q\simeq1/{t^{1+b\eta_0}}.$

\noindent{\bf 2. The range $\e_0\leq\eta\leq\eta_0+1$.} Now that we have chosen the function $q(t)$, we need to choose $\xi(t)$ and we do it here. this time the term that will 
help us will be $T_1=\dot\xi S_*(t,x)$. We have indeed, still taking for granted that $\dot\xi\geq0$: $T_1\gtrsim\dot\xi x\gtrsim\sqrt t\dot\xi,$ and $\vert T_4\vert\lesssim\xi/{2t},$
while all he other terms are estimated by $1/t$. Let us, for instance, examine $T_9$. We have  
$$
\int_\RR K_*(x-y)\biggl\vert\bigl(\Gamma(\frac{x}{\sqrt t})-\Gamma(\frac{y}{\sqrt t})\bigl)\bigl(e^{-bx/\sqrt t}-e^{-by/\sqrt t}\bigl)\biggl\vert\lesssim\int_{\RR}K_*(x-y)\frac{(x-y)^2}tdy\lesssim\frac1t.
$$
And so, it is sufficient to impose $\sqrt t\dot\xi\geq(1+\xi)/t$,
that is, $\xi(t)=O(1+1/{\sqrt t})$. Let us therefore choose $q(t)$ and $\xi(t)$ in this manner.

\noindent{\bf 3. The range $0\leq\eta\leq\e_0$.} This time we have $\Gamma=\Gamma'\equiv0$, so that only $T_1$ and $T_4$ matter. In this range, however, we have $\partial_xS_*(t,x)\geq0$; as $\dot\xi\geq0$ we still have the super-solution property.

\noindent Note that our computations may not be, strictly speaking, valid up to $x=0$. However they certainly hold for $x$ larger than a suitably large constant, which will be sufficient for the sequel. To sum up,  there is $x_0>0$ and $t_0>0$ such that the function $\bar  v^+(t,x)$ given by \eqref{e2*.5.801} solves
\begin{equation}
\label{e2*.5.802}
\biggl(\partial_t+\II_*+\frac3{2t}(\partial_x-1)\biggl)\bar  v^+(t,x)\geq0\ \hbox{for $t\geq t_0$ and $x\geq x_0$.}
\end{equation}
In other words, $\bar  v^+(t,x)$ is a super-solution to \eqref{e2*.5.602}.
\subsubsection{Construction of $\bar  v^-$} 
\noindent As $\bar  v^+(t,x)$ becomes negative for negative $x$, we need another ingredient. What we have gained is that we know with $O(1)$ precision where the nontrivial level sets of $u(t,x)$ are, that is, they are at finite distance from 0. 
In particular,  a sign of this fact is equation \eqref{e2*.5.601}, that we did not know before proving Theorem \ref{t4.4.4}, which opens new possibilities. For small $\omega>0$, let us repeat the following estimate, which is by now usual:
$
\II_*e^{\omega x}\gtrsim-\omega^2e^{\omega x}.
$
This implies
$$
\biggl(\partial_t+\II_*+\frac3{2t}(\partial_x-1)+g'(u)\biggl)e^{\omega x}\geq \biggl(q_*+O(\omega^2)+O(\di\frac1t)\biggl)e^{\omega x}.
$$
So, there is $t_1>0$ such that 
$\bar  v^-(x)=e^{\omega x},
$
is a super-solution to \eqref{e2*.5.602} on $\{t\geq t_1,\ x\leq-1\}$.
\subsubsection{Construction of the super-solution}
\noindent As equation \eqref{e2*.5.802} is linear, we may work with any multiples of $\bar  v^+$ or $\bar  v^-$. One may wonder how one can glue $\bar  v^+$ and $\bar  v^-$ together so as to make a global super-solution, the remark that will
make everything work is that $\bar  v^+$ is a super-solution not only of \eqref{e2*.5.802} on $[t_0,+\infty)\times[x_0,+\infty)$, but also of the equation without the term $g'(u)v$, in other words a translation invariant equation.
Thus, for every $a>0$, the function $\bar  v^+(t,x+a)$ is also a super-solution to \eqref{e2*.5.802} on $[t_0,+\infty)\times[x_0-a,+\infty)$. Moreover, we have by construction:
$
\bar  v^+(t,x)\simeq (x+a)e^{-x^2/4d_*t}+O(1)$ for $t\geq t_0$ and $x\geq x_0-a$,
the $O(1)$ term being independent of $a$.  There is $\alpha_\infty>0$ such that 
$
\di\lim_{t\to+\infty}\bar  v^+(t,x+a)=\alpha_\infty(x+a)
$
uniformly on every interval of the form $[x_0-a,t^\delta]$ with $\delta\in(0,1/2)$.  Consider the function $a\alpha_\infty\bar  v^-(x)$; it is larger than $\alpha_\infty(x+a)$ on $\RR_+$ and smaller on an interval $(-x_a,0]$ with $\lim_{a\to+\infty}x-a=-\infty$. Consequently, for $t>0$ large enough, there is a function $\bar x(t)$ tending to 0 as $t\to+\infty$ such that $\bar  v^+(t,x)$ is below $a\alpha_\infty\bar  v^-(x)$ on $\bar x(t)-1,\bar x(t)]$, and above for $x\geq x(t)$. Thus, the function 
\begin{equation}
\label{e2*.5.809}
\bar  v(t,x)=
\left\{
\begin{array}{rll}
&\bar  v^+(t,x)\ \hbox{if $x\geq \bar x(t)$}\\
&a\alpha_\infty\bar  v^-(x)\ \hbox{if $x\leq \bar x(t).$}
\end{array}
\right.
\end{equation}
It remains to pick any compactly supported, nonnegative, nonzero function $\bar  v_0(x)$; a large enough  multiple of $\bar  v(t_0,.)$  will dominate the initial datum $v_0$. 

\noindent {\sc Proof of Theorem \ref{t2.4.6}.} All the previous development  argument proves the global boundedness of $v(t,x)$, hence the boundedness of $u_x(t,x)$ on every half line containing $\RR_+$. As the equation is essentially linear on $\RR_+$, it also proves bounds for all spatial derivatives of $v$, henceforth all spatial derivatives of $u$ on every positive half line.

\noindent The function $u_x(t,x)$ is then easily bounded on every negative half line: indeed, from equation \eqref{e2*.5.6010} we have
\begin{equation}
\label{e2*.5.6011}
\biggl(\partial_t+\mathcal{J}-c_*\partial_x+q_*\biggl)\vert u_x\vert\leq0\ \ \hbox{for $x\leq0$}.
\end{equation}
Let $M$ be an upper bound for $u_x$ for $t\geq0$ and $x\geq-1$. Then 
$
\bar  u(t)=M+e^{-q*t}
$
is a super-solution to \eqref{e2*.5.6011} that dominates $\vert u_x\vert$ for $t=0$ and all $x\leq-1$, as well as for $t\geq 0$ and $x\in[-1,0]$. Therefore it dominates $u$ for $t\geq0$ and $x\leq0$.  The successive bounds are proved by induction. \hfill$\Box$
\section{Large time convergence}\label{s4.4}
\noindent We now have all the elements for the precise study of the solution $u(t,x)$ of the Cauchy Problem \eqref{e4.1.1}, starting from a nonnegative, nontrivial, compactly supported initial datum $u_0(x)$.   
\begin{theorem}
\label{t2.4.1}
There exists $x_\infty\in\RR$ such that  $\di\lim_{t\to+\infty}\biggl(u(t,x)-\varphi_{c_*}\bigl(x-c_*t-\frac3{2\lambda_*}\mathrm{ln}~t+x_\infty\bigl)\biggl)=0,$
uniformly in $x\in\RR_+$. A similar statement holds on $\RR_-$, with a possibly different $x_\infty$.
\end{theorem}
Once Theorem \ref{t2.4.1} is proved, the access to Theorem \ref{t2.1.1} is easy, as $x_\infty(\theta)+x_\infty=\varphi_{c_*}^{-1}(\theta)$.
\subsection{Analysis in the region $x\sim t^\gamma$}\label{s4.5.1}
\noindent Given how we have proved the weak version of the logarithmic delay, Theorem \ref{t4.4.4}, it is clear that the first objective that we should pursue is to obtain a refinment of Theorem \ref{t2.4.5} on the slope of $v(t,x)$ 
in the sub-diffusive region. This is the object of the next section, where we prove that we have an actual equivalent of $v(t,t^\gamma)$ for $\gamma\in(0,1/2]$. This is where the main part of the effort is cast. 
\begin{theorem}
\label{t2.4.3}
There exists $p_\infty>0$ such that, for all $\gamma\in(0,1/2)$ we have
\begin{equation}
\label{e2*.4.81}
v(t,t^\gamma)={p_\infty}/{t^{3/2-\gamma}}\bigl(1+o_{t\to+\infty}(1)\bigl).
\end{equation}
\end{theorem}
Once it is proved, retrieval of 
information to $x=O(1)$ will essentially the line of the preceding section \ref{s2.5.3}.

\noindent The proof of Theorem \ref{t2.4.3} is similar in spirit to that of Theorem \ref{t2.4.5}, the difference being that, instead of letting the process run by itself from $t=0$ onwards, we interrupt it at arbitrary large times in order to estimate it more precisely by solutions of the linear equation. The validity of the barriers that we construct will rely on slope estimates performed at various times, that we detail next. With this in hand, we may proceed to the construction of the barriers, leading to the proof of Theorem \ref{t2.4.3}.

\noindent The sought for sub and super solutions, that we denote by $v_s(t,x)$ for the moment, will be defined for $s>1$ and $t\geq s$. We want to look for them under the form
\begin{equation}
\label{e2*.4.68}
 v_s^\pm(t,x)=\xi_s^\pm(t)e^{-(t-s)\II_*}{\cal S}_{*\e}^\pm v(s,.)\pm\di\frac1{t^{3/2-\beta}}\cos(\frac{x}{t^\alpha})\un_{(-1,\frac{3\pi t^\alpha}2]}(x),
\end{equation}
where ${\cal S}_{*\e}^\pm v(s,.)$ is the odd extension around $\pm\e\sqrt s$ of  a modification of $v(s,.)$.
The important imposed feature is here that $\xi_s^\pm(s)=O(1)$,  and $\dot\xi_s^+\geq0$, $\dot\xi^-_s(t)\leq0$.  In other words, we update the observation at time $s$ by replacing $v_0^*$ by $v^*(s,.)$ but the perturbations $\xi_s^+$ and $q_s^+$ keep the memory of the past.
 In all that follows, the time $s$ will be assumed to be as large as needed. 
 
 \noindent We will arrange that $v_s^+(s,x)\geq v(s,x)$ if $x\geq0$. So, if we manage to prove 
\begin{equation}
\label{e2*.4.87}
\partial_t v_s^++\II_*v_s^+\geq0,\ \  \partial_t v_s^-+\II_* v_s^-\leq0\ \ \ \ \ \ (t>s,\ x>0),
\end{equation}
and 
\begin{equation}
\label{e2*.4.870}
v(t,x)\leq  v_s^+(t,x),\ \ \ v(t,x)\geq v_s^-(t,x)\ \ \ \ \ \hbox{for $t\geq s$ and $x\in[-1,0]$,}
\end{equation}
we will have upper and lower bounds for $v(t,x)$ in the range $t\geq s$.

\noindent We will essentially need to
prove that $v^\pm$ dominates $v$ at the boundary of our domain, in other words, that $e^{(t-s)\II_*}$ is not too big in comparison to the perturbation function at that 
place. The other property is, from the experience of the proof of Lemma \ref{l2.4.1}, is an estimate of the slope of $e^{-(t-s)\II_*}v^*(s,.)$ from below, which will in turn allow
us to choose the modulation function $\xi_s^\pm$. So, our main building block will be the set of estimates that have been stated separately in Theorem \ref{t4.3.3000} above.

%
%
\subsubsection{Barrier property  and convergence}\label{s2.6.2}
\noindent If now the functions $v^\pm_s(t,x)$, defined for $t\geq s$ by \eqref{e2*.4.68}, and if the functions $\xi^\pm_s(t)$ satisfy, as in the proof of the rough bounds:
\begin{equation}
\label{e2*.6.10}
\dot\xi_s^\pm(t)\approx1/{(t+1)^{3\alpha-2\beta}},
\end{equation}
with $\alpha>1/3$, $\beta<1$ small enough so that $3\alpha-2\beta>1$, we wish to  reproducef the computations in Lemma \ref{l2.4.1}, and  arrive at the integro-differential inequalities 
\eqref{e2*.4.87}. 

\noindent However, the ordering \eqref{e2*.4.870} between $v$ and $v^\pm_s$ for $x\in[-1,0]$ is not necessarily true: as it is certainly true for $t\gg s$, given the asymptotics of $e^{-(t-s)\II_*}$, there is still a large time frame, of the order, say, $s^\delta$, where this is not so obvious. One could indeed  imagine a scenario where the cosine perturbation would not be large enough to counterbalance $e^{-(t-s)\II_*}v^*(s,x)$ for $x$ bounded and $t\leq s+s^\delta$; this is not 
forbidden by the asymptotics concerning the heat semigroup with well spread data. This is where we need a modification of $v(s,.)$ in an $\e\sqrt s$ neighbourhood of 0,  in order to neutralise possible spurious effects of the semigroup generated by $\II_*$.

\noindent The modification of the definition of $v_s^\pm$ is thus conducted as follows. Choose $\e>0$ small, and $s>0$ such that $\e\gg1/s$.  A first ingredient is the
\begin{lemma}
\label{l2.6.100}
There is a function $w_s(\eta)$, locally bounded on $\RR_+$, such that
\begin{equation}
\label{e2*.6.500}
(\eta-s^{\delta-1/2})\un_{[0,{s^\delta]}}(\eta)\lesssim w_s(\eta)\lesssim(\eta+s^{\delta-1/2})e^{-A\eta},
\end{equation} 
and such that 
\begin{equation}
\label{e2*.6.600}
v_s(x)=\frac1sw_s(\frac{x}{\sqrt s}).
\end{equation}
 \end{lemma}
\noindent{\sc Proof.} The left handside of \eqref{e2*.6.500} is straightforward in view of Lemma \ref{l2.4.1} for the large time estimates of $v(t,x)$, and Theorem \ref{t2.4.2} for the behaviour of the heat kernel. As for the right handside, the exponential factor needs a justification. It stems from two points. The first one is   Lemma \ref{l2.4.1}, which ensures that $v(s,.)$ has the correct order of magnitude ${x}/{s^{3/2}}$ at a distance from 0 larger than $t^\alpha$, and less than a large multiple of $\sqrt s$. The second is that Proposition \ref{p2.4.1} is applicable to $v$, which ensures that it decays like an exponential for $x\gg\sqrt s$. This proves the lemma. \hfill$\Box$
 
\noindent At this point, we are guaranteed  that, for large enough $s>0$, identities \eqref{e2*.6.5} and \eqref{e2*.6.6} hold for $v(s,.)$. Let us now define the operators $\mathcal{S}_{\e,*}^\pm$, and let us start with 
$\mathcal{S}_{\e,*}^+$. What we look for is a function that coincides with $v(s,.)$ almost down to $x=0$, that is, down to $x=\e\sqrt s$, that well exceeds $v$ around $x=0$, and that we may oddly extend at our 
leisure so that, when we solve the linear Cauchy Problem starting from the odd extension, we retrieve something that is still above $v(t,x)$ for $t\geq s$ and $x$ close to 0. A few minutes of deep thought should be enough to
convince us that this is not an unfeasable task. Let us recall here the notation \eqref{e2*.6.600}, we know that $w(\xi)$ is close to $s^\delta$ in the vicinity of $\xi=0$ and controlled from above and below by positive multiples of 
$\xi$, for $\xi\geq0$. Let $K$ and $1/{K}$ be these multiples. We define $\overline w_s(\xi)$ as
$$
\overline w_s(\xi)=\left\{
\begin{array}{rll}
&w_s(\xi)\ \hbox{if $\xi\geq\e$}\\
&w_s(\e)+(\xi+\e)/K\ \hbox{if $\xi\leq\e$}.
\end{array}
\right.
$$
As $w_s(\e)\geq\e/{K}$, and as $v(s,\xi)$ is no more than $1/{s^{3/2}}$ in an $O(1)$ vicinity of $x=0$, we  have $\overline v_s(\xi)\geq v(s,\xi)$ for all $\xi\in(-\e/2,+\infty)$. Taking the odd extension of
$\overline w_s(\xi)$ around $\xi=-\e$ provides the sought for $\mathcal{S}_{\e,*}^+v(.,x)$. The quantity $\mathcal{S}_{\e,*}^-v(.,x)$ is obtained in a similar fashion, we just force it to be 0 at $x=s^\delta$, $0<\delta\ll1$, instead of 
$-\e\sqrt s$.

\noindent We then use Corollary \ref{c4.3.1} to infer that
$
e^{-(t-s)\II_*}\mathcal{S}_{*,\e}^-v(s,x)\leq0\leq e^{-(t-s)\II_*}\mathcal{S}_{*,\e}^+v(s,x)$ for $t\geq s$ and $x\in[-1,0]$.
On the other hand,  Lemma \ref{l2.4.1}, and more precisely inequality \eqref{e2*.4.45},
  estimates $v(t,x)$ in the vicinity of 0. It means, in particular, that $v(t,x)$ is dominated by $1/{t^{3/2-\beta}}$, for any small $\beta>0$. We therefore modify slightly the definition of $\xi^\pm_s$ by setting:
  \begin{equation}
\label{e2*.6.8}
 v_s^\pm(t,x)=\xi_s^\pm(t)e^{-(t-s)\II_*}\mathcal{S}_{*,\e}^\pm v(s,x)\pm\di\frac1{t^{3/2-2\beta}}\cos(\frac{x}{t^\alpha})\un_{(-1,\frac{3\pi t^\alpha}2]}(x),
\end{equation}
with
$
\dot\xi_s^\pm(t)\approx1/{(t+1)^{3(\alpha-\beta)}},
$
and omit the dependence in $\e$ in order to avoid an inflation of indices.
The real number  $\beta$ is small enough so that $3(\alpha-\beta)>1$. 
We have, by construction: $v_s^\pm(t,x)\leq v(t,x)\leq v^\pm_s(t,x)$ for $t\geq s$ and $x\in[-1,0]$. 
 Then, the  definition \eqref{e2*.6.8} of $v_s^\pm$ implies:
$$
v_s^-(t,x)\leq v(t,x)\leq v_s^+(t,x)\ \ \ \hbox{for $t\geq s$ and $x\in[-1,0]$,}
$$
and we finally have: $v_s^-(t,x)\leq v(t,x)\leq v_s^+(t,x)$ for $t\geq s$ and $x>0$.
It remains to analyse the large time behaviour of $v_s^\pm(t,x)$. The behaviour of the cosine perturbation being clear, let us analyse the main term $\xi_s^\pm(t)e^{-(t-s)\II_*}\mathcal{S}_{*,\e}^\pm v(s,x)$.
Recall that $\xi_s^\pm(t)$ is chosen, according to \eqref{e2*.4.600}, as $\dot\xi^\pm_s(t)=1/t^{3\alpha-2\beta}$. 
 This entails
$$
\xi_s^\pm(t)=1+O\bigl(1/{s^{3\alpha-2\beta-1}}\bigl).
$$
On the other hand Theorem \ref{t2.4.2} that $e^{-(t-s)\II_*}$ can, at no real cost, be replaced by $e^{(t-s)d_*\partial_{xx}}$. Let us set $w_{s,\e}^\pm(\eta)=\mathcal{S}_{*,\e}^\pm v(s,\eta\sqrt s)$, and notice that
\begin{equation}
\label{e2*.6.12}
\int_0^{+\infty}\eta w^\pm_{s,\e}(\eta)~d\eta=\int_0^{+\infty}\eta w_{s}(\eta)~d\eta+O(\e).
\end{equation}
We denote
$
p_s=\di\frac1{2\sqrt\pi d_*^{3/2}}\int_0^{+\infty}\eta w_s(\eta)~d\eta.
$
From equations \eqref{e2*.5.106} and\eqref{e2*.6.12},
we have, in the triple limit $x/\sqrt t\to0$, $x/t^\alpha\to+\infty$ and $t\to+\infty$:
\begin{equation}
\label{e2*.6.9}
v_s^\pm(t,x)\sim\frac{\bigl(1+o_{s\to+\infty}(1)\bigl)x}{2\sqrt\pi d_*^{3/2}t^{3/2}}\int_0^{+\infty}\eta w_{s,\e}^\pm(\eta)~d\eta,
\end{equation}
Theorem \ref{t2.4.3}  now becomes an elementary topology exercise. We have, for all $\gamma\in(\beta,1/2)$:
$$
\bigl(1+o_{\substack{s\to+\infty\\
\e\to0}}(1)\bigl)p_s\leq\liminf t^{3/2-\gamma}v(t,t^\gamma)\leq\limsup_{t\to+\infty}t^{3/2-\gamma}v(t,t^\gamma)\leq\bigl(1+o_{\substack{s\to+\infty\\
\e\to0}}(1)\bigl)p_s.
$$
In other words, all asymptotic  values (that is, all limits of subsequences) of $t^{3/2-\gamma}v(t,t^\gamma)$ lie in an interval centred at $p_s$ and of width $1+o_{\substack{s\to+\infty\\
\e\to0}}(1)$. As the function $s\mapsto p_s$ is bounded due to Theorem \ref{t2.4.5}, we may let $s\to+\infty$ and $\e\to0$, with $s\e\to+\infty$, and discover that all limiting values of 
 $t^{3/2-\gamma}v(t,t^\gamma)$ coincide. To end the proof of the theorem, it remains to show that $\beta$ can be chosen as small as needed, so that $\gamma$ can in the end be chosen arbitrarily small, and $\beta<\gamma$ .
\subsection{Convergence to the shifted wave}
\noindent At this stage, there is no additional idea, other than those presented in Section \ref{s2.5.3} for the large time bound of the level sets. It is, however, a good time to round up things and give an informal derivation of the full quantity by which we should shift a travelling wave in order to make it an asymptotic representation of $u$. Once this is understood, we will explain the (trivial) modifications of the arguments of Section \ref{s2.5.3} that have to be made in order to obtain the full Theorem \ref{t2.4.1}, however, an informal derivation of the shift. The philosophy is quite simple: it is a $C^0$ matching argument. Select $\gamma\in(0,1/2)$; if $u(t,x)$ is the solution of the initial problem \eqref{e4.2.8} or \eqref{e4.2.9}, the just proved Theorem \ref{t2.4.3} yields
\begin{equation}
\label{e4.5.340}
u(t,t^\gamma)\sim p_\infty/t^{3/2-\gamma}.
\end{equation}
If we believe for one moment that $u(t,x)$ should converge to a travelling wave profile $\varphi_{c_*}(x-\sigma_\infty(t))$, then it should certainly match the expression \eqref{e4.5.340} as $t\to+\infty$ and $x=t^\gamma$. Assume for convenience that $\varphi_{c_*}=xe^{-\lambda_*x}+O(e^{-(\lambda_*+\delta)x})$ as $x\to+\infty$; matching this expression with \eqref{e4.5.340} at $x=t^\gamma$ yields (the $=$ sign below should of course be taken with a pinch of salt):
$$
\bigl(t^\gamma-\sigma_\infty(t)\bigl)e^{-\lambda_*(t^\gamma-\sigma_\infty(t))}=p_\infty/t^{3/2-\gamma}e^{-\lambda_*t^\gamma}.
$$
Taking the $\mathrm{ln}$ of both sides yields
$$\sigma_\infty(t)=\di-\frac3{2\lambda_*}\biggl(\mathrm{ln}~t-\lambda_*\mathrm{ln}~p_\infty\biggl)+o_{t\to+\infty}(1).
$$

\noindent In order to justify this rigorously, we operate this time in the reference frame moving to the left like $3\mathrm{ln}~t{/2\lambda_*}$. We have, from Theorem
\ref{t2.4.3}: again
$
v(t,t^\gamma)\sim_{t\to+\infty}\alpha_\infty t^\gamma.$
Note indeed that operating in the reference frame moving to the left like $3\mathrm{ln}~t{/2\lambda_*}$ entails a multiplication by $t^{3/2}$ in the asymptotic expression of $u(t,x)$ as $x\to+\infty$. So, one has
 to remove the $t^{-3/2}$ factor in the
the new asymptotic expression of $v(t,x)$. Adding a logarithmic correction does not  alter it at $x=t^\gamma$.

\noindent For any $\e>0$, define the two translations $\sigma_{\infty,\e}^\pm$ such that
$
\varphi_*\bigl(t^\gamma-\sigma_{\infty,\e}^\mp(t)\bigl)=(\alpha_\infty\mp\e)t^\gamma e^{-\lambda_*t^\gamma}.
$
For $x\in[t^\gamma-1,t^\gamma]$ we have
$
\varphi_*\bigl(t^\gamma-\sigma_{\infty,\e}^-(t)\bigl)\leq u(t,x)\leq\varphi_*\bigl(t^\gamma-\sigma_{\infty,\e}^+(t)\bigl).
$
Arguing as in Section \ref{s2.5.3} we obtain
$$
\limsup_{t\to+\infty}\biggl(u(t,x)-\varphi_*\bigl(t^\gamma+\sigma_{\infty,\e}^-(t)\bigl)\biggl)\geq0,\ \ \liminf_{t\to+\infty}\biggl(u(t,x)-\varphi_*\bigl(t^\gamma-\sigma_{\infty,\e}^+(t)\bigl)\biggl)\leq0.
$$
As $\e>0$ is arbitrarily small, this puts an end to the proof of Theorem \ref{t2.4.1}.

\section{Bibliographical elements, comments, open questions}\label{s4.8}
\subsection*{The spreading speed (Section \ref{s4.1})}
\noindent This is a rather standard issue for this sort of spatially homogeneous models,  and we do not claim any novelty in what we are presenting. Spreading was figured out by Aronson \cite{A}. Actually, Aronson's work concerns the SIR model woth nonlocal contaminations \eqref{e2*.1.5}, which has a lot to do with our basic model \eqref{e4.1.1}; at least we hope that, at this stage, the reader is convinced of this. 
His results have, by the way, much to do with the classical Aronson-Weinberger \cite{AW} theorems for the 
Fisher-KPP equation.  

\noindent Equation \eqref{e2*.1.6} can further be  reduced to an integral equation on $v(t,x):=K*u(t,x)$:
\begin{equation}
\label{e2*.8.2}
v(t,x)=f(t,x)+\int_0^t\int_{\RR}e^{-\alpha(t-s)}K(x-y)(1-e^{-v(s,y)})~dy,
 \end{equation}
 where $f(t,x)$ depends on the data.  This is an important remark, as models including more effects may not always be amenable to a reduction of the form \eqref{e2*.1.6}, but are amenable to a reduction to \eqref{e2*.8.2}.
 For instance, there is the subsequent study  of Diekmann \cite{Diek}, who puts into the model the fact that an infected individual does not have the same infecting power during the course of the infection. He reduces it to an integral equation of the type \eqref{e2*.8.2}, where the exponential $e^{-\alpha s}$ is replaced by a more general integrable function of time. These results were proved independently and at the same time by Thieme \cite{Thi}, on a somehow different model, but still with a formulation of the type \eqref{e2*.8.2}. This has generated an important line of research on such integral equations, leading in particular to the beautiful abstract theory of monotone systems. I have limited myself to quoting Hirsch \cite{Hir}; the rest is simply  too abundant to be cited exhaustively here. 

\noindent While the entire book concerns the spatial dimension 1, it is appropriate to give at least some hints of what the speed of spreading becomes in higher space dimensions. To describe the general question,  consider a   model of the form
\begin{equation}
\label{e4.7.2}
u_t+Au=f(x,u),\quad t>0,\ x\in\RR^N.
\end{equation}
 In ecology, $u(t,x)$  models the density of population. The operator $A$ accounts for  population dispersal, a common choice being   the Laplacian. Due to nonlocal effects, an integral operator of the form
\begin{equation}
\label{4.7.3}
Au(t,x)=\int_{\RR^N}K(x,x')(u(t,x)-u(t,x'))~dx'.
\end{equation}
 is often a relevant choice, the fractional Laplacian being one among many.  If the diffusion is given by the Laplacian or a  rapidly decaying kernel $K$, it is called {\it short range}. When given by a slowly decaying kernel or   the fractional Laplacian, it is called {\it long range}.
The expression 
\begin{equation}
\label{e4.7.4}
f(x,u)=\mu(x)u-u^2
\end{equation} 
is relevant in many cases.  The term $\mu(x)u$ is a growth rate, while $-u^2$
is a saturation term that prevents the density from growing  unlimited.  When $\mu\equiv1$ and $A=-\Delta$, Model (\ref{e4.7.2}) is the Fisher-KPP equation. Assuming  $\mu(x)$ in (\ref{e4.7.4}) to be variable leads quickly to fascinating unsolved questions. Another important class of nonlinearities is that of the bistable type: if $f(x,u)=f(u)$ we have $f'(0)<0$ and $f'(1)>0$ (so that 1 and 0 are stable points of the ODE $\dot u=f(u)$).  In ecology, this corresponds to a situation when the medium is unfavourable when the population density drops below a certain threshold. 

\noindent  The basic question is the following: when the initial density  is a function with bounded support, the unstable  (or less stable) state 0 will be invaded (by the state $u\equiv 1$ if $\mu(x)\equiv 1$), and 
a transition will form between the region $\{u(t,x)\sim0\}$, and the region $\{u(t,x)\sim1\}$. The "speed of spreading" question is how fast this transition moves in each direction. In other words the issue is to find, for every direction $e$, a quantity $w_*(e)$ such that we have, uniformly on each compact in  $x$:
$$
\lim_{t\to+\infty}u(t,twe+x)=0\ \hbox{for all $w>w_*(e)$},\quad \liminf_{t\to+\infty}u(t,twe+x)>0\ \hbox{for all $w<w_*(e)$}.
$$

\noindent KPP's  result is generalised to higher dimensions by Aronson-Weinberger \cite{AW}, that is, $R_e(t)=2t+o(t)$. 
If $\mu$ to be periodic,  $A=-\Delta$, in any space dimension $N$,  Freidlin-G\"artner \cite{FrGa}, using the Feynman-Kac formula,
propose the following beautiful expression for $R_e(t)$:
\begin{equation}
\label{e33}
\frac{R_e(t)}t=\inf_{e'.e>0,\vert e'\vert=1}\frac{c_*(e')}{e.e'}.
\end{equation}
Here,  $c_*(e')$ is the smallest linear wave speed in the direction $e'$, that is, the least $c>0$ so that the linearised equation around  0, $v_t-\Delta v=\mu(x)v$, has solutions of the form $\phi(x)e^{-\lambda(x.e'-ct)}$. Several alternative proofs have been proposed,
one, due to Weinberger \cite{Wein} using abstract dynamical systems arguments, and being valid on periodic networks.  The results of \cite{FrGa} extend to the Fisher-KPP equation with a stationary ergodic growth rate $\mu(x,\omega)$. Here,  stationarity means that 
 $\mu (x+y,\omega)= \mu (x,\pi_{y}\omega)$ for all  $(x,y)$, where $(\pi_{y})_{y}$ is a group of measure-preserving transformations. 
 
 \noindent When propagation is expected to occur at asymptotically constant speed, it is natural to scale $t$ and $x$ in the same way, and so  to introduce the change of variables $(t/\e,x/\e)$, with $\e>0$, then to send $\e$ to 0. This 
leads to a homogenisation problem whose limit is an eikonal equation  (Evans-Souganidis \cite{ES1}), thus underlining a deep link with Hamilton-Jacobi equations. The  theory has since then been extended to second order geometric movements, with numerous important contributions. 

\subsection*{Logarithmic behaviour (Sections \ref{s4.3} and \ref{s4.4})}

\noindent   The Fisher-KPP equation 
\begin{equation}
\label{e2*.8.1}
u_t-u_{xx}=f(u)\ (t>0,x\in\RR)\ \ \ \ 
u(0,x)=1-H(x),
\end{equation}
with $f(u)\leq f'(0)u$, and where $H$ is the Heaviside function, has been the subject of fundamental studies. The first one, with $f(u)=u(1-u)^2$, is the pioneering paper of Kolmogorov-Petrovskii-Piskunov  \cite{KPP}. It is one of the works that one can read many times and still learn something at each reading session. My personal opinion is that it has not been emphasised enough in Kolmogorov's biographies. The term $f(u)=u-u^2$ plays an important role in probability theory. While it is not the primary intent of this book to dwell on the probabilistic side of the models studied here, it would be difficult to ignore this aspect, just because the first sharp asymptotics for \eqref{e2*.8.1} were issued by Bramson in a large paper \cite{Br1},  and a masterful memoir \cite{Br2}.  If $X(t)$ is, as usual, the rightmost point $x$ so that $u(t,x)=1/2$, the work \cite{Br1} proves that 
\begin{equation}
\label{e2*.8.200}
X(t)=2t-\frac32\mathrm{ln}~t+O(1),
\end{equation}
while \cite{Br2} shows that the $O(1)$ term is asymptotically constant. These  results are based on various probabilistic interpretations of Model \eqref{e2*.8.1} {\it via} the Brownian Motion. There is a remarkable bridge, drawn by McKean \cite{McK}, between $u(t,x)$ and the rightmost particle $\Omega_t$ in the Branching Brownian Motion (BBM): we have
\begin{equation}
\label{e2*.8.3}
u(t,x)=\mathbb{P}\bigl(\{\Omega_t\geq x\}\bigl).
\end{equation}
This is the path used to retrieve \eqref{e2*.8.200} in \cite{Br1}, while \cite{Br2} relies on a representation of $u(t,x)$ in terms of the classical Brownian Motion $(B_t)_{t>0}$, through the Feynman-Kac formula: 
it turns out indeed that the solution $u(t,x)$ of the equation
$
u_t-\di\frac12u_{xx}+c(t,x)=0\ (t>0,x\in\RR),$ with $u(0,x)=u_0(x),
$
has the beautiful representation
\begin{equation}
\label{e4.8.1}
u(t,x)=\mathbb{E}\biggl(u_0(x+B_t)\mathrm{exp}~\bigl(-\int_0^tc(s,x+B_s)~ds\bigl)\biggl).
\end{equation}
See for instance {\O}ksendal \cite{Ok} for a proof of \eqref{e4.8.1}. In the case of \eqref{e2*.8.1}, with $f(u)=u-u^2$, one takes $c(t,x)=1-u(t,x)$, something that one would love to identify to $1-H(x)$.
  This seemingly innocent idea is put to work in a large part of  \cite{Br2}, which unleashes an impressive aresenal of  methods and computations (I believe of independent interest) for the study of the Brownian Motion. 
 A more recent work of Roberts \cite{Rob} retrieves Bramson's results by a detailed understanding of the law of $\Omega_t$.  All this raises the question of what nonlinearities, other than $u-u^2$,  lead to a representation of $u(t,x)$ of the McKean type \eqref{e2*.8.3}. A  study of An, Henderson and Ryzhik \cite{AHR3} assesses what kind of nonlinearites are amenable to a Mc Kean type representation with the Branching Brownian Motion. It turns out that, surprisingly,  even $f(u)=u-u^n$ cannot, except for $n=2$. This comes from   stubborn facts in elementary algebra. For even more elaborate studies of the properties of the Branching Brownian Motion, which keep pouring at the time of writing of this book, I refer, for instance, to 
  A\"idekon-J. Berestycki-Brunet-Shi \cite{ABBS}.

\noindent The above history (thanks to my colleagues and friends from Probability who have patiently endeavoured to teach it to me), together with the independent interest of understanding things within the PDEs realm, legitimates an
  approach alternative to the probabilistic one. And so, less than thirty five years after \cite{Br1}-\cite{Br2}, a short and flexible PDE proof  of \eqref{e2*.8.200} was proposed by Hamel, Nolen, Ryzhik and the author \cite{HNRR}. It was later converted in a full proof of the Bramson theorem by Nolen, Ryzhik and the author in \cite{NRR1}. The reader will, incidentally, notice that the idea of studying the Feynman-Kac formula \eqref{e4.8.1} with $c(t,x)=1-u(t,x)$  is not at all far from solving Dirichlet problems with the heat equation, as is insistently done in this monograph. This is of course something that occurred to me some time after \cite{HNRR} and \cite{NRR1} were issued. On the other hand, developping this independent approach has not only enabled to understand 
  nonlocal models that are not so accessible to probabilistic representations, as presented in this monograph or in the selection of problems, it has also opened the door to further understanding of the Branching Brownian Motion, through
   unconventional links with equation \eqref{e2*.8.1}; an example was worked out by Mytnik, Ryzhik and the author of the present book \cite{MRR}.

\noindent A remarkable fact is that Model \eqref{e4.1.1}, with $f(u)=u-u^2$, also comes from a McKean type representation. The Branching Brownian Motion is this time replaced by a Branching Random Walk.  In slightly more explicit terms, a particle jumps, splits, and the offsprings reproduce the ancestors' behaviour. Splitting and jumping times are random, given by Poisson type  distributions.  The distribution of jumps is given by the density $K$. We may then introduce the position $\Omega_t$ of the rightmost particle, and the function $u(t,x)$ defined by \eqref{e2*.8.3} is a solution of our initial model \eqref{e4.1.1}. As the particles may split into more than two offsprings, nonlinearities of the form $f(u)=u-u^2P(u)$, where 
$P$ is a power series with positive coefficients such that $P(1)=1$, may be allowed. However, the theory would not apply for all $P$'s, similarly to the Branching Brownian Motion.  The reader is referred to Graham \cite{Gr1}
for an extremely clear description of the bridge between the random walk and $u(t,x)$. The paper \cite{Gr1} covers kernels $K$ that can be much more wild than those considered in the book, and explores the extreme assumptions on $K$ such that the results displayed here still work. All in all,  if $X(t)$ is the position of the rightmost point such that $u(t,x)=1/2$, an expansion of $X(t)$ up to $o_{t\to+\infty}(1)$ terms is  proved by  A\"idekon \cite{Ai}, following $O_{t\to+\infty}(1)$ type estimates in Addario Berry-Reed \cite{ABR}.   The result is given in the more general framework of discrete time branching random walks.

\noindent Interestingly, while not having any comparison principle, the Berestycki-Nadin-Perthame-Ryzhik model for competition \cite{BNPR}, that reads
\begin{equation}
\label{e4.8.4}
u_t-u_{xx}=u(1-\phi*u),
\end{equation}
 displays a logarithmic delay in the propagation. We present several versions of this in the problems below, the first result in this direction being due to Bouin-Henderson-Ryzhik  \cite{BoHR}. However it is not known whether the next term is a constant or not; it may well be an oscillating function, depending on possible instabilities of the basic travelling wave. A probabilistic counterpart is studied, for instance, in Addario Berry-J. Berestycki-Penington \cite {ABPen}.    

\subsection*{ Open questions}
\noindent   What happens when the initial datum for \eqref{e4.1.1}, instead of being compactly supported, has some exponential decay at infinity, that is
$u_0(x)\lesssim e^{-r\vert x\vert},
$
for some (possibly quite small) $r>0$? This question looks benign and innocent, it is in fact truly nontrivial. If $r>0$ is small, we have just seen that \eqref{e4.1.1} has a lot of travelling waves with different speeds, 
all decaying faster than $e^{-r\vert x\vert}$ at infinity. Hamel and Nadirashvili \cite{HN}, in a important contribution, show for the Fisher-KPP equation with standard diffusion, that these waves will joyfully mix in order to form a big and intricate set of eternal solutions. In $N$ space dimensions, it contains a set homeomorphic to the set of bounded measures of the unit sphere. By the way, the classification of the eternal solutions is not yet complete, as degeneracies in the vicinity of the bottom speed complicate the task a lot. It is to be noted that the question is just as wide open in the case where the diffusion is Gaussian.

\noindent As soon as one leaves the realm of models for which a comparison principle is available, one quickly enters the domain of more or less open questions, and most often more open than understood. If we think, for instance, of models in ecology, it is quite legitimate to take competition into account. A prototype model is \eqref{e4.8.4}, whose behaviours is studied in Problems \ref{P2.103}-\ref{P2.105} (and these problems are not speculative ones). When spatial behaviour is coupled with diffusion between traits, 
one may (at least I firmly believe so) carry out a complete study of the Berestycki-Chapuisat type models, such as  \eqref{e1.3.8}, with no really new ideas. Such is, however, no more the case when competition between traits occur. An example of this stituation is given by the equation
\begin{equation}
\label{e4.8.5}
u_t-\Delta u+\alpha y^2u=\biggl(1-\int_\RR K(z)u(t,x,z)~dz\biggl)u,\quad (t>0,\ (x,y)\in\RR^2),
\end{equation}
which introduces a nonlocal competition term in \eqref{e1.3.8}. The work \cite{BJS}, by Berestycki-Jin-Silvestre, proves the existence of a propagation speed, and provides a study of travelling waves. The second derivatives in $x$ could be replaced by the operator $\mathcal{J}$. A related model is treated by Alfaro, Coville and Raoul \cite{AlCoRa}. In this class of models, one may also put the model for the proliferation of the cane toad (an unpleasant batracian that invaded a large part of Australia in, roughly, the second half of the XXth century), as treated by Bouin, Henderson and Ryzhik \cite{BHR2}. In all these models, the newest issues are less the front location up to $O(1)$ terms, as it can be figured out with more or less the same ideas as displayed in this chapter (and a good deal of additional technical work), than what happens in the zone ahead of the diffusive zone. What controls the behaviour of the solutions is the (global) stability of the travelling waves in spaces of functions having a special behaviour at infinity, and very little is known in general on this subject. The reason is of course the failure of the maximum principle.

\noindent Another fascinating line of research concerns models which have derivatives of higher order, a basic instance being the fourth order equation
\begin{equation}
\label{e4.8.6000}
u_t+u_{xxxx}-u_{xx}=f(u),\quad (t>0,\ x\in\RR)
\end{equation}
Here, $f$ is of the Fisher-KPP type. An interesting remark is that this model can easily be obtained as the truncation to higher order of our favourite nonlocal problem \eqref{e1.1.1} (exercise: work it out). This innocent 
procedure carries a high price, that is, the loss of the comparison principle. A highly inspiring paper of Ebert and Van Saarloos \cite{EVS} (see also \cite{VS} for many examples advocating for the physical relevance of higher order Fisher-KPP models)
proposes a full expansion of the front position, in the formal style. It pioneers the detailed study of the standard Fisher-KPP equation. Despite the surprising form of the proofs, and despite the fact that the findings seem to be contradicted by other formal 
studies (J. Berestycki-Brunet-Derrida \cite{BBD}) and rigorous studies (Graham \cite{Gr1}), it is a source of information whose value is difficult to match. Needless to say, any rigorous confirmation of the proposed expansions for the higher order
 problems would be an achievement. A good starting point would be to treat the fourth order term as an annoying, but benign perturbation   and see how Theorem \ref{t2.1.1} is modified.
\section{Problems}\label{s4.9}
\begin{problem}\label{P4.8.0}
Extend the theory as much as possible when $K$ is smooth but not compactly supported anymore. It should, however, be assumed that $K$ has finite second moment. 

\noindent When $K$ decays exponentially fast at infinity, a probabilistic alternative approach, which does not, however, capture the whole shift, can be found in  Boutillon \cite{Bou}. 
Among the places that deserve care are the comparison principle, and the study of the heat kernel.
\end{problem}
\begin{problem}\label{P4.8.15}
Let $u(t,x)$ be the solution of \eqref{e4.1.1} starting from a nonnegative compactly supported initial datum. Show that $u(t,x)$ converges exponentially fast to 1 in every set of the form $\vert x\vert\leq ct$, with $c<c_*$.
\end{problem}

\begin{problem}\label {P4.8.1} In the basic model \eqref{e2*.1.5}, assume that $K$ is an approximation of the identity. The nonlinearity is assumed to be of the form $\e^2f(u)$, with $f(u)\leq f'(0)u$.

\noindent Find out the large time behaviour of the solution $u_\e(t,x)$ of \eqref{e4.1.1}, modulo a rescaling in time, with your favourite compactly supported initial datum $u_0$. The challenge is here to show that it has some uniformity in $\e$, in other words, and that the limits commute. In particular, to what extent is the expansion of $u_\e$ for large time uniform in $\e$? 

\noindent A hint of what can be expected may be find in Pauthier \cite{Pau}. There, he investigates some uniform dynamical properties of the propagation directed by a line of fast diffusion, when the exchanges between the line and the ambient space become localised.
\end{problem}
\begin{problem}\label{P4.8.3}
In the context of Problem \ref{P4.8.1}, rescale $t$ like $1/\e$ and  investigate the regularity of the solution for small $\e$. For $\e=0$ there is a gradient bound, and for $\e>0$ Section \ref{s2.4} yields and $\e$-dependent gradient bound.
Is there a  gradient bound that is uniform in time and uniform in $\e$? 
\end{problem}
\begin{problem}\label{P8.4.200}
An important application of Theorem \ref{t2.1.1} is the complete  study of models arising from epidemiology. I guess that it is a good time to
 recall, from Chapters \ref{Intro} and \ref{Cauchy}, the SIR model with nonlocal interactions: if $x\in\RR$, let $S(t,x,)$ denote the fraction of susceptible individuals, 
 and  $I(t,x)$ the fraction of infected individuals.  The model   reads
\begin{equation}
\label{e2*.1.5}
\left\{
\begin{array}{rll}
\partial_tI+\alpha I=&SK*I\quad(t>0,\ x\in\RR)\\
\partial_tS=&-SK*I
\end{array}
\right.
\end{equation}
The quantity $\beta>0$ is still the mass of $K$. We assume  $S(0,x,y)\equiv S_0>0$, and that  $I(0,x)=I_0(x)$ compactly supported.   The cumulative numbers of infected individuals $u(t,x)$ solves 
\begin{equation}
\label{e2*.1.6}
\partial_t u+\alpha u=S_0(1-e^{-K*u(t,.)})+I_0(x).
\end{equation}
\noindent Start from the result displayed in Problem \ref{P2.108}, or, otherwise said, Kendall'a {\it pandemic threshold theorem}: Let $R_0={S_0\beta}/\alpha$, assume it is $>1$. Then
 the stationary equation
\begin{equation}
\label{e2*.2.25}
\alpha u=S_0\bigl(1-e^{-\beta K*u}\bigl)+I_0(x),
\end{equation}
has a unique solution $u_\infty[I_0](x)$. Moreover we have
$\di\lim_{t\to+\infty}u(t,x)=u_\infty[I_0]$ locally uniformly in $x\in\RR$. We denote by $u_*$ the unique positive constant solution of \eqref{e2*.2.25} with $I_0\equiv0$ and we have $u_\infty[I_0]\geq u_*$ with 
$\di\lim_{\vert x\vert\to+\infty}u_\infty[I_0](x)=u_*$.  Let us still denote by $c_*$ the critical travelling wave speed.
\begin{itemize}
\item [---] (Spreading at critical speed). Show that 
$$
\hbox{for all $0<c<c_*$,}\ \lim_{t\to+\infty}\sup_{\vert x\vert\geq ct}u(t,x)=u_*,\ \hbox{and for all $c>c_*$,} \  \lim_{t\to+\infty}\sup_{\vert x\vert\geq ct}u(t,x)=0.
$$
The presence of $I_0$ should be seen as a minor inconvenience. One can ignore it for the lower bound, simply by starting from $t=1$ and putting a well chosen datum below $u(1,x)$. For the upper bound, argue outside the support of $I_0$.
\item [---] Run with speed $c_*$ to the right, so that \eqref{e2*.1.6} becomes
\begin{equation}
\label{e4.7.20}
\partial_t u-c_*u_x+\alpha u=S_0(1-e^{-K*u(t,.)})+I_0(x+c_*t).
\end{equation}
Show the analogue of Theorem \ref{t2.4.3}, that is:
there exists $p_\infty>0$ such that, for all $\gamma\in(0,1/2)$ we have
$u(t,t^\gamma)={p_\infty e^{-\lambda_*t^\gamma}}{t^{3/2-\gamma}}\bigl(1+o_{t\to+\infty}(1)\bigl).$
Here, the attention of the reader is drawn upon the fact that, is is tempting to expand \eqref{e4.7.20} into an equation of the form
\begin{equation}
\label{e4.7.21}
\partial_tu-c_*u_x+S_0\beta\bigl(K*u/\beta-u\bigl)=\alpha(R_0-1)u+h(K*u)+I_0(x+ct).
\end{equation} 
The function $h$ is nonpositive, quadratic in the vicinity of 0 and decreasing. Therefore, equation \ref{e4.7.21} does not have a comparison principle. 
Yet, equation \eqref{e4.7.20} does have one.
\item[---] (Intermezzo) Let $a(x)$ and $b(x)$ be  continuous positive functions such that $a(x),b(x)\leq S_0\beta$. Show that $\bar u(x):=e^{-\lambda_*x}$ solves
\begin{equation}
\label{e4.7.22}
a(x)(\bar u-K*\bar u/\beta)-c\bar u'-\bigl(b(x)-\alpha)\bar u\geq0.
\end{equation}
\item[---] We set $f(u)=1-e^{-\beta u}$. Let $u_1(t,x)$ and $u_2(t,x)$ be two solutions of \eqref{e4.7.20}. Show that the difference $v(t,x)=u_1(t,x)-u_2(t,x)$ solves
\begin{equation}
\label{e4.7.23}
\partial_tv+c_*v_x+S_0f'(\bigl(a(t,x)\bigl)\bigl(K*v/\beta-v\bigl)=\bigl(S_0f'\bigl(a(t,x)\bigl)-\alpha)v,
\end{equation}
where the quantity $a(t,x)$ is to be expressed with a suitable Taylor integral remainder.
\end{itemize}
\item [---] From the Intermezzo, figure out a super-solution to \eqref{e4.7.23} of the form $e^{-\lambda_*x}\cos(x/t^\alpha)/t^A$ on an interval of the form $(-t^{\gamma},t^\gamma)$, $\gamma>0$ small.
\item[---] Show the existence of a solution $\phi_{c_*}$ of the travelling wave equation \eqref{e3.8.21} and $x_\infty\in\RR$  such that $u(t,x)=\phi_{c_*}\bigl(x+\di\frac3{2\lambda_*}\mathrm{ln}~t+x_\infty\bigl)+o_{t\to+\infty}(1)$.
\item[---] Adapting the arguments of Section \ref{s4.10}, show that $u_x(t,x)=\phi_{c_*}'\bigl(x+\di\frac3{2\lambda_*}\mathrm{ln}~t+x_\infty\bigl)+o_{t\to+\infty}(1)$.  
\item[---] Reverting to the initial $t$ and $x$ variables, that is, $x$ not translated by $c_*t$, conclude that 
$$I(t,x)=-c_*\phi_{c_*}'\bigl(x+\di\frac3{2\lambda_*}\mathrm{ln}~t+x_\infty\bigl)+o_{t\to+\infty}(1).
$$
\end{problem}
\begin{problem}\label{P4.8.2}
Still in the   SIR model  with nonlocal contaminations \ref{e2*.1.5}, let us assume that the infected diffuse, with  diffusion coefficient $d$. This yields the system
\begin{equation}
\label{e2*.9.1}
\left\{
\begin{array}{rll}
\partial_tI-d\partial_{xx} I+\alpha I=&SK*I\quad(t>0,\ x\in\RR)\\
\partial_tS=&-SK*I
\end{array}
\right.
\end{equation}
Work out the whole large time study, and show that the presence of the coefficient $d$ does not alter the results: only the quantities change. Show, in particular, that the limit $d\to0$ is nonsingular. Estimate the distance between the solution of \eqref{e2*.9.1} and that with $d=0$ uniformly in time.
\end{problem}
\begin{problem}\label{P4.8.4}
Conversely to Problem \ref{P4.8.2}, we assume $d>0$ in \eqref{e2*.9.1}, but we now assume that $K$ is an approximation of the identity, similar to that of Problem \ref{P4.8.1}. Investigate the large time dynamics of the system, trying to keep is as uniform as possible in $\e$.
What happens if $d=0$?
\end{problem}
\begin{problem}\label{P$.8.14} This problem, and the next one, are devoted to proving the logarithmic delay up to $O(1)$ terms in the nonlocal model for competing species. 
Consider, in the present problem, the original model \eqref{e4.8.4}.
\begin{itemize}
\item[---] Show that the minimal wave speed is $c_*=2$.
\item[---] In the reference frame moving with speed $c_*$, 
produce, in the spirit of Section \ref{s4.3}, a pair of sub/super-solutions leading to the estimate
$$
{xe^{-x}}/{t^{3/2}}\lesssim u(t,x)\lesssim{xe^{-x}}/{t^{3/2}},\ \ \ \  t^\delta\leq x\leq t^{1/2-\delta}.
$$
\item[---] Recall, from Problem \ref{P2.103}, the existence of $M>0$ such that $0\leq u(t,x)\leq M$. Produce a concave function $f_M$ such that $u(t,x)\bigl(1-\phi*u(t,x)\bigl)\leq f\bigl(u(t,x)\bigl)$, with $f'(0)=1$ and $f(0)=f(2M)=0$. 
\item[---] If $\varphi_M$ solves
\begin{equation}
\label{e4.8.6}
-\varphi_M''-c_*\varphi_M'=f_M(\varphi_M)\ (x\in\RR),\ \ \ \ \varphi_M(-\infty)=2M,\ \varphi_M(+\infty)=0,
\end{equation}
show the existence of $K_M>0$ such that, back in the original reference frame, we have $u(t,x)\leq\varphi_M(x-2t-3\mathrm{ln}~t/2-K_M)+o_{t\to+\infty}(1)$.
\item[---] Consider $\e>0$ small and $(t_0,x_0)\in[2,+\infty)\times\RR$ such that $u(t_0,x_0)=\e$. Show, using the fact that $u$ is globally bounded, the estimate
$$
\e\geq\delta\int_\RR e^{-y^2/4}u(t_0-1,y)~dy.
$$
Deduce from this fact the (nonoptimal) estimate $\phi*u(t_0,x_0)\lesssim\sqrt\e$. This entails the estimate $1-\phi*u(t,x)\geq1-u(t,x)/\e^{1/3}$, for a sufficiently small $\e$.
\item[---] If $\varphi_\e$ solves \eqref{e4.8.6}, with this time $M$ replaced by $\e^{1/3}$, show the existence of $K'>0$ such that,
still in the original reference frame, we have $u(t,x)\leq\varphi_\e(x-2t-3\mathrm{ln}~t/2-K')+o_{t\to+\infty}(1)$.
\end{itemize}
All this proves that, if $\theta>0$ is a small constant, then we have
$
X_\theta(t)=c_*t-3\mathrm{ln}~t/2+O(1).
$
\end{problem}
\begin{problem}\label{P4.8.11}
\noindent Investigate the large time behaviour of a solution to the Cauchy Problem for \eqref{e4.1.1} that is trapped between two different travelling waves of the same speed $c$. An easy start is what happens if
the initial datum is a perturbation of a travelling wave with speed $c>c_*$ that is compactly supported in $\RR_+$. 

\noindent The same question, with $c=c_*$ is  more interesting. One may ask if the results of an insightful work of Gallay \cite{Gal}, for the Fisher-KPP equation, that gives the optimal rate of convergence to the wave, are valid in this setting. One may also consult another proof of this result by Faye and Holzer \cite{FH}. In any case, the heat kernel estimates developped in Section \ref{s2.4} should be useful.

\noindent In the  setting \eqref{e4.1.1}, one may use the work of Bages, Martinez and the author \cite{BMR} for the homogeneous Fisher-KPP equation as a guideline, one may suspect that the results are true, although one should probably be careful in adapting the result. Once again the job is expected to be more difficult for the critical velocity. 

\end{problem}
\begin{problem}
\label{P4.50}
Consider the time-dependent version of \eqref{e1.1.1}, recalled below:
\begin{equation}
\label{e4.8.50}
u_t-K*u+u=f(t,u)\ \ (t>0,x\in\RR)
\end{equation}
with, this time, and to simplify things, an initial datum $u_0$ such that $\mathrm{supp}~u_0\subset\RR_-$ and $\di\liminf_{x\to-\infty}u_0(x)>0$. Let $\varphi_{c_*}(t,x)$ be the periodic travelling wave, decaying as $xe^{-\lambda_*x}$ as $x\to+\infty$. Show that
$u(t,x)=\varphi_{c_*}\bigl(t,x-c_*t-\di\frac3{2\lambda_*}\mathrm{ln}~t+x_\infty)+o_{t\to+\infty}(1)$.
\end{problem}
\begin{problem}
\label{P4.20}
Consider the original Berestycki-Chapuisat model, that is, Model \eqref{e1.3.8} but with standard diffusion in all variables:
\begin{equation}
\label{e4.8.25}
u_t-\Delta u+\alpha y^2u=f(u)\ \ \bigl(t>0,(x,y)\in\RR^2)\bigl).
\end{equation}
The function $f$ is a standard Fisher-KPP term, and we start from our favourite smooth compactly supported initial datum $u_0(x,y)$.
\begin{itemize}
\item[---] [Warm-up] For the original one-dimensional Fisher-KPP equation $u_t-u_{xx}=f(u)$, the proof of Theorem \ref{t2.1.1} can be to a large extent simplified. Do it as much as you can; for the front position up to $O_{t\to+\infty}(1)$ terms a possible reference is \cite{BFRZ}. For the location up to $o_{t\to+\infty}(1)$ terms, Theorem \ref{t2.4.4} is not needed anymore, but this is not the only possible simplification.
\item[---] Investigate the large time behaviour of the initially compctly supported solutions $v(t,x,y)$ of the Dirichlet problem
$$
v_t-\Delta v-\bigl(c_*-2\lambda_*\bigl)v_x+\bigl(\alpha y^2-\lambda_*^2-f'(0)\bigl)v=0\ \bigl(t>0,(x,y)\in\RR^2\bigl)\quad
v(t,0,y)=0.
$$
Notice that  an odd extension in $x$ provides the solution.
\item[---] Deduce that  $u(t,x,y)$ behaves, as $t\to+\infty$, as $u(t,x,y)=\di\varphi_{c_*}\bigl(x-c_*t+\frac3{2\lambda_*}\mathrm{ln}~t+x_\infty,y\bigl)+o_{t\to+\infty}(1),$
the convergence being uniform in $y\in\RR$.
\end{itemize}
\end{problem}
\begin{problem}
\label{P4.21}
We put the longitudinal diffusion back in, that is, we come back to Model \eqref{e1.3.8}. Adapt all the steps leading to the proof of Theorem \ref{t2.1.1} for the basic model \eqref{e1.1.1}, and derive a sharp asymptotic estimate for $u(t,x)$.
\end{problem}
\begin{problem}
\label{p4.8.11}
Consider a large $a>0$, and the following problem:
\begin{equation}
\label{e4.8.20}
\left\{
\begin{array}{rll}
\partial_tu+\JJ u=&f(u)\quad (t>0,\ -a\leq x\leq a)\\
u(t,x)=&0\quad (t>0,\ -a-1\leq x<a)
\end{array}
\right.
\end{equation}
starting from $u_0$ smooth, small, supported in $[-1,1]$.
  Show that $u(t,x)$ converges, as $t\to+\infty$, to the solution $v_a$ of Problem \ref{P2.115}. The convergence is uniform in $x$.
Then,  study the stabilisation of $u(t,x)$ to $v_a$ in the double limit $t\to+\infty$ and $a\to+\infty$. One should first see, as an intermediate asymptotics, a behaviour very close to two counter-propagating waves at speed $c_*$ with the logarithmic delay. How long can one observe this?
\end{problem}
\begin{problem}\label{P4.8.5}
Instead of taking a smooth kernel in Model \eqref{e4.1.1}, assume that we have
$$
K(z)\sim_{z\to0}k/{\vert z\vert^\gamma},
$$
with $0\leq\gamma<3$. Also assume $K$ to be compactly supported.
 If $\gamma<1$, it is a more or lass harmless modification of what has already been displayed, as $K$ still has a finite mass. This is not true anymore for $\gamma\geq 1$, and, for $\gamma>1$, the convolution looks like the fractional Laplacian studied in the next section, apart from the fact that $K$ is compactly supported. For $\gamma\in(1,2)$, the integral
$$
\int_{\RR}K(x-y)\bigl(u(x)-u(y)\bigl)~dy
$$
makes sense as soon as $u$ is $C^1$. For $\gamma>2$, it should be understood in the distributional sense.

\noindent How does the study of Problem \eqref{e4.1.1} gets modified? The crucial ingredient should be the study of the heat kernel. For $\gamma<1$, the singularity should just be a cumbersome, harmless inconvenient. For $\gamma>1$, regularising effects should take place: see, for instance, Caffarelli-Silvestre \cite{CS} and the references therein. You are not obliged to digest the whole theory: it is a huge piece of mathematics. However, the study of the heat kernel should be an important part of the resolution of this problem.
\end{problem}
\begin{problem}\label{P4.8.10}
Work out completely the study of the solution $v(t,x)$ of 
$$
v_t+\frac3{2\lambda_*t}(v_x-v)+\mathcal{I}_*v+h(x,v)=0.
$$
 in the diffusive  area, that is ${x/}{\sqrt t}\in(t^{-\delta},t^\delta)$, and how it matches with the area of larger $x$. 

\noindent As an application, propose an asymptotic expansion of the shift in terms of half powers of $1/t$. For the standard Fisher-KPP equation, a formal asymptotic expansion is proposed in the fascinating work of Ebert and Van Saarloos; the first term of the expansion is rigorously confirmed by Nolen, Ryzhik and the author \cite{NRR2}, and the second term is worked out (in contradiction to \cite{EVS}) by Graham in \cite{Gr1}). Another formal asymptotic expansion, still for the standard Fisher-KPP equation, is proposed by J. Berestycki, Brunet and Derrida \cite{BBD}. Find out whether this expansion is still valid.
\end{problem}

\chapter{Sharp ZFK spreading}\label{ZFK_short_range}
\section{Introduction}\label{s3.1}
\noindent This introductory part is not only the occasion to set the problem and ask the precise question that is at stake, it  will also be a good occasion to take a fresh view of what we have achieved in the preceding chapter, and what remains to be done when the seemingly innocent Fisher-KPP assumption  for the nonlinearity is removed.
\subsection{Model and question}
We investigate the solutions of the integral equation
\begin{equation}
\label{e3.7.2}
u_t+\mathcal{J}u=f(u),\ \ t>0,x\in\RR,
\end{equation}
with $\mathcal{J}u(x)=u(x)-K*u(x)$. The assumptions on the kernel $K$ are kept unchanged, it is nonnegative, not identically 0, and compactly supported; given the expression chosen for $\mathcal{J}$ its mass is 1. Instead of being compactly supported, the initial datum $u_0(x)$ will be assumed to satisfy
\begin{equation}
\label{e3.7.3}
\lim_{x\to-\infty}u_0(x)=1, \ \ \ u_0(x)=0\ \hbox{if $x\geq0$}.
\end{equation}
We apologise to the reader for this apparent change of paradigm. In fact, this is not such a drastic change, our choice simply assumes that invasion has still occurred in the left part of the landscape, and that we only need to concentrate on what happens to the right. The interested reader may check, with a pen and paper, that we are not cheating. Moreover, once we walk into this sort of assumption at $-\infty$, an important relaxation can be made for the behaviour at $-\infty$, this is the object of one of the problems below.

\noindent The function $f$ is smooth, and satisfies  $f(0)=f(1)=0$. It is positive on $(0,1)$ and satisfies $f'(0)>0$ and $f'(1)<0$.  As we wish, in this chapter, to highlight the differences in the behaviour of the Cauchy problem with this sort of nonlinearity, we do not aim at having the most general assumptions. Therefore we make the following additional hypothesis, which is really a commodity one: there is $\theta_0\in(0,1)$ such that
\begin{equation}
\label{e5.1.140}
f(u)\geq f'(0)u\ \hbox{for $u\in[0,\theta_0]$}.
\end{equation}
A typical case will be  a function $f$ such that $f'(0)$ is a small positive constant, while its mass is big. The class of functions satisfying this set of assumptions will be called ZFK nonlinearities, from the physicists Zeldovich and Frank-Kamenestkii, for reasons that are explained at the beginning of the book.

\noindent The question is, once again, we ask what happens when $t\to+\infty$.  Recall that, for a given $\theta>0$, we look for an asymptotic expansion, as $t\to+\infty$, of the quantity
\begin{equation}
\label{e3.1.10}
X_\theta(t)=\sup\{x\in\RR:\ u(t,x)=\theta\}.
\end{equation}
We already have a hint that an invasion at a positive speed will occur. Indeed, let $g(u)$ be any concave function satisfying $g(u)\leq f(u)$ and $g'(0)>0$. Such a function is not really hard to find, and we have $u(t,x)\geq v(t,x)$, where $v$ is the solution of the Cauchy Problem for \eqref{e3.7.2} with $g(u)$ instead of $f(u)$, and the same initial datum. Then, if $X^g_\theta(t)$ is the corresponding level set, and $c_K^g$ the bottom speed defined in Chapter \ref{short_range} (we had denoted it by $c_*$, and the reason for this change of notations will be made clear soon) we have
$$
X_\theta(t)\geq X_\theta^g(t)\geq c_K^gt+O(\mathrm{ln}~t).
$$
Our goal is to make this estimate much more precise. If $c_K$ is the Fisher-KPP speed corresponding to $f$, we are going to see that at least two sorts of behaviours are possible.
We will see that there is a class of nonlinearities $f$ for which we have $
\di \liminf_{t\to+\infty}{X_\theta(t)}/t>c_K.
$
In this case we will have, for some $c_*>c_K$, $\omega>0$ and $x_\infty\in\RR$:
\begin{equation}
\label{e3.1.301}
X_\theta(t)=c_*t+x_\infty(\theta)+O(e^{-\omega t}).
\end{equation}
Conversely, we will also meet a class of nonlinearities $f$ for which we have
$
\di\limsup_{t\to+\infty}{X_\theta(t)}/t\leq c_K.
$
Setting this time $c_*=c_K$, we will see that there is $\mu_*>0$ such that
\begin{equation}
\label{e3.1.300}
X_\theta(t)=c_*t-\mu_*\mathrm{ln}~t+x_\infty(\theta)+o_{t\to+\infty}(1).
\end{equation}
In cases that we will encounter, we will have $\mu_*=3/{2\lambda_*}$, where $\lambda_*$ is the decay exponent at infinity of the bottom wave.  The goal of this chapter is to explain  what causes these apparently very different  behaviours.
\section{Large time behaviour for the Cauchy Problem when $c_*>c_K$}\label{s3.4}
\noindent The result is once again very different from what we have already seen in the classical Fisher-KPP case. We have now that 
$\varphi_{c_*}$ behaves like $e^{-\lambda_+(c_*)x}$ as $x\to+\infty$, as a very important difference with the Fisher-KPP setting. And, as a matter of fact, the behaviour of the Cauchy problem is also 
very different. 
\begin{theorem}
\label{t3.7.3}
Under the assumptions \eqref{e3.7.3} for $u_0$, there exists $x_\infty\in\RR$ and $\omega>0$, universal, such that 
\begin{equation}
\label{e3.7.19}
\sup_{x\in\RR}\bigl\vert u(t,x)-\varphi_{c_*}(x-c_*t+x_\infty)\bigl\vert=O(e^{-\omega t}).
\end{equation}
\end{theorem}
\noindent In order to prove the theorem for all initial datum, what better start than to try to understand what happens when the solution initiates from something that is close to a wave. This is our first step,
through which we will go {\it via}  a rather classical stability argument. The least standard part will be the computation of the exponential of the linearised operator; if the result is itself standard the lack of dissipativity 
prevents the use of the classical tools of functional analysis, and forces us to resort to push the basic tools to their limits. A first consequence is that we will have an in-depth understanding of the underlying processes, and a second consequence will be that,
 if convergence occurs, it will be exponential. 
 
 \noindent Once this is understood, what is at stake is to prove that the set of travelling waves is approached for large times, or, in other words, that a travelling wave is approached for an infinite sequence of times.
 Stability will set in, and imply the theorem. For this step, the first task is to prove that all nontrivial things happen in the reference frame moving exactly like $c_*t$. 
 A final argument will show that all asymptotic profiles are travelling wave profiles.
 
 \noindent The starting point is simple: let us do the change of reference frames $x'=x-c_*t$ and immediately set $x:=x'$. The unknown function is still
denoted $u(t,x)$, and the  new
equation is 
\begin{equation}
\label{e3.7.35}
u_t+\mathcal{J}u-c_*u_x=f(u) \ \ \ (t>0,\ x\in\RR).
\end{equation}
We are now going to work on \eqref{e3.7.35}, for which $\varphi_{c_*}$ is simply a rest point connecting 1 at $-\infty$ to 0 at $+\infty$. 

\noindent In order to assess whether an initial datum close to a wave in, say, $B_{w_r,1}$, will entail convergence to a
(possibly nontrivial) translate of the wave, the standard idea is to linearise around $\varphi_{c_*}$ and to consider the linear equation
\begin{equation}
\label{e3.7.1000}
v_t+\mathcal{J}v-c_*\partial_xv-f'(\varphi_{c_*})v=0,\ \ \ v(0,x)=v_0(x)\in B_{w_r,1}.
\end{equation}
Recall (see \eqref{e3.7.43}) that the linearised operator around $\varphi_{c_*}$ is denoted by $\mathcal{M}_{c_*}$, so that an important effort will be devoted  to computing $e^{-t\mathcal{M}_{c_*}}$. The section will be organised in two parts: first
we study \eqref{e3.7.1000}, then we use this knowledge to study the full nonlinear equation \eqref{e3.7.35}.
  \subsection{Linear stability}
\noindent Fix $r\in(\lambda_-(c_*),\lambda_+(c_*))$. The main result in this section is 
\begin{theorem}
\label{t3.7.5}
There is a continuous linear functional $e_*$ on $B_{w_r,0}$ and $\beta_*>0$ such that, for all $u_0\in B_{w_r,0}$, we have
\begin{equation}
\label{e3.7.48}
\Vert e^{-t\mathcal{M}_{c_*}}u_0-<\!e_*,u_0\!>\varphi_{c_*}'\Vert_{w_r,1}\lesssim e^{-\beta_*t}\Vert u_0\Vert_{w_r,1}.
\end{equation}
\end{theorem}

\noindent The path that we will adopt for the proof of this result will be quite pedestrian, and will rely on computing and inverting Laplace transforms. Thus it is appropriate to recall the bases in a few lines here. Consider a function $f(t)$ defined for $t>0$ and extended by 0 on $\RR_-$. The Laplace transform of $f$, denoted by $\hat f(\lambda)$ in this section, is
\begin{equation}
\label{e3.7.50}
\hat f(\lambda)=\int_{-\infty}^{+\infty}e^{\lambda t}f(t)~dt.
\end{equation}
For this definition to have a meaning, we need to assume $f$ to behave in a reasonable manner at infinity, let  us assume that it is bounded. Then \eqref{e3.7.50} is defined as soon as ${\mathrm Re}~\lambda<0$, moreover it defines a holomorphic function. The question we ask is how to reconstruct $f$ in terms of $\hat f$. The first idea that comes to the mind is Fourier inversion: as $\hat f $ is defined on all of $\{{\mathrm Re}~\lambda<0\}$, let us pick any $\beta>0$; for all $\lambda\in\RR$ we have, as an obvious consequence of \eqref{e3.7.50}:
$
\hat f(-\beta-i\lambda)=\di\int_{-\infty}^{+\infty}e^{-i\lambda t}\bigl(e^{-\beta t}f(t)\bigl)~dt,
$
so that we have
\begin{equation}
\label{e3.7.51}
f(t)=\frac{e^{\beta t}}{2\pi}\int_{-\infty}^{+\infty}e^{i\lambda t}\hat f(-\beta-i\lambda)~d\lambda=\frac1{2i\pi}\int_{\mathrm{Re}~\lambda=-\beta}e^{-\lambda t}\hat f(\lambda)~d\lambda.
\end{equation}
A small caveat is  that the integration path goes downwards, from $i\infty$ to $-i\infty$. Formula \eqref{e3.7.51} should, by the way, be taken with a pinch of salt, as $\hat f$ may not be the fastest decreasing function at infinity: remember that $f$ has brutally been put to 0 on $\RR_-$, something that one should pay somewhere. What is realistic, however, is to assume $\hat f(-\beta+i.)\in L^2(\RR)$, for all $\beta>0$. A reasonable way to ensure that condition - this is, in fact, what we will realise - is
\begin{equation}
\label{e3.7.53}
\vert\hat f(\lambda)\vert\lesssim1/{\vert\lambda\vert}\ \hbox{as $\vert\lambda\vert\to+\infty$.}
\end{equation}
In such a configuration, \eqref{e3.7.51} should be replaced, from Plancherel's formula, by
\begin{equation}
\label{e3.7.520}
f(t)=\frac1{2i\pi}\lim_{A\to+\infty}\biggl(\int_{\substack{\mathrm{Re}~\lambda=-\beta\\\mathrm{Im}~\lambda\in[-A,A]}}e^{-\lambda t}\hat f(\lambda)~d\lambda\biggl)\ \ \hbox{in the $L^2$ sense}.
\end{equation}
If now $\hat f$ can holomorphically extended from $\RR_-\times\RR$ to $(-\infty,\beta_0)\times\RR$, for some $\beta_0>0$, with Estimate \eqref{e3.7.53} being conserved in the extension, the Cauchy theorem applied to $\hat f$ on the rectangle whose vertices are $\pm A\pm i\beta$, with $\beta\in(0,\beta_0)$ implies
$$
\int_{\substack{\mathrm{Re}~\lambda=-\beta\\\mathrm{Im}~\lambda\in[-A,A]}}e^{-\lambda t}\hat f(\lambda)~d\lambda=\int_{\substack{\mathrm{Re}~\lambda=\beta\\\mathrm{Im}~\lambda\in[-A,A]}}e^{-\lambda t}\hat f(\lambda)~d\lambda+O(\frac1A).
$$
This yields
$$
f(t)=\frac1{2i\pi}\lim_{A\to+\infty }\biggl(\int_{\substack{\mathrm{Re}~\lambda=-\beta\\\mathrm{Im}~\lambda\in[-A,A]}}e^{-\lambda t}\hat f(\lambda)~d\lambda\biggl)\ \hbox{in the $L^2$ sense},
$$
so that, by the Cauchy-Schwarz inequality we have:
$$
\Vert e^{\beta t}f\Vert_{L^2(\RR)}\lesssim\Vert \hat f(\beta+.)\Vert_{L^2(\RR)}.
$$
The $L^2$ norm in the right handside is  finite because of the preservation of Estimate \eqref{e3.7.53}. If $f$ can be shown to have more regularity, for instance $f$ globally Lipschitz on $\RR+$, interpolation proves that 
$e^{\beta't}f(t)$ is bounded, for all $\beta'\in(0,\beta)$.

\noindent Therefore we have a roadmap to the proof of Theorem \ref{t3.7.5}. As $\mathcal{M}_{c_*}\varphi_{c_*}'=0$, and as $\varphi_{c_*}'$ has, this time, the maximal decay, it is in the weigted space $B_{w_r,0}$ . Thus, it
cannot be pushed outside the domain of definition of $\mathcal{M}_{c_*}$. A special attention will have, therefore, to be given to the eigenvalue 0. However, if we want to follow the path laid out a few lines ago, the first thing we should do is 
to prove that the Laplace transform of $e^{-t\mathcal{M}_{c_*}}$ is well defined. 
\begin{proposition}
\label{p5.7.2}
Consider $v_0\in B_{w_r,1}$. Then $t\mapsto \Vert e^{-t\mathcal{M}_{c_*}}v_0\Vert_{w_r,1}$ is bounded on $\RR_+$.
\end{proposition}
\noindent{\sc Proof.} It is sufficient to prove an inequality of the form
\begin{equation}
\label{e3.7.52}
e^{-t\mathcal{M}_{c_*}}v_0(x)\lesssim\frac1{1+e^{-rx}}.
\end{equation}
Then, \eqref{e3.7.52} will also apply to $-v_0$, showing the uniform boundedness of $\Vert e^{-t\mathcal{M}_{c_*}}v_0\Vert_{w_r,0}$. As the function 
$\mathcal{M}_{c_*}v_0=-\di\partial_t\bigl(e^{-t\mathcal{M}_{c_*}}v_0\bigl)\biggl\vert_{t=0}
$
will also be in $B_{w_r,0}$, estimate \eqref{e3.7.52} applied to $\mathcal{M}_{c_*}v_0$ shows that $\Vert\partial_t\mathcal{M}_{c_*}v_0\Vert_{w_r,0}$ is bounded. This implies the boundedness of 
$t\mapsto e^{-t\mathcal{M}_{c_*}}v_0$ in $B_{w_r,1}$.

\noindent Estimate \eqref{e3.7.52} is proved by putting a supersolution $\bar v(t,x)$ of the form
\begin{equation}
\label{e5.7.53}
\bar v(t,x)=-\xi(t)\varphi_{c_*}'(x)+q(t)\inf(1,e^{-rx}),
\end{equation}
where, and this should be no surprise for the reader that has accepted to follow us to this point, we ask $\dot\xi\geq0$, while $q(t)$ should decay exponentially. We have, denoting by $H(x)$ the Heaviside function:
$$
(\partial_t+\mathcal{M}_{c_*}\overline v=\dot q(t)\inf(1,e^{-rx})+q(t)H(x)\mathcal{M}_{c_*}e^{-rx}-f'(\varphi_{c_*})\bigl(1-H(x)\bigl)q(t)-\dot\xi\varphi_c'.
$$
We have, for $x\geq0$:
$$
e^{rx}\mathcal{M}_{c_*}e^{-rx}=2\di\int_0^{+\infty}\bigl(1-\cosh{ry}\bigl)K(y)~dy+cr-f'(\varphi_{c_*})
=-D_{c_*}(r)+f'(0)-f'(\varphi_{c_*}).
$$
One is then  invited to introduce $M>0$ so that 
\begin{equation}
\label{e3.7.47}
f'\bigl((\varphi_{c_*}(x)\bigl)\leq \frac{f'(1)}2\ \hbox{for $x\leq-M$},\
\hbox{and}\ 
 \vert f'\bigl((\varphi_{c_*}(x)\bigl)-f'(0)\vert\leq-\frac{D_{c_*}(r)}2\ \hbox{for $x\geq M$,}
\end{equation}
so that, in particular, $\mathcal{M}_{c_*}e^{-rx}>0$ on $[M,+\infty)$. This number $M$ will appear over and over again until the rest of the section.

\noindent If $x\leq M$,  a sufficient condition for $(\partial_t+\mathcal{M}_{c_*})\overline v\geq 0$ is, because $\dot\xi\geq0$ and $-\varphi_{c_*}'\geq0$:
\begin{equation}
\label{e5.7.5300}
\dot q-{f'(1)}q/2\geq0.
\end{equation}
If $x\geq M$, a sufficient condition is
\begin{equation}
\label{e3.7.54}
\dot q-{D_{c_*(r)}}q/2\geq0.
\end{equation}
Thus, $q(t)=q_0e^{-\frac12\inf(-f'(1),-D_{c_*}(r))t}$ will do, for any $q_0>0$.
Finally, if $x\in(-M,M)$, it is sufficient to have
\begin{equation}
\label{e3.7.55}
\dot \xi(t)\geq\frac{\vert\dot q(t)\vert+q(t)H(x)\vert \mathcal{M}_{c_*}e^{-rx}\vert+\vert f'(\varphi_{c_*})\vert q(t)}{\min_{-M\leq x\leq M}\bigl(-\varphi_{c_*}'(x)\bigl)},
\end{equation}
which defines a bounded increasing function. Fix $\xi(0)$ and $q_0$ so that $\bar v(0,x)\geq v_0(x)$, we have constructed our upper barrier. \hfill$\Box$

\noindent For $v_0\in B_{w_r,1}$ we set $v(t,x)=e^{-t\mathcal{M}_{c_*}}v_0(x),$
extending $v$ by 0 for negative times we may define its Laplace transform $\hat v(\lambda,x)$, which satisfies
\begin{equation}
\label{e3.7.57}
\mathcal{M}_{c_*}\hat v-\lambda \hat v=v_0.
\end{equation}
Our task is now to extend $\hat v$ to a strip of the right complex half plane, and the first thing we need to worry about is $\lambda=0$. Trying to understand why $\varphi_{c_*}$ decays like the maximal exponent opened up our eyes, almost incidentally, to some realities of the operator $\mathcal{M}_{c_*}$, one of them being an easy decomposition lemma for 
$B_{w_r,0}$. Let $N(\mathcal{M}_{c_*})$ be the null space of $\mathcal{M}_{c_*}$, and $R(\mathcal{M}_{c_*})$ its range. 
\begin{lemma}
\label{l3.7.3}
We have $B_{w_r,0}=N(\mathcal{M}_{c_*})\oplus R(\mathcal{M}_{c_*})$, the sum being algebraic and topological. Additionally, there is a continuous linear functional $e_*$ on $B_{w_r,0}$ and a continuous projection $\pi_*$ on $R(\mathcal{M}_{c_*})$ such that,
for all $u\in B_{w_r,0}$ we have
\begin{equation}
\label{e3.7.49}
u=<\! e_*,u\!>\varphi_{c_*}'+\pi_*u.
\end{equation}
In addition we have $\pi_*\varphi_{c_*}'=0$, and $\pi_*\mathcal{M}_{c_*}=\mathcal{M}_{c_*}\pi_*$.
\end{lemma}
\noindent{\sc Proof.} Define $\mathcal{M}_c^0$ as \eqref{e3.7.430}. As $\mathcal{M}_{c_*}=\bigl(\mathcal{I}_{B_{w_r,0}}+\mathcal{K}_{c_*}(\mathcal{M}_{c_*}^0)^{-1}\bigl)\mathcal{M}_{c_*}^0$, with $\mathcal{K}_{c_*}(\mathcal{M}_{c_*}^0)^{-1}$ compact in $B_{w_r,0}$, the Fredholm property implies that the result will be at hand as soon as we have proved that $N(\mathcal{M}_{c_*})=N(\mathcal{M}_{c_*}^2)=\RR\varphi_{c_*}'$.  Let us prove $N(\mathcal{M}_{c_*})=\RR\varphi_{c_*}'$; if $u\in B_{w_r,1}$ solves $\mathcal{M}_{c_*}u=0$ we take the inspiration from the preceding section to claim the existence of $k>0$ such that $u\leq -k\varphi_{c_*}'$. The only difference is that, as $u$ is allowed to decay less fast that $-\varphi_{c_*}'$ at $+\infty$, we have to say something. The inspiration once again comes from the preceding section. Let $M>0$ such that \eqref{e3.7.47} is true. Pick $k_M>0$ be such that $u\leq -k_M\varphi_{c_*}'$ on $(-\infty,M]$, the case of negative infinity being investigated in the preceding section. 
On $(M,+\infty)$ we already know that, for all $r'\leq r$ we have $u(x)\leq-k_M\varphi_{c_*}'+me^{-r'x}$, $m>0$ possibly very large. Let us pick any $r'\in(\lambda_-(c_*),r)$, close enough to $r$ so that $\mathcal{M}_{c_*}e^{-r'x}>0$ on $[M,+\infty)$.
And let $m_0$ be the least $m$ such that the inequality holds.
Assume $m_0>0$, then there has to be an interior  contact point between $u$ and $-k_M\varphi_{c_*}'+m_0e^{-r'x}$. However we have, by construction:
$$
\mathcal{M}_{c_*}\bigl(u+k_M\varphi_{c_*}'-m_0e^{-r'x}\bigl)\leq 0.
$$
This is a contradiction, entailing $m_0=0$. From then on, the Linear Sliding Argument applies.

\noindent If now $u\in B_{w_r,2}$ solves $\mathcal{M}_{c_*}^2u=0$, then there is $k\in\RR$ such that $\mathcal{M}_{c_*}u=k\varphi_{c_*}'$, and we may assume $k\geq0$, so that $\mathcal{M}_{c_*}u\leq0$. And, from then on, the whole previous argument applies. The remaining properties of $\pi_*$ are elementary. \hfill$\Box$

\noindent Now that we have the main actors of the play, namely the Laplace transform $\hat v(\lambda,.)$, the linear functional $e_*$ and the projection $\pi_*$, we may start the extension. The equation satisfied by $\hat v$ is
\begin{equation}
\label{e3.7.56}
\mathcal{M}_{c_*}\hat v-\lambda\hat v=v_0,\ \ \ \hat v(\lambda,.)\in B_{w_r,1}.
\end{equation}
Notice that, although $\hat v(\lambda,.)$ exists, we still need to show that it has the correct decay in $\lambda$, even for $\mathrm{Re}~\lambda<0$. Let us, however, take the problems as they arise and let us use our newly acquired knowledge
 on the spectral decomposition of $\mathcal{M}_{c_*}$. The operator $\pi_*\mathcal{M}_{c_*}$ being continuously invertible from $\pi_*B_{w_r,1}$ to $R(\mathcal{M}_{c_*})$, any perturbation of the form 
 $\pi_*\mathcal{M}_{c_*}-\lambda\mathcal{I}_{R(\mathcal{M}_{c_*})}$, $\vert\lambda\vert$ small, is also continuously invertible between the same spaces, with holomorphic dependence in $\lambda$. Projection of \eqref{e3.7.56} onto 
 $N(\mathcal{M}_{c_*}$  and $R(\mathcal{M}_{c_*}$ readily gives a unique solution:
 \begin{equation}
 \label{e5.7.57}
 \hat v(\lambda,.)=\di\frac{<\!e_*,v_0\!>}{\lambda}\varphi_{c_*}'+\bigl(\pi_*\mathcal{M}_{c_*}-\lambda\mathcal{I}_{R(\mathcal{M}_{c_*})}\bigl)^{-1}\pi_*v_0.
 \end{equation}
 This extends $\hat v(\lambda,.)$ meromorphically in a vicinity of $\lambda=0$, and allows us to proceed to the extension of $\hat v$ past the imaginary axis.
 \begin{proposition}
 \label{p3.7.4}
 There is $\beta_0>0$ such that the Laplace transform $\hat v(\lambda,.)$ can be holomorphically extended to $\{\mathrm{Re}~\lambda<\beta_0\}\backslash\{0\}$, and meromorphically in $\{\mathrm{Re}~\lambda<\beta_0\}$. Moreover, for $\vert\lambda\vert$ large inside the new domain, we have
 \begin{equation}
 \label{e3.7.58}
 \Vert \hat v(\lambda,.)\Vert_{w_r,0}\lesssim\frac1{\vert\lambda\vert}.
 \end{equation}
 \end{proposition}
 It is important to note that \eqref{e3.7.58} only holds in the $B_{w_r,0}$ norm and not in $B_{w_r,1}$. This is sufficient to conclude to the estimate of $v(t,.)$ in $B_{w_r,0}$. A $B_{w_r,1}$ estimate follows as in the proof of Proposition \ref{p5.7.2}.
 
 \noindent The consequence of Proposition \ref{p3.7.4} is an easy proof of Theorem \ref{t3.7.5}. We indeed pick any $\beta\in(0,\beta_0)$, the number $\beta_0$ being given by the proposition. Consider $v_0\in B_{w_r,1}$;
 from the warm-up introduction we have, for every $x\in\RR$:
 $$
 e^{-t\mathcal{M}_{c_*}}=\frac1{2i\pi}\lim_{A\to+\infty}\biggl(\int_{\substack{\mathrm{Re}~\lambda=-\beta\\\mathrm{Im}~\lambda\in[-A,A]}}e^{-\lambda t}\hat v(\lambda,x)~d\lambda\biggl)\ \ \hbox{in the $L^2(\RR_t)$ sense}.
 $$
 Consider a large $A>0$, and $\rho\in(0,A)$, from the proposition we have
 $$
 \begin{array}{rll}
\di\frac1{2i\pi} \di\int_{\substack{\mathrm{Re}~\lambda=-\beta\\\mathrm{Im}~\lambda\in[-A,A]}}e^{-\lambda t}\hat v(\lambda,x)~d\lambda=&\di\frac1{2i\pi}\int_{\substack{\mathrm{Re}~\lambda=\beta\\\mathrm{Im}~\lambda\in[-A,A]}}e^{-\lambda t}\hat v(\lambda,x)~d\lambda+\frac1{2i\pi}\int_{\vert\lambda\vert=\rho}e^{-\lambda t}\hat v(\lambda,x)~d\lambda,
 \end{array}
 $$
 the integral on the circle of radius $\rho$ being counter-clockwise. From \eqref{e5.7.57} and the residue theorem, this integral is exactly $-<\!e_*,v_0\!>\varphi_{c_*}'$, while the first integral is estimated, as in the warm-up introduction, as $e^{-\beta't}$, $\beta'\in(0,\beta)$. So, let us turn to the proof of Proposition \ref{p3.7.4}. We will first prove a Fredholm property for $\mathcal{M}_{c_*}-\mathcal{I}_{B_{w_r,1}}$, which will allow us to concentrate on eigenvalues. Then we will show that the nonzero eigenvalues are far away from the imaginary axis, as well as the resolvent estimate.
 \subsubsection*{The Fredholm property}
 \noindent It should be no surprise that the operator $\mathcal{M}_{c_*}^0$, defined by \eqref{e3.7.430}, reappears. Recall that we have defined
 $\mathcal{K}_{c_*}=\mathcal{M}-\mathcal{M}_{c_*}^0$; Consider 
\begin{equation}
\label{e3.7.59}
 \beta_0=\min\bigl(-{D_{c_*}(r)}/4,-{f'(1)}/4\bigl),
\end{equation}
We are going to prove that $\mathcal{M}_{c_*}^0-\lambda\mathcal{I}$ is continuously invertible from $B_{w_r,1}$ to $B_{w_r,0}$ as soon as $\mathrm{Re}\lambda\leq\beta_0$,  readily implying that $\mathcal{M}_{c_*}$ has the Fredholm property.

\noindent To carry out this program, consider $\beta\in(0,\beta_0)$, $\lambda=\lambda_1+i\lambda_2$ with $\lambda_1\leq\beta$, $f\in B_{w_r,1}$ and $v(t,x)=e^{-t(\mathcal{M}_{c_*}^0-\lambda\mathcal{I})}f(x)$. From (a slight adaptation of) Lemma \ref{l3.7.10}, we have
$$
\begin{array}{rll}
\vert v(t,x)\vert\leq &e^{\lambda_1t}e^{-t\mathcal{M}_{c_*}^0}\vert f(x)\vert\\
\leq&e^{-(\beta_0-\beta)t}\Vert f\Vert_{w_r,0}\ \ \hbox{from (a slight adaptation of) Lemma \ref{l3.7.10}.}
\end{array}
$$
We may therefore consider
\begin{equation}
\label{e5.7.5301}
v(x)=\int_0^{+\infty}e^{-t(\mathcal{M}_{c_*}^0-\lambda\mathcal{I})}f(x)~dt,
\end{equation}
as an excellent candidate for a solution to $(\mathcal{M}_{c_*}^0-\lambda\mathcal{I})v=f$. Notice, by the way, that the above computation shows that the only solution to $(\mathcal{M}_{c_*}^0-\lambda\mathcal{I})v=0$ is zero: indeed, such a solution also solves $v(x)=e^{-t(\mathcal{M}_{c_*}^0-\lambda\mathcal{I})}v=0$, which proves $v\equiv0$ by letting $t$ go to infinity. Let us therefore go back to $v(x)$ given by \eqref{e5.7.5301}. As $f\in B_{w_r,1}$ we have, in the most classical sense:
$$
\mathcal{M}_{c_*}^0v(x)=\di\int_0^{+\infty}\mathcal{M}_{c_*}^0e^{-t(\mathcal{M}_{c_*}^0-\lambda\mathcal{I})}f(x)~dt
=-\di\int_0^{+\infty}\partial_t\bigl(e^{-t(\mathcal{M}_{c_*}^0-\lambda\mathcal{I})}f(x)\bigl)~dt
=f(x)-\lambda v(x).
$$
Have we finished and may we move on? Not quite, there is a small issue. We indeed need $f$ to be in $B_{w_r,0}$; we dodged it when we used the sub/super solution method to invert ${\mathcal{M}^0}_{c_*}$, something not possible now. So, we argue as follows:
approximate $f\in B_{w_r,0}$ by a sequence $(f_\e)_\e$ in $B_{w_r,1}$, and let $v_\e$ be given by \eqref{e5.7.5301} with $f=f_\e$. As $e^{-t({\mathcal{M}^0}_{c_*}-\lambda\mathcal{I})}$ decays exponentially in $B_{w_r,0}$, the sequence $(v_\e)$ is a Cauchy sequence, that converges in $B_{w_r,0}$ to some function $v\in B_{w_r,0}$. There is, however, more: as $v_\e$ solves the eigenvalue problem, the equation
$$
cv_\e'=-\mathcal{J}v_\e+\biggl(\gamma(x)f'(1)+\bigl(1-\gamma(x)\bigl)f'(0)\biggl)v_\e+f_\e
$$
shows that $(v_\e')_\e$ converges in $B_{w_r,0}$. Therefore $v\in B_{w_r,1}$ and solves the eigenvalue problem. And now, we may move to the next item.
 \subsubsection*{Location of nonzero eigenvalues.}
 \noindent The Fredholm property of $\mathcal{M}_{c_*}-\lambda\mathcal{I}$, for $\mathrm{Re}~\lambda<\beta_0$, allows to concentrate only on the possible eigenvalues of $\mathcal{M}_{c_*}$. We saw that $\lambda=0$ was one of them, and that it was semi-simple, that is, $B_{w_r,0}=N(\mathcal{M}_{c_*})\oplus R(\mathcal{M}_{c_*})$. The goal of this section is to prove that, in the half plane $\{\mathrm{Re}\lambda<\beta\}$, for some $\beta\in(0,\beta_0)$, there is no other eigenvalue. This last item will be proved in the last place, as a culmination of this section. The argument will not, however, be uniform in the size of $\lambda$, so that we must, as is usual, make a separate study for the large eigenvalues. This is the object of the following
 \begin{lemma}
 \label{l5.7.10}
Let $\beta_0$ be defined by \eqref{e3.7.59}. There exist $\Lambda>0$ and $\beta_*\in(0,\beta_0)$ such that, if $\mathrm{Im}~\lambda\notin(-\Lambda,\Lambda)$ and  $\mathrm{Re}~\lambda<\beta_*$, then
$\mathcal{M}_{c_*}-\lambda\mathcal{I}$ is invertible. Moreover, if $v_0$ in $B_{w_r,1}$; then we have
\begin{equation}
\label{e5.7.58}
 \Vert(\mathcal{M}_{c_*}-\lambda\mathcal{I})^{-1}v_0\Vert_{w_r,0}\lesssim{\Vert v_0\Vert_{w_r,1}}/{\vert\lambda\vert}.
 \end{equation}
 \end{lemma}
 \noindent{\sc Proof.}  Consider such a $\lambda$, we define a new operator from $B_{w_r,1}$ to $B_{w_r,0}$:
 \begin{equation}
 \label{e3.7.60}
  \mathcal{Q}_*(\lambda)=-cv'-\bigl(\lambda+f'(\varphi_{c_*})-1\bigl)v
 \end{equation}
 We claim that $  \mathcal{Q}_*(\lambda)$ is invertible for $\vert\lambda\vert$ large enough, and that 
 \begin{equation}
 \label{e3.7.61}
\Vert\mathcal{Q}_*(\lambda)^{-1}\Vert_{\mathcal{L}(B_{w_r,1},B_{w_r,0})}\lesssim1/{\vert\lambda\vert}.
 \end{equation}
 There is no mistake in \eqref{e3.7.61}, we take a function in $B_{w_r,1}$ and estimate its pre-image as $\mathcal{Q}_*(\lambda)$ by $1/{\vert\lambda\vert}$. The estimate would be wrong if we wanted to estimate the $B_{w_r,1}$ norm of the pre-image.
 
 \noindent Once \eqref{e3.7.61} is proved, the proof of the lemma follows easily, as $f$ is assumed 
 write the equation $(\mathcal{M}_{c_*}-\lambda\mathcal{I})v=f$ under the form
 $$
 \mathcal{Q}_*(\lambda)v=K*v+v_0.
 $$
 A solution $v\in B_{w_r,0}$ is found by expanding $(\mathcal{I}_{B_{w_r,0}}-\mathcal{Q}_*(\lambda)^{-1}K*.)^{-1}$ in a Neumann series, which is convergent if $\vert\lambda\vert$ is larger than some $\Lambda>0$. Note that, if $v\in B_{w_r,0}$, then $K*v\in B_{w_r,1}$. Once $v\in B_{w_r,0}$ is found, inspection of the expression \eqref{e3.7.60} of $\mathcal{Q}_*(\lambda)$ shows that, in fact, we have $v\in B_{w_r,1}$. Uniqueness should, at this point, be a simple exercise.
 
 \noindent So, let us prove the invertibility of $\mathcal{Q}_*(\lambda)$ for $\vert\lambda\vert$ large enough. An easy, but important preliminary step, is to notice that the equation $\mathcal{Q}_*(\lambda)v=0$ has only the zero solution. Indeed, a potential candidate $v(x)$ would behave like $e^{-\frac{f'(0)-1+\mathrm{Re}~\lambda}{c_*}x}$. However, the definition of $r$ gives $D_{c_*}(r)<0$, that is,
 $
 2\di\int_0^1K(y)\bigl(\cosh(r y)-1\bigl)~dy-c_*r+f'(0)<0.
 $
Define
$\beta_*=\min\biggl(\beta_0, 2\di\int_0^1K(y)\cosh(r y)~dy\biggl);
$
if $\mathrm{Re}~\lambda<\beta_*$, we have
 $
 c_*r>f'(0)-1+\mathrm{Re}~\lambda,
 $
  in other words $v\notin B_{w_r,0}$ unless $v\equiv0$. This computation allows us to write that, if $v_0\in B_{w_r,1}$, the equation $\mathcal{Q}_*(\lambda)v=v_0$ amounts to
$$
v(x)=\frac1{c_*}\int_x^{+\infty}\mathrm{exp}~\biggl(\frac1{c_*}\di\int_x^y\bigl(1-\lambda-f'(\varphi_{c_*})\bigl)~dz~v_0(y)~dy\biggl), 
$$
and estimate \eqref{e5.7.58} is obtained by integration by parts, where 
$\di\mathrm{exp}~\biggl(-\frac1{c_*}\di\int_x^yf'(\varphi_{c_*})~dz\biggl)~v_0(y)$ is differentiated, and  $\mathrm{exp}~\bigl(1-\lambda\bigl)(x-y)/{c_*}$ is integrated. This brings out the sought for $1/\lambda$, at the expense of 
sacrificing a derivative for $v_0$. \hfill$\Box$

\noindent As a consequence of this last lemma, one obtains a positive $\beta$ such that the spectrum of $\mathcal{M}_{c_*}$ inside the half plane $\{\mathrm{Re}~\lambda\leq\beta\}$ is composed of the sole eigenvalue $\lambda=0$, which is,
let us recall it, semi-simple. Indeed, each element of the resolvent set of $\mathcal{M}_{c_*}$ is isolated, while there is no spectrum outside a horizontal strip $\{\vert\mathrm{Im}~\lambda\vert\geq\Lambda\}$. The $L^2$ estimate on the resolvent is the last ingredient for the validity of Theorem \ref{t3.7.5}.
\subsection{Global exponential stability}\label{s5.2.2}
\noindent As the reader will notice, there is a parallel between this section and the preceding one, especially because some arguments will have a flavour very similar to those pertaining to the analysis of the linear semigroup. In particular, there is a strict parallel between multiplication in the linear setting and translation in the full nonlinear context. We will first apply Theorem \ref{t3.7.5} to the Cauchy Problem for Equation \ref{e3.7.35} with an initial datum $u_0$ close to a wave in $B_{w_r,1}$, and prove that the solution converges exponentially fast, in $B_{w_r,1}$, to a translate of the wave. The best part of the section will then be devoted to proving that a solution that is initially far from the set of travelling waves will approach it in $B_{w_r,1}$ for large times. This will imply the theorem.
\subsubsection*{Local  nonlinear stability}
\noindent The precise statement is the following. Note that the normalisation of $\varphi_{c_*}$ in the statement is unconsequential, as the picture can be translated by any quantity. 
\begin{theorem}
\label{t3.7.6}
Consider an intial datum $u_0$ to \eqref{e3.7.35}, and let $\varphi_{c_*}$ be the solution of \eqref{e3.7.1} with $c=c_*$ and $\varphi_{c_*}(0)=1/2$. Let $u(t,x)$ be the solution of \eqref{e3.7.35} with initial datum $u_0$.
There are $\omega>0$ and $\e_0>0$ such that, for all $\e\in[0,\e_0]$, if $\Vert u_0-\varphi_{c_*}\Vert_{w_r,1}\leq\e$, there is $x_\e=O(\e)$ such that 
\begin{equation}
\label{e3.7.65}
\Vert u(t,.)-\varphi_{c_*}(.+x_\e)\Vert_{w_r,1}=O(e^{-\omega t}).
\end{equation}
\end{theorem}
\noindent {\sc Proof.} If $x_\e$ the expected asymptotic translation, we may translate the whole picture by $x_\e$ and 
 look for a solution $u(t,x)$ under the form
$u(t,x)=\varphi_{c_*}(x)-\xi(t)\varphi_{c_*}'(x)+v(t,x),
$
where, in order to fix $\xi(t)$, we ask $v(t,x)$ to be in $R(\mathcal{M}_{c_*})$ for all time; that is
$<\!e_*,v(t,x)\!>=0$, for all $t\geq0$.
%
Write 
$$f(u)=f(\varphi_{c_*})+f'(\varphi_{c_*})\bigl(-\xi(t)\varphi_{c_*}'+v\bigl)v+g(x,\xi,v),
$$
with
$
g(x,\xi,v)=O(\xi^2+v^2),\ g'(v)=O(\vert\xi\vert+\vert v\vert)$ as $(\xi,v)\to0$,
so that \eqref{e3.7.65} becomes
\begin{equation}
\label{e3.7.69}
-\dot\xi\varphi_{c_*}'+v_t+\mathcal{M}_{c_*}v=g(x,\xi,v).
\end{equation}
Projecting \eqref{e3.7.69} successively on $-\varphi_{c_*}'$ and on $R(\mathcal{M}_{c_*})$ yields
\begin{equation}
\label{e3.7.690}
\dot\xi=-<\!e_*,g(.,\xi,v)\!>,\ \ \ \ v_t+\mathcal{M}_{c_*}v=\pi_*g(x,\xi,v).
\end{equation}
We expect $\xi(t)$ and $\Vert v(t,x)\Vert_{w_r,1}$ to decay exponentially fast, \eqref{e3.7.690} may be merged into the following system for the unknown $\bigl(\xi(t),v(t,x)\bigl)$:
$$
\xi(t)=\di\int_s^{+\infty}<\!e_*,g\bigl(.,\xi(\tau),v(\tau,.)\bigl)\!>d\tau,\ \ \ \
v(t,.)=e^{-t\mathcal{M}_{c_*}}v_0+\di\int_0^te^{-(t-s)\mathcal{M}_{c_*}}\pi_*g\bigl(.,\xi(s),v(s,.)\bigl)~ds.
$$
As $g$ is quadratic in $\xi$ and $v$, this system can easily be solved. Consider $\beta>0$ such that $\Vert e^{-t\mathcal{M}_{c_*}}\Vert_{\mathcal{L}(R(\mathcal{M}_{c_*})}=O(e^{-\beta t})$, and pick $\omega<\beta$. As 
$g$ is quadratic both in $\xi$ and $v$, the Implicit Functions Theorem applies to the system, in the space of all functions 
$(\xi,v)\in C(\RR_+)\times C(\RR_+,B_{w_r,1})$ such that  $\di\sup_{t\geq0}e^{\omega t}\bigl(\vert\xi(t)\vert+\Vert v(t,.)\Vert_{w_r,1}\bigl)<+\infty).$
The relation $\xi(0)=\di\int_0^{+\infty}<\!e_*,g\bigl(.,\xi(\tau),v(\tau,.)\bigl)\!>d\tau$ allows for the reconstruction of $x_\e$. \hfill$\Box$

\subsubsection*{Compactness}
\noindent We are going to prove here that the solution is nontrivial in the reference frame moving exactly like $c_*t$, with no unbounded correction of any sort. The reader is welcome to check 
why the argument that will follow fails in the Fisher-KPP case. 
\begin{theorem}
\label{t3.7.4}
For every $r\in\bigl(\lambda_-(c_*),\lambda_+(c_*)\bigl)$, the family $(u(t,.))_{t>0}$ is precompact in $B_{w_r,0}$.
\end{theorem}
This result suggests that we will have, at some point, to prove some regularity properties for $u(t,.)$. We have already seen in the Fisher-KPP case that such properties are not automatic, 
we will see that, here too, we have to say something. Also, precompactness in the uniform norm is something one has to fight for, even if one has regularisation. 

\noindent The first step, however, is to show that $u(t,.)$ is uniformly bounded in $B_{w_r,0}$. This is the object of the following proposition, that parallels exactly Proposition \ref{p5.7.2}.
\begin{proposition}
\label{p3.7.5}
There exist two constants $\xi^-<\xi^+$, and a function $q(t)\lesssim e^{-\omega t}$ with $\omega>0$, such that
\begin{equation}
\label{e3.7.72}
\varphi_{c_*}(x-\xi^-)-\frac{q(t)}{1+e^{-rx}}\leq u(t,x)\leq\varphi_{c_*}(x-\xi^+)+\frac{q(t)}{1+e^{-rx}}.
\end{equation} 
\end{proposition}
\noindent{\sc Proof.} The idea is to trap the solution $u(t,x)$ to \eqref{e3.7.35} between a pair of sub and super solutions of the form
$$
\underline u(t,x)=\varphi_{c_*}\bigl(x-\xi^-(t)\bigl)-q(t)\inf(1,e^{-r(x-\xi^-(t))})\ \ \ \ 
\overline u(t,x)=\varphi_{c_*}\bigl(x-\xi^+(t)\bigl)+q(t)\inf(1,e^{-r(x-\xi^+(t))});
$$
the positive function $q(t)$ and the functions $\xi^{\pm}(t)$ being adjusted to the best of our interests. Let us construct a supersolution, the principle for a subsolution being exactly the same. We have, setting $\zeta=x-\xi^+(t)$:
$$
\mathcal{N}_{c_*}\overline u(t,x)=-\dot\xi^+(t)\varphi_{c_*}(\zeta)+\dot q(t)\inf(1,e^{-r\zeta})+\bigl(\mathcal{J}-c_*\partial_\zeta\bigl)\inf(1,e^{-r\zeta})-f\bigl(\overline u(t,x)\bigl)+f\bigl(\varphi_{c_*}(\zeta)\bigl).
$$
We will impose $\dot\zeta^+\geq0$. Re-introduce $M>0$ such that \eqref{e3.7.47}  holds, and distinguish as is by now standard, the following theree zones:
\begin{itemize}
\item[---] Zone $\{\zeta\leq -M\}$: we have, because $\dot\zeta^+\geq0$:
$
\mathcal{N}_{c_*}\overline u(t,x)\geq\dot q(t)-{f'(1)}q(t)/2,
$
so that it suffices to choose $q(t)\geq q_0e^{tf'(1)/2}$.
\item[---] Zone $\{\zeta\geq M\}$: we have this time, by a simple application of the mean value theorem to $f\bigl(\overline u(t,x)\bigl)-f\bigl(\varphi_{c_*}(\zeta)\bigl)$:
$ 
\mathcal{N}_{c_*}\overline u(t,x)\geq\bigl(\dot q(t)-{D_{c_*}(r)}q(t)/2\bigl)e^{-rx},
$
so that it suffices to choose $q(t)\geq q_0e^{tD_{c_*}(r)/2}$.
\item [---] Zone $\{-M\leq\zeta\leq M\}$: we now choose $q(t)=q_0e^{-\omega t}$, $\omega={q_0}\min\bigl(-f'(1)/2,-D_{c_*}(r)/2\bigl)$. Then we have
$\mathcal{N}_{c_*}\overline u(t,x)+\dot\xi^+(t)\varphi_{c_*}\gtrsim -e^{-\omega t}.$
This entails the choice of $\xi^+(t)$, an increasing function tending exponentially to its   limit.
\end{itemize}
As $\xi^+(t)$ and $q(t)$ are chosen to tend exponentially in time, respectively to a positive limit and to 0, it only remains to choose $\xi^+(0)$ large enough and $q_0$ small enough to ensure
that $u_0(x)\leq\overline u(0,x)$, thus proving the proposition. \hfill $\Box$

\noindent If the diffusion was given by a second order elliptic operator, this would immediately yield the same sort of bounds for the derivatives. Here, and this is one of the main differences, the best we can hope for is 
preservation of smoothness. This is the goal of the next step.
\begin{proposition}
\label{p3.7.6}
For every $k\in{\mathbb N}$, we have $\di\sup_{t\geq0}\Vert u(t,.)\Vert_{B_{w_r,k}}<+\infty$.
\end{proposition}
\noindent{\sc Proof.} The case $k=0$ is treated in Proposition \ref{p3.7.5}, so that one may start investigating the case $k=1$. Let us examine the solution $u(t,x)$ in the steady reference frame, that is, we leave the reference frame moving
like $c_*t$ for one moment. Set $v(t,x)=u(t,x)$, the equation for $v(t,x)$, obtained by differentiating the equation for $u$ with respect to $x$, is written as
$v_t+\bigl(1-f'(u)\bigl)v=K'*u.
$
This implies
\begin{equation}
\begin{array}{rll}
\label{e3.7.76}
v(t,x)=&e^{-t+\int_0^tf'(u(s,x))~ds}v_0'(x)+\di\int_0^te^{-(t-s)+\int_s^tf'(u(\tau,x))~d\tau}K'*u(s,.)~ds\\
:=&v_1(t,x)+v_2(t,x)
\end{array}
\end{equation}
As all nontrivial phenomena happen in the range $x\sim c_*t$, it is a safe bet to organise things around this position. Let, once again, $M>0$ be defined by \eqref{e3.7.47}. From Proposition \ref{p3.7.5} we may enlarge it so that there is $t_0>0$ such that
\begin{equation}
\label{e5.7.76}
\begin{array}{rll}
f'\bigl(u(t,x)\bigl)\leq&\di{f'(1)}/2\ \ \hbox{for $t\geq t_0$ and $x\leq c_*t-M$,}\\
f'\bigl(u(t,x)\bigl)-f'(0)\leq&-\di{D_{c_*}(r)}/2\ \ \hbox{for $t\geq t_0$ and $x\geq c_*t+M$.}e^{-r(x-c_*t)}
\end{array}
\end{equation}
\begin{itemize}
\item [---] Zone $\{x\geq c_*t-M\}$. The term $v_1(t,x)$ is the simplest and allows a good understanding of the mechanism at work: we have 
$$
\begin{array}{rll}
\vert v_1(t,x)\vert\leq&\Vert v_0\Vert_{w_r,1}\mathrm{exp}~\biggl(-rx-t+\di\int_0^tf'\bigl(u(s,x)\bigl)~ds\biggl)\\
\leq&e^M\Vert v_0\Vert_{w_r,1}\mathrm{exp}~\biggl[-t\biggl(-D_{c_*}(r)-\di\frac1t\int_0^t\bigl[f'\bigl(u(s,x)\bigl)-f'(0)\bigl]ds\biggl)\biggl]\\
\leq&e^{M+t_0\Vert f'\Vert_\infty}\Vert v_0\Vert_{w_r,1}e^{-r(x-c_*t)}e^{tD_{c_*}(r)/2}.
\end{array}
$$
Therefore $v_1(t,x)$ decays exponentially fast in time in this zone. This gives the correct intuition to treat the Duhamel term $v_2(t,x)$. Indeed we have
$$
\big\vert K'*u(s,x)\big\vert\lesssim e^{-r(x-c_*s)},
$$
so that the integrand in $v_2(t,x)$ decays like $e^{-r(x-c_*t)}e^{tD_{c_*}(r)/2}.$ Integrating in $x$ and taking the estimate for $v_1$ into account yields
\begin{equation}
\label{e3.7.77}
\vert \partial_xu(t,x)\vert\lesssim e^{-r(x-c_*t)}\ \hbox{for $x-c_*t\geq -M$.}
\end{equation}
\item [---] Zone $\{x-c_*t\leq-M\}$. We may revert to the usual frame moving like $c_*t$, and directly move to $\vert v(t,x)$. From the first inequality of \eqref{e5.7.76}  an equation for $v(t,x)$ is
\begin{equation}
\label{e3.7.78}
\begin{array}{rll}
\vert v\vert_t+\mathcal{J}\vert v\vert -c_*\vert v\vert_x+\di\frac{\vert f'(1)\vert}2\vert v\vert\leq&0\ \ (x\leq-M)\\
\vert v(t,x)\vert\ \hbox{bounded on}& \RR_+\times[-M,-M+1].
\end{array}
\end{equation}
The maximum principle entails 
$\vert v(t,x)\vert\leq\sup_{t>0}\Vert v(t,.)\Vert_{L^\infty([-M,-M-1])},
$ which ends the proof of the $B_{w_r,1}$ bound.
\end{itemize}
The cases $k\geq 2$ are treated as the case $k=1$, by induction and successive differentiations in $x$. \hfill$\Box$

\noindent We may extract a little more information from equation \eqref{e3.7.78}, once again with the aid of the maximum principle. Consider $\rho>0$ such that 
$
\bigl(\mathcal{J}-c_*\partial_x-{f'(1)}/4\bigl){e^{\rho x}}\geq-{f'(1)}e^{\rho x}/8,
$
and $\gamma\in(0,-{f'(1)}/8)$. Propositions \ref{p3.7.5} and \ref{p3.7.6} entail a corollary that will be the starting point of our final effort to prove the exponential stability of the travelling waves.
\begin{corollary}
\label{c3.7.1}
For all $k\geq1$ we have $\bigl\vert\partial_x^ku(t,x)\big\vert\lesssim\inf\bigl(e^{\rho x}+e^{-\gamma t},e^{-\omega t}.e^{-rx}\bigl).$
\end{corollary}
\noindent{\sc Proof.} The only point left to prove is the first argument in the infimum, but it follows readily from inspection of inequation \eqref{e3.7.78}, and differentiations of the equation for $u$. \hfill$\Box$
\subsubsection*{Final argument}
From Corollary \ref{c3.7.1}, the set $\{u(t,.),\ t>0\}$ is precompact in $B_{w_r,k}$ for all $k\in\mathbb{N}$. Let $\omega(u_0)$ be the set of all limits of sequences $\bigl(u(t_n,.)\bigl)_n$, it is a compact set of $B_{w_r,k}$ for all $k$, which is in addition
comprised in the interval $\bigl[\varphi_{c_*}(.-\xi^-),\varphi_{c_*}(.-\xi^+)\bigl]$. It is also invariant with respect to the Cauchy Problem \eqref{e3.7.35}. Finally, as the just mentionned Cauchy Problem does not have a regularising effect, it is not irrelevant to stress the fact that the $B_{w_r,k}$ distances in $\omega(u_0)$ are somehow equivalent. Here is the form that will be useful to us.
\begin{lemma}
\label{l3.7.6}
There are constants $C_1,\gamma>0$ such that, for all $(\psi_1,\psi_2)$ in $\omega(u_0)$ we have
\begin{equation}
\label{e3.7.81}
\Vert\psi_1-\psi_2\Vert_{w_r,1}\leq C_1\Vert\psi_1-\psi_2\Vert_\infty^\gamma.
\end{equation}
\end{lemma}
\noindent{\sc Proof.} Consider two such $\psi_1$ and $\psi_2$, we may assume, without loss of generality, that $\Vert\psi_1-\psi_2\Vert_\infty\leq1$.  Let us show \eqref{e3.7.81} for $k=0$, the inequality for  $k=1$ being deduced from $k=0$ by interpolation, and the fact that there is a common $C^2$ bound for all element of $\omega(u_0)$.  So, it all boils down to examining what happens for $x\geq0$. Pick any small $\delta>0$, if $x\leq\delta\mathrm{ln}~\bigl(\Vert\psi_1-\psi_2\Vert_\infty\bigl)^{-1}$ we have
$$
e^{rx}\bigl\vert\psi_1(x)-\psi_2(x)\bigl\vert\leq\frac{\bigl\vert\psi_1(x)-\psi_2(x)\bigl\vert}{\Vert\psi_1-\psi_2\Vert_\infty^{r\delta}}\leq\Vert\psi_1-\psi_2\Vert_\infty^{1-r\delta}.
$$
For $x\geq\delta\mathrm{ln}~\bigl(\Vert\psi_1-\psi_2\Vert_\infty\bigl)^{-1}$, recall that $\psi_1$ and $\psi_2$  are trapped between $\varphi(.-\xi^-)$ and $\varphi(.-\xi^+)$,  so that
$$
e^{rx}\bigl\vert\psi_1(x)-\psi_2(x)\bigl\vert\lesssim e^{-(\lambda_+(c_*)-r)x}
\leq\Vert\psi_1-\psi_2\Vert_\infty^{\delta((\lambda_+(c_*)-r)}.
$$
And so, inequality \eqref{e3.7.81} is valid as soon as $1-\delta r>0$ and $\gamma=\inf\bigl(1-\delta r,\delta((\lambda_+(c_*)-r)\bigl)$. \hfill$\Box$

\noindent Consider the smallest $\xi$ such that $\omega(u_0)\leq\varphi_{c_*}(.+\xi)$. We baptise it $\xi_m$, it exists and is finite, as  it lies between $\xi^-$ and $\xi^+$. 
Set
\begin{equation}
\label{e3.7.82}
d=\inf_{\psi\in\omega(u_0)}\Vert\varphi_{c_*}(.+\xi_m)-\psi\Vert_\infty.
\end{equation}
If $d=0$, Lemma \ref{l3.7.6} and the compactness of $\omega(u_0)$ in $B_{w_r,1}$ implies $\varphi_{c_*}(.+\xi_m)\in\omega(u_0)$, which in turns implies, from the definition of $\omega(u_0)$ and Theorem \ref{t3.7.6}, that $u(t,.)$ converges exponentially in time towards $\varphi_{c_*}(.+\xi_m)$. So, we have to prove that the case $d>0$ cannot happen. Suppose, therefore, that this is the situation.

\noindent We first claim that, for all $a>0$, there exists $d_a>0$ such that, for all $\psi\in\omega(u_0)$, we have
$\di \inf_{x\in[-a,a]}\bigl(\varphi_{c_*}(x-\xi_m)-\psi(x)\bigl)\geq d_a.
$
Suppose indeed that this is not the case, the compactness of $\omega(u_0)$ in $B_{w_r,1}$ entails the existence of $a_0>0$, $\psi_0\in\omega(u_0)$ such that 
$\psi_0$ and $\varphi_{c_*}(.-\xi_m)$ have a contact point that we denote $x_0$.
For commodity, let us denote by $\mathcal{S}(t)\psi$ the solution of the Cauchy Problem \eqref{e3.7.35} emanating from $\psi$. The point is that, if $\psi\in\omega(u_0)$, then $\mathcal{S}(t)\psi$ is uniquely defined for negative $t$. The function $v(t,x)=\varphi_{c_*}(x-\xi_m)-\mathcal{S}(t)\psi(x)$ is a nonnegative function that solves an equation of the form
$$
v_t+\mathcal{J}v-c_*v_x+a(t,x)v=0,\ \ -1\leq t\leq0,\ x\in\RR.
$$
At $t=0$, we have $\mathcal{J}v(0,x_0)=0$, which implies $v(0,.)\equiv0$ from Proposition \ref{p2.2.3}. This is in contradiction with $d>0$.

\noindent At this stage, we may reintroduce the solution $u(t,x)$ to the Cauchy Problem \eqref{e3.7.35}, and reintroduce $M>0$ defined by \eqref{e5.7.76}, with $x-c_*t$ replaced by $x$, as we are already in the moving frame. From the definition of $d_M$, there is
$t_0$ large enough so that, for all $t\geq t_0$ we have
\begin{equation}
\label{e3.7.83}
\inf_{x\in[-M,M]}\bigl(\varphi_{c_*}(x-\xi_m)-u(t,x)\bigl)\geq {d_M}/2.
\end{equation}
This entails the existence of $\delta_0>0$ such that, for all $\delta\in[0,\delta_0]$ we have $\varphi_{c_*}\bigl(x-(\xi_m-\delta)\bigl)\geq u(t,x)$ for $t\geq t_0$ and $x\in[-M,M]$. Let us examine what happens outside this interval; pick $\delta\in[0,\delta_0]$ and
set
$
v(t,x)=\biggl(u(t,x)-\varphi_{c_*}\bigl(x-(\xi_m-\delta)\bigl)\biggl)^+.
$
For $x\leq-M$, we have
$$
\partial_tv+\mathcal{J}v-c_*v_x-{f'(1)}v/2\leq0,\ \ \ v(t,x)=0\ \hbox{for $t\geq t_0$ and $x\in[-M,-M+1]$}.
$$
This implies, by the maximum principle, the existence of $q_0>0$ such that $v(t,x)\geq -q_0e^{tf'(1)/2}$. For $x\geq M$ we have
$$
\partial_tv+\mathcal{J}v-c_*v_x-\bigl(f'(0)-{D_{c_*}(r)}/2\bigl)v\leq0,\ \ \ v(t,x)=0\ \hbox{for $t\geq t_0$ and $x\in[M-1,M]$}.
$$
Again by the maximum principle, and  even if it means enlarging $q_0$, we have  $v(t,x)\geq -q_0e^{tD_{c_*}/2}$.
Therefore we have 
$\di\liminf_{t\to+\infty}v(t,x)=0,
$ uniformly in $x\in\RR$. This implies that, for all $\psi\in\omega(u_0)$ we
have
$
\varphi_{c_*}\bigl(x-(\xi_m-\delta)\bigl)\geq\psi(x).
$
This contradicts the minimality of $\xi_m$, meaning that the situation $d=0$ is the only possible. This ends the proof of Theorem \ref{t3.7.3}.

\section{The large time behaviour when $c_*=c_K$}\label{s5.3}
\noindent In this section, we will not only assume $c=c_K$, but that the wave with bottom speed $c_K$ has the minimal decay, that is, $\varphi_{c_K}(x)=(x+k)e^{-\lambda_Kx}
+O(e^{-(\lambda_*+\delta)x})$ as $x\to+\infty$. We saw in Chapter \ref{TW}, Theorem \ref{t3.7.100}, that such waves exist and are abundant. In particular, for commodity reasons, we will use $f(u)=f'(0)u+g(u)$, with $g(0)=g'(0)=0$, $g>0$ on $(0,\theta)$, $g<0$ on $(\theta,+\infty)$ with, say, $g(u)\lesssim -u^2$ for $u$ large. We will also arrange that $f(1)=0$ and $f'(1)<0$. This is again not contradictory with the existence of waves with minimal decay.  

\noindent In order to keep some unity in the notations, we will from now on denote by $c_*$ the bottom speed, having in mind that it is equal to $c_K$. It is not out of place to recall the notations that we will use. Let us write $\MM_{c_*}=\JJ-c_*\partial_x-f'(\varphi_{c_*})$, we have
\begin{equation}\label{e5.3.1}
\MM_{c_*}=\JJ-c_*\partial_x-f'(0)-g'(\varphi_{c_*})=\LL_*-g'(\varphi_{c_*}).
\end{equation}
\noindent Our assumptions will lead to a result similar to Theorem \ref{t2.4.1} in Chapter \ref{short_range}.
\begin{theorem}
\label{t5.2.1}
There exists $x_\infty\in\RR$ such that $\di\lim_{t\to+\infty}\biggl(u(t,x)-\varphi_{c_*}\bigl(x-c_*t-\frac3{2\lambda_*}\mathrm{ln}~t+x_\infty\bigl)\biggl)=0,$
uniformly in $x\in\RR_+$.
\end{theorem}
We will see that, similarly to what happens for Fisher-KPP nonlinearities, the dynamics of the solution in the diffusive zone ahead of $x\sim c_*t$ is what drives the whole dynamics. 

\noindent There is, however, an important technical issue. Indeed, the Fisher-KPP assumption allowed an easy retrieval of the information from this diffusive zone to the region moving like $c_*t-3\mathrm{ln}~t/2\lambda_*$, as the convergence to the wave with the proper shift resulted from a change of unknown and a simple spectral consideration on a linear operator. This is not going to work here, as the assumption we are making on $f$ is exactly orthogonal to the Fisher-KPP assumption. In order to compensate for this, we will have to fully use the assumption that we have a wave of minimal decay, as well as a not totally standard analysis of the linear operator $\MM_{c_*}$. This analysis will have the merit of highlighting what general assumptions are really at work in the convergence process, as well as the crucial role played by the diffusive region ahead of $x\sim c_*t$.
\subsection{The operator $\MM_{c_*}$ outside its comfort zone}
\noindent Of interest to us will be the solutions of the eigenvalue problem
\begin{equation}
\label{e3.5.2}
\MM_{c_*}u=\mu u\ (x\in\RR), 
\end{equation}
with $\mu>0$. What justifies the title of the section is that we are investigating properties of  the operator $\MM_{c_*}$ that do not result from a well charted application of statndard functional analysis results, like the Krein-Rutman theorem, that was at the basis of almost all the developments of the preceding section.
\begin{theorem}
\label{t5.3.2}
There is $\mu_0>0$ such that Problem \eqref{e3.5.2}, together with the normalisation condition that $u$ has unit mass over $\RR_-$, has a $C^1$ branch of solution $(u_\mu)_{0\leq\mu\leq\mu_0}$, in the norm of uniform convergence. Moreover, if $\lambda_\pm(\mu)$ are the decay exponents at $+\infty$ given by \eqref{e2.2.27}, there are two constants $\alpha^\pm_\mu$ such that 
\begin{equation}
\label{e5.3.5}
u_\mu(x)=\alpha^+_\mu e^{-\lambda_+(\mu)x}+\alpha_\mu^-e^{-\lambda_-(\mu)x}+O\bigl(e^{-(\mathrm{Re}~\lambda_\pm(\mu)+\delta)x}\bigl).
\end{equation}
\end{theorem}
Before delving into the proof of this result, we propose a short pause to show the reader that this result is perfectly logical. Replace indeed $\JJ$ by $-\partial_{xx}$, and assume that $f'$ is constant, and equal to $f'(0)$ on $[0,\theta_0]$. Let $\varphi_{c_*}(x)$ be the wave with bottom speed such that $\varphi_{c_*}(0)=\theta_0$. A standard application of the Krein-Rutman theorem shows us that, for small $\mu>0$, the problem
$$
\bigl(-\partial_{xx}-c_*\partial_x-f'(\varphi_{c_*})\bigl)u=\mu u\ (x<0),\quad u(0)=1
$$
has a unique solution $u_\mu(x)$.  Indeed, we have $\MM_{c_*}\varphi_{c_*}'=0$, so that, by the Krein-Rutman theorem, the bottom Dirichlet eigenvalue of $\MM_{c_*}$ on $\RR_-$ is positive. Then, solving the linear equation $\bigl(-\partial_{xx}-c_*\partial_x-f'(0)-\mu\bigl)u=0$ on $\RR_+$, with $u(0)=1$ and $u'(0)=u_\mu'(0)$ yields a solution that behaves exactly as in \eqref{e5.3.5}, thus yielding an overly easy proof of the theorem. The trouble here is that we do not have here any Cauchy theory to help us; getting around this inconvenience is what is going to occupy us now.

\noindent The first ingredient is, similarly to the above argument, the well-posedness of the  boundary value problem
\begin{equation}
\label{e5.3.6}
\MM_{c_*}u=\mu u\ \ \bigl(x\in(-\infty,M)\bigl),\quad u(x)=u_M(x)\ \ (M< x< M+1),
\end{equation}
where $M>0$ is a (possibly large) constant, while $u_M$ is a continuous function of $[M,M+1]$.
\begin{proposition}\label{p5.3.2}
For $M>0$ there is $\mu_M>0$ such, for all $\mu\in[0,\mu_M]$,  problem \eqref{e5.3.6} has a unique solution, denoted by $\TT_M(\mu) v_M$. Moreover, for all $\mu\in[0,\mu_M]$, the operator $\TT_M$ is compact from $C\bigl([M,M+1]\bigl)$ to $C\bigl((-\infty,M]\bigl)$.
\end{proposition}
\noindent{\sc Proof.} It suffices to prove the proposition for $\mu=0$, which will imply the result for every small $\mu$. Then, we notice that, for $\gamma>0$ large enough, the function $\gamma\varphi_{c_*}'$ is a super-solution to \eqref{e5.3.6} that is above $v_M$ on $[M,M+1]$, while $-\gamma\varphi_{c_*}'$  is a sub-solution to \eqref{e5.3.6} that is above $v_M$. This entails the existence of a solution $\TT_Mv_M$ that lies between $-\gamma\varphi_{c_*}'$ and $\gamma\varphi_{c_*}'$. $C^1$ regularity of $\TT_Mv_M$ is obvious, moreover, the decay of $\varphi_{c_*}'$ to 0 at $-\infty$ entails the sought for compactness. Uniqueness is proved by the classical linear sliding method. \hfill$\Box$

\noindent The other important ingredient will be the representation of the solutions of
\begin{equation}
\label{e3.5.3}
\LL_*u=h,\quad x\in\RR,
\end{equation}
where $h\in L^1(\RR)$. This last requirement for $h$ is in contrast with the beginning of the proof of the exponential behaviour \eqref{e2.3.3} and \eqref{e2.3.4} of the waves in Theorem \ref{t2.3.1}, where we were faced with the (relatively minor) inconvenience that the right handside was not in $L^1$. However, the function $h$ will have no reason to be positive, so that the beautiful dodging argument devised in Lemma \ref{l2.3.1} will not be useful to us anymore. 

\noindent Recall the definition \eqref{e2.3.6} of $\DD_c$:
\begin{equation}
\label{e5.3.4}
\DD_c(\xi)=2\int_{\RR_+}K(z)\bigl(1-\cos(z\xi)\bigl)~dz-ic\xi-f'(0).
\end{equation}
The main information, that we do not deem as entirely obvious, is the location of the zeroes of $\DD_{c_*}$. Once they are known, lightheartedly shifting the line of integration becomes possible.
\begin{proposition}\label{p5.3.1}
Consider $c\geq c_*$, and let $\xi=\omega+i\lambda$, with $\lambda<\lambda_-(c_*)$. Then $\DD_{c}(\xi)\neq0$.
\end{proposition}
\noindent{\sc Proof.} If $\DD_{c}(\omega+i\lambda)=0$, then $u(x)=e^{ix\xi}$ solves $\LL_cu=0$, that is, $u$ is a complex linear wave. 



\noindent  Consider $c>c_*$ and let us make the homotopy 
$K_\e(x)=\di\frac\beta{\e^3} K(\di\frac{x}\e),
$
 for $\e$ between 0 and 1.  Let us denote by $\DD_{\e,c}$ the function $\DD_c$ of \eqref{e5.3.4} with $K$ replaced by $K_\e$. 
 We adjust the constant $\beta>0$, more properly baptised $\beta_\e$,   so that, for all $\e>0$, the function $\DD_{\e,c}(\lambda)$ has $\lambda_-(c)$ as its least real zero.  Note that, from the construction of $c_*$ in Chapter \ref{Cauchy}, the family $(\beta_\e)_\e$ is bounded and bounded away from 0: we have indeed $\DD_{0,c}(\xi)=d\beta \xi^2-ic\xi-f'(0)$, where $d$ is the second moment of $K$.  
Notice that the kernel $K_\e$ is proposed, albeit with a different power of $\e$, in many problems of this book;  the reader should therefore not be surprised to see that it is put to good use at some point. Let $c_\e>0$ such that The sequence $\bigl(\DD_{\e,c}\bigl)_{0\leq\e\leq1}$ is a sequence of entire functions such that $\DD_{1,c}=\DD_c$ and whose power series are nonsingular with respect to $\e$, hence uniformly convergent with respect to $\e$, and with all derivatives, on every compact set in $\xi$.
 
\noindent For large $A>0$ , let $R_A$ be the rectangle $\{0\leq\omega\leq A, \ 0\leq\lambda\leq\lambda_-(c)\}$. We are going to show that no zero of $\DD_{\e,c}$ else than $\lambda_-(c)$ may enter $R_A$. Still because of the asymptotics of $\DD_{\e,c}$ as $\e\to0$, we may find a small disk $D_0$ around $\lambda_-(c)$ in which $\DD_{\e,c}$ has no other zero than $\lambda_-(c)$. Rouch\'e's theorem will imply that $\DD_{\e,c}$ has the same number of zeros as $\DD_(0,\e)$ inside $R_A\cup D_0$, that is, only one. We have just seen that no zero can enter from the top of $R_A$, that is, the segment $[-A,A]\times\{\lambda_-(c)\}$. Trivially, no zero can enter from the bottom or the left side, so we are left to examine whether the equation $\DD_{\e,c}(\xi)=0$ has solutions of the form $\xi=A+i\lambda$, $0\leq\lambda\leq\lambda_-(c)$. 

\noindent Assume first that $\e A\leq\pi/2$, and $u(x)=e^{i\xi x}$. Let $\overline u(x)$ coincide with $\mathrm{Re}~u(x)$ on $[0,\pi/2\omega]$ and be zero everywhere else. As $\e A\leq\pi/2$, it solves $\LL_c\underline u\leq0$ from Proposition \ref{p2.2.2}. A suitable translation of $\underline u$, still denoted $\underline u$, is below $e^{-\lambda_-(c)x}$ with a contact point. This contradicts the strong maximum principle. Assume now that $\pi/2\leq\e A\leq\sqrt A$. Equating $\mathrm{Re}\bigl(\DD_{\e,c}(\xi)\bigl)$ to 0 yields the equation
$$
\int_\RR K_\e(\e x)\biggl(1-\cos(\e Ax)\cosh\bigl(\e\lambda_-(c)x\bigl)\biggl)~dx=\e^2\bigl(-c\lambda_-(c)+f'(0)\bigl),
$$
so that, taking into account the that $\e A\geq\pi/2$:
$
\di\int_\RR K(x)\bigl(1-\cos(\e Ax)\bigl)~dx=O(\e^2).
$
However, if $\e A$ is in an interval of the form $[2k\pi,(2k+1/2)\pi]$ or $[(2k+3/2)\pi,2(k+1)\pi]$, we have $1-\cos(\e Ax)\geq 1/2$ on an interval of the size at least $\pi/(3\e A)$. Otherwise, we have $1-\cos(\e Ax)\geq 1$. Consequently, there is $a>0$ independent of $\e$ such that $\di \int_\RR K(x)\bigl(1-\cos(\e Ax)\bigl)~dx\geq a$. This is impossible if $A$ is large enough. Finally, assume $\e\geq1/\sqrt A$. 
Equating $\mathrm{Im}\bigl(\DD_{\e,c}(\xi)\bigl)$ to 0 yields the equation
$$
\e A=\int_\RR K_\e(x)\sin(\e Ax)\sinh\bigl(\e\lambda_-(c)x\bigl)~dx=\lambda_-(c)\int_\RR K_\e(x)\sin(\e Ax)~dx+O(\e).
$$
As $\e A\geq\sqrt A$, and as the Riemann-Lebesgue lemma asserts that $\di\lim_{\omega\to+\infty}\int_\RR K(x)\sin(\omega x)~dx=0$, we once again meet an impossibility. This proves the result.
\hfill$\Box$

\noindent{\sc Proof of Theorem \ref{t5.3.2}.} It consists, essentially, in putting the equation $\MM_{c_*}u=\mu u$ in an equivalent form. We start from 
$\LL_*u=g'(\varphi_{c_*})u+\mu u,
$
 and we pick $M>0$ large, we will freeze it in due time. Write
 $$
 u=\bigl(\un_{(-\infty, M-1]}+\un_{[M,M+1]}+\un_{[M+1,+\infty)}\bigl)u\\
 =\bigl(I+\TT_M)\un_{[M,M+1]}\bigl)+\un_{[M+1,+\infty)}u.
 $$
 Taking the Fourier transform of both sides of the equation, inverting, and shifting the line of integration yields
\begin{equation}
\label{e5.3.10}
\begin{array}{rll}
2\pi u(x)=&\EE_-\biggl[g'(\varphi_{c_*})\biggl(\bigl(I+\TT_M)\un_{[M,M+1]}\bigl)+\un_{[M+1,+\infty)}\biggl)u\biggl](x)\\
&+\EE_+\biggl[\lambda,g'(\varphi_{c_*})\biggl(\bigl(I+\TT_M)\un_{[M,M+1]}\bigl)+\un_{[M+1,+\infty)}\biggl)u\biggl](x),
\end{array}
\end{equation} 
the integral operators $\EE_\pm$ being given by \eqref{e2.3.8} and \eqref{e2.3.12}, the expression being valid as long as the function $\xi\mapsto \DD_c(\xi+i\lambda)$ does not vanish on $\RR$. From proposition \ref{p5.3.2}, we may choose $\lambda=\mathrm{Re}~\lambda_+(\mu)=\mathrm{Re}~\lambda_-(\mu)$, to alleviate the notation we denote this quantity $\bar\lambda(\mu)$. Let $\e_M>0$ denote the supremum of $g'(\varphi_{c_*})$ over $[M,+\infty)$, we have
$$
\biggl\Vert\EE_-\biggl[g'(\varphi_{c_*})\un_{[M+1,+\infty)}u\biggl]+\EE_+\biggl[\lambda,g'(\varphi_{c_*})\un_{[M+1,+\infty)}u\biggl]\biggl\Vert_{L^\infty([M+1,+\infty)}\lesssim\e_M\Vert u\Vert_{L^\infty([M+1,+\infty))}.
$$
If $M$ is large enough, this entails the existence of a continuous operator $\KK_M(\mu)$, with norm $\lesssim\e_M$, such that 
$$
\un_{[M+1,+\infty)}u=(\frac{I}{2\pi}+\KK_M(\mu))\biggl(\EE_-\biggl[g'(\varphi_{c_*})\bigl(I+\TT_M)\un_{[M,M+1]}\bigl)u\biggl]+
\EE_+\biggl[\lambda,g'(\varphi_{c_*})\bigl(I+\TT_M)\un_{[M,M+1]}\bigl)u\biggl]\biggl).
$$
Inserting the above inequality into \eqref{e5.3.10}, and considering it for $x\in[M,M+1]$ yields an integral equation of the form $u=\QQ_M(\mu)u$, the operator $\QQ_M$ being compact from $C([M,M+1)$ into itself. In order to show that there are solutions, it suffices to show uniqueness, and, as $\mu$ is small, it suffices to show that the problem
 $
 u-\QQ_M(0)u=0$, with the normalisation $\di\int_{\RR_-}\TT_Mu=0,
 $
  has zero as its unique solution. However, this  is equivalent to studying the equation $\MM_{c_*}u=0$, with the normalisation condition $\di\int_{\RR_-}u=0$. From \eqref{e5.3.10} we infer that $u$ decays at most like $xe^{-\lambda_*x}$ as $x\to+\infty$, so that the linear sliding argument yields $u\equiv0$. The behaviour \eqref{e5.3.5} is guaranteed by a last examination of \eqref{e5.3.10}, together with Lemma \ref{l2.3.2}. \hfill$\Box$
  
  \noindent Two remarks are in order. The first one is that  the asymptotic expansion \eqref{e5.3.5}, while factually true, will not be the most convenient for the purpose that we will pursue in the next section, namely to devise barriers to the left of the diffusive zone. These barriers will be devised with the aid of 
 $u_{1/t^\delta}$, $t$ large, so that we really have to tailor an asymptotic expansion that will be valid in a smaller zone, but that will bear some uniformity in $\mu$. 
 
 \noindent This leads us to the second remark. In the statement of Proposition \ref{p5.3.1}, we have been rather elusive about how the branch $(u_\mu)_{0\leq\mu\leq\mu_0}$ approaches the limit $\mu=0$. Let us declare that $u_{\mu=0}=-\varphi_{c_*}'$.
 Of course there is convergence  in the $C^1$ topology, there is no mistake about it. However we need to know a little more, especially as we have to deal with uniform expansions at $x=+\infty$. Examination of equation \eqref{e5.3.10} shows that we may bootstrap the $C^1$ convergence into the following: for all $\e>0$ there is $C_\e>0$ (certain to blow up as $\e\to0$) such that 
 \begin{equation}
 \label{e5.3.21}
 \bigl\vert u_\mu(x)+\varphi_{c_*}(x)\bigl\vert\leq C_\e\mu e^{-(\lambda_*-\e)x},\ \ \ x>0.
 \end{equation}

\noindent So, let us revert once again to the argument that led to the integral equation \eqref{e2.3.310}
 for the behaviour of $\varphi_{c_*}$ at $+\infty$. As there is no conceptual difference between applying it to $\varphi_{c_*}$ or $u_\mu$, we just extrapolate it, and obtain the following 
 expressions for $\alpha_\mu^\pm$ in \eqref{e5.3.5}:
$$
\alpha_\mu^\pm=\frac1{D_{c_*}'\bigl(\lambda_\pm(\mu)\bigl)}\int_\RR e^{\lambda_\pm(\mu)y}g'\bigl(\varphi_{c_*}(y)\bigl)u_\mu(y)~dy.
$$
 The exponential remainder is uniform in $\mu$. At this stage, it is worth pointing out that the signs are consistent with those in \eqref{e2.3.310}, as $g$ has now become $-g$.
 
 \noindent Let us rewrite the  formula \eqref{e2.2.27}, that even a very serious reader would be forgiven to have forgotten at this stage, as 
 $\lambda_\pm(\mu)=\lambda_*\pm ia\sqrt\mu+b\mu+O(\mu^{3/2})$, with $a>0$,$ b>0$. And so, we have
 $$
 D_{c_*}'\bigl(\lambda_\pm(\mu)\bigl)=D_{c_*}''(\lambda_*)\bigl(\lambda_\pm(\mu)-\lambda_*\bigl)+O\bigl(\lambda_\pm(\mu)-\lambda_*\bigl)^2
 =\pm iaD_{c_*}''(\lambda_*)\sqrt\mu+O(\mu).
 $$
Recall that, as $\mu\to0$, the limit of $u_\mu$ is $-\varphi_{c_*}'$. So  we have, from \eqref{e5.3.21} and the fact that $g'(0)=0$:
 $$
 \begin{array}{rll}
\di \int_\RR e^{\lambda_\pm(\mu)y}g'\bigl(\varphi_{c_*}(y)\bigl)u_\mu(y)~dy=&-\di\int_\RR e^{\lambda_*y}g'\bigl(\varphi_{c_*}(y)\bigl)\varphi_{c_*}'(y)~dy+O(\mu)\\
=&\di\lambda_*\int_\RR e^{\lambda_*y}g\bigl(\varphi_{c_*}(y)\bigl)~dy+O(\mu):=\alpha_*+O(\mu).
\end{array} 
$$
Reverting to the full formula \eqref{e5.3.5}, we have
$$
\begin{array}{rll}
u_\mu(x)=&\di\frac{\alpha_*e^{-\lambda_*x}}{ia\sqrt\mu D_{c_*}''(\lambda_*)}\bigl(e^{(ia\sqrt\mu+O(\mu))x}-e^{(-ia\sqrt\mu +O(\mu))x}\bigl)+O(\sqrt\mu)e^{-\lambda_*x}\\
=&\di\frac{2\alpha_*e^{-\lambda_*x}\sin(a\sqrt\mu x)}{a\sqrt\mu D_{c_*}''(\lambda_*)}+O\bigl(\sqrt\mu e^{-\lambda_*x}\bigl).
\end{array}
$$
Recall that, from formula \eqref{e2.3.3100}, we have
$$
\varphi_{c_*}(x)=-\frac{e^{-\lambda_*x}}{D_{c_*}''(\lambda_*)}\int_{\RR} e^{\lambda_* y}(x-y)g(\varphi(y))~dy+O(e^{-(\lambda_*+\e_0/2)x})=\frac{\alpha_*}{\lambda_*D_{c_*}''(\lambda_*)}e^{-\lambda_*x}\bigl(x+O_{x\to+\infty}(1)\bigl).
$$
As our basic assumption here is that $\varphi_{c_*}$ decays to 0 as $xe^{-\lambda_*x}$, we have $\alpha_*>0$, something that was not entirely obvious from the mere inspection of its expression. This leads us to the estimate that we were anxious to discover: for 
every $\nu\in(0,1/2)$, there is $C_\nu>0$ (also certain to blow up as $\nu\to0$) such that
\begin{equation}
\label{e5.3.210}
\bigl\vert u_\mu(x)-\bigl(2\alpha_*/D_{c_*}''(\lambda_*)+O(\sqrt\mu)xe^{-\lambda_*x}\bigl)\bigl\vert\leq C_\nu e^{-\lambda_*x},\ \ \ \ 0\leq x\leq\mu^{-1/2+\nu}.
\end{equation}
We are not totally done yet with the study of $u_\mu$, as we need to say something about its derivatives as $x\to+\infty$. This is the object of the last proposition of this section.
\begin{proposition}
\label{p5.3.3}
For all $\nu>0$ there is $C_\nu>0$ such that, for all sufficiently small $\mu$ we have, for $0\leq x\leq\mu^{-1/2+\nu}$:
$\bigl\vert \partial_xu_\mu(x)\bigl\vert\leq C_\nu u_\mu(x)$ and $\bigl\vert \partial_xu_\mu(x)\bigl\vert\leq C_\nu x u_\mu(x).$
\end{proposition}
\noindent{\sc Proof.} The $x$-derivative $\partial_xu_\mu$ satisfies 
$\bigl(\mathcal{M}_{c_*}-\mu\bigl)\partial_xu_\mu=g''(\varphi_{c_*})\varphi_{c_*}'u_\mu$. As the coefficient $g''(\varphi_{c_*})\varphi_{c_*}'$ decays to 0 exponentially fast as $x\to+\infty$,
  the argument leading to the estimate \eqref{e5.3.210} applies word by word to $\partial_xu_\mu$. For $\partial_\mu u_\mu$ there is a small catch. Indeed the equation is this time
  $\bigl(\mathcal{M}_{c_*}-\mu\bigl)\partial_xu_\mu=u_\mu$, so that the integrals of the type $\di\int_\RR e^{\lambda_\pm(\mu)y}u_\mu(y)~dy$ supposed to yield the component of $\partial_\mu u_\mu$ on 
  the characteristic exponentials do not converge anymore. This is, however, a minor inconvenience. Examination of \eqref{e2.3.12} in Chapter \ref{TW} shows that the coefficients $\alpha_\mu^\pm$ are limits, when they exist, of integrals between $-\infty$ and $x$, in other words we should rather consider the coefficients
  $$
  \alpha_\mu^\pm(x)=\frac1{D_{c_*}'\bigl(\lambda_\pm(\mu)\bigl)}\int_{-\infty}^x e^{\lambda_\pm(\mu)y}\biggl(\mu+g'\bigl(\varphi_{c_*}(y)\bigl)\biggl)u_\mu(y)~dy,
  $$
  and the extra factor $x$ in the estimate of $\partial_\mu u_\mu$ comes from the integral $\di\int_{-\infty}^x e^{\lambda_\pm(\mu)y}u_\mu(y)~dy$. \hfill$\Box$
  
\subsection{How boundary data on expanding intervals drive the convergence}
\noindent Let is now believe that the reference frame moving like $c_*t-3\mathrm{ln}~t/2\lambda_*$ is the correct one. In this reference frame, the initial problem \eqref{e3.7.2} is reformulated into
\begin{equation}
\label{e5.3.100}
u_t+\JJ u-\bigl(c_*-\frac3{2\lambda_*}\mathrm{ln}~t~\partial_x\bigl)u=f(u),
\end{equation}
and of interest to us in this section will be the following issue. Assume that \eqref{e5.3.100} is solved, instead of the whole line, in an expanding interval $(-\infty,Y(t))$, where $Y(t)$ goes to infinity as $t$ grows, and $u$ is made to roughly coincide with $\varphi_{c_*}(x)$ on $\bigl(Y(t),Y(t)+1\bigl)$. Will $u(t,x)$ converge to $\varphi_{c_*}$, uniformly on $\bigl(-\infty,Y(t)\bigl)$? We will see that the answer is yes, and that it very much hinges on the construction of $u_\mu$ of the preceding section, and on its precise asymptotic behaviour at infinity. 

\noindent We  denote by $\varphi_{c_*}$ the only wave solution whose asymptotic behaviour at infinity coincides with \eqref{e5.3.11}, that is:
 $\varphi_{c_*}(x)\sim_{x\to+\infty} (x+k_*)e^{-\lambda_*x}.
$
Here comes the main result of this section:.
\begin{theorem}
\label{t5.3.3}
Pick $Y(t)=t^\gamma-3\mathrm{ln}~t/2\lambda_*:=Y_\gamma(t)$, with $\gamma\in(0,1/2)$ and let $u(t,x)$ be the solution of \eqref{e5.3.100} with Dirichlet data
\begin{equation}
\label{e5.3.11}
u(t,x)=xe^{-\lambda_*x},\quad Y_\gamma(t)\leq x\leq Y_\gamma(t)+1, 
\end{equation}
and such that $\di\lim_{x\to-\infty}u_0(x)=0$.
Then we have $\di\lim_{t\to+\infty}\sup_{x\leq Y_\gamma(t)}(1+e^{\lambda_*x}/x)\bigl\vert u(t,x)-\varphi_{c_*}(x)\bigl\vert=0$.
\end{theorem}

\noindent{\sc Proof.} The strategy looks at first sight quite different from that leading to the convergence proof to $\varphi_{c_*}$ when $c_*>c_K$. Indeed, we need to prove that the dynamics of the solution is slaved to what happens around $Y_\gamma(t)$. In other words, we have lost the degrees of freedom ensured by the translation. On the other hand,  the existence of waves with speed $c_*=c_K$, for nonlinearities slightly above or below $f$, is a welcome additional degree of freedom. We are going to play with it in order to compare the solution with translates of waves.

\noindent So, for small  $\e>0$, let $\varphi_{c_*}^\e$ be the wave given by Corollary \ref{c3.4.1}, and that has the behaviour \eqref{e5.3.11} at infinity. We always may assume that $u_0$ sits below a translate of $\varphi_{c_*}^\e$, something we could not really assume for $\varphi_{c_*}$. This allows us to define, for all $s>0$,  e the smallest $y$ such that 
\begin{equation}
\label{e5.3.22}
u(s,x)\leq\varphi_{c_*}^\e(x-y)\ \ \ \ \ \hbox{for $x\leq Y_\gamma(s)$.}
\end{equation}
Call it $y_s$. We have $u(t,x)\leq\varphi_{c_*}^\e(x-y)$ for $t\geq s$ and $x\leq Y_\gamma(t)$, so that $s\mapsto y_s$ is a nonincreasing function. Let $y_\infty$ its limit as $s\to+\infty$, we wish to prove that $y_\infty=0$. Let us assume the contrary, and consider a small $\kappa>0$ (at least small enough to ensure the validity of the subsequent considerations), and let $t_\kappa$ be the first $s>0$ such that $y_s\leq y_\infty+\kappa$. From \eqref{e5.3.22} and the behaviour \eqref{e5.3.210} of $u_\mu$, we may pick a small $q_\kappa^0>0$ such that 
$
u(t_\kappa,x)\leq \varphi_{c_*}^\e\bigl(x-(y_\infty-2\kappa)\bigl)+q_\kappa^0u_{t_\kappa^{-2\gamma}}(x).
$
We hope very much that an inequality of this sort will persist for $t\geq t_\kappa$, let us set
$$
\overline u(t,x)=\varphi_{c_*}^\e\bigl(x-(y_\infty-2\kappa)\bigl)+q(t)u_{t^{-2\gamma}}(x),
$$
and compute, using the fact that $\mathcal{M}_{c_*}u_{t^{-2\gamma}}=t^{-2\gamma} u_{t^{-2\gamma}}$:
$$
\mathcal{N}_{c_*}\overline u=\frac3{2\lambda_*t}(\varphi_{c_*}^\e)'+\biggl(\dot q+\bigl(\frac{q}{t^{2\gamma}}+\frac3{2\lambda_*t}\frac{\partial_xu_{t^{-2\gamma}}}{u_{t^{-2\gamma}}}-\frac{2\gamma}{t^{1+2\gamma}}\frac{\partial_\mu u_{t^{-2\gamma}}}{{u_{t^{-2\gamma}}}}\bigl)q\biggl)u_{t^{-2\gamma}}+O(q^2u_{t^{-2\gamma}}^2).
$$
Here, $\varphi_{c_*}^\e$ is evaluated at $x-(y_\infty-2\kappa)$, something we have omitted in order to keep the size of the identity under control. From Proposition \ref{p5.3.3}, we have $\partial_xu_{t^{-2\gamma}}/u_{t^{-2\gamma}}=O(1)$ on $[0,Y_\gamma(t)]$, while we have $\partial_\mu u_{t^{-2\gamma}}/u_{t^{-2\gamma}}=O(x)=O(t^{2\gamma})$ on the same interval. Also, we have $(\varphi_{c_*}^\e)'/u_{t^{-2\gamma}}=O(1)$. Thus, in order to have $\mathcal{N}_{c_*}\overline u\geq0$, it suffices, at least if $t_\kappa$ is large enough, to take
\begin{equation}
\label{e5.3.50}
q(t)=\frac1{t^{1-4\gamma}}.
\end{equation}
Now, as we have $\overline u(t,x)\geq u(t,x)$ on $[Y_\gamma(t),Y_{\gamma}(t)+1]$ by assumption, we have $u(t,x)\leq \overline u(t,x)$ for $t\geq t_\kappa$. And so, there exists $t_\kappa'>t_\kappa$ such that $u(t,x)\leq \varphi_{c_*}(x-y_\infty-\kappa)$, thus contradicting the minimality of $y_\infty$. 

\noindent In the same way, if $z_s$ is the smallest $z$ such that $u(s,x)\geq\varphi_{c_*}^{-\e}(x+z)$, we have $\di\lim_{s\to+\infty}z_s=0$. This is enough to prove the proposition. \hfill$\Box$
\subsection{Convergence in the diffusive zone and final argument}
\noindent In a way that is quite similar to what happens with KPP nonlinearities, what drives the whole process is the dynamics at the entrance of the diffusive zone, that is, $x\sim t^\gamma$ with $\gamma<1/2$. And,
 indeed, similarly to Theorem \ref{t2.4.1}
in Section \ref{s4.5.1}. Write the full equation \eqref{e3.7.2} in the reference frame moving withspeed $c_*$:
\begin{equation}
\label{e5.3.200}
u_t+\mathcal{J}u-c_*u_x-f'(0)u=g(u);
\end{equation}
the main result is the
\begin{theorem}
\label{t5.3.4}
There is $p_\infty>0$ such that, for all $\gamma\in(0,1/2)$ we have:
$
u(t,t^\gamma)\sim\di\frac{p_\infty}{t^{3/2-\gamma}}e^{-\lambda_*t^{\gamma}}.
$
\end{theorem}
As we will see, Theorem \ref{t5.3.4} bears a lot of common points with  Theorem \ref{t2.4.1} in the preceding chapter \ref{short_range}. The main difference is that, in the latter, we used the KPP assumption $f'(u)\leq f'(0)$. This assumption is of course invalid here, and what we use instead is the existence of a wave with bottom speed and minimal decay. We will see that it plays an important role in the whole process.

\noindent Let us also point out that, once Theorem \ref{t5.3.4} is proved, the main theorem of the section, that is, Theorem \ref{t5.2.1} becomes immediate in view of the preparations already done. Indeed, it follows from Theorem \ref{t5.3.3} with $Y(t)=Y_\gamma(t)$ with $\gamma$ small enough, and comparison with solutions of \eqref{e5.3.100}-\eqref{e5.3.11}.

\noindent{\sc Proof of Theorem \ref{e5.3.4}.} the scheme of the proof is quite similar to the proof of Theorem \ref{t2.4.1}: (i) identify the correct reference frame up to $O(1)$ terms, that is, find out the unbounded terms in the expansion of $X_\theta(t)$, (ii) gradient estimates in the reference frame moving with speed $c_*$, (iii) final convergence argument. The difference is that what happens away from
the diffusive zone is not entirely unrelated to what happens at its beginning, or inside. And so, throughout the proof of the theorem, there will be an interplay between the solution $u(t,x)$ of \eqref{e5.3.200} and its reduced form $v(t,x)=e^{\lambda_*x}u(t,x)$, that solves
\begin{equation}
\label{e5.3.201}
v_t+\mathcal{I}_*v=e^{\lambda_*x}g(e^{-\lambda_*x}v).
\end{equation}
The nonlinear term is indeed small for large $x$, as $g$ is at least quadratic at the origin. 

\noindent So, the first step is to prove the weaker estimate, valid for all $\delta\in(0,1/2)$:
\begin{equation}
\label{e5.3.202}
{e^{-\lambda t^\delta}}/{t^{3/2-\delta}}\lesssim u(t,x)\lesssim{e^{-\lambda t^\delta}}/{t^{3/2-\delta}}.
\end{equation}
While interesting in itself and yielding, eventually, the location of the level sets $X_\theta(t)$, its main utility was to unveil the structure of the reduced unknown function $v(t,x)$ in the diffusive zone. Recall that we may, indeed, write it under the form
$v(t,x)=\di\frac1tw_t(\frac{x}{\sqrt t}),\quad x\geq t^\delta,
$
where $w_t(\eta)$ is a locally bounded function such that $\eta\lesssim w_t(\eta)\lesssim \eta$ for $\eta>0$. This is exactly what we aim at proving now.
For this, we will trap $u(t,x)$ between a sub and a super solution bearing a lot of common points with those devised in Theorem \ref{t2.4.1}, the main difference being that the wave $\varphi_{c_*}$ will come into play, in order to make up for the assumption $g(u)\geq0$ for small $u$, that goes against us at the beginning of the diffusive zone. Let us start with the super-solution.  Consider $\e_0$ and $\bar x_0$ such that $u_0\leq\varphi_{c_*}^\e(x-\bar x_0)$; we may as well assume that $\bar x_0=0$, by translating $u_0$ by the amount $\bar x_0$. Consider, also, a smooth, nonnegative,  even function $\rho(x)$, supported in $[-1,1]$, $a>0$.  For $M>0$, let $\rho_{*,M}$ be the odd extension of $\rho(.-M)$. Taking $\Vert\rho\Vert_\infty$  suitably small allows us to assume 
\begin{equation}
\label{e5.3.205}
e^{-t\II_*}\rho_{*,M}(x)\leq x/4t^{3/2},\quad x\geq 1.
\end{equation}

\noindent Pick, once and for all, a small number $\gamma>0$. We organise the constants $\alpha$, $\beta$ and $\delta$ as follows:  we choose $\beta<\gamma$ and $\alpha>1/3$, close to $1/3$. For $x\geq(1+t)^\delta$ we define
$$\bar v^+(t,x)=\xi^+(t)e^{-t\II_*}\rho_{*,M}(x)+\frac1{(1+t)^{3/2-\beta}}\cos~\frac{x}{(1+t)^\alpha}\un_{[0,3\pi (1+t)^\alpha/2}(x),
$$
and its counterpart $\bar u^+(t,x)=e^{-\lambda_*x}\bar v^+(t,x)$. 
The construction is similar to that in Chap \ref{short_range}, Section \ref{s2.5.2}. We always may assume, even if it 
means enlarging $M$, that $e^{-t\II_*}\rho_{*,M}(x)\geq0$ for $x\geq (1+t)^\delta$ with $0<\delta<\beta$. For $x\leq t^{2\gamma}$, consider 
\begin{equation}
\label{e5.3.2040}
\bar u^-(t,x)=\varphi_{c_*}^{\e_0}\bigl(x+\frac3{2\lambda_*}\mathrm{ln}~t\bigl)+q(t)u_{(1+t)^{-2\gamma}}(x),
\end{equation}
with $q(t)\propto(1+t)^{1-4\gamma}$, as constructed in the preceding section. Both $\bar u^+$ and $\bar u^-$ are super-solutions to \eqref{e5.3.200} in their domain of definition. 
The asymptotic behaviour of $u_\mu$ described in Proposition \ref{p5.3.3} and the normalisation of $\varphi_{c_*}^\e$ entail 
$u^-(t,x)\geq\bigl(1+o_{x\to+\infty}(1)\bigl)e^{-\lambda_*x}/t^{3/2}.
$ 
Therefore, from estimate \eqref{e5.3.205}, the graph of $x\mapsto \bar u^-(t,x)$ intersects that of $x\mapsto\bar u^+(t,x)$ at (at least) two points; call $\bar x^+(t)$ the furthest
 to the right. The next one is at least $(1+t)^{\delta}$ far from $\bar x^+(t)$, thus at distance at least 1. Notice that, at least for times of order 1, we may play with the initial datum to ensure this property. Subsequently, we define $\bar u(t,x)$ as $\bar u^+(t,x)$ if $x>\bar x^+(t)$ and $\bar u^-(t,x)$ if $x<x^+(t)$; the function $\bar u$ is a super-solution that may be put above $u_0$ at $t=0$ by translation.
 
\noindent A sub-solution is devised  more easily, as it is enough to define
$\underline u(t,x)=\di e^{-\lambda_*x}\bigl(e^{-t\II_*}\rho_{*,-M}\bigl)_+,
$
where, this time, the graph of $\rho_*$ has been translated by the amount $-M$. Choosing $M$ large enough is sufficient to ensure $\underline u(0,x)\leq u_0(x)$. All in all we have $\underline u(t,x)\leq u(t,x)\leq\overline u(t,x)$, which is sufficient to prove  
\eqref{e5.3.202}. As a consequence, given the heat kernel like behaviour of $e^{-t\II_*}$, the expression \eqref{e5.3.202} for $v(t,x)$ is valid.

\medskip
\noindent Let us turn to gradient estimates, the objective being to estimate $v_x$. Let us first examine $w=u_x$,  in the initial reference frame, that is, the equation is written without the advection terms.
The Duhamel formula for $w$ was given in equation \eqref{e3.7.76}, let us recall it:
$$w(t,x)=e^{-t+\int_0^tf'(u(s,x))~ds}w_0'(x)+\di\int_0^te^{-(t-s)+\int_s^tf'(u(\tau,x))~d\tau}K'*u(s,.)~ds\:=w_1(t,x)+w_2(t,x),
$$
and the only term of real interest is $w_2(t,x)$. We first use the classical argument that $u(t,x)$ only spends a finite time in the region $\{f'(u)\geq f'(0)\}$, that is: 
$\di\int_s^tf'\bigl(u(\tau,x)\bigl)~d\tau\leq f'(0)+O(1).
$
Then, we apply our knowledge of $u(t,x)$ to infer $u(s,x)\lesssim \bigl(1+(x-c_*s)^+\bigl)e^{-\lambda_*(x-c_*s)}/s^{3/2}$. Elementary algebra then yields:
$$
w_2(t,x)\lesssim \bigl(1+(x-c_*t)^+\bigl)\int_0^te^{(t-s)\bigl(f'(0)-1-\lambda_*c_*\bigl)(t-s)}~ds,
$$
Recall $f'(0)-c_*\lambda_*-1=-2\di\int_0^1K(x)\cosh x~dx$, so that
$w_2(t,x)\lesssim \bigl(1+(x-c_*t)_+\bigl)e^{-\lambda_*(x-c_*t)};
$
 of course, $-w_2(t,x)$ may be estimated in the same fashion. For $x$ far to the left of $x-c_*t$, the zero order term $f'(u)$ is $<0$, so that the maximum principle implies an upper bound for $w$. At this stage, one only needs (i) go back to the reference frame moving like $c_*t$, and (ii)
 to reproduce the proof of Theorem \ref{t4.4.4} to prove gradient estimates for $v(t,x)=e^{\lambda_*x}u(t,x)$.

\medskip
\noindent To finish the proof of the theorem, it is now possible to repeat the argument in the proof of Theorem \ref{t2.4.1} because all its ingredients are there. For all $s>0$, we define $\bar u^+_s(t,x)$ for all $t\geq s$ as
$$
\bar u^+_s(t,x)=\xi_s^+(t)e^{-(t-s)\II_*}\mathcal{S}_{*\e}^+v_*(s,x)+\frac1{(1+t)^{3/2-\beta}}\cos~\frac{x}{(1+t)^\alpha}\un_{[0,3\pi (1+t)^\alpha/2]}(x),
$$ 
with $\xi_s^+(t)$ defined as in \eqref{e2*.6.10}. We also define $\bar u^-(t,x)$ as in \eqref{e5.3.2040}, the wave $\varphi_{c_*}^{\e_0}$ being suitably translated so that the graph of $x\mapsto\underline u^-(t,x)$ as a rightmost intersection with that of $\bar u^+(t,x)$ located at a point of the order $t^\gamma$. A sub-solution is defined as 
$
\underline u(t,x)=\di\biggl(e^{-(t-s)\II_*}\mathcal{S}_{*\e}^-v_*(s,x)\biggl)_-;
$
by construction we have $\underline u(s,x)\leq\bar u(s,x)$. Thus we have $\underline u(t,x)\leq u(t,x)\leq\bar u(t,x)$ for $t\geq s$, and the arguent now runs just as in the proof of Theorem \ref{t2.4.1}\hfill$\Box$
\section{Conclusion}
\noindent We have seen that our basic model \eqref{e1.1.1} has very different behaviours according to whether the nonlinearity $f$ is of the KPP type or the ZFK type in the case $c_*>c_K$. However, there are deep analogies between both the result and the path towards its proof that we have chosen. Indeed, in both cases, there is uniform convergence to the wave with bottom speed, the difference being in the unbounded terms of the reference frame. To prove the result, one makes the same three steps: (i) identifying the correct reference frame (Theorem \ref{t4.4.4} in the KPP case and Proposition \ref{p3.7.5} in the ZFK case with $c_*>c_K$), (ii) compactness estimates (in both cases, gradient estimates; Theorem \ref{t4.4.4} in the KPP case and Proposition \ref{p3.7.6} in the ZFK case with $c_*=c_K$), (iii) proof of convergence through, in both cases, a Harnack type argument, the argument being a little more well-hidden in the KPP case.

\section{Bibliographical elements, open questions}\label{s3.8}

\subsection*{The large time behaviour when $c>c_K$ (Section \ref{s3.4})}
The original idea for the proof of convergence to travelling waves is due to Fife-McLeod \cite{FML}, given for the diffusion equation
\begin{equation}
\label{e5.7.1000} 
\partial_tu-\partial_{xx}u=f(u),
\end{equation}
with a bistable $f$. It is hard to overstate the importance of this work.  As we have seen for the Fisher-KPP model, finding out the reference frame where the solution will be nontrivial is not an easy task. However, for this sort of problems with bistable nonlinearities, the discovery of Fife and Mc Leod is that one can use translates of the wave, possibly corrected by exponential perturbations, to locate this reference frame. The sub and super solutions that enable to trap the solution between two waves is reproduced from them, and this idea deserves a very careful analysis. In my opinion it is, according to the principles of the current Newspeak, a ground-breaking work. The impact of \cite{FML} goes far beyond bistable reaction-diffusion equations; in particular, reading Section \ref{s4.3} of Chapter \ref{short_range} in the light of this work should convince the reader that there is a deep filiation between two apparently disconnected methods of investigation.

\noindent To prove the convergence to the wave with bottom speed for a solution emanating from an initial datum whose support in $\RR_+$ is compact, I have chosen to follow the Fife-Mc Leod scheme. As said before, the sub  and super solution argument is adapted from \cite{FML}.  Convergence of a subsequence of the form $(u(t_n,.))_n$ was proved in \cite{FML} by compactness, and a Lyapunov function argument. Compactness follows from parabolic regularisation and the trapping between two wave-like sub and super-solutions, the first ingredient being unavailable here. What is, however, possible to prove is the preservation of the gradient bound, I have exposed it in  detail. The underlying idea is that $\varphi_{c_*}(c-c_*t)$ does not spend too much time in the region $\{f'(u)\geq0\}$; it was already noticed  by several authors; the reader may for instance consult  Coville, Davila, Martinez \cite{CDM}, or Shen, Shen \cite{ShSh}.
 The existence of a Lyapunov functional in our context is probably true, but I have chosen not to develop this aspect here. Instead,  I have  presented an alternative argument, borrowed from \cite{Roq-monot}, tailored to cases when looking for such a functional may look especially not promising.  
 
\noindent The weigthed spaces needed for the study of the linearised operator were introduced, I believe, by Sattinger \cite{Sat} in the context of reaction-diffusion systems of the second order, and have become an important tool in the study of the stability of travelling waves. Many arguments, but not all, carry over to nonlocal equations. As a matter of fact, these questions have been extensively studied at the beginning of this century, and I do not pretend to present new material here. The reader can consult Mallet-Paret \cite{MP} for a study of the Fredholm property for this sort of nonlocal operators. As for the the final stability argument, I have developped it from an exercise of Chapter 5 in the important reference book of Henry \cite{Hy}.  For more information on the stability of reaction-diffusion waves in heterogeneous environments, one may read the review of Xin \cite{Xin-review}.

\noindent While the convergence to a wave is, as far as I know, not written anywhere for ZFK nonlinearities, the convergence to a wave for nonlocal equations with bistable nonlinearity was proved by Chen   \cite{Chen}. His argument should be adaptable here, at some technical cost. The point in \cite{Chen} is to avoid compactness arguments, as regularisation is not true here. The argument that he uses consists in refining the sub and super solution argument so as to make them get closer and closer to each other, by making use of a Harnack type inequality. Such a property, contrary to regularisation, is true in some form.  As the Harnack inequality is, in its most basic form, a combination of strong maximum principle and regularity, the author devises a clever integral formulation that allows him to go forward. One of the main points in the proof that I have chosen to develop is to show that, in fact, one may resort to regularity as well as a more classical Harnack type argument. This leads to a proof that looks slightly less technical at some points, and that also lays out the many common points that the nonlocal setting shares with the case of a diffusion given by the Laplacian.
\subsection*{Large time behaviour when $c=c_K$ (Section \ref{s5.3})}
\noindent This case, while quite  important, has been much less studied than the case $c_*>c_K$, and the only results concerning the Cauchy Problem I am aware of concern the case with Gaussian diffusion, that is, $\JJ$ is replaced by $\partial_{xx}$. 
The earliest significant contribution that I know of is a very nice paper by 
 Giletti \cite{Gil}. What is interesting about this work is that it makes minimal assumptions on the nonlinearity, and takes as a working assumptions the position of the bottom speed with respect to the KPP speed, and the asymptotic behaviour of the wave at infinity. This classification is {\it peu ou prou} what I have reproduced here. Giletti's work identifies three possibilities for the level set $X_\theta(t)$: either $X_\theta(t)=c_*t+O_{t\to+\infty}(1)$ in the case
$c_*>c_K$, $X_\theta(t)=c_*t-3\mathrm{ln}~t/2\lambda_*+O_{t\to+\infty}(1)$ in the case $c_*=c_K$ and the wave has minimal decay, and (the newest part, which unveils a five-legged sheep/pushmi-pullyu) 
$X_\theta(t)=c_*t-\mathrm{ln}~t/2\lambda_*+O_{t\to+\infty}(\mathrm{ln}~t)$  in the case $c_*=c_K$ and the wave has maximal decay. I am, however, exceedingly unhappy that the paper lists my work with Hamel, Nolen and Ryzhik \cite{HNRR} as Reference [13].

\noindent A  detailed analysis of the convergence of the solutions of the Cauchy Problem to the /pushmi-pullyu front (as well as a new proof of all the known cases) is provided by An, Henderson and Ryzhik \cite{AHR1}, \cite{AHR2}. While it is not possible to give a reasonable account of these 120$+$ page-long papers, let me at least explain what they do.  As a travelling wave $\varphi_c$ is strictly monotone one may consider, for any $\lambda$ between 0 and 1, the function 
$$
\eta_c(\lambda)=\varphi_c'\bigl(\varphi_c^{-1}(\lambda)\bigl).
$$
A line of ideas, dating back to the original paper of KPP, consists in noticing that the function 
$$
w(t,x)=u_x(t,x)-\eta_c\bigl(u(t,x)\bigl)
$$
remains of constant sign for all $t>0$, if it is initially of constant sign. In geometric terms, this says that the slope of the solution and that of the travelling wave at a level set are ordered, if they are initially so. This has led to a number of beautiful extensions,
due in particular to Matano and his collaborators, on the "convergence in shape" of $u(t,x)$, that is, on the approach of the solution to the family of travelling waves.  In particular, the poetic concept of "propagating terraces", that is, solutions described by combinations of waves of increasing speeds, thus presenting themselves as terraces, has been thoroughly investigated. The strategy in these works is (a generalisation of) the KPP approach, as well as an extensive use  of the lap number decrease principle (something that was also in germ in the KPP paper).  As an illustration that  I do not  not deem complete, I quote the interesting contribution of Ducrot-Giletti-Matano \cite{DGM} that introduces the concept for KPP type nonlinearities, as well as the study \cite{Pola} of Polacik, which explores the concept for very general nonlinearities. This last approach is very general, but remains qualitative in nature. In particular, it does not pin down a particular translate of the wave to which to which $u(t,x)$ should converge. Especially, Bramson's result looks impossible to attain with this sort of arguments. The main reason is that the equation for $w$ looks intractable. 
The nontrivial contribution of An, Henderson and Ryzhik  is to go beyond this impression and to estimate $w$, using a host of Nash-type inequalities, energy functionals, and an unexpected auxiliary Burgers-type equation.  

\noindent  Whether the An-Henderson-Ryzhik approach carries over in the nonlocal case is unclear, because of the particular algebra generated by the second order diffusion. Trying to understand whether one can fish something out of it would be quite interesting but, for the moment being, Theorem \ref{t5.2.1} is the first to say something about the nonlocal setting.

\subsection*{Open questions}
\noindent A first question concerns the asymptotics of the level sets $X_\theta(t)$ beyond the constant terms, in the case $c=c_K$ and the waves have minimal decay. While, in the case of Fisher-KPP nonlinearities, I have put, perhaps wrongly, this question in the category of problems (see Problem \ref{P4.8.10}), I think - perhaps wrongly, too - that it falls in the range of open questions in the ZFK case.

\noindent I have also left convergence to pushmi-pullyu fronts (or five legged sheep, according to the reader's preferences between Martin du Gard and Lofting) when the diffusion is nonlocal as an interesting open question, as the interaction between the diffusive zone and the region where the solution is nontrivial is certain to be even more subtle than what has been presented here. 

\noindent By Murphy's law (or, if the reader prefers, the principle of conservation of trouble), there is no reason why this transition phenomenon, which appears beautifully in the context of the simple equations $u_t-u_{xx}=f(u)$, or $u_t+\JJ u=f(u)$, should not be present, if not amplified, in all the models for directed propagation that have been presented. I have proposed a rather complete study of the Berestycki-Chapuisat model, but the transition mechanism between KPP and ZFK remains to be elucidated. The same is true for the "road-field" models of various kinds presented in Section \ref{s1.4} of Chapter \ref{Intro}, but in this last class of models one may prefer, as a matter of taste, first understand completely the slow-fast transition issues. 

\noindent Another line of questions is, for a frankly  ZFK nonlinearity, to try to understand the dynamics of a solution emanating from an initial datum that decays like $e^{-rx}$ at infinity, $r>0$ small. It is just as fascinating and open as in the KPP case. I suspect that the setting where $c_*>c_K$ is less difficult, but not much less. Tthe complexity of the problem, at least when the diffusion is $-\partial_{xx}$, can be read from Hamel-Sire \cite{HS}. Adding a nonlocal ingredient is certainly the source of additional difficulties, as regularisation is not present here. Working without it will certainly help the understanding of new features that should have an impact on more standard diffusion equations.

\section{Problems}\label{s5.700}
\begin{problem}
\label{P5.32}
We have always assumed the initial datum $u_0$ to \eqref{e5.7.1} to be $\leq 1$. Show that this assumption can be removed.
\end{problem}
\begin{problem}\label{P5.23} (Stretching the assumptions on the initial datum) Consider the problem
\begin{equation}
\label{e5.7.1}
u_t+\JJ u=f(u)\quad (t>0,\ x\in\RR)
\end{equation}
with $f\geq0$, $f(0)=f(1)=0$, $f'(0)>0$, and $f>0$ on $(0,1)$, and $c_*>c_K$. 
\begin{itemize}
\item[---] If  $u_0$ is the initial datum of \eqref{e3.7.2}, show that the assumption \eqref{e3.7.3} at $-\infty$ can be modified into
$ \di\liminf_{x\to-\infty}u_0(x)>0,$
and that Theorem \ref{t3.7.3} still holds. Show that this assumption can be relaxed even more, as $u_0$ may be allowed to vanish on a fairly big set. How big?
\item[---] Assume that $u_0$ is compactly supported. Show the convergence to two counter-propagating waves, one with speed $c_*$, the other with speed $-c_*$.
\end{itemize}
\end{problem}
\begin{problem}
\label{P5.40}
Study the large time behaviour of $u(t,x)$ when there is $r>0$ small such that  
\begin{equation}
\label{e5.7.10}
e^{-rx}\lesssim u_0(x)\lesssim e^{-rx},
\end{equation} 
\end{problem}
\begin{problem}
\label{P5.41}
Assume $c_*>c_K$, and consider an initial datum for \eqref{e5.7.1} satisfying \eqref{e5.7.10}. Study the solution $u(t,x)$ of \eqref{e5.7.1} in the double limit $r\to\lambda_+(c_*)$ and $t\to+\infty$.
\end{problem}
\begin{problem}
\label{P5.1}
In \eqref{e5.7.1}, assume $f'(0)=f''(0)=0$ and $f>0$ on $(0,1)$, so that we obviously have $c_*>c_K$. Consider an initial datum $u_0$ such that $\di\liminf_{x\to-\infty}u_0(x)=1$, and $supp~u_0\cap\RR_+$ compact. Show that $u(t,x)$ converges, exponentially fast, to a travelling wave of speed $c_*$. How far can one relax the assumption  $supp~u_0\cap\RR_+$ compact? 
\end{problem}
\begin{problem}
\label{P5.15}
Consider the bistable nonlinearity $f(u)=u(1-u)(u-\theta)$, with $0<\theta<1$. The goal of the problem is to provide a proof of existence, and uniqueness (up to translations) of travelling waves for Model \eqref{e1.1.1} with the nnlinearity $f$, that uses a homotopy different from that of Bates et al \cite{BFRW}. Recall that we are solving the problem
\begin{equation}
\label{e5.7.11}
\mathcal{J}\varphi-c\varphi'=f(\varphi)\ (x\in\RR),\quad \varphi(-\infty)=1,\varphi(+\infty)=0.
\end{equation}
Consider the family of functions $(f_t)_{0\leq t\leq\theta}$ given by $f_t(u)=u(u-t)(1-u)$. Problem \eqref{e5.7.11}, with $f(u)=f_0(u)$, has a wave with bottom speed that will be called $c_*^0$, and the corresponding wave  $\varphi_0$.
\begin{itemize}
\item[---] Show that, for $t\in(0,\theta)$, Problem \eqref{e5.7.11} with $f=f_t$  has at most a solution $(c_*^t,\varphi_t)$.
\item[---] For $t\in(0,\theta)$, show an analogue of Lemma \ref{l3.7.3}, with the spaces $B_{w_r}$ replaced by the space of uniformly continuous functions on $\RR$. Deduce that, for all $t_0\in(0,\theta)$, Problem \eqref{e5.7.11} with $f=f_t$ has a solution for $t$ close to $t_0$ as soon as the problem with $f=f_{t_0}$ has one.
\item[---] Show the same result for $t_0=0$. One simply has to be careful at the asymptotic behaviour of the bottom wave with $f=f_0$.
\item[---] Pick $t_0\in[0,\theta]$ and let $(t_n)_n$ a sequence tending to $t_0$. Assume that Problem \eqref{e5.7.11} with $t=t_n$ has a solution $(c_*^n,\varphi_n)$. 
\begin{itemize}
\item [--] Prove that the sequence $(c_*^n)_n$ is bounded. {\rm Hint:} compare a solution of the Cauchy problem with $f=t_{t_n}$ and a solution of the same Cauchy problem with $f=f_0$.
\item[--] Show that the sequence $(c_*^n,\varphi_{t_n})_n$ converges to a solution $(c_*^0,\varphi_0)$.
\end{itemize}
\item[---] Show that the set of $t$'s such that  Problem \eqref{e5.7.11} with $f=f_t$ has a solution is open and closed in $[0,\theta]$. Conclude.
\end{itemize}
\end{problem}
\begin{problem}
\label{P5.43}
In the case $c_*>c_K$, use the stability theorem \ref{t3.7.6} to provide a proof with no computation for Proposition \ref{p3.7.5}. By the way, explain why the strategy of this proposition to get hold of the correct reference framce fails miserably in the case $c_*=c_K$.
\end{problem}
\begin{problem}
\label{P5.42}
Assume $c_*=c_K$, and consider an initial datum for \eqref{e5.7.1} satisfying \eqref{e5.7.10} again. Study the solution $u(t,x)$ of \eqref{e5.7.1} in the double limit $r\to\lambda_*$ and $t\to+\infty$.
\end{problem}
\begin{problem}
\label{P5.5}
In Proposition \ref{p5.3.1}, find asymptotic expansions of $u_\mu+\varphi_{c_*}'$ in terms of $\mu$ valid for $x$ in the whole $\RR_+$.
\end{problem}
\begin{problem}
\label{P5.30}
Study the complex zeroes of the function $\DD_c$. Show, in particular, the existence of an infinite number of zeroes of the form $\xi=i\lambda+\rho(\lambda)$, the function $\rho$ behaving logarithmically as $\lambda\to+\infty$. Derive an asymptotic expansion of $\rho$. It could be useful to study the case where $K=\un_{[-1,1]}$.
\end{problem}
\begin{problem}
\label{P5.31}
Write down all the details of the proof of Theorem \ref{t5.3.4}.
\end{problem}
\begin{problem}
\label{P5.29}
Study the large time behaviour of the discrete version of the ZFK problem 
\begin{equation}
\label{e5.7.30}
\dot u_i=2u_i-u_{i+1}-u_{i-1}+f(u_i).
\end{equation}
In particular, show that there are cases when the bottom speed is the Fisher-KPP speed, while it is strictly larger in other instances. In this last case, study the global stability of the waves.

\noindent  The case of a Fisher-KPP bottom speed is rather involved. The front location up to $O(1)$ terms is derived in a recent preprint  \cite{BFRZ} by Besse, Faye, Zhang and the author of the present book. While the location up to $o_{t\to+\infty}(1)$ terms seems to follow the general scheme set up in Chapter \ref{short_range}, its derivation, which is a work in progress, is expected to pose significant technical issues. When $f$ is of the ZFK type, things are manageable when the bottom speed exceeds the KPP speed. In the opposite case, things have not been written. I expect that no further ideas than those dispayed in Section \ref{s5.3} to be needed, but I do not exclude, either, that some innovation could be needed.
\end{problem}
\begin{problem}
\label{P5.52}
\noindent Let $(f_\e)_{\e>0}$ be a decreasing family of ZFK reaction terms such that $f_0\equiv0$ on an interval of the form $[0,\theta]$, while being $>0$ on $(\theta,1)$. Consider the Cauchy Problem \eqref{e3.7.2} with an initial datum $u_0(x)$  that coincides with $e^{-rx}$ on $[0,+\infty)$. Assume $r>0$ to be small.

\noindent Study the large time behaviour of the solution $u(t,x)$ in the double limit $\e\to0$ and $t\to+\infty$. One should observe a transition between two speeds of propagation, the moment at which this transition happens, and how long it lasts, being of particular interest.
\end{problem}
\begin{problem}
\label{P5.33}
Assume the initial datum $u_0$ in Problem \eqref{e5.7.1} has all its derivatives bounded, with, of course, $\mathrm{supp}~u_0\cap\RR_+$ compact. Show that all the partial derivatives of $u(t,x)$ are bounded.\end{problem}
\begin{problem}
\label{P5.34}
Consider the kernel
$K_\e(x)=\di\frac\beta{\e^3} K(\di\frac{x}\e),
$
and the corresponding problem
$$
u_t+u-K_\e*u=f(u).
$$
The nonlinearity is either of the KPP, or ZFK type. Examine how the travelling waves organise themselves as $\e\to0$, and examine the Cauchy Problem in the double limit $\e\to0$ and $t\to+\infty$.
\end{problem}
\begin{problem}
\label{P5.45}
Figure out a kernel $K_\e$ such that, in the limit $\e\to0$, Problem \eqref{e5.7.1} will become \eqref{e5.7.30}. Examine the organisation of the travelling waves as $\e\to0$ and, if you have sufficient stamina left, study the Cauchy Problem in the double limit $\e\to0$ and $t\to+\infty$. \end{problem}
\begin{problem}
\label{P5.17}
\noindent Consider the Berestycki-Chapuisat model with nonlocal longitudinal diffusion
\begin{equation}
\label{e5.7.15}
\begin{array}{rll}
\partial_tu+\JJ u-\partial_{yy}u+\alpha^2yu=&f(u)\\
u(0,x,y)=&u_0(x,y)\ \hbox{with}\ 0\leq u_0(x,y)\lesssim\mathrm{exp}(-\delta y^2),
\end{array}
\end{equation}
for some small constants  $\e_i>0$ and $\delta>0$. We assume as usual $\mu_1(\alpha)<f'(0)$, and also that the mass of $f$ is so big that we have $c_*>c_K$. The goal of the problem is to show the exponential convergence of $u(t,x,y)$ to a wave with bottom speed.
\begin{itemize}
\item[---] Comment upon the assumptions on $u_0$. To what extent can they be relaxed further in order to grant convergence to the bottom wave? 
\item[---] Consider $r\in\bigl(\lambda_-(c_*),\lambda_+(c_*)\bigl)$ and $\delta\leq1/8$. Set $w_{r,\delta}(x,y)=e^{-\delta y^2}/(1+e^{rx})$ and let $B_{w_{r,\delta}}(\RR^2)$ denote the set of all functions $f(x,y)$ on $\RR^2$ such that $f/w_{r,\delta}$ is bounded and uniformly continuous. Denote by $\MM_{c_*}$ the operator as
$$
\MM_{c_*}u=\JJ u-c_*u_x-\partial_{yy}u+\\bigl(\alpha y^2-f'(\varphi_{c_*}\bigl)u.
$$
Show the existence of a linear form $e*$ acting on  $B_{w_{r,\delta}}(\RR^2)$ and $\omega>0$ such that, if $u_0\in B_{w_{r,\delta}}(\RR^2)$we have: 
$$\Vert e^{-t\MM_{c_*}}-<\!e_*,u_0\!>\varphi_{c_*}\Vert_{w_{r,\delta}}=O(e^{-\omega t}).
$$
\item[---] Adapt the arguments of Section \ref{s5.2.2} to prove the analogue of Theorem \ref{t3.7.3}, that is: there is $x_\infty\in\RR$ such that $u(t,x,y)=\varphi_{c_*}(x-c_*t+x_\infty,y)+O(e^{-\omega t})$.
\end{itemize}
\end{problem}

\chapter{Final remarks}\label{Final}
\noindent   A complete and self-contained study of reaction-diffusion front propagation with nonlocal diffusion in homogeneous one-dimensional media has been given. I have the hope -- and firm belief -- that the ideas displayed will be useful for the study of problems in several space dimensions, but where the propagation is forced, by some mechanism, to occur in a single direction. In doing so, I have kept outside the scope of the book three important directions, that await a treatment as complete as the one given here to homogeneous models in one space dimension. These three directions, listed below, are of course not limitative, as new models keep pouring in as this book is being written.
\section{Is linear spreading the only possible propagation mode?} 
\noindent Far from that, of course.  It suffices to depart from the realm of compactly supported kernels to obtain accelerated modes of propagation, provided that the tail of the kernel does not decay too fast at infinity.  If, for instance, we take $K(z)$ of the order $\vert z\vert^{-(1+2\alpha)}$, with $\alpha>0$, so as to make it integrable at infinity, the $\gamma$ level set of $u(t,x)$ will move like $e^{f'(0)t/(1+2\alpha)}$ as $t\to+\infty$. This fact has long been identified by physicists or theoretical ecologists, see for instance - and without attempt to find the earliest reference -  Kot-Lewis-Van de Driessche \cite{KLV}. Rigorous studies of spreading speed date back to Cabr\'e and the author \cite{CR} when the diffusion is given by the fractional Laplacian, and Garnier \cite{Gar} for smooth kernels of the type studied here. A very serious look at these issues is given in the work of A.-C. Chalmin \cite{ACC}, including a study of what happens in the models of front propagation directed by a line of fractional diffusion. Acceleration by long range diffusion  is not, and by far again, the only acceleration mechanism. I mention the beautiful work of Bouin, Calvez and Nadin \cite{BCN} on accelerating solutions of kinetic-like equations, where the acceleration mechanism comes from the fact that the velocities can be chosen in an unbounded set. Coming back to models of the type \eqref{e1.1.1} with slowly decaying kernels, the full asymptotic behaviour of the solutions is still to be understood, although many partial results exist. The (almost finalised) ongoing work (Chalmin and the author \cite{ChR}) answers the question for homogeneous media.
\section{Nonlocal 1D models in heterogeneous media} 
\noindent There is not a long way to depart from the initial model \eqref{e1.1.1} to be in an almost uncharted territory. Consider, for instance, the equation
\begin{equation}
\label{e5.4.1}
u_t+\mathcal{J}u=f(x,u),\quad(t>0,\ x\in\RR),
\end{equation}
the function $f$ satisfying, uniformly in $x$, the KPP or ZFK assumptions. Even an innocent looking nonlinearity of the form $f(x,u)=\mu(x)u-u^2$, $\mu>0$, is problematic as soon as $\mu$ does not have some structure.
The situation can even be complicated by assuming $f$ to be of the KPP or the ZFK type, depending on $x$, with all the degeneracies that it may entail. One could also make  $\mathcal{J}$ inhomogeneous by setting
$$
\JJ u=\di\int_\RR K(x,y)\biggl(u(x)-u(y)\biggl)~dy.
$$
 In this last case, a reasonable assumption is $K_1(x-y)\leq K(x,y)\leq K_2(x-y)$, where the kernels $K_i$ are the honest even, compactly supported kernes that have been studied in the book. One could, however, be a little more perverse and assume that the function $y\mapsto K(x,y)$ is compactly supported in $\vert x-y\vert$ for some values of $x$, while decaying like $\vert x-y\vert^{-\gamma}$, for some $\gamma>1$. Many interesting partial results on exist, and I will not even try to give an account of them. 

\noindent When the diffusion is Gaussian, a lot is known. The book of Berestycki-Hamel \cite{BH-book}, which has been preceded by a wealth of works of these two authors, gives a general theory of how things organise. 
In particular, they identify the correct notion that should replace travelling waves: the suitable concept is that of {\it transition wave}, in other words, an eternal solution of the equation that connects -- perhaps in quite a distorted fashion - the rest points, and such that, loosely speaking, the width of the front is uniformly bounded in time. Many of these objects have been constructed, a remarkable fact being that they also appear in homogeneous problems. Hamel-Nadirashvili \cite{HN} construct a wealth of these. Nevertheless, fundamental questions concerning their attractivity are still wide open. 

\noindent Transition fronts, in the spirit of those constructed for Gaussian diffusion, have been constructed for integral diffusions in inhomogeneous media, see for instance Shen \cite{Sh}. A complete theory for nonlocal equations is still to be devised, and I expect inhomogeneous nonlocalities to generate even wilder behaviours than the ones known for Gaussian diffusion.

\noindent Even a seemingly benign  time dependence in the coefficients has the potential to sow chaos in the logarithmic harmony that we have depicted. Consider once again the simplest looking model
\begin{equation}
\label{e6.1.1} 
u_t-u_{xx}=\mu(t)u-u^2.
\end{equation}
When $\mu(t)$ is periodic, everything is known, either with the local diffusion operator $-\partial_{xx}$ or the nonlocal diffusion $\JJ$ (see the problems scattered in the various chapters). 

\noindent However, when $\mu(t)$ is time increasing and tends slolwly to 1, that is, for instance, $\mu(t)=1-t^{-\gamma}\bigl(1+o_{t\to+\infty}(1)\bigl)$, a forthcoming work of Rossi, Ryzhik and the author \cite{RRR} gives an asymptotic expansion of $X_\theta(t)-c_*t$, and show that it is much larger than logarithmic, in fact it is slightly smaller than $t^{1/3}$.  This raises the following interesting question: does there exist a universal bound for the discrepancy $X_\theta(t)-c_*t$ (if $u(t)$ tends to 1) or, more generally, between $X_\theta(t)$ and the $\theta$-level set of the solution of the heat equation $u_t-u_{xx}=\mu(t)u$? The work of Fang-Zeitouni \cite{FZ} on the branching Brownian motion, completed (and interpreted for the Fisher-KPP equation) by Nolen, Ryzhik and the author \cite{NRR3} support a conjecture of the form $\vert X_\theta(t)-c_*t\vert\lesssim t^{1/3}$. This fascinating question, as well as how it extends to the nonlocal setting, is very much open.
\section{Multidimensional models} 
\noindent As far as the spreading velocity is concerned, and still with the Gaussian diffusion, an important step has been made by Berestycki-Nadin \cite{BN}. In a very general multi-dimensional medium, they identify a lower and an upper asymptotic velocity,
given by the solutions of very degenerate eikonal equations (whose study is in itself a fascinating open question) that coincide with the known expressions. Whether this beautiful theory holds for integral diffusion of the type studied in this book
is not known at the moment. 

\noindent The sharp asymptotics aspects remain quite far from a general understanding.   One could of course try to generalise the main results of chapters \ref{short_range} and of the present one to models that would look like the original model of this book, namely Model \eqref{e1.1.1}, but with this time $x\in\RR^N$ and the kernel $K$ radially symmetric. The nonlinearity is either of the Fisher-KPP type or of the ZFK type, and the initial datum is compactly supported. When the diffusion is given by the standard Laplacian, that is, when we are interested in
\begin{equation}
\label{e5.4.2}
u_t-\Delta u=f(u)\ \ (t>0,\ x\in\RR^N), \quad u(0,x)\ \hbox{compactly supported}
\end{equation}
the large time behaviour is understood. If $(r,\Theta)$ are the polar coordinates of $\RR^N$, the $\gamma$-level set of $u(t,x)$ is an almost spherical surface. In the Fisher-KPP case it has (Rossi, Roussier and
 the author of these lines \cite{RRR1}) the form
$r=c_*t-\di\frac{N+2}{2\lambda_*}\mathrm{ln~}t+s_\infty(\Theta),
$
while, in the ZFK case (independent proofs by Nara \cite{Nara} and Roussier and the author of these lines \cite{RR2}) it has the form
$r=c_*t-\di\frac{N-1}{2\lambda_*}\mathrm{ln~}t+s_\infty(\Theta).
$
I should say that, in these two works, the function $f$ is bistable rather than ZFK, but I would be quite surprised if the result did not extend easily to ZFK nonlinearities. In all cases, the function $s_\infty$ is Lipschitz, and whether this regularity can be improved is something unknown to me. Another interesting issue
is the existence of nonspherical initial data giving raise to a constant function $s_\infty$; this is also unknown to me - the only thing I know is that such a behaviour would definitely be of another five-legged sheep type. 

\noindent Generalising the above results to the nonlocal case does not seem to the hardest undertaking. However, if the kernel is not assumed to be radially symmetric anymore, a radically different behaviour is to be expected. The spreading speed is indeed given by the Freidlin-G\"artner formula  \eqref{e4.7.4}; this is a recent result of Finkelshtein-Kondratiev-Tkachov \cite{FKT}. What happens, by the way, beyond the Freidlin-G\"artner term is an outstanding question that is far from being understood, to the exception of a beautiful preprint of Shabani \cite{Sha} which locates, up to $O(1)$ term, the level sets of the solution in a Fisher-KPP equation in a periodic medium. The sharp asymptotic behaviour is an even more open question as one departs from the simple geometries 
depicted so far. Curvy domains, such as those stidies in Berestycki-Bouhours-Chapuisat \cite{BBC}, or Fisher-KPP models in non-euclidian spaces, such as the hyperbolic plane, as envisaged by Matano, Punzo and Tesei \cite{MPT} open the door to fascinating considerations.

\end{document}